\documentclass[11pt,letterpaper,twoside,reqno,nosumlimits]{amsart}

\allowdisplaybreaks

\usepackage[usenames,dvipsnames]{xcolor}
\usepackage{fancyhdr}
\usepackage{amsmath,amsfonts,amsbsy,amsgen,amscd,amssymb,amsthm}
\usepackage{url}

\usepackage{mathtools}

\usepackage[font=small,margin=0.25in,labelfont={bf},labelsep={colon}]{caption}

\usepackage{subcaption}

\usepackage{tikz}
\usepackage{microtype}
\usepackage{enumitem}

\definecolor{dark-gray}{gray}{0.3}
\definecolor{dkgray}{rgb}{.4,.4,.4}
\definecolor{dkblue}{rgb}{0,0,.5}
\definecolor{medblue}{rgb}{0,0,.75}
\definecolor{rust}{rgb}{0.5,0.1,0.1}

\usepackage[colorlinks=true]{hyperref}

\hypersetup{urlcolor=rust}
\hypersetup{citecolor=dkblue}
\hypersetup{linkcolor=dkblue}

\usepackage{setspace}

\usepackage{graphicx}
\usepackage{booktabs,longtable,tabu} 
\usepackage{multirow} 

\usepackage{float}

\usepackage[full]{textcomp}

\usepackage[scaled=.98,sups,osf]{XCharter}
\usepackage[scaled=1.04,varqu,varl]{inconsolata}
\usepackage[type1]{cabin}
\usepackage[charter,vvarbb,scaled=1.07]{newtxmath}
\usepackage[cal=boondoxo]{mathalfa}
\usepackage[T1]{fontenc}

\usepackage{bm} 

\graphicspath{{figures/}}

\newtheorem{theorem}{Theorem}[section]

\theoremstyle{definition}

\newtheorem{remark}[theorem]{Remark}

\numberwithin{equation}{section} 
\numberwithin{figure}{section}
\numberwithin{table}{section}

\floatstyle{plaintop}
\newfloat{recipe}{thp}{lor}
\floatname{recipe}{Recipe}
\numberwithin{recipe}{section}

\providecommand{\mathbold}[1]{\bm{#1}}  
                                
\newcommand{\econst}{\mathrm{e}}

\newcommand{\diff}[1]{\mathrm{d}{#1}}
\newcommand{\idiff}[1]{\, \diff{#1}}

\newcommand{\vct}[1]{\mathbold{#1}}
\newcommand{\mtx}[1]{\mathbold{#1}}

\newcommand{\triplenorm}[1]{{\left\vert\kern-0.25ex\left\vert\kern-0.25ex\left\vert #1
    \right\vert\kern-0.25ex\right\vert\kern-0.25ex\right\vert}}

\usepackage[numbers]{natbib}

\newcommand{\om}{\omega}
\newcommand{\vom}{\vct{\omega}}
\newcommand{\vu}{\vct{u}}
\newcommand{\vpsi}{\vct{\psi}}

\begin{document}

\title[Potential singularity of Euler with degenerate viscosity]{Potential singularity formation of incompressible axisymmetric Euler equations with degenerate viscosity coefficients}

\author[T. Y. Hou and D. Huang]{Thomas Y. Hou and De Huang}

\date{\today}
 
\begin{abstract}
In this paper, we present strong numerical evidences that the incompressible axisymmetric Euler equations with degenerate viscosity coefficients and smooth initial data of finite energy develop a potential finite-time locally self-similar singularity at the origin. An important feature of this potential singularity is that the solution develops a two-scale traveling wave that travels towards the origin. The two-scale feature is characterized by the scaling property that the center of the traveling wave is located at a ring of radius $O((T-t)^{1/2})$ surrounding the symmetry axis while the thickness of the ring collapses at a rate $O(T-t)$. The driving mechanism for this potential singularity is due to an antisymmetric vortex dipole that generates a strong shearing layer in both the radial and axial velocity fields.  Without the viscous regularization, the $3$D Euler equations develop a sharp front and some shearing instability in the far field. On the other hand, the Navier--Stokes equations with a constant viscosity coefficient regularize the two-scale solution structure and do not develop a finite-time singularity for the same initial data. 
\end{abstract}

\maketitle

\section{Introduction}

The three-dimensional ($3$D) incompressible Euler equations in fluid dynamics describe the motion of inviscid incompressible flows. 
Despite their wide range of applications, the question regarding the global regularity of the $3$D Euler equations with smooth initial data of finite energy has remained open \cite{majda2002vorticity}. 
The main difficulty associated with the global regularity of the $3$D Euler equations is the presence of vortex stretching, which is absent in the corresponding $2$D problem. 
In 2014, Luo and Hou \cite{luo2014potentially,luo2014toward} presented strong numerical evidences that the $3$D axisyemmtric Euler equations with smooth initial data and boundary develop potential finite-time singular solutions at the boundary. The presence of the boundary and the symmetry properties of the initial data seem to play a crucial role in generating a sustainable finite-time singularity of the $3$D Euler equations.

In this paper, we present strong numerical evidence that the incompressible axisymmetric Euler equations with smooth degenerate viscosity coefficients and smooth initial data of finite energy seem to develop a two-scale locally self-similar singularity. 
Unlike the Hou--Luo blowup scenario, the potential singularity for the Euler equations with degenerate viscosity coefficients occurs at the origin. Without the viscous regularization, the $3$D Euler equations develop an additional small scale characterizing the thickness of the sharp front. The degenerate viscosity coefficients are designed to select a stable locally self-similar two-scale solution structure and stabilize the shearing induced instability in the far field.

We also study the Navier--Stokes equations with a constant viscosity coefficient using the same initial data. Our study shows that the Navier--Stokes equations will regularize the two-scale solution structure and destroy the strong nonlinear alignment of the vortex stretching term. Moreover, we will present some preliminary numerical results indicating that the $3$D Euler equations seem to develop a three-scale solution structure. The rapid collapse of the thickness of the sharp front makes it extremely difficult to resolve the potential Euler singularity numerically.

\subsection{Major features of the potential blowup and the blowup mechanism}

One of the important features of the potential blowup solution is that it develops a two-scale traveling solution approaching the origin. We denote by $u^\theta$ and $\omega^\theta$ the angular velocity and angular vorticity, respectively, and define $u_1 = u^\theta/r$ and $\omega_1 = \omega^\theta/r$ with $r = \sqrt{x^2+y^2}$. Let $(R(t),Z(t))$ be the location where $u_1$ achieves its global maximum in the $rz$-plane. The traveling wave is centered at a ring  with radius $R(t)$ surrounding the symmetry axis $r=0$ and the thickness of the ring is roughly of order $Z(t)$. 
The two-scale traveling wave solution is characterized by the scaling property that $R(t)= O((T-t)^{1/2})$ and $Z(t) = O(T-t)$. 
Another important feature is that the odd symmetry (in $z$) of the initial data of $\om_1$ induces a vortex dipole and an antisymmetric local convective circulation. This convective circulation is the cornerstone of our blowup scenario, as it has the desirable property of pushing the solution near $z=0$ towards the symmetry axis $r=0$.

An important guiding principle for constructing our initial data is to enforce a strong nonlinear alignment of vortex stretching. First of all, the vortex dipole induces a negative radial velocity $u^r$ near $z=0$, i.e. $u^r = - r \psi_{1,z} < 0$, which implies $ \psi_{1,z} >0$. Moreover, $\psi_{1,z}(R(t),z,t)$ is a monotonically decreasing function of $z$ and is relatively flat near $z=0$. Through the vortex stretching term $2 \psi_{1,z} u_1$ (see \eqref{eq:as_NSE_1_a}), the large value of $\psi_{1,z}$ near $z=0$ induces a traveling wave for $u_1$ that approaches $z=0$ rapidly. The strong nonlinear alignment in vortex stretching overcomes the stabilizing effect of advection in the upward $z$ direction (see e.g. \cite{hou2008dynamic,lei2009stabilizing}). The oddness of $u_1$ in $z$ then generates a large positive gradient $u_{1z}$, which contributes positively to the rapidly growth of $\omega_1$ through the vortex stretching term $2 u_1 u_{1z}$ (see \eqref{eq:as_NSE_1_b}). The rapid growth of $\omega_1$ in turn feeds back to the rapid growth of $\psi_{1,z}$. The whole coupling mechanism described above forms a positive feedback loop.

Moreover, we observe that the $2$D velocity field $(u^r(t),u^z(t))$ in the $rz$-plane forms a closed circle right above $(R(t),Z(t))$. The corresponding streamline is then trapped in the circle region in the $rz$-plane. This local circle structure of the $2$D velocity field is critical in stabilizing the blowup process, as it keeps the bulk parts of the $u_1,\om_1$ profiles traveling towards the origin instead of being pushed upward. The strong shear layer in $u^r$ and $u^z$ generates a sharp front for $u_1$ in both $r$ and $z$ directions.

Another important feature of our initial data is that it generates a local hyperbolic flow in the $rz$-plane.
Due to the odd symmetry of $u_1$, $u_1$ is almost zero in the region near $z=0$. The strong upward transport near $r=0$ makes $u_1$ really small in a no-spinning region between the sharp front of $u_1$ and the symmetry axis $r=0$. 
Within this no-spinning region, the angular velocity $u^\theta = r u_1$ is almost zero, which implies that there is almost no spinning around the symmetry axis. The flow effectively travels upward along the vertical direction inside this no-spinning region. Outside this no-spinning region, $u_1$ becomes very large and the flow spins rapidly around the symmetry axis. Moreover, the streamlines induced by the velocity field travel upward along the vertical direction and then move outward along the radial direction. The local blowup solution resembles the structure of a tornado. 
For this reason, we also call the potential singularity of the Euler equations with variable viscosity coefficients ``a tornado singularity''.

\subsection{Asymptotic scaling analysis }
To confirm that the potential singular solution develops a locally self-similar blowup, we perform numerical fitting of the growth rates for several physical quantities. Our study shows that the maximum of the vorticity vector grows like $O((T-t)^{-3/2})$, and $R(t)=O((T-t)^{1/2})$, $Z(t)=O(T-t)$. The fact that $\| \nabla {\vu} \|_{L^\infty} \geq \| \nabla {\vom} \|_{L^\infty} \geq C_0(T-t)^{-3/2}$ gives that $\int_0^T \| \nabla {\vu(t)} \|_{L^\infty} dt = \infty$, which implies that the solution could potentially develop a finite-time singularity \cite{majda2002vorticity}. 

We have also performed an asymptotic scaling analysis to study the scaling properties of potential locally self-similar blowup solution. By balancing the scales in various terms, we show that $u_1$ and $\psi_{1,z}$ must blow up with the rate $O(1/(T-t))$ if there is a locally self-similar blowup. Due to the conservation of total circulation $\Gamma = r^2 u_1$ and the degeneracy of the viscosity coefficients, we show that $\Gamma$ remains $O(1)$ at $(R(t),Z(t))$. This property and the scaling property that $u_1 = O(1/(T-t))$ imply that $R(t) = O((T-t)^{1/2})$. Moreover, the balance between the vortex stretching term and the degenerate viscosity term suggests that $Z(t) = O(T-t)$. Similarly, we can show that $\omega_1$ blows up like $O(1/(T-t)^2)$. In terms of the original physical variables, the vorticity vector blows up like $O(1/(T-t)^{3/2})$ and the velocity field blows up like $O(1/(T-t)^{1/2})$. The results obtained by our scaling analysis are consistent with our numerical fitting of the blowup rates for various quantities. 

\subsection{Computational challenges}
The two-scale nature of the potential singular solution presents considerable challenges in obtaining a well-resolved numerical solution for the Euler equations with variable viscosity coefficients. 
To resolve this potential two-scale singular solution, we design an adaptive mesh by constructing two adaptive mesh maps for $r$ and $z$ explicitly. More specifically, we construct our mapping densities in the near field (phase $1$: resolving the $Z(t)$ scale), the intermediate field (phase $2$: resolving the $R(t)$ scale), and the far field (phase $3$: resolving the $O(1)$ scale) with a transition phase in between. We then allocate a fixed number of grid points in each phase and update the mesh maps dynamically according to some criteria. 
This adaptive mesh strategy achieves a highly adaptive mesh with the smallest mesh size of order $O(10^{-10})$. 
Our adaptive mesh strategy is more complicated than the one presented in \cite{luo2014potentially,luo2014toward} since we have a two-scale traveling wave. 

We use a $2$nd-order finite difference method to discretize the spatial derivatives, and a $2$nd-order explicit Runge--Kutta method to discretize in time.
We choose an adaptive time-step size according to the standard time-stepping stability constraint with the smallest time-step size of order $O(10^{-11})$. 
We adopt the $2$nd-order B-spline based Galerkin method developed in \cite{luo2014potentially,luo2014toward} to solve the Poisson equation for the stream function. We also design a $2$nd-order filtering scheme to control some mild oscillations in the tail region. The overall method is $2$nd-order accurate. We have performed resolution study to confirm that our method indeed gives $2$nd-order accuracy in the maximum norm.

\subsection{Comparison with results obtained in two subsequent papers}

Inspired by the work presented in this manuscript, the first author of this paper investigated potential singular behavior of the $3$D Euler and Navier--Stokes equations using a different but relatively simple initial condition in two subsequent papers \cite{Hou-euler-2021,Hou-nse-2021}. Although the solution presented in this paper and the solutions obtained in \cite{Hou-euler-2021,Hou-nse-2021} share many similar properties, there are some important differences between the two blowup scenarios. One important difference is that the potential Euler singularity considered in  \cite{Hou-euler-2021} is essentially a one-scale traveling wave singularity instead of  a two-scale traveling wave singularity considered in this paper. More importantly, the scaling properties of the potential Euler singularity presented in \cite{Hou-euler-2021} are compatible with those of the Navier--Stokes equations. In \cite{Hou-nse-2021}, it is shown that the maximum vorticity of the Navier--Stokes equations using a relatively large constant viscosity coefficient $\nu = 5\times10^{-3}$ increases by a factor of $10^7$ relatively to its initial maximum vorticity. In comparison, the maximum vorticity of the Navier--Stokes equations with $\nu = 10^{-5}$ using the initial condition considered in this paper has increased less than $2$. Another important feature of the initial condition considered in  \cite{Hou-euler-2021,Hou-nse-2021} is that it decays rapidly in the far field. As a result, its solution does not suffer from the shearing instability that we observe in this paper and there is no need to use a degenerate viscosity or a low pass filter to stabilize the shearing instability in the far field as we did in this paper.

\subsection{Review of related works in the literature}

One of the best known results for the $3$D Euler equations is the Beale--Kato--Majda non blowup criteria \cite{majda2002vorticity}, which states that the smooth solution of the $3$D Euler equations ceases to be regular at some finite time $T$ if and only if $\int_0^T \| \nabla {\vom} \|_{L^\infty} dt = \infty$. In \cite{cfm1996}, Constantin--Fefferman--Majda showed that the local geometric regularity of the vorticity vector near the region of maximum vorticity could lead to the dynamic depletion of vortex stretching, thus preventing a potential finite-time singularity (see also\cite{dhy2005}). An exciting recent development is the work by Elgindi \cite{elgindi2019finite} (see also \cite{elgindi2021stable}) who proved that the $3$D axisymmetric Euler equations develop a finite-time singularity for a class of $C^{1,\alpha}$ initial velocity with no swirl. There have been a number of very interesting results inspired by the Hou--Lou blowup scenario \cite{luo2014potentially,luo2014toward}, see e.g. \cite{kiselev2014small,choi2014on,kryz2016,chen2019finite2,chen2021finite}
and the excellent survey article \cite{kiselev2018,Elg22}.
 
There has been a number of previous attempts to search for potential Euler singularities numerically. In \cite{gs1991}, Grauer--Sideris presented numerical evidences that the $3$D axisymmetric Euler equations with a smooth initial condition could develop a potential finite-time singularity away from the symmetry axis. In \cite{es1994} , E--Shu studied the potential development of finite-time singularities in the $2$D Boussinesq equations with initial data completely analogous to those of \cite{gs1991,ps1992} 
and found no evidence for singular solutions. Another well-known work on potential Euler singularity is the two anti-parallel vortex tube computation by Kerr in \cite{kerr1993}. In \cite{hl2006}, Hou--Li repeated the computation of \cite{kerr1993} with higher resolutions and only observed double exponential growth of the maximum vorticity in time.
In \cite{bp1994}, Boratav--Pelz presented numerical evidence that the $3$D Euler equations with Kida's high-symmetry initial data would develop a finite-time singularity. In \cite{hl2008}, Hou--Li also repeated the computation of \cite{bp1994} and found
that the singularity reported in \cite{bp1994} is likely an artifact due to under-resolution. 
In \cite{brenner2016euler}, Brenner--Hormoz--Pumir considered an iterative mechanism for potential singularity formation of the $3$D Euler equations.
We refer to an excellent review paper \cite{gibbon2008} for more discussion on potential Euler singularities.

The rest of the paper is organized as follows. In Section \ref{sec:setup}, we describe the setup of the problem in some detail, including the analytic construction of our initial data and the variable viscosity coefficients. In Section \ref{sec:first_sign}, we report the major findings of our numerical results, including the first sign of singularity and the main features of the potentially singular solution. 
In Sections \ref{sec:original_NSE} and \ref{sec:Euler}, we present a comparison study with the standard Navier--Stokes equations with constant viscosity coefficient and the inviscid Euler equations, respectively. 
In Section \ref{sec:scaling_study}, we investigate the asymptotic blowup scaling both numerically and by asymptotic scaling analysis. Some concluding remarks are made in Section \ref{sec:conclusion}.

\section{Description of the Problem}\label{sec:setup}
We study the $3$D incompressible Euler equations with variable viscosity coefficients:
\begin{equation}\label{eq:NSE_vc}
\begin{split}
\vu_t + \vu\cdot \nabla\vu &= -\nabla p + \nabla\cdot(\nu \nabla \vu),\\
\nabla \cdot \vu &= 0,
\end{split}
\end{equation}
where $\vu = (u^x,u^y,u^z)^T:\mathbb{R}^3\mapsto\mathbb{R}^3$ is the $3$D velocity vector, $p:\mathbb{R}^3\mapsto \mathbb{R}$ is the scalar pressure, $\nabla = (\partial_x,\partial_y,\partial_z)^T$ is the gradient operator in $\mathbb{R}^3$, and $\nu:\mathbb{R}^3\mapsto \mathbb{R}^{3\times 3}$ is the variable viscosity tensor.  In the inviscid case (i.e. $\nu = \mtx{0}$), \eqref{eq:NSE_vc} reduce to the $3$D Euler equations. By taking the curl on both sides, the equations can be rewritten in the equivalent vorticity form 
\begin{equation*}\label{eq:vorticity_form}
\vom_t + \vu\cdot \nabla\vom = \vom\cdot \nabla \vu + \nabla \times (\nabla\cdot(\nu \nabla \vu)), 
\end{equation*}
where $\vom = \nabla\times \vu$ is the $3$D vorticity vector. The velocity $\vu$ is related to the vorticity via the vector-valued stream function $\vpsi$:
\[-\Delta \vpsi = \vom,\quad \vu = \nabla\times \vpsi.\]

\subsection{Axisymmetric Euler equations} We will study the potential singularity formulation for the axisymmetric Euler equations. In the axisymmetric scenario, it is convenient to rewrite equations~\eqref{eq:NSE_vc} in cylindrical coordinates. Consider the change of variables 
\[x = r\cos\theta,\quad y = r\sin\theta,\quad z = z,\]
and the decomposition 
\[\vct{v}(r,z) = v^r(r,z)\vct{e}_r + v^\theta(r,z)\vct{e}_\theta + v^z(r,z)\vct{e}_z,\]
\[\vct{e}_r = \frac{1}{r}(x,y,0)^T,\quad \vct{e}_\theta = \frac{1}{r}(-y,x,0)^T,\quad \vct{e}_z = (0,0,1)^T,\]
for radially symmetric vector-valued functions $\vct{v}(r,z)$. The Euler equations ~\eqref{eq:NSE_vc} with variable viscosity coefficients can be rewritten in the axisymmetric form:
\begin{subequations}
\begin{align}
u^\theta_t+u^ru^\theta_r+u^zu^\theta_z &=-\frac{1}{r}u^ru^\theta + f_u,\label{eq:as_NSE_a}\\
 \om^\theta_t+u^r\om^\theta_r+u^z\om^\theta_z &=\frac{2}{r}u^\theta u^\theta_z+\frac{1}{r}u^r\om^\theta + f_\om,\label{eq:as_NSE_b}\\
 -\left(\Delta-\frac{1}{r^2}\right)\psi^\theta &= \om^\theta,\label{eq:as_NSE_c}\\
 u^r=-\psi^\theta_z,\quad u^z&=\frac{1}{r}(r\psi^\theta)_r, \label{eq:as_NSE_d}
\end{align}
\label{eq:axisymmetric_NSE}
\end{subequations}
where $f_u$ and $f_\om$ are the viscosity terms. The incompressibility condition 
\[\frac{1}{r}(ru^r)_r + u^z_z = 0\]
is automatically satisfied owing to equations~\eqref{eq:as_NSE_d}. Equations~\eqref{eq:axisymmetric_NSE}, together with appropriate initial and boundary conditions, completely determine the evolution of the solution flow.

To determine the viscosity terms, we choose the variable viscosity tensor such that 
\[\nu = \text{diag}(\nu^r,\nu^r,\nu^z)\]
in the cylindrical coordinates, where $\nu^r = \nu^r(r,z),\nu^z = \nu^z(r,z)$ are functions of $(r,z)$. We remark that this is equivalent to choosing $\nu = \text{diag}(\nu^x,\nu^y,\nu^z)$ with $\nu^x=\nu^y=\nu^r$ in the Euclidean coordinates $(x,y,z) = (r\cos\theta,r\sin\theta,z)$. In order for $\nu$ to be a smooth function in the primitive coordinates $(x,y,z)$, we require that $\nu^r(r,z),\nu^z(r,z)$ are even functions of $r$ with respect to $r=0$. With this choice of the variable viscosity coefficients, $f_u$ and $f_\om$ have the expressions
\begin{align*}
f_u &= \nu^r\left(u^\theta_{rr} + \frac{1}{r} u^\theta_{r} -\frac{1}{r^2}u^\theta\right) + \nu^z u^\theta_{zz} + \nu^r_ru^\theta_r + \nu^z_zu^\theta_z,\\
f_\om &= \nu^r\left(\om^\theta_{rr} + \frac{1}{r} \om^\theta_{r} -\frac{1}{r^2}\om^\theta\right) + \nu^z \om^\theta_{zz} + \nu^r_r\om^\theta_r + \nu^z_z\om^\theta_z \\
 &\ \quad + \nu^r_z\left(u^r_{rr} + \frac{1}{r}u^r_r - \frac{1}{r^2}u^r\right) + \nu^z_zu^r_{zz} - \nu^r_r\left(u^z_{rr} + \frac{1}{r}u^z_r\right) - \nu^z_ru^z_{zz} \\
 &\ \quad + \nu^r_{rz}u^r_r + \nu^z_{zz}u^r_z  - \nu^r_{rr}u^z_r - \nu^z_{rz}u^z_z.
\end{align*}

The axisymmetric Euler equations~\eqref{eq:axisymmetric_NSE} have a formal singularity at $r = 0$, which sometimes is inconvenient to work with. To remove the singularity, Hou--Li \cite{hou2008dynamic} introduced the variables 
\[u_1 = u^\theta/r,\quad \om_1=\om^\theta/r,\quad \psi_1 = \psi^\theta/r \]
and transform equations~\eqref{eq:axisymmetric_NSE} into the form 
\begin{subequations}\label{eq:axisymmetric_NSE_1}
\begin{align}
u_{1,t}+u^ru_{1,r}+u^zu_{1,z} &=2u_1\psi_{1,z} + f_{u_1},\label{eq:as_NSE_1_a}\\
\om_{1,t}+u^r\om_{1,r}+u^z\om_{1,z} &=2u_1u_{1,z} + f_{\om_1},\label{eq:as_NSE_1_b}\\
 -\left(\partial_r^2+\frac{3}{r}\partial_r+\partial_z^2\right)\psi_1 &= \om_1,\label{eq:as_NSE_1_c}\\
 u^r=-r\psi_{1,z},\quad u^z&=2\psi_1 + r\psi_{1,r}, \label{eq:as_NSE_1_d}
\end{align}
\end{subequations}
where the viscosity terms $f_{u_1},f_{\om_1}$ are given by
\begin{subequations}\label{eq:viscosity}
\begin{align}
f_{u_1} = \frac{1}{r}f_u &= \nu^r\left(u_{1,rr} + \frac{3}{r}u_{1,r}\right) + \nu^zu_{1,zz} + \frac{1}{r}\nu^r_r u_1 + \nu^r_ru_{1,r} + \nu^z_zu_{1,z}, \label{eq:viscosity_u1}\\
f_{\om_1} = \frac{1}{r}f_\om &= \nu^r\left(\om_{1,rr} + \frac{3}{r}\om_{1,r}\right) + \nu^z\om_{1,zz} + \frac{1}{r}\nu^r_r \om_1 + \nu^r_r\om_{1,r} + \nu^z_z\om_{1,z}  \nonumber\\
&\ \quad + \frac{1}{r}\left(\nu^r_z\left(u^r_{rr} + \frac{1}{r}u^r_r - \frac{1}{r^2}u^r\right) + \nu^z_zu^r_{zz} - \nu^r_r\left(u^z_{rr} + \frac{1}{r}u^z_r\right) - \nu^z_ru^z_{zz}\right) \nonumber \\
 &\ \quad + \frac{1}{r}\Big(\nu^r_{rz}u^r_r + \nu^z_{zz}u^r_z  - \nu^r_{rr}u^z_r - \nu^z_{rz}u^z_z\Big). \label{eq:viscosity_w1}
\end{align}
\end{subequations} 
Note that $f_u,f_\om$ are even functions of $r$, thus $f_{u_1},f_{\om_1}$ are well defined as long as $f_u,f_\om$ are smooth. Therefore, $u_1$, $\om_1$, and $\psi_1$ are well defined as long as the corresponding solution to equations~\eqref{eq:axisymmetric_NSE} remains smooth. 

\subsection{Settings of the solution}\label{sec:settings} We will solve the transformed equations~\eqref{eq:axisymmetric_NSE_1} in the cylinder region
\[\mathcal{D} = \{(r,z): 0\leq r\leq 1\}.\]
In particular, we will enforce the following properties of the solution:
\begin{subequations}\label{eq:BC}
\begin{enumerate}
\item \label{periodicity} $u_1,\om_1,\psi_1$ are periodic in $z$ with period $1$:
\begin{equation}\label{eq:periodicity}
\begin{split}
u_1(r,z,t) &= u_1(r,z+1,t),\quad \om_1(r,z,t) = \om_1(r,z+1,t),\\
\text{and}\quad \psi_1(r,z,t) &= \psi_1(r,z+1,t)\quad \text{for all $(r,z)\in \mathcal{D}$}.
\end{split}
\end{equation}
\item \label{odd_z} $u_1,\om_1,\psi_1$ are odd functions of $z$ at $z=0$:
\begin{equation}\label{eq:odd_z}
\begin{split}
u_1(r,z,t) &= -u_1(r,-z,t),\quad \om_1(r,z,t) = -\om_1(r,-z,t),\\
\text{and}\quad \psi_1(r,z,t) &= -\psi_1(r,-z,t)\quad \text{for all $(r,z)\in \mathcal{D}$}.
\end{split}
\end{equation}
\item \label{even_r} The smoothness of the solution in the Cartesian coordinates implies that $u^\theta,\om^\theta,\psi^\theta$ must be an odd function of $r$ \cite{liuwang2006}. Consequently, $u_1,\om_1,\psi_1$ must be even functions of $r$ at $r=0$, which imposes the pole conditions:
\begin{equation}\label{eq:even_r}
u_{1,r}(0,z) = \om_{1,r}(0,z) = \psi_{1,r}(0,z) = 0.
\end{equation}
\item \label{no-flow} The velocity satisfies a no-flow boundary condition on the solid boundary $r=1$:
\begin{equation}\label{eq:no-flow}
\psi_1(1,z,t) = 0\quad \text{for all $z$}.
\end{equation}
\item \label{no-slip} Due to the existence of non-degenerate viscosity at $r=1$, the tangent flows on the solid boundary should satisfy a no-slip boundary condition:
\begin{equation}\label{eq:no-slip}
u^\theta(1,z,t) = u^z(1,z,t) = 0,\quad \text{for all $z$}.
\end{equation}
In view of \eqref{eq:as_NSE_1_d} and \eqref{eq:no-flow}, this further leads to $\psi_{1,r}(1,z,t) = 0$. Therefore, the no-slip boundary in terms of the new variables $u_1,\om_1,\psi_1$ reads
\begin{equation}\label{eq:no-slip1}
u_1(1,z,t) = 0,\quad \om_1(1,z,t) = -\psi_{1,rr}(1,z,t),\quad \text{for all $z$}.
\end{equation}
\end{enumerate}
\end{subequations}
In fact, equations~\eqref{eq:axisymmetric_NSE_1} automatically preserve properties \eqref{periodicity}--\eqref{no-flow} of the solution for all time $t\geq0$ if the initial data satisfy these properties and if the variable viscosity coefficients $v^r,v^z$ satisfy 
\begin{itemize}
\item[(i)] $v^r,v^z$ are periodic in $z$ with period 1.
\item[(ii)] $v^r,v^z$ are even functions of $z$ at $z=0$.
\item[(iii)] $v^r,v^z$ are even functions of $r$ at $r=0$.
\end{itemize}
The no-slip boundary condition \eqref{no-slip} will be enforced numerically.
By the periodicity and the odd symmetry of the solution, we only need to solve equations~\eqref{eq:axisymmetric_NSE_1} in the half-period domain 
\[\mathcal{D}_1 = \{(r,z): 0\leq r\leq 1, 0\leq z\leq 1/2\}.\] 
Note that the properties \eqref{eq:periodicity} and \eqref{eq:odd_z} together imply that $u_1,\om_1,\psi_1$ are also odd functions of $z$ at $z=1/2$. The boundaries of $\mathcal{D}_1$ then behave like ``impermeable walls'' since 
\[u^r = -r \psi_{1,z} = 0 \quad \text{on $r=1$}\quad \text{and}\quad u^z = 2\psi_1+r\psi_{1,r} = 0\quad \text{on $z=0,1/2$}.\]

\subsection{Initial data}\label{sec:initial_data} 
\begin{figure}[!ht]
\centering
    \includegraphics[width=0.4\textwidth]{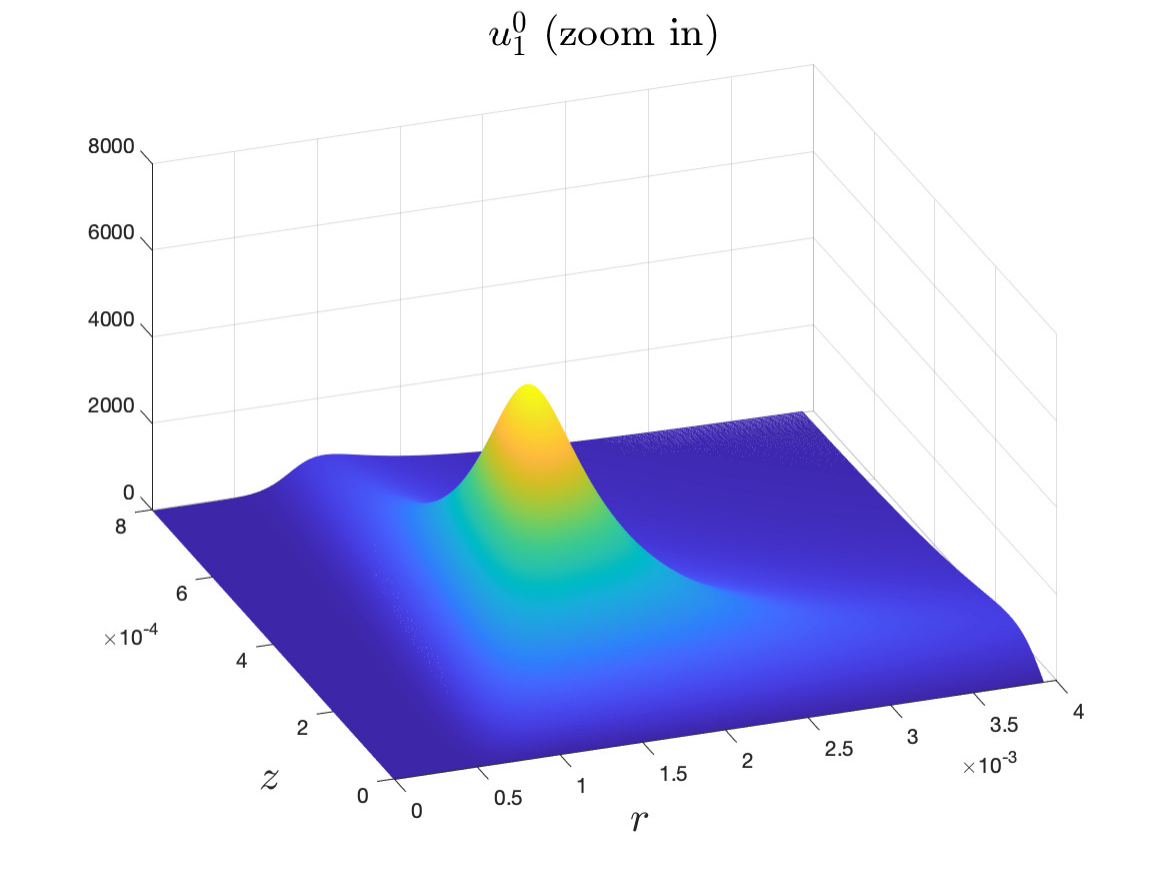}
    \includegraphics[width=0.4\textwidth]{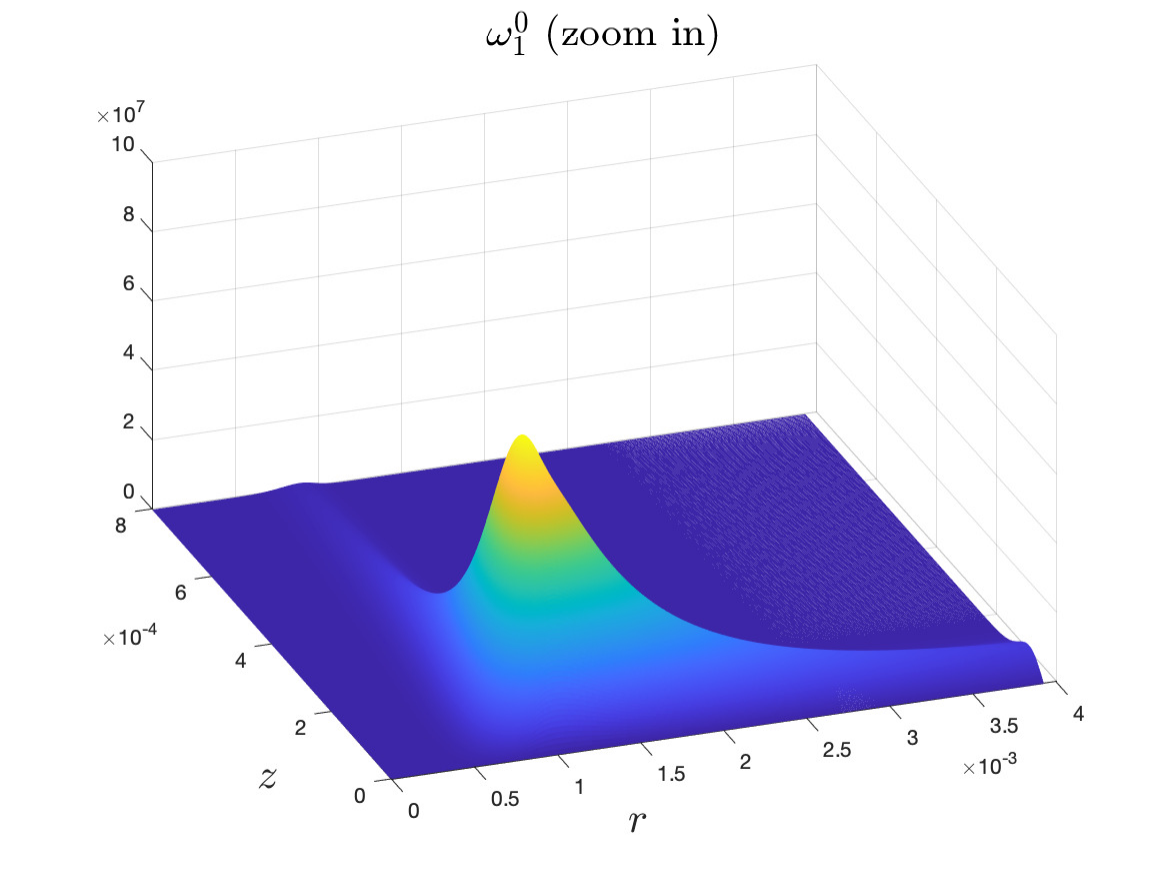}
    \includegraphics[width=0.4\textwidth]{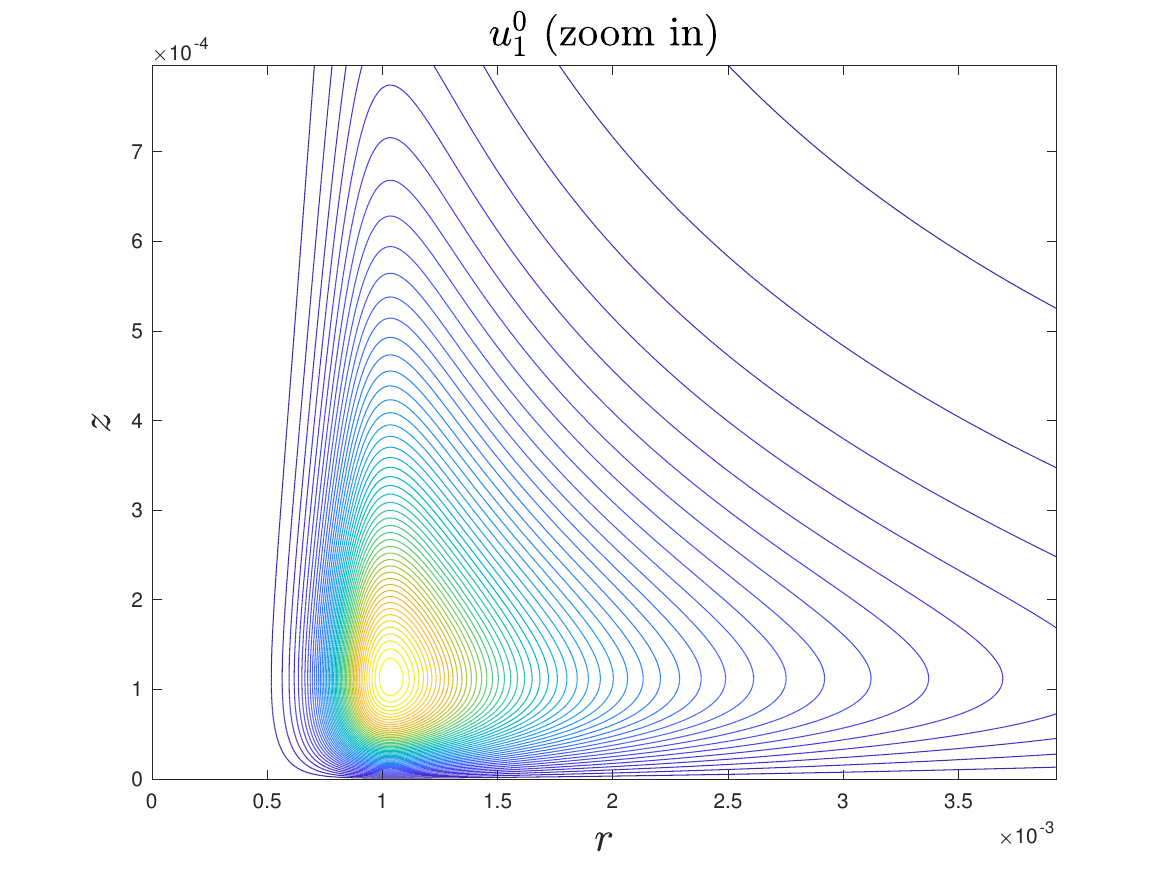}
    \includegraphics[width=0.4\textwidth]{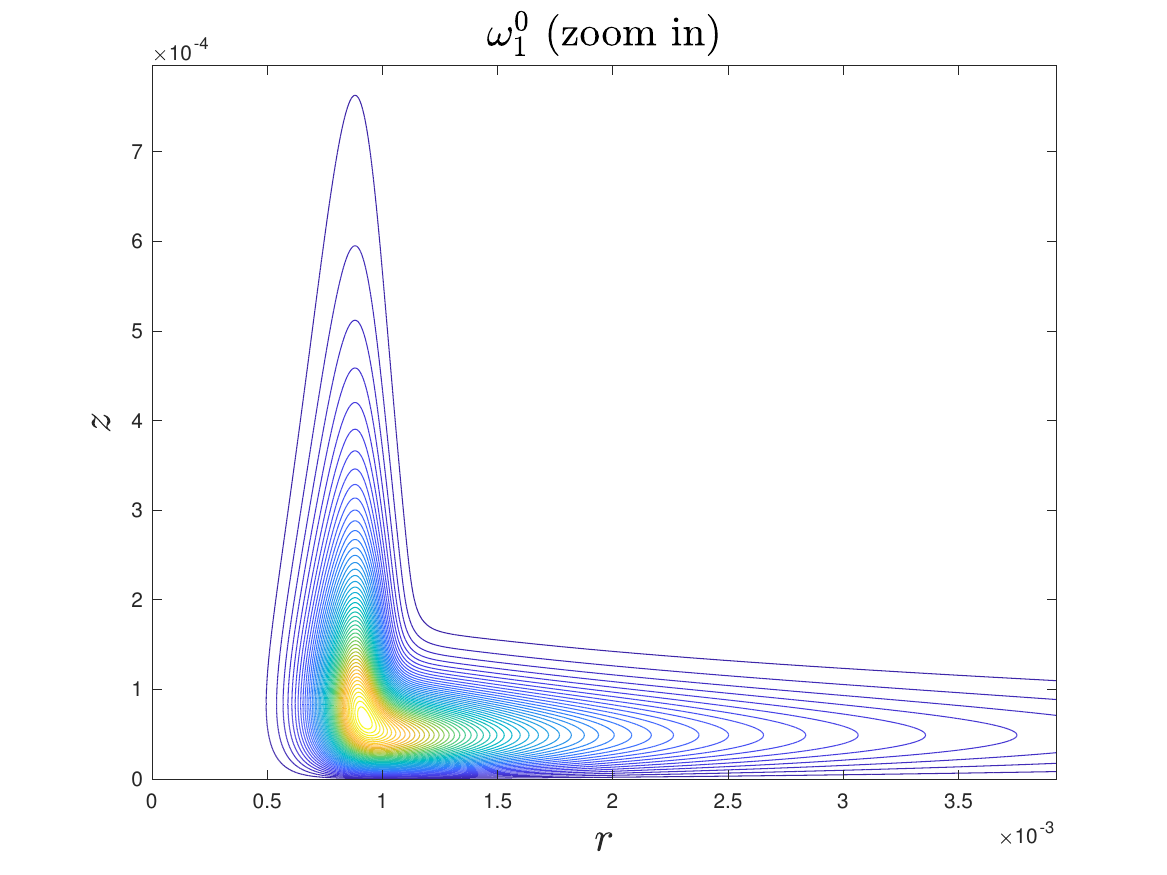}
    \caption[Initial data]{Profiles (first row) and contours (second row) of the initial data $u_1^0$ and $\om_1^0$ on a zoom-in domain $(r,z)\in[0,4\times 10^{-3}]\times [0,8\times 10^{-4}]$.}
    \label{fig:initial_data}
    \vspace{-0.05in}
\end{figure}

We construct the initial data based on our empirical insights and understanding of the potential blowup scenario that we shall explain later. The initial data are given by 
\begin{equation}\label{eq:initial_data}
u_1^0(r,z) = m_u^{(1)}\frac{u_1^{(1)}(r,z)}{\|u_1^{(1)}\|_{L^\infty}} + m_u^{(2)}u_1^{(2)}(r,z),\quad \om_1^0(r,z) = m_\om^{(1)}\frac{\om_1^{(1)}(r,z)}{\|\om_1^{(1)}\|_{L^\infty}} + m_\om^{(2)}\om_1^{(2)}(r,z),
\end{equation}
where 
\begin{align*}
u_1^{(1)} &= \frac{\sin(2\pi z)}{1+(\sin(\pi z)/a_{z1})^2+ (\sin(\pi z)/a_{z2})^4}\cdot \frac{r^8(1-r^2)}{1+(r/a_{r1})^{10} + (r/a_{r2})^{14}},\\
u_1^{(2)} &= \sin(2\pi z)\cdot r^2(1-r^2),\\ 
\om_1^{(1)} &= g(r,z)\cdot \frac{\sin(2\pi z)}{1+(\sin(\pi z)/b_{z1})^2+ (\sin(\pi z)/b_{z2})^4}\cdot \frac{r^8(1-r^2)}{1+(r/b_{r1})^{10} + (r/b_{r2})^{14}},\\
\text{and}\quad \om_1^{(2)} &= \sin(2\pi z)\cdot r^2(1-r^2).
\end{align*}
The parameters are chosen as follows:
\begin{align*}
&m_u^{(1)} = 7.6\times 10^3,\quad m_u^{(2)} = 50,\quad m_\om^{(1)} = 8.6\times 10^7,\quad m_\om^{(2)}= 50,\\
&a_{z1} = (1.2\times 10^{-4})\pi,\quad a_{z2} = (2.5\times 10^{-4})\pi,\quad a_{r1} = 9\times 10^{-4},\quad a_{r2} = 5\times 10^{-3},\\
&b_{z1} = (1\times 10^{-4})\pi,\quad b_{z2} = (1.5\times 10^{-4})\pi,\quad b_{r1} = 9\times 10^{-4},\quad b_{r2} = 3\times 10^{-3}.
\end{align*}
The function $g(r,z)$ is defined through a soft-cutoff function, and it forces the profile of $\om_1^0$ to have a smooth ``corner'' shape. Define the soft-cutoff function
\begin{equation}\label{eq:cutoff}
f_{sc}(x;a,b) = \frac{\econst^{(x-a)/b}}{\econst^{(x-a)/b} + \econst^{-(x-a)/b}}.
\end{equation}
Then $g(r,z)$ is given by the formula
\begin{align*}
g(r,z) &= \big(1-f_{sc}(\sin(\pi z)/\pi;0.7b_{z1},0.5b_{z1})\cdot f_{sc}(r; b_{r1}+0.5b_{z1}, b_{z1})\big)\\
&\quad \times \big(1-f_{sc}(-\sin(\pi z)/\pi;0.7b_{z1},0.5b_{z1})\cdot f_{sc}(r; b_{r1}+0.5b_{z1}, b_{z1})\big).
\end{align*}
Moreover, the initial stream function $\psi_1^0$ is obtained from $\om_1^0$ via the Poisson equation
\[-\left(\partial_r^2+\frac{3}{r}\partial_r+\partial_z^2\right)\psi_1^0(r,z) = \om_1^0(r,z)\quad \text{for $(r,z)\in \mathcal{D}_1$},\]
subject to the homogeneous boundary conditions 
\[\psi_{1,r}^0(0,z) =\psi_1^0(1,z)=\psi_1^0(r,0) = \psi_1^0(r,1/2)= 0.\]
It is not hard to check that the initial data $u_1^0,\om_1^0,\psi_1^0$ satisfy all the conditions \eqref{periodicity}--\eqref{no-flow}. Figure~\ref{fig:initial_data} shows the profiles and contours of the initial data $u_1^0$ and $\om_1^0$.

As we will see in Section \ref{sec:first_sign}, the solution to the initial-boundary problem \eqref{eq:axisymmetric_NSE_1}--\eqref{eq:initial_data} develops a potential finite singularity that is focusing at the origin $(r,z)= (0,0)$ and has two separated spatial scales. The sustainability and stability of this two-scale singularity crucially rely on a coupling blowup mechanism that will be discussed in details in Section \ref{sec:mechanism}. Our construction of the initial data serves to trigger this blowup mechanism, owing to the following principles:
\begin{itemize}
\item $\om_1^0$ is chosen to be an odd function of $z$ at $z=0$, so that the angular vorticity $\om^\theta = r\om_1$ has a dipole structure on the whole period $\{(r,z): r\in[0,1],z\in[-1/2,1/2]\}$ that induces a strong inward radial flow near the symmetry plane $z=0$ (see Figure \eqref{fig:dipole}). This flow structure has the desirable property that pushes the blowing up part of the solution towards the symmetry axis.
\item $u_1^0$ is also chosen to be odd in $z$ at $z=0$, so that the derivative $u_{1,z}^0$ is non-trivial and positive between the maximum point of $u_1^0$ and $z=0$. As a result, the forcing term $2u_1u_{1,z}$ in the $\om_1$ equation \eqref{eq:as_NSE_1_b} is positive and large near $z=0$, which contributes to the rapid growth of $\om_1$. 
\item For $u_1$ and $\om_1$ to have a good initial alignment, we manipulate the initial data so that the maximum point of $\om_1^0$ is slightly below the maximum point of $u_1^0$, where $u_{1,z}^0$ large and has the proper sign. 
\item We scale the magnitude of $u_1^0$ and $\om_1^0$ so that $\|\om_1^0\|_{L^\infty} \approx \|u_1^0\|_{L^\infty}^2$. This is because we have observed the blowup scaling property that $\|\om_1(t)\|_{L^\infty} \sim \|u_1(t)\|_{L^\infty}^2$ in our computations. The reason behind this blowup scaling property will be made clear in Section \ref{sec:scaling_study}. 
\end{itemize} 

These principles are critical for the solution to trigger a positive feed back mechanism that leads to a sustainable focusing blowup. We will have more discussions on the understanding of this mechanism in Section \ref{sec:mechanism}. We remark that the potential blowup is robust under relatively small perturbation in the initial data. In fact, we have observed similar two-scale blowup phenomena from a family of initial data that satisfy the above properties. 

\subsection{Variable viscosity coefficients} 
In our main cases of computation, we choose the variable viscosity coefficients $\nu^r,\nu^z$ to be the sum of a space-dependent part and a time-dependent part:
\begin{subequations}\label{eq:viscosity_coefficient}
\begin{align}
\nu^r(r,z,t) &= \frac{10r^2}{1+10^8r^2} + \frac{10^2(\sin(\pi z)/\pi)^2}{1+10^{11}(\sin(\pi z)/\pi)^2} + \frac{2.5\times 10^{-2}}{\|\om^\theta(t)\|_{L^\infty}}, \label{eq:nu_r}\\
\nu^z(r,z,t) &= \frac{10^{-1}r^2}{1+10^8r^2} + \frac{10^4(\sin(\pi z)/\pi)^2}{1+10^{11}(\sin(\pi z)/\pi)^2} + \frac{2.5\times 10^{-2}}{\|\om^\theta(t)\|_{L^\infty}}. \label{eq:nu_z}
\end{align}
\end{subequations}
We remark that that the space-dependent parts of $\nu^r,\nu^z$ are very small (below $10^{-7}$) on the whole domain and are of order $O(r^2) + O(z^2)$ for $r\leq 10^{-4}$ and $z\leq 10^{-5}$. Since the quantity $\|\om^\theta(t)\|_{L^\infty}$ is growing rapidly in our computation, the time-dependent part in $\nu^r,\nu^z$ is also very small (below $4\times 10^{-7}$ initially) and is decreasing rapidly in time. In fact, the time-dependent part is non-essential for the potential singularity formation in our scenario; it only serves to regularize the solution in the very early stage of our computation and will quickly be dominated by the space-dependent parts (see Figure \ref{fig:SDPvsTDP} in Section \ref{sec:original_NSE}). We can even remove the time-dependent part of $\nu^r,\nu^z$ after the solution enters a stable phase, and the phenomena we observe would remain almost the same. 

In all, we only have an extremely weak viscosity effect with smooth degenerate coefficients. Nevertheless, the viscosity plays an important role in the development of the singularity in our blowup scenario. On the one hand, the non-trivial viscosity in the far field prevents the shearing induced instability from disturbing the locally self-similar solution and the nonlinear alignment of vortex stretching in the near field. On the other hand, we will see in Section \ref{sec:original_NSE} that the degeneracy of the viscosity is crucial for a two-scale singularity to survive in a shrinking domain near the origin $(r,z) = (0,0)$. Furthermore, we will argue in Section \ref{sec:scaling_study} that the order of degeneracy of $\nu^r,\nu^z$ contributes to the formation of a potential locally self-similar blowup. 

For the sake of comparison, we will also study the solution of the $3$D Euler equations and of the original Navier--Stokes equations with constant viscosity coefficient on the entire domain. We will compare the numerical results from different choices of viscosity coefficients to study the effect of viscosity in the potential singularity near the symmetry axis $r=0$. As a preview, the geometric structure of the solution without viscosity quickly becomes too singular to resolve, while there is no blowup observed for the solution of the Navier--Stokes equations with a constant viscosity coefficient. This verifies the criticality of degeneracy in the viscosity coefficients. 

\subsection{The numerical methods} To numerically compute the potential singularity formation of the equations \eqref{eq:axisymmetric_NSE_1}, we have designed a composite algorithm that is well-tailored to the solution in our blowup scenario. The detailed descriptions of our numerical methods will be presented in Appendix \ref{apdx:numerical_methods}. Here we briefly introduce the main ingredients of our overall algorithm. 

\subsubsection{Adaptive mesh} We will see in Section \ref{sec:first_sign} that the profile of the solution quickly shrinks in space and develops complex geometric structures, which makes it extremely challenging to numerically compute the solution accurately. In order to overcome this difficulty, we design in Appendix \ref{apdx:adaptive_mesh} a special adaptive mesh strategy to resolve the singularity formation near the origin $(r,z) = (0,0)$. More precisely, we construct a pair of mapping functions 
\[r = r(\rho),\quad  z = z(\eta),\quad (\rho,\eta)\in[0,1]\times[0,1]\]
that maps the square $[0,1]\times[0,1]$ bijectively to the computational domain $\mathcal{D}_1=\{(r,z): 0\leq r\leq 1, 0\leq z\leq 1/2\}$. These mapping functions are dynamically adaptive to the complex multi-scale structure of the solution, which is crucial to the accurate computation of the potential singularity. Given a uniform mesh of size $n\times m$ on the $\rho\eta$-domain, the adaptive mesh covering the physical domain is produced as
\[r_i = r(ih_\rho),\quad h_\rho = 1/n;\quad z_j = z(jh_\eta),\quad h_\eta = 1/m.\]
The precise definition and construction of the mesh mapping functions are described in Appendix \ref{apdx:adaptive_mesh_construction}.

\subsubsection{B-spline based Galerkin Poisson solver} One crucial step in our computation is to solve the Poisson problem (2.3c) accurately. The Poisson solver we use should be compatible to the our adaptive mesh setting. Moreover, the finite dimensional system of this solver needs to be easy to construct from the mesh, as the mesh is updated frequently in our computation. For these reasons, we choose to implement the Galerkin finite element method based on a tensorization of B-spline functions, following the framework of Luo--Hou \cite{luo2014toward} who used this method for computing the potential singularity formation of the $3$D Euler equations on the solid boundary. The description of the Poisson solver will be given in Appendix \ref{apdx:poisson_solver}.

\subsubsection{Numerical regularization} The potential blowup solution we compute develops a long thin tail structure, stretching from the sharp front to the far field. This tail structure can develop some shearing-induced instability in the late stage of the computation, which may disturb the blowup mechanism. Therefore, we have chosen to apply numerical regularization to stabilize the solution, especially in the tail part. In particular, a low-pass filtering operator with respect to the $\rho\eta$-coordinates is introduced in Appendix \ref{apdx:regularization} for our regularization purpose.

\subsubsection{Overall algorithm} We use $2$nd-order centered difference schemes for the discretization in space and a $2$nd-order explicit Runge--Kutta method for marching the solution in time. We apply the low-pass filtering to mildly regularize the update of the solution in every time step. We have chosen to employ $2$nd-order methods since the low-pass filtering scheme we use introduces a $2$nd-order error of size $O(h_\rho^2+h_\eta^2)$. The resulting overall algorithm for solving the initial-boundary value problem \eqref{eq:axisymmetric_NSE_1}--\eqref{eq:initial_data} is formally $2$nd-order accurate in space and in time, which will be verified in Section \ref{sec:resolution_study}.

\section{Numerical Results: Features of Singularity}\label{sec:first_sign}
We have numerically solved the initial-boundary value problem \eqref{eq:axisymmetric_NSE_1}--\eqref{eq:initial_data} on the half-period cylinder $\mathcal{D}_1=\{(r,z):0\leq r\leq 1, 0\leq z\leq 1/2\}$ with meshes of size $(n,m) = (256p,128p)$ for $p = 2, 3, \dots, 8$. In particular, we have performed computations in $3$ cases:
\begin{enumerate}
\item[] Case $1$: 
$\nu^r,\nu^z$ given by \eqref{eq:viscosity_coefficient}. 
\item[] Case $2$:  
$\nu^r = \nu^z = \mu$ for some constant $\mu$.
\item[] Case $3$:  
$\nu^r = \nu^z = 0$.
\end{enumerate}
We will focus our discussions on the potential blowup phenomena in Case $1$. Our results suggest that the solution will develop a singularity on the symmetry axis $r=0$ in finite time, and we will provide ample evidences to support this finding. We first present, in this section, the major features of the potential finite-time singularity revealed by our computation. In Section \ref{sec:performance_study}, we carry out a careful resolution study of the numerical solutions. Then we further quantitatively investigate the properties of the potential singularity and analyze the potential blowup scaling properties in Section \ref{sec:scaling_study}.

Cases $2$ and $3$ are mainly for comparison purposes. Case $2$ compares the solution of the original Navier--Stokes equations and the solution of the Euler equations with degenerate viscosity coefficients using the same initial data. The results in Case $2$ show that the degeneracy of the viscosity coefficients near $(r,z)=(0,0)$ is critical for a sustainable singularity that approaches the symmetry axis $r=0$. The corresponding numerical results and discussions are presented in Section~\ref{sec:original_NSE}.

In Cases $3$, we study the evolution of the solution to the $3$D Euler equations from the same initial data. We will see that the solutions in Case $1$ and Case $3$ evolve almost in the same way during the early stage of the computation. However, the Euler solution quickly develops some oscillations due to under-resolution and shearing instability that prevent us from pushing the computation to the stable phase of the solution. Note that we did not apply any numerical regularization in Case $3$ to suppress the instabilities. Based on our preliminary results, we conjecture that the Euler solution will develop a similar or even more singular behavior in a later stage. However, our current adaptive mesh strategy does not allow us to resolve the potential Euler singularity to reach a convincing conclusion.

\subsection{Profile evolution}\label{sec:profile_evolution}
In this subsection, we investigate how the profiles of the solution evolve in time. We will use the numerical results in Case $1$ computed on the adaptive mesh of size $(n,m) = (1024,512)$. We have computed the numerical solution in this case up to time $t=1.76\times 10^{-4}$ when it is still well resolved. We cannot guarantee the reliability of our computation in Case $1$ beyond this time due to the loss of resolution, which will be discussed in Section \ref{sec:performance_study}. The computation roughly consists of three phases: a warm-up phase ($t\in[0,1.6\times 10^{-4}]$), a stable phase ($t\in(1.6\times 10^{-4},1.75\times 10^{-4}]$), and a phase afterwards ($t>1.75\times 10^{-4}$). In the warm-up phase, the solution evolves from the smooth initial data into a special structure. In the stable phase, the solution maintains a certain geometric structure and blows up stably. Beyond the stable phase, the solution starts to exhibit some unstable features that may arise from under-resolution, and the tail part of the solution also generates some shearing induced oscillations that are hard to resolve.

\begin{remark}
To have a better understanding of the solution behavior during the time interval $[0,1.76\times 10^{-4}]$, we first discuss the characteristic time and length scales of our problem. Since the solution of the Euler equations with degenerate viscosity coefficients develops a potential focusing singularity,  the characteristic length scale of the solution will decrease rapidly in time, and the maximal magnitude of the velocity will grow in time. For our initial data, the characteristic length scale is $10^{-4}$ and $\|\vct{u}\|_{L^\infty} \sim 10$, so the characteristic time scale is about $10^{-4}/\|\vct{u}\|_{L^\infty}\sim 10^{-5}$. At the time instant $t = 1.76\times10^{-4}$, the characteristic length scale drops to $10^{-6}$ and $\|\vct{u}\|_{L^\infty} \sim 50$, so the characteristic time scale is about $10^{-6}/\|\vct{u}\|_{L^\infty}\sim 2\times 10^{-8}$. For the record, an extrapolation of our numerical fitting of the potentially singular solution implies that the potential blow-up time is around $1.79\times 10^{-4}$. Thus, $1.76\times 10^{-4}$ is quite close to the potential blowup time. 
\end{remark}

\begin{figure}[!ht]
\centering
    \includegraphics[width=0.32\textwidth]{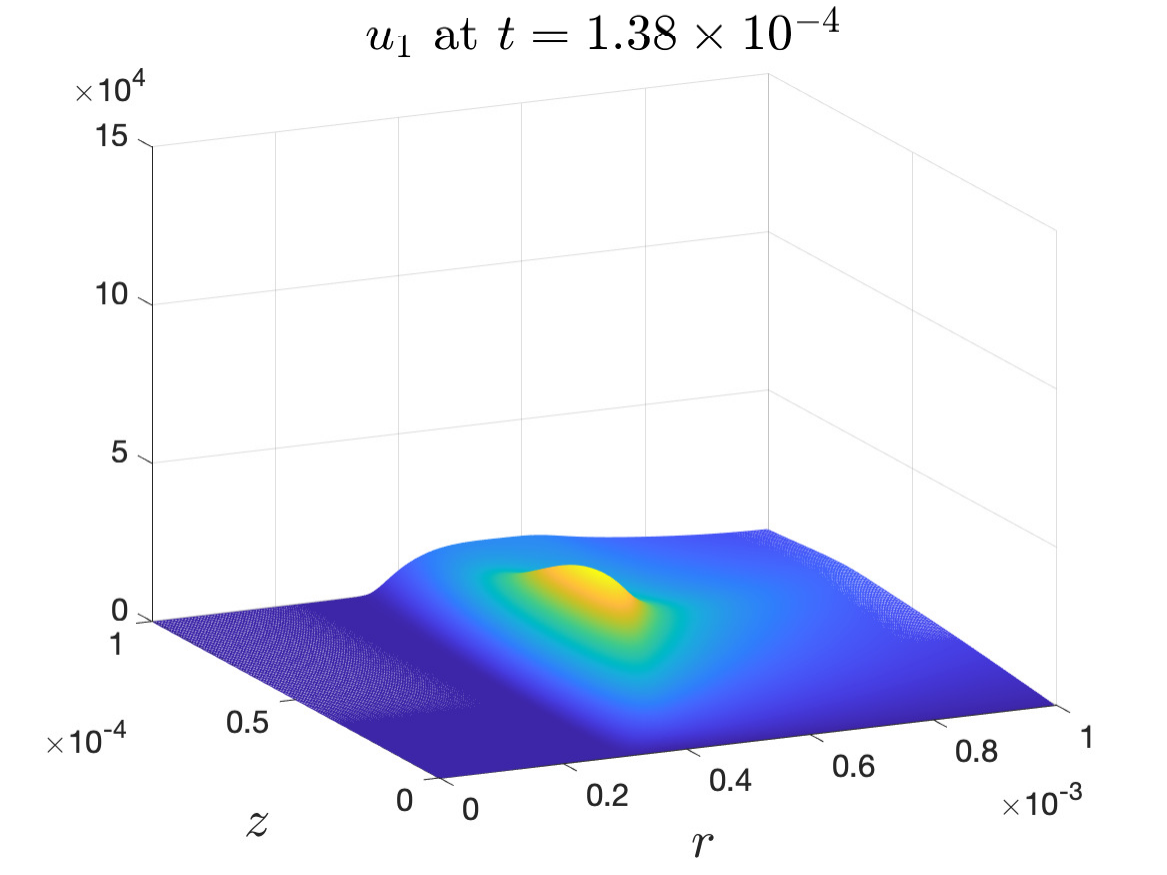}
    \includegraphics[width=0.32\textwidth]{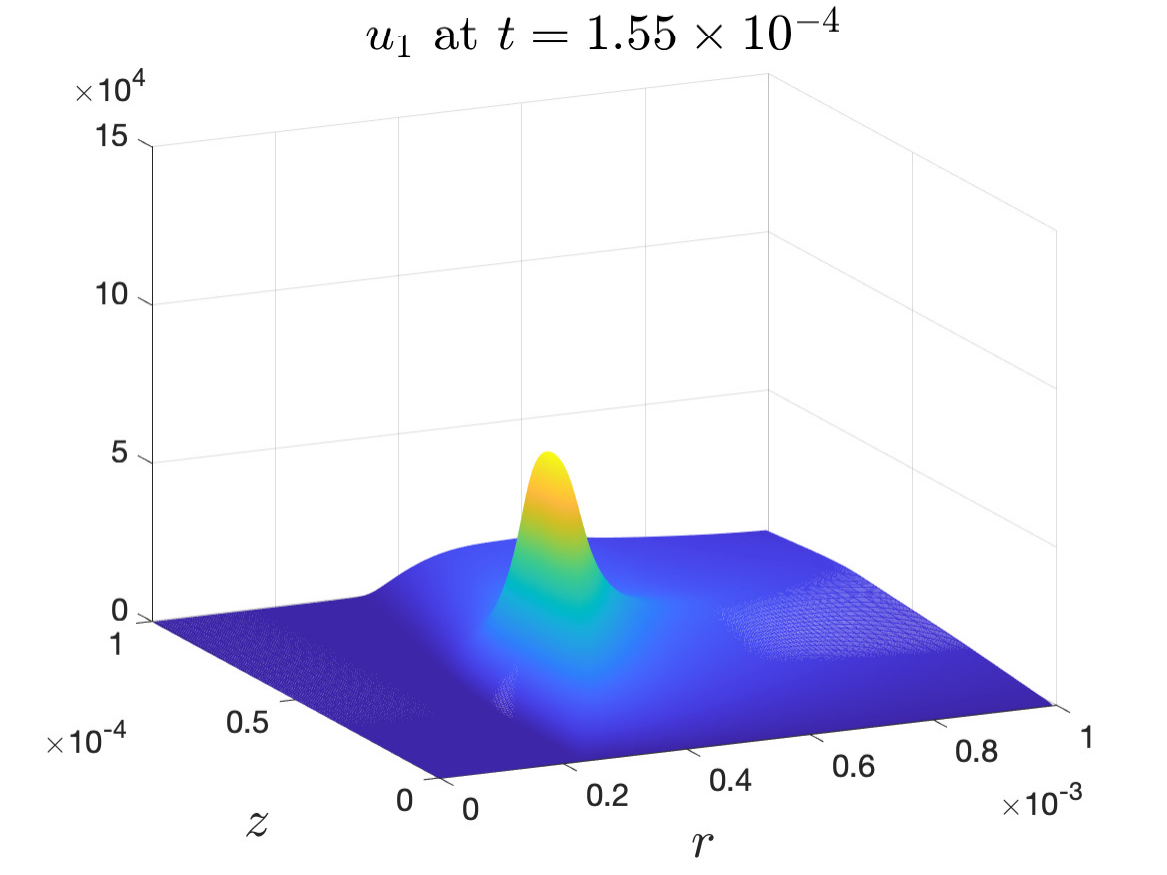}
    \includegraphics[width=0.32\textwidth]{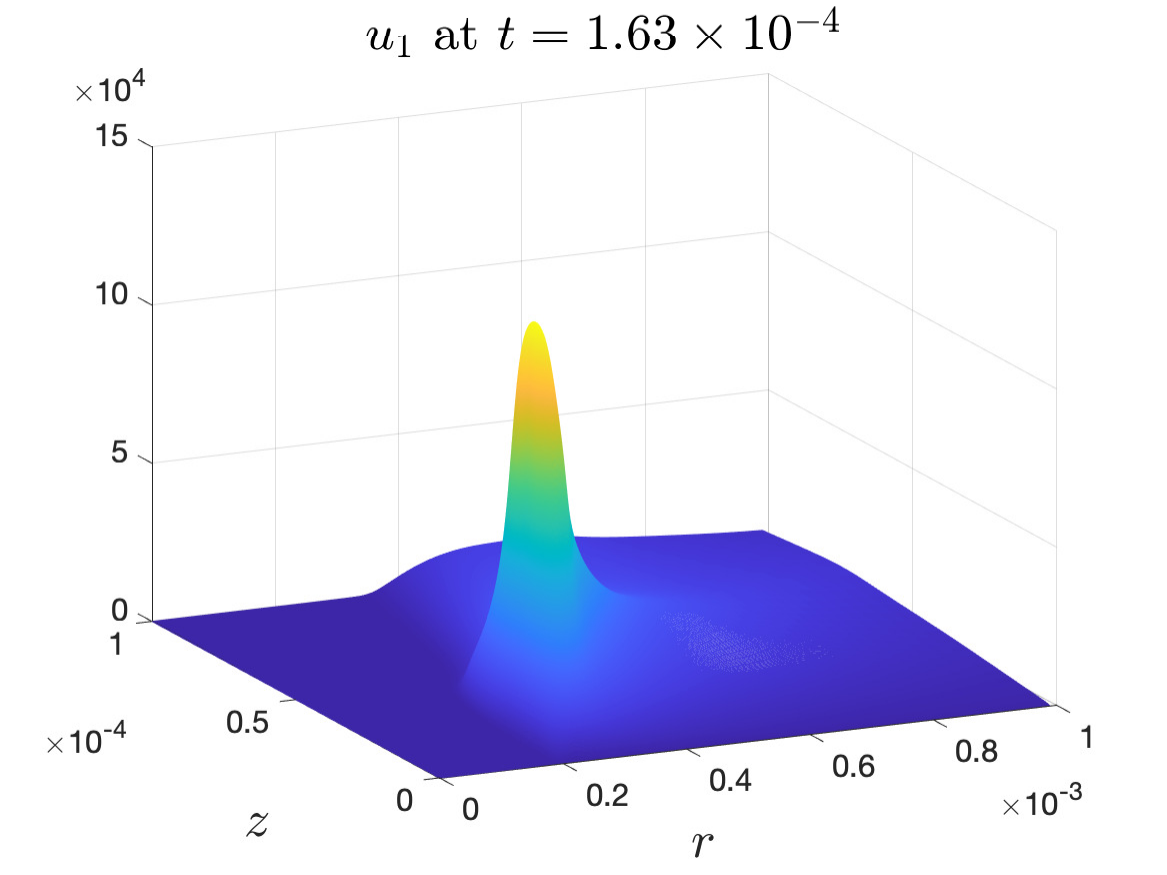}
    \includegraphics[width=0.32\textwidth]{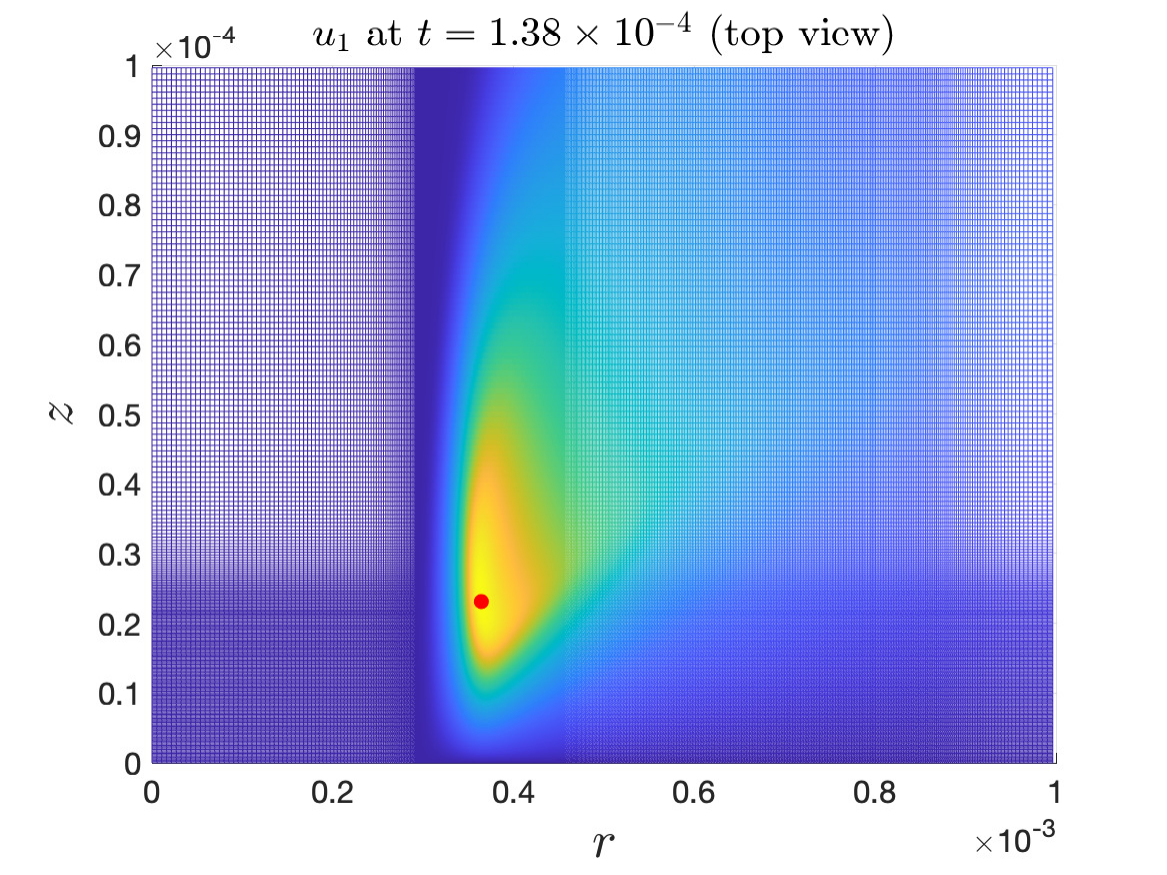}
    \includegraphics[width=0.32\textwidth]{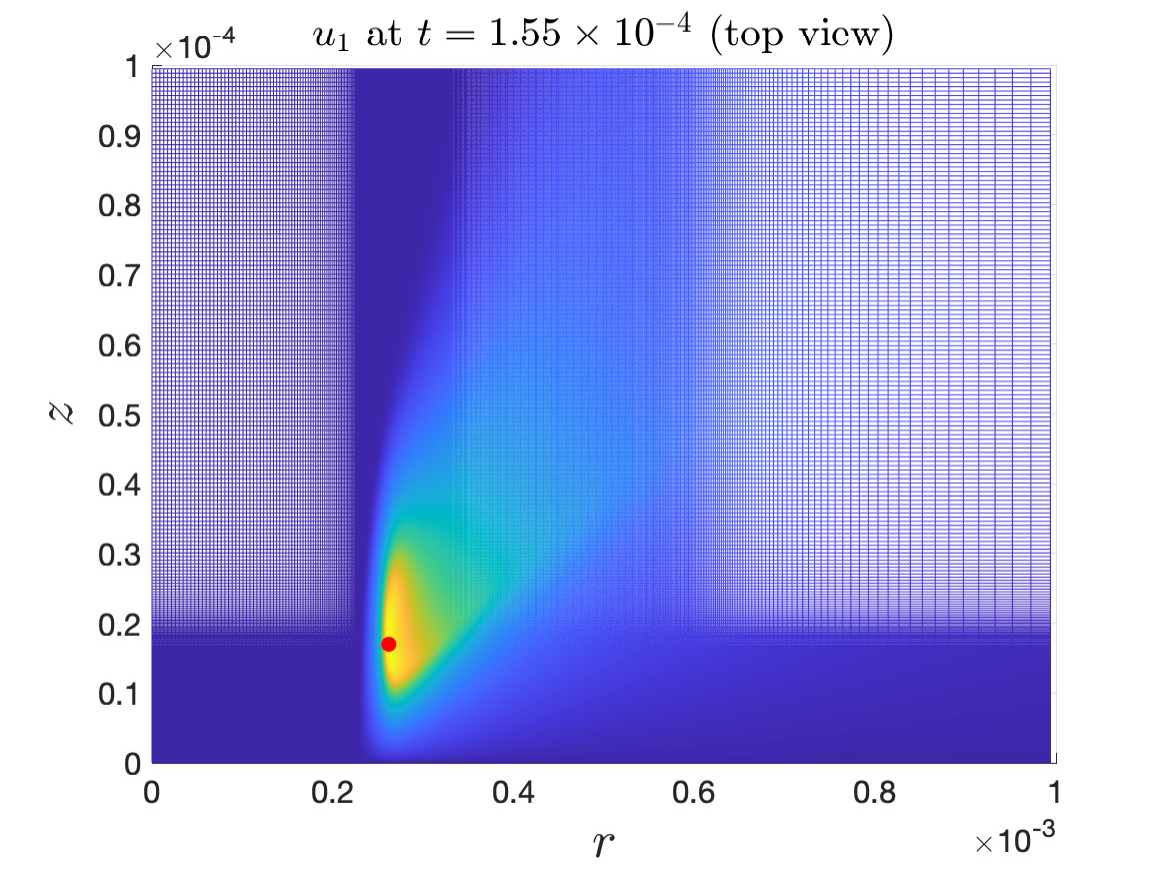}
    \includegraphics[width=0.32\textwidth]{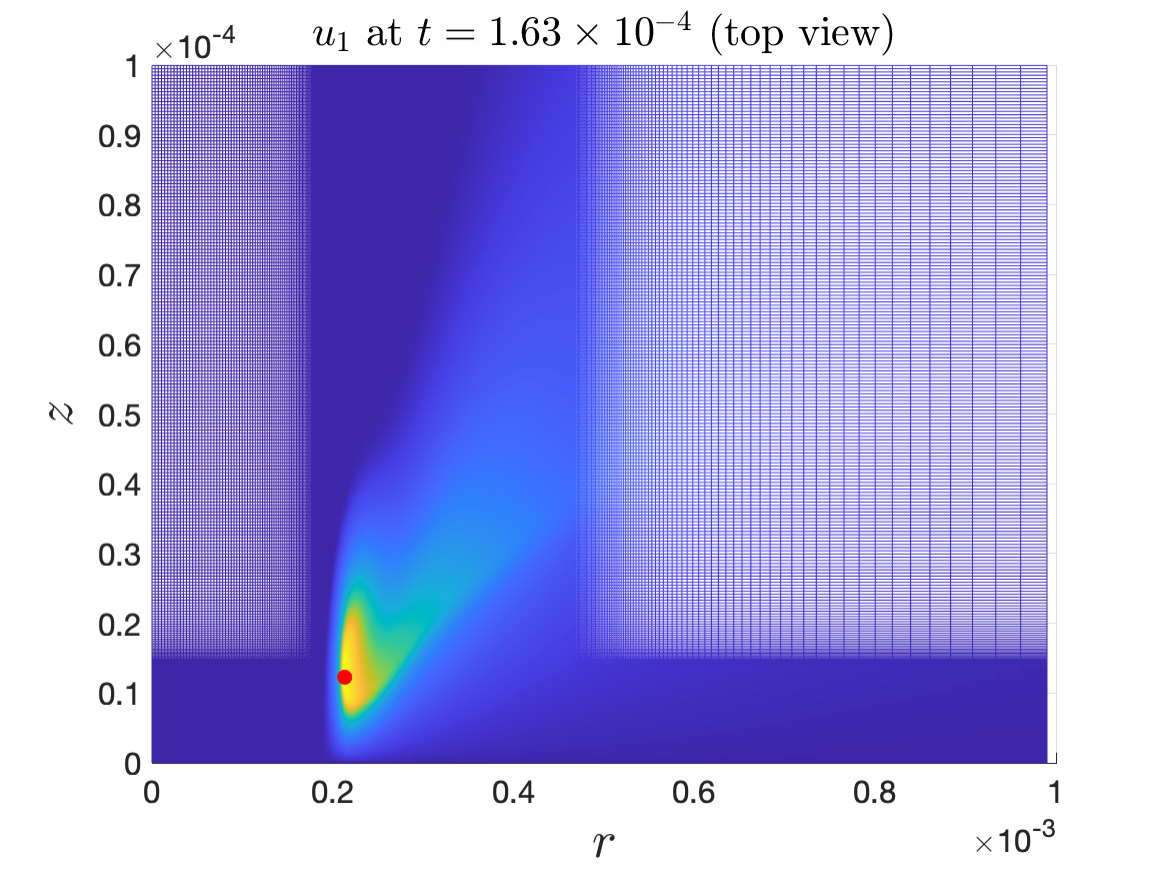}
    \includegraphics[width=0.32\textwidth]{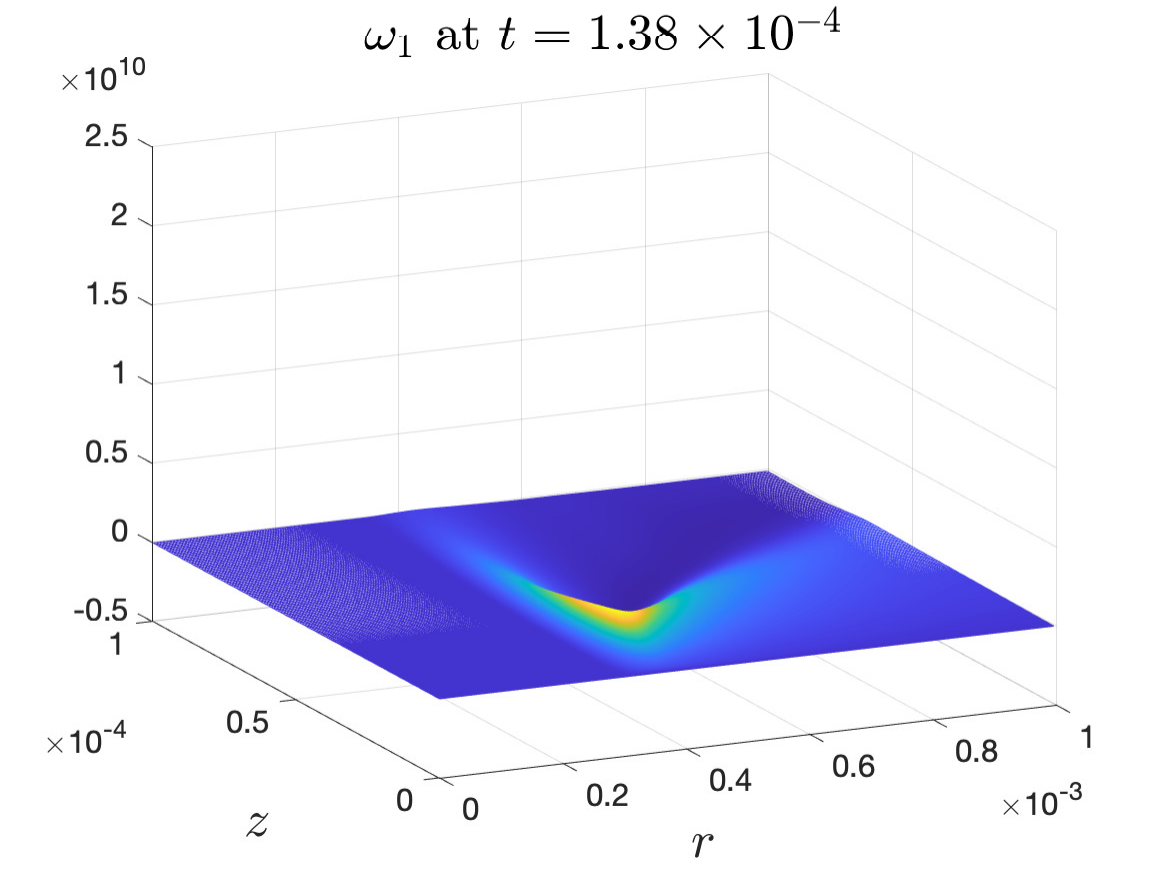}
    \includegraphics[width=0.32\textwidth]{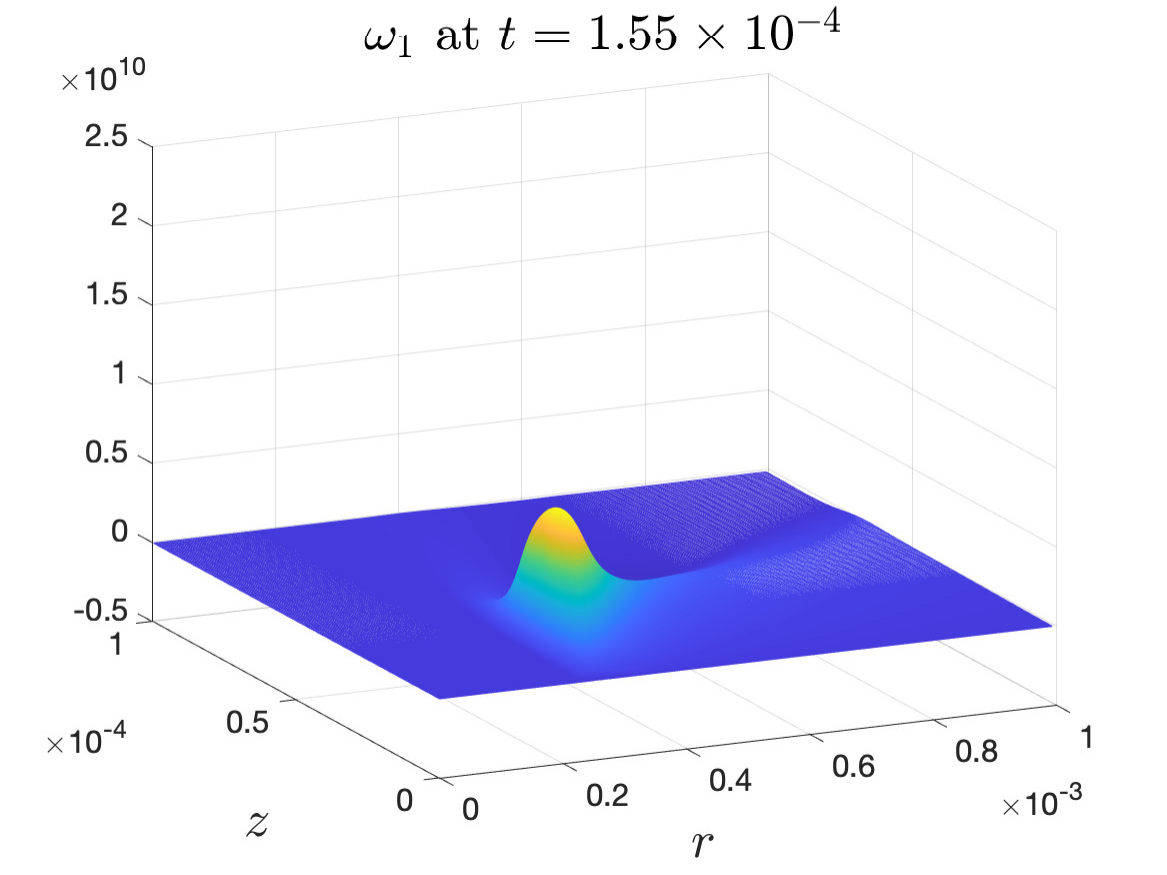}
    \includegraphics[width=0.32\textwidth]{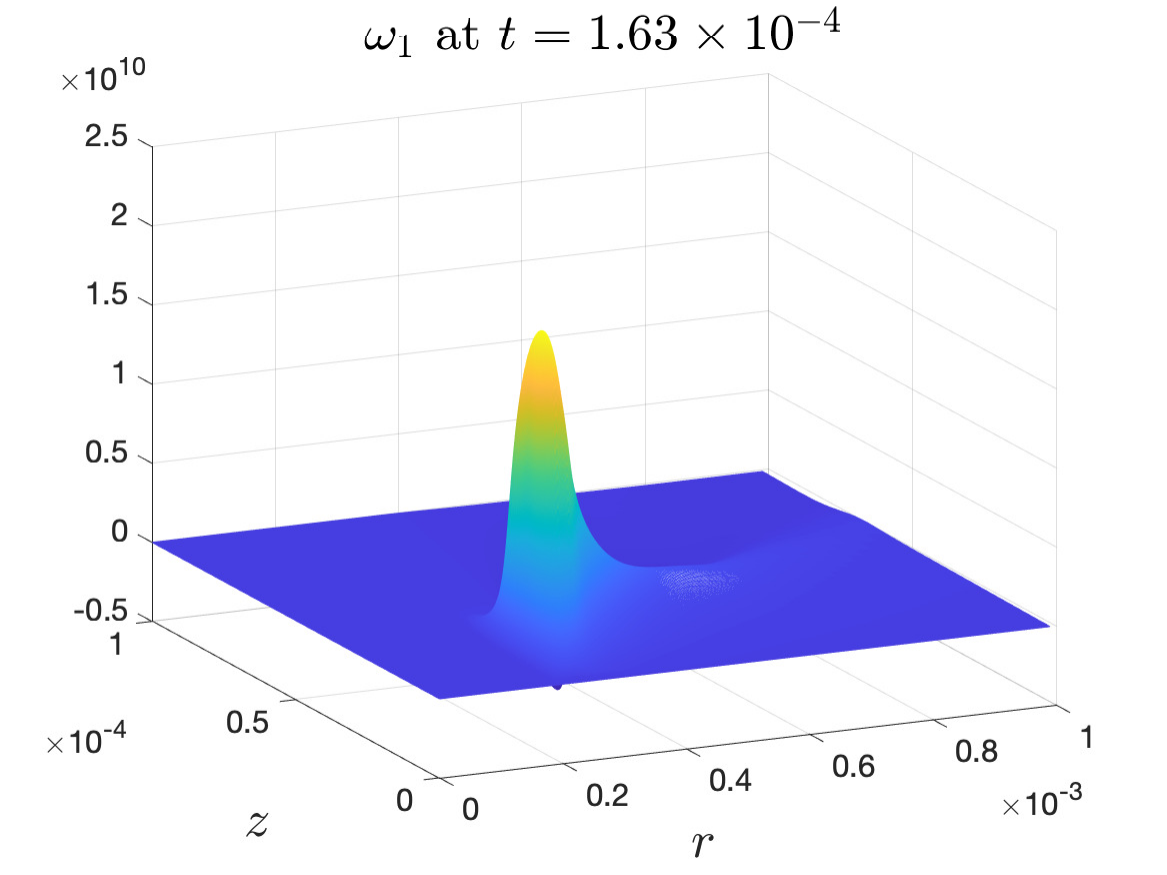}
    \includegraphics[width=0.32\textwidth]{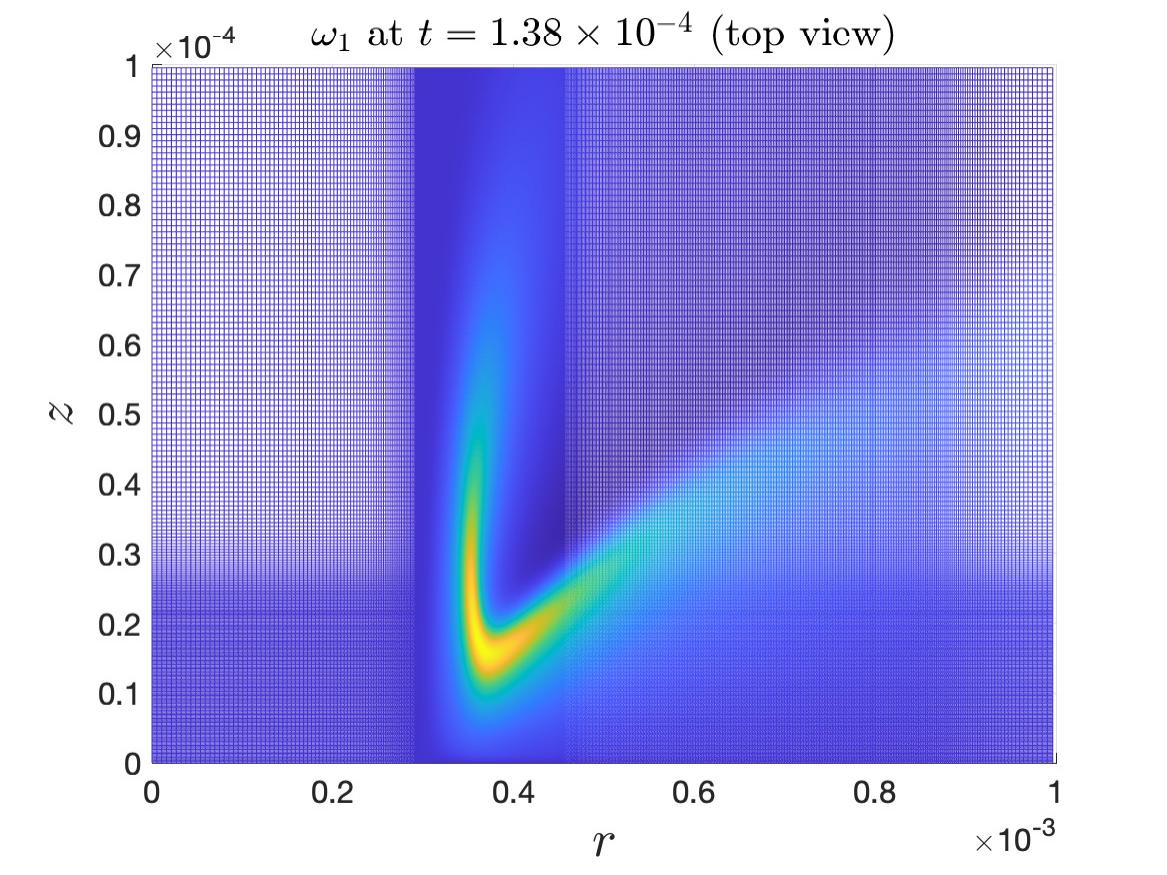}
    \includegraphics[width=0.32\textwidth]{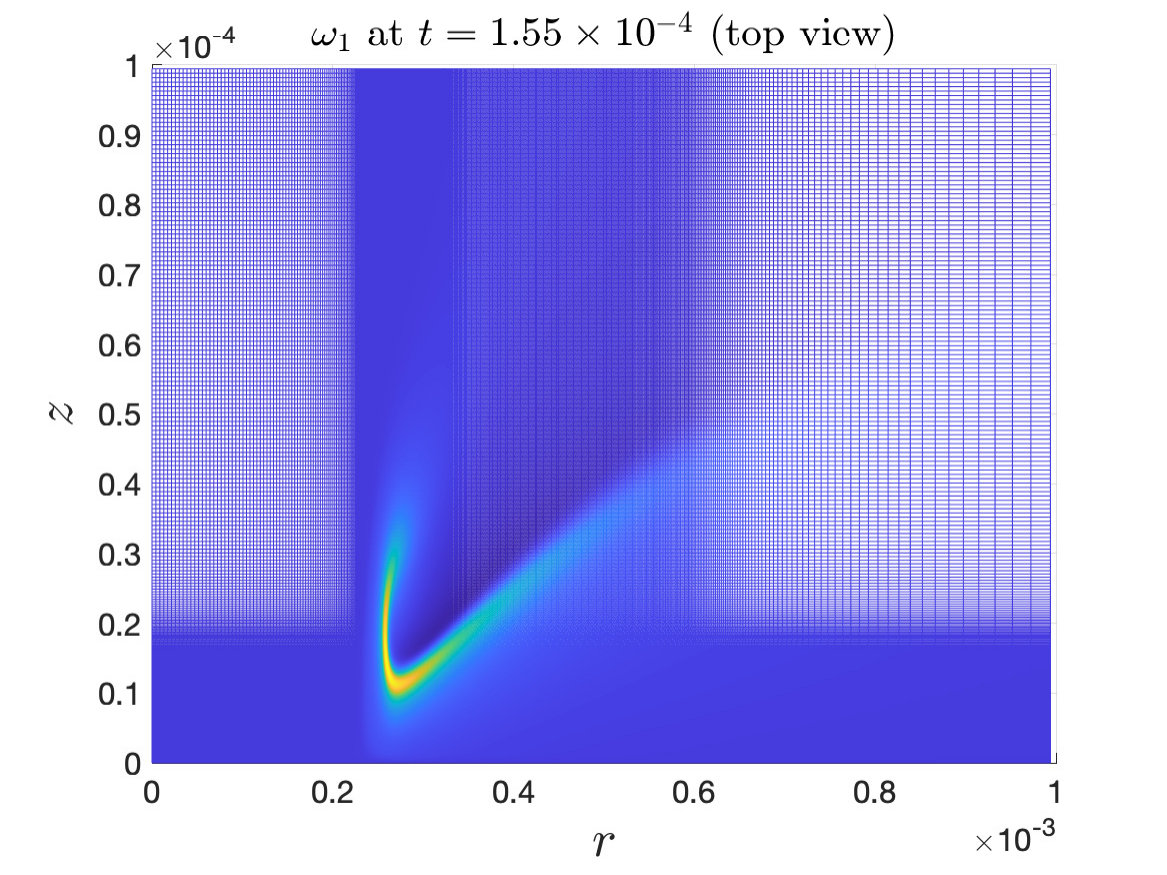}
    \includegraphics[width=0.32\textwidth]{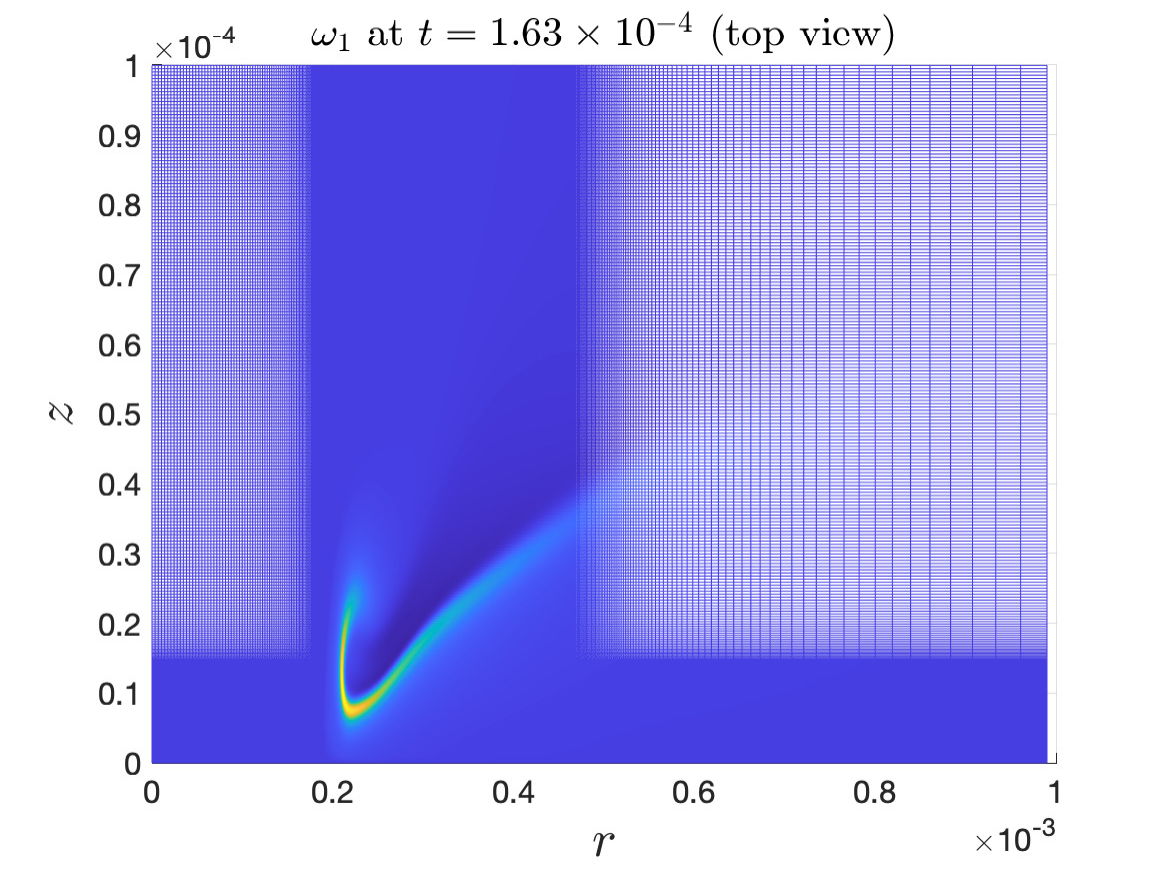}
    \caption[Profile evolution]{The evolution of the profiles of $u_1$ (row $1$ and $2$) and $\om_1$ (row $3$ and $4$) in Case $1$. Line $1$ and $3$ are the profiles of $u_1,\om_1$ at three different times;  Line $2$ and $4$ are the corresponding top-views. The red dot is the location of the maximum point of $u_1$.}  
    \label{fig:profile_evolution}
\end{figure}

\begin{figure}[!ht]
\centering
    \includegraphics[width=0.4\textwidth]{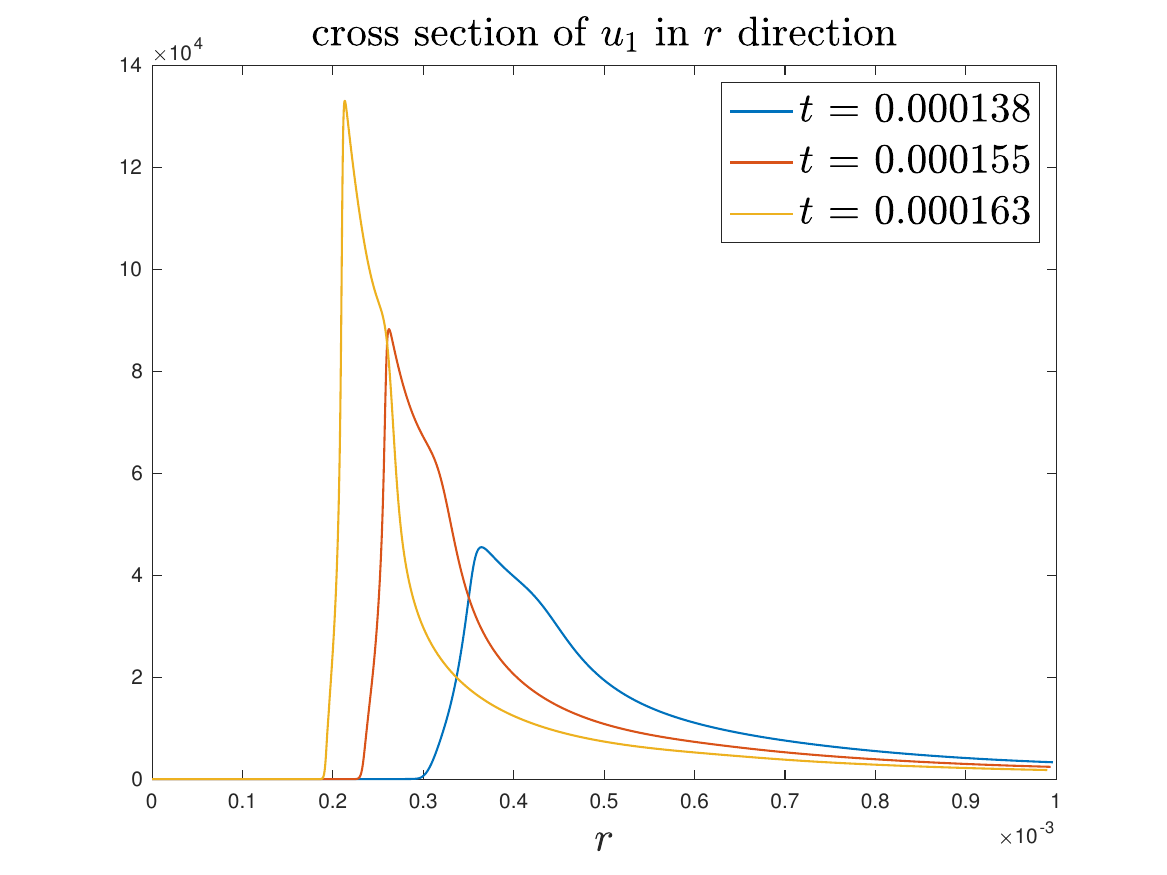}
    \includegraphics[width=0.4\textwidth]{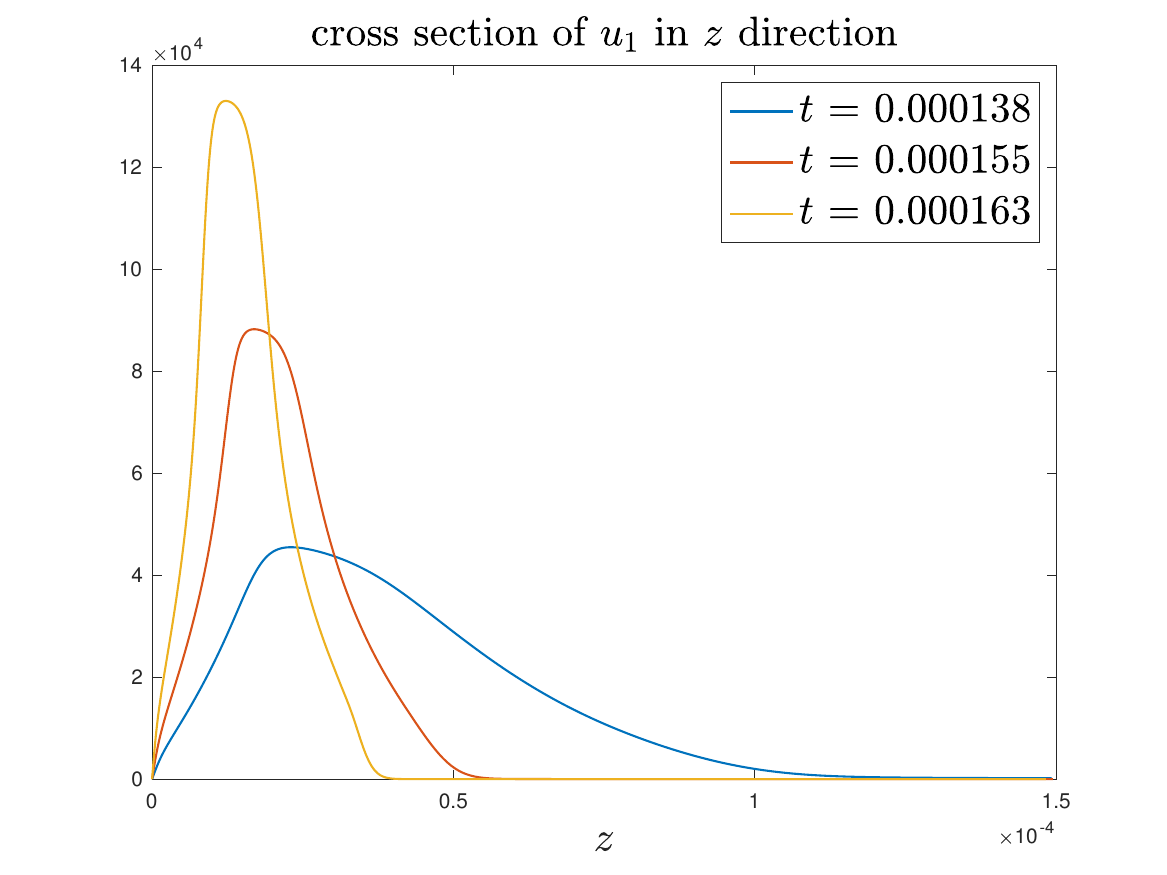}
    \caption[Cross section]{Cross sections of $u_1$ in both directions at different times.}  
    \label{fig:cross_section}
\end{figure}

Figure \ref{fig:profile_evolution} illustrates the evolution of $u_1,\om_1$ in the late warm-up phase by showing the solution profiles at $3$ different times $t = 1.38\times 10^{-4}, 1.55\times 10^{-4}, 1.63\times 10^{-4}$. We can see that the magnitudes of $u_1,\om_1$ grow in time. The ``peak'' parts of the profiles travel towards the symmetry axis $r=0$ and shrink in space. The profile of $u_1$ develops sharp gradients around the peak; in particular, it develops a sharp front in the $r$ direction. This is clearer if we look at the cross-sections of $u_1$ in both directions (Figure \ref{fig:cross_section}). Moreover, $\om_1$ develops a thin curved structure. Between the sharp front and the symmetry axis $r=0$, there is a no-spinning region where $u_1,\om_1$ are almost $0$. On the outer side, both $u_1$ and $\om_1$ form a long tail part propagating towards the far field. 
 
Let $(R(t),Z(t))$ denote the maximum point of $u_1(r,z,t)$. We will use this notation throughout the paper. Figure \ref{fig:cross_section} shows the cross sections of $u_1$ going through the point $(R(t),Z(t))$ in both directions. That is, we plot $u_1(r,Z(t),t)$ versus $r$ and $u_1(R(t),z,t)$ versus $z$, respectively. Again, it is clear that $u_1$ develops sharp gradients in both directions. In the $r$ direction, $u_1$ forms a sharp front and a no-spinning region between the sharp front and $r=0$. In the $z$ direction, the profile of $u_1$ seems to develop a compact support that is shrinking towards $z=0$. 

\subsection{Two scales}\label{sec:two_scale} 
\begin{figure}[!ht]
\centering 
    \includegraphics[width=0.4\textwidth]{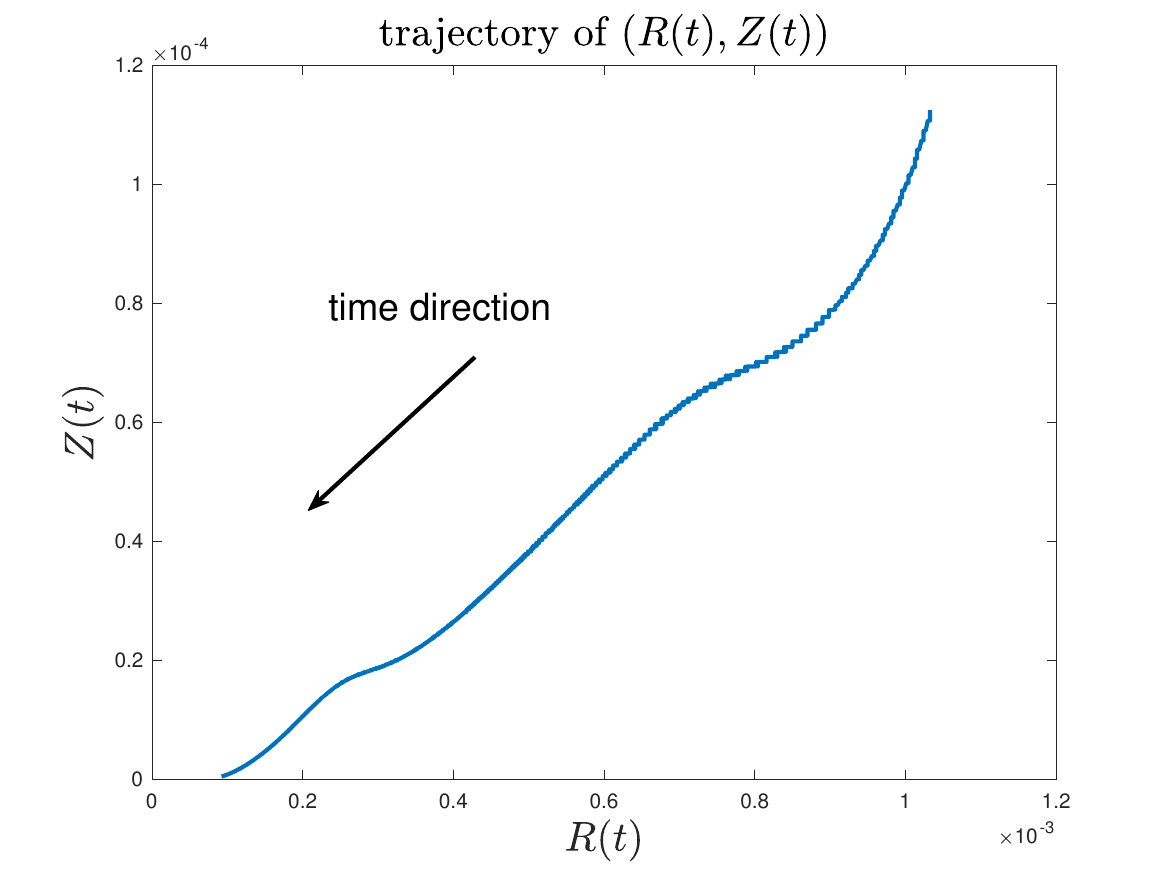}
    \includegraphics[width=0.4\textwidth]{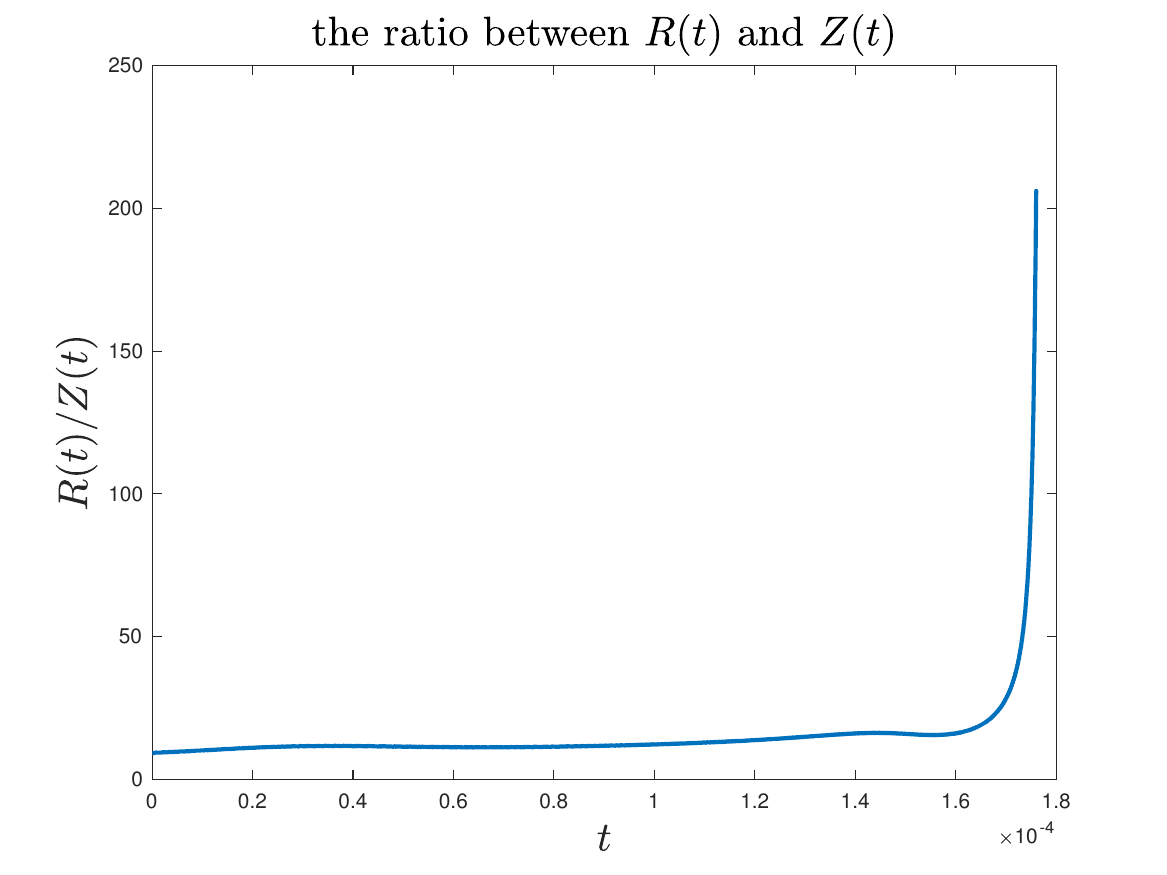}
    \includegraphics[width=0.40\textwidth]{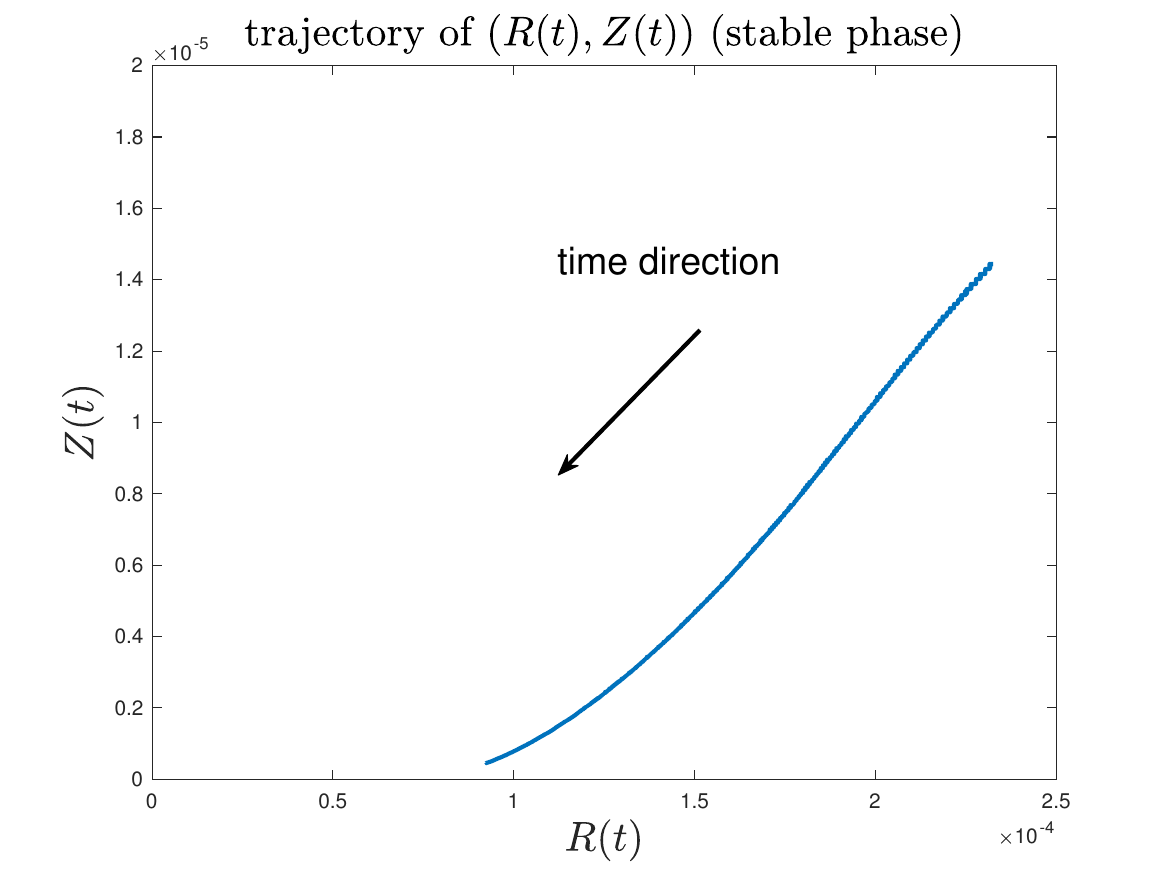}
    \includegraphics[width=0.40\textwidth]{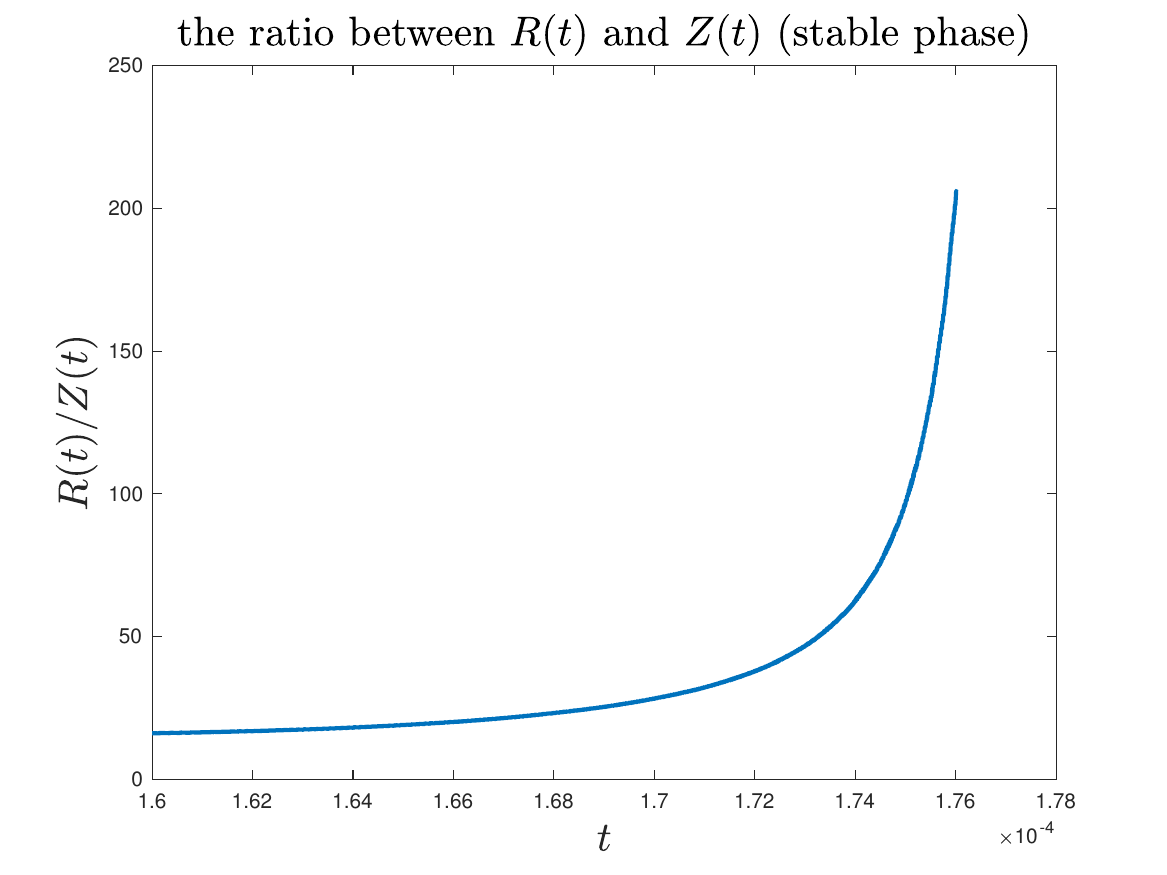}
    \caption[Trajectory]{The trajectory of $(R(t),Z(t)$ and the ratio $R(t)/Z(t)$ as a function of time for $t\in[0,1.76\times 10^{-4}]$. First row: the whole computation. Second row: the stable phase.}  
     \label{fig:trajectory} 
        \vspace{-0.05in}
\end{figure}

\begin{figure}[!ht]
\centering
    \includegraphics[width=0.40\textwidth]{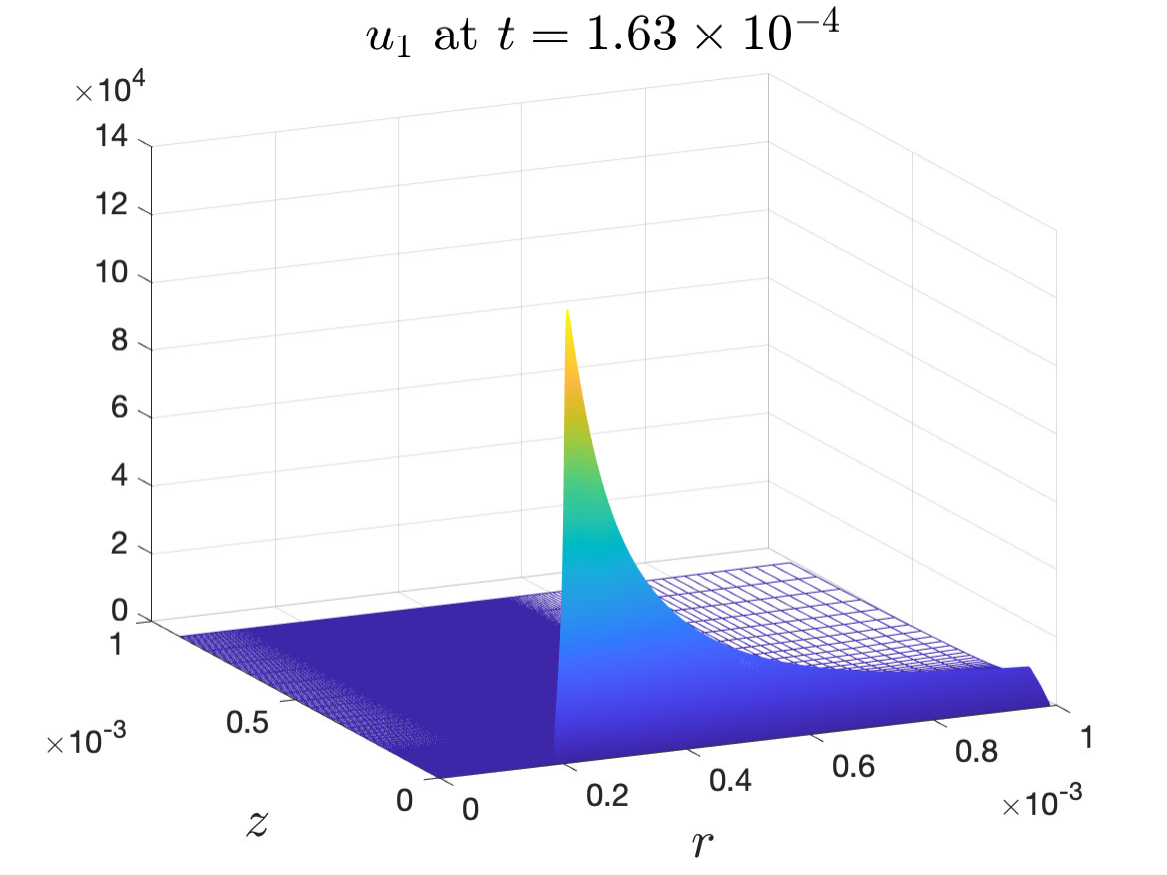}
    \includegraphics[width=0.40\textwidth]{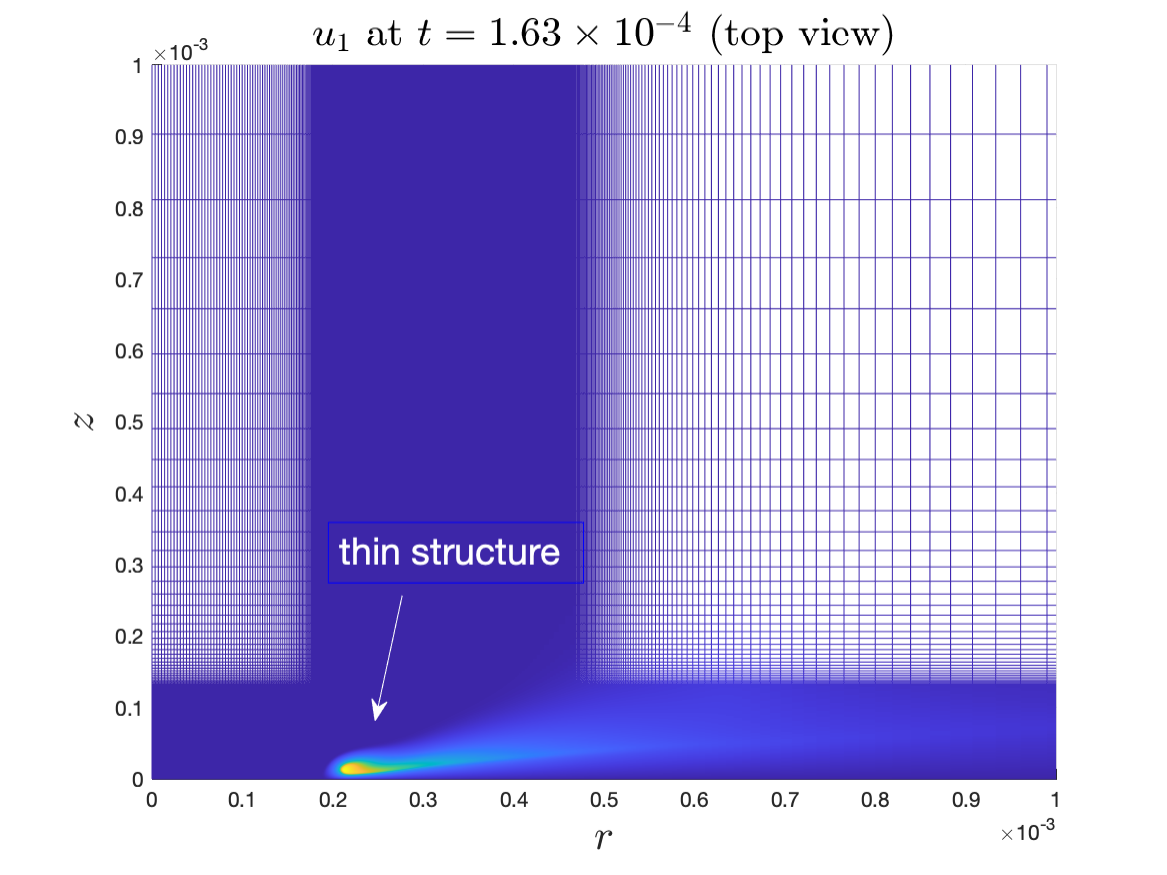}
    \includegraphics[width=0.40\textwidth]{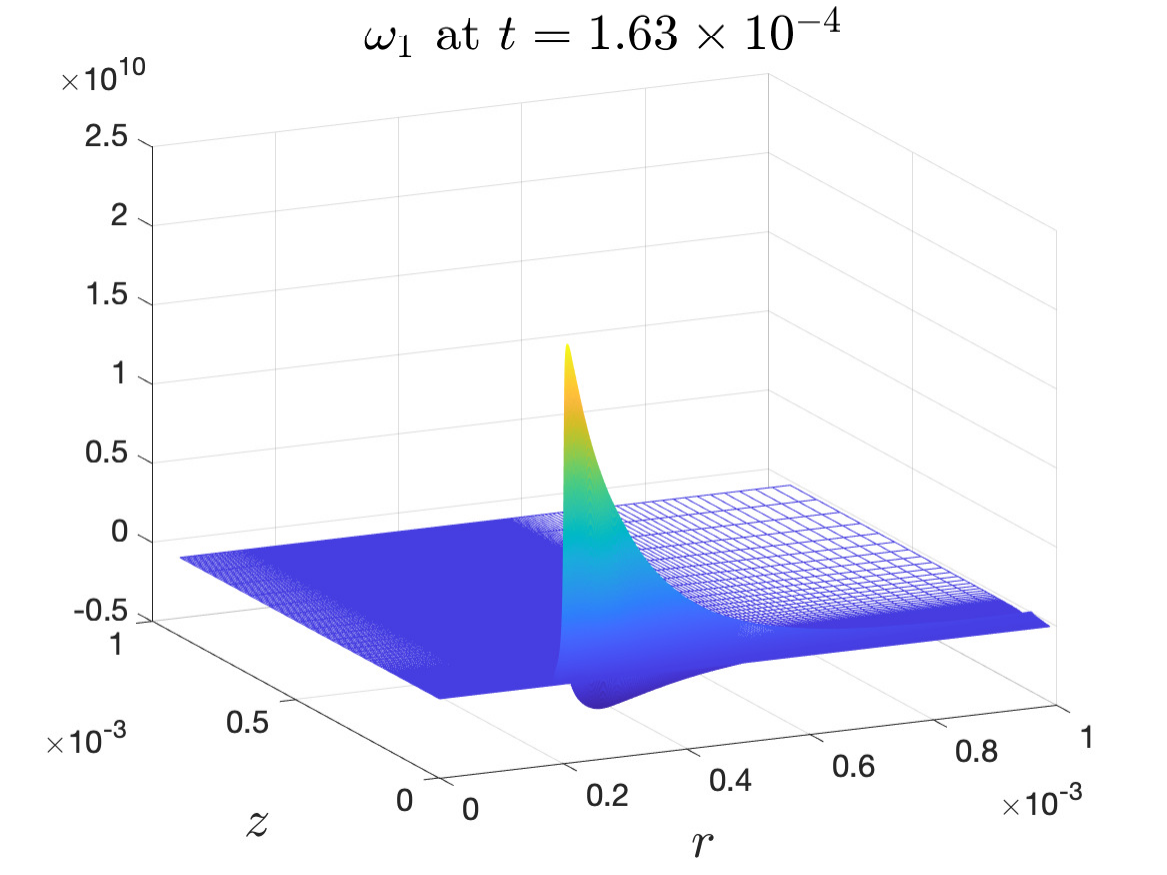}
    \includegraphics[width=0.40\textwidth]{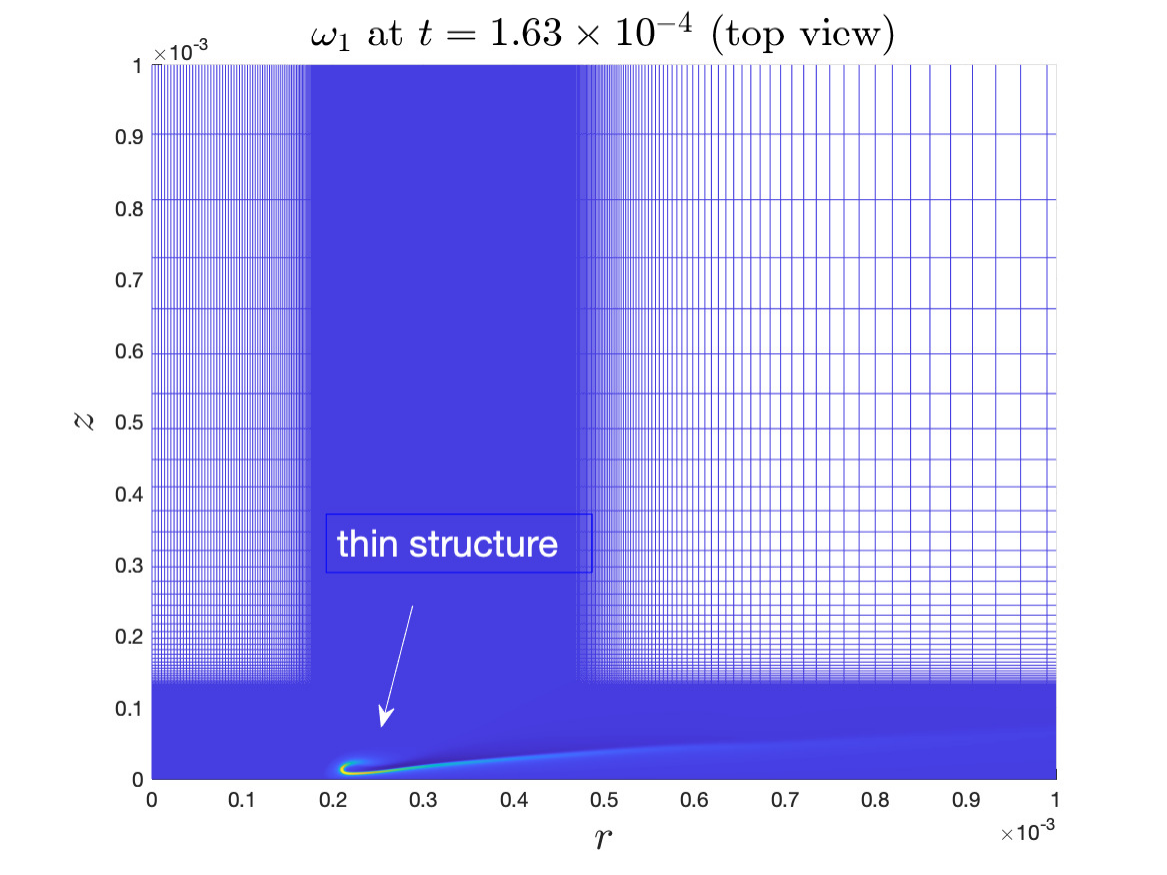} 
    \caption[Profile in square domain]{Profiles and level sets of $u_1$ (first row) and $\om_1$ (second row) at time $t= 1.63\times 10^{-4}$ in the square domain $\{(r,z):0\leq r\leq 10^{-3},0\leq z\leq 10^{-3}\}$.}  
     \label{fig:profile_square_domain}
        \vspace{-0.05in}
\end{figure}

Figure \ref{fig:trajectory} (first column) shows the trajectory of the maximum point $(R(t),Z(t))$ of $u_1(r,z,t)$. We can see that $(R(t),Z(t))$ moves towards the origin $(r,z)=(0,0)$, but with different rates in the two directions. This trajectory tends to become parallel to the horizontal axis $z=0$ in the stable phase, which means that $Z(t)$ may approach $0$ much faster than $R(t)$. As shown in the second column of Figure \ref{fig:trajectory}, the ratio $R(t)/Z(t)$ grows rapidly in time, especially in the stable phase. These evidences imply that there are two separate spatial scales in the solution. We can see this more clearly if we plot the solution profiles in a square domain in the $rz$-plane. For example, Figure \ref{fig:profile_square_domain} shows the profiles and level sets of $u_1,\om_1$ at time $t=1.63\times 10^{-4}$ in a square domain $\{(r,z):0\leq r\leq 10^{-3},0\leq z\leq 10^{-3}\}$. The profiles have a sharp front in the $r$ direction and are extremely thin in the $z$ direction, which corresponds to the scale of $Z(t)$ (the smaller scale). The long spreading tails of the profiles and the distance between the sharp front and the symmetry axis $r=0$ correspond to the scale of $R(t)$.

\begin{figure}[!ht]
\centering
    	\includegraphics[width=0.4\textwidth]{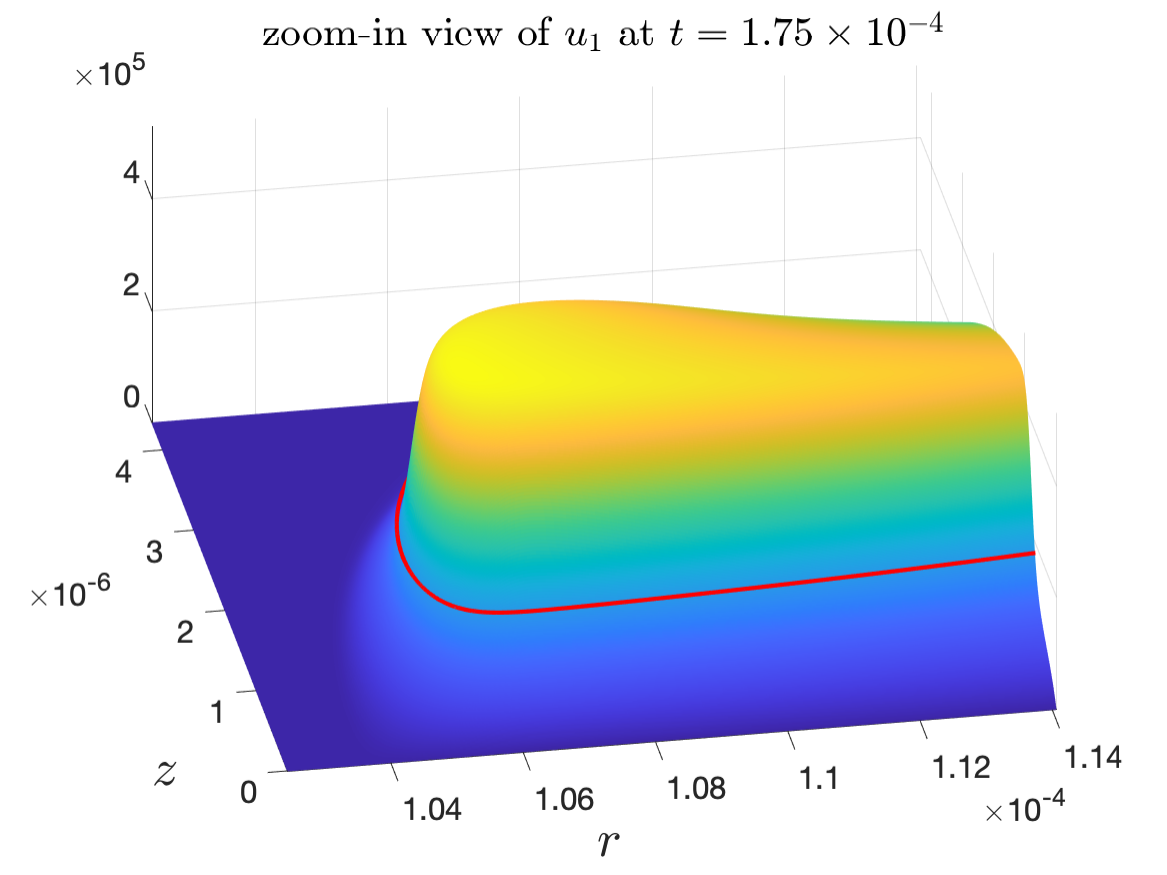}
    	\includegraphics[width=0.4\textwidth]{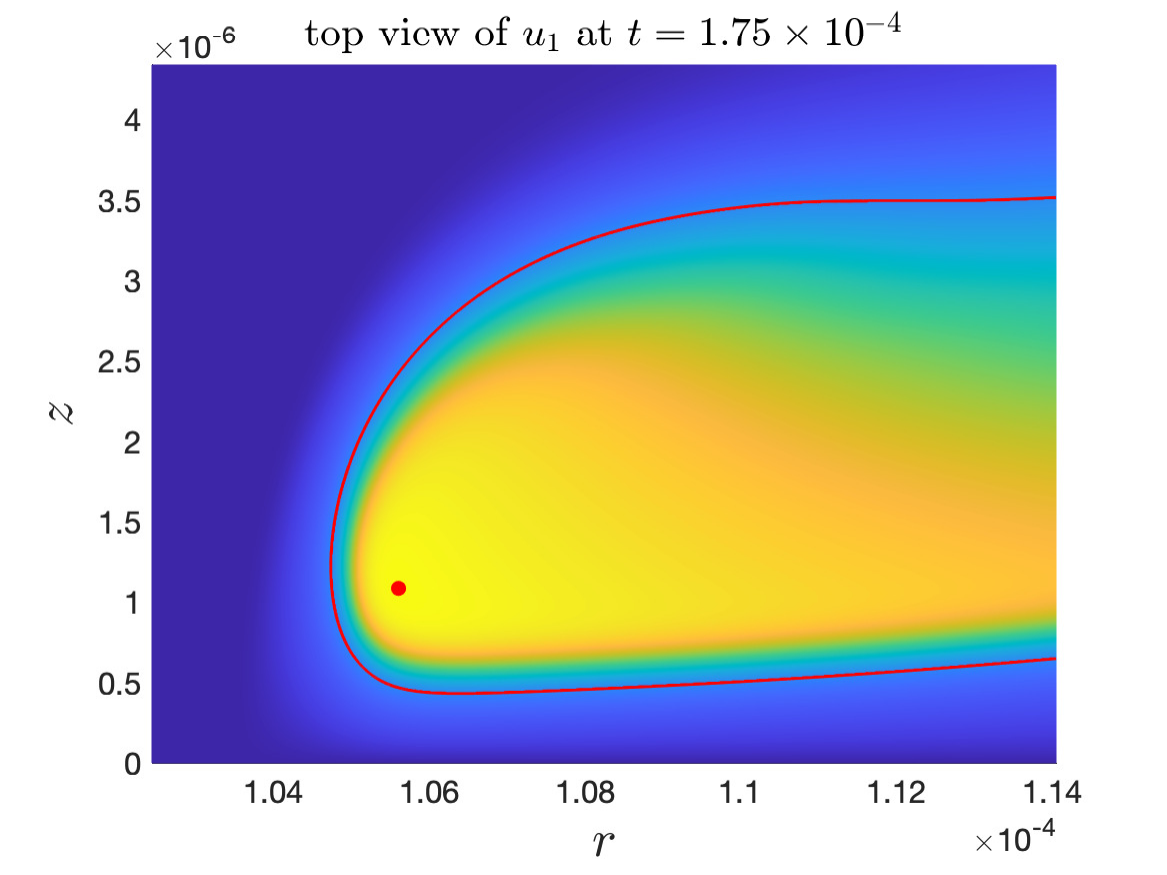}
    	\includegraphics[width=0.4\textwidth]{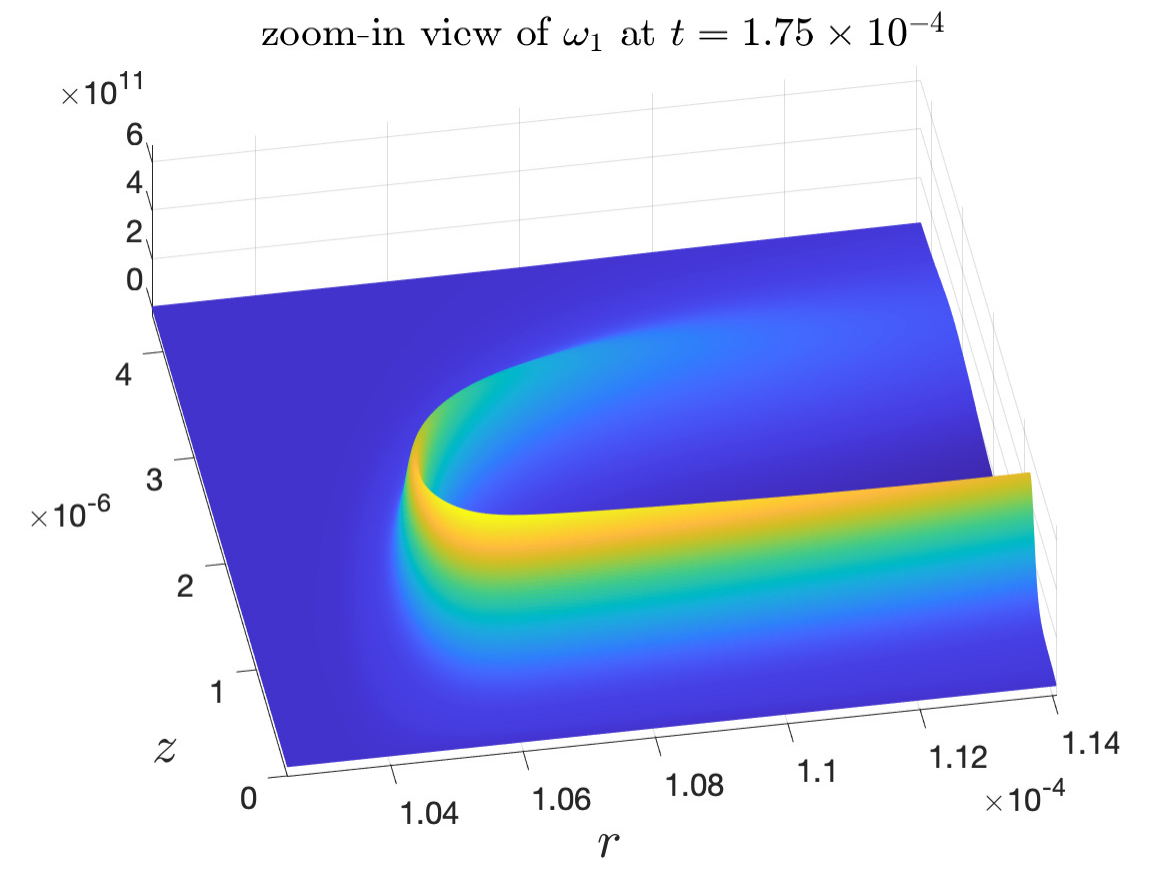}
    	\includegraphics[width=0.4\textwidth]{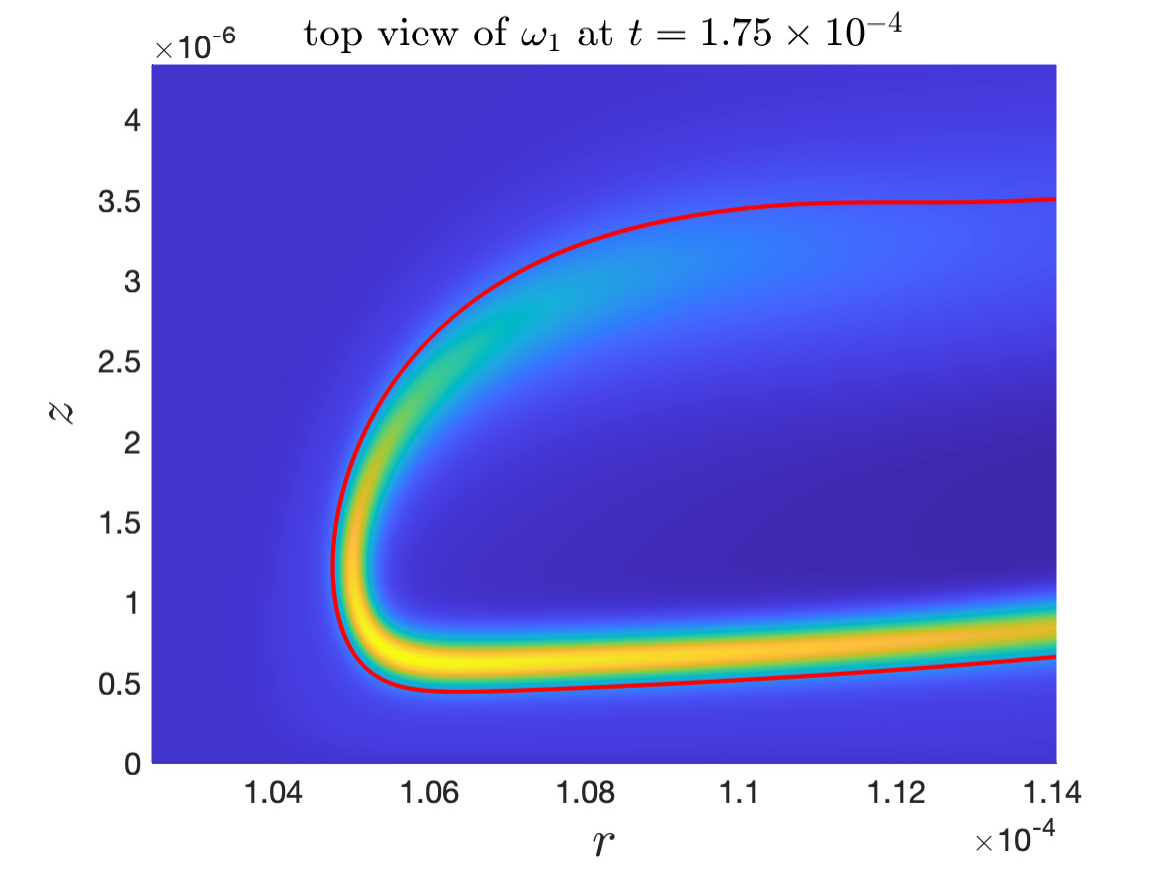}
    \caption[Zoom-in profile in square domain]{Zoom-in views of $u_1,\om_1$ at time $t = 1.75\times 10^{-4}$. First row: profile and tow view of $u_1$. Second row: profile and top view of $\om_1$. The red curve (in all figures above) is the level set of $u_1$ for the value $0.3\|u_1\|_{L^\infty}$, and the red point is the maximum point of $u_1$.}  
     \label{fig:zoomin_profile}
        \vspace{-0.05in}
\end{figure}

If we zoom into a neighborhood of the smaller scale around the point $(R(t),Z(t))$, we can see that the smooth profiles of $u_1,\om_1$ are locally isotropic. Figure \ref{fig:zoomin_profile} shows the local isotropic profiles of $u_1,\om_1$ near the sharp front at a later time $t = 1.75\times 10^{-4}$. These profiles are very smooth with respect to the smaller scale. In fact, such local structures have been developed ever since the solution enters the stable phase ($t\geq 1.6\times 10^{-4}$), and they remain stable afterwards. We will further investigate this in Section \ref{sec:evidence_of_self-similar}.


It is curious that the contours of $u_1$ and $\om_1$ seem to have the same shape. The thin structure of $\omega_1$ behaves like a regularized $1$D delta function supported along the ``boundary'' of the bulk part of $u_1$, which is roughly indicated by the red curve. In fact, we will see in Section \ref{sec:scaling_study} that this phenomena is an evidence of the existence of a two-scale, locally self-similar blowup. 

We remark that the numerical solutions computed in Case $3$ have almost the same features as described above in Case $1$. What varies most is how long these features can maintain stable in time. 

\subsection{Rapid growth}\label{sec:rapid_growth} 
\begin{figure}[!ht]
\centering
    \includegraphics[width=0.32\textwidth]{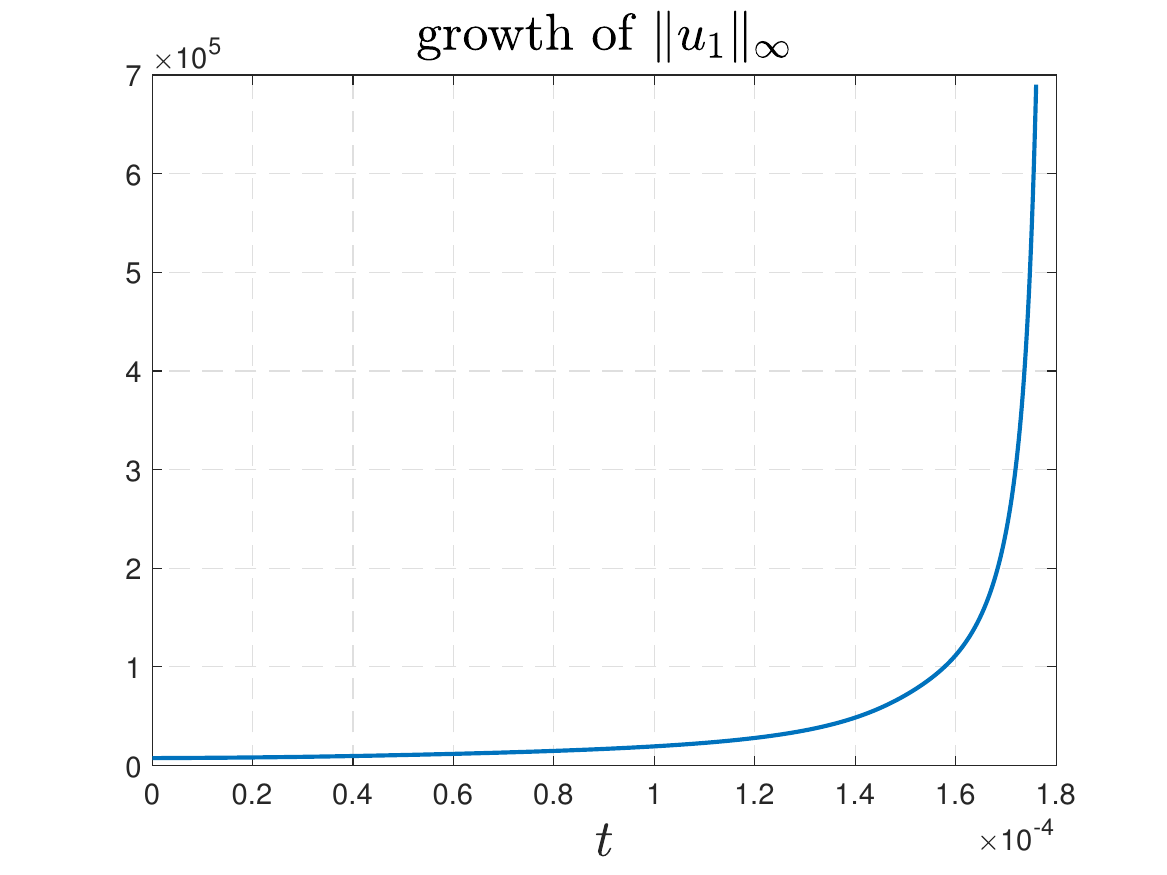}
    \includegraphics[width=0.32\textwidth]{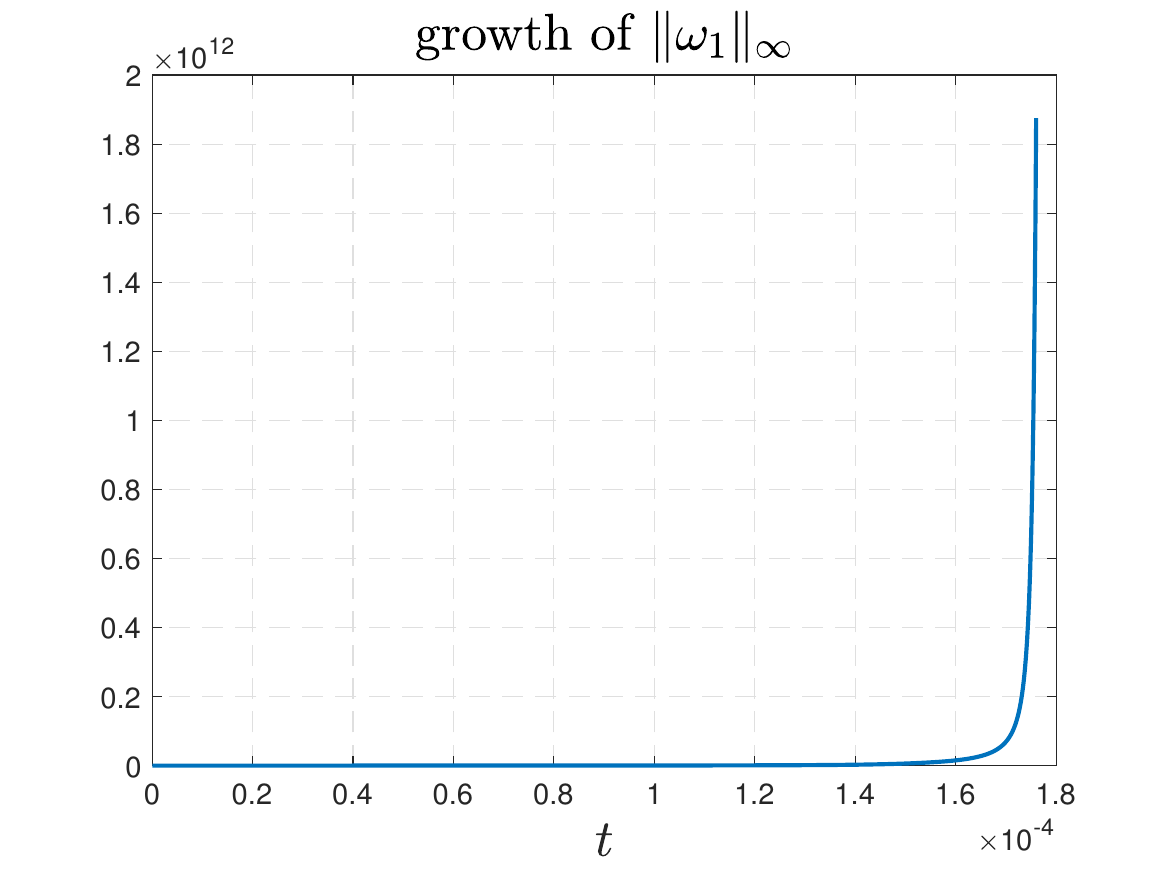} 
    \includegraphics[width=0.32\textwidth]{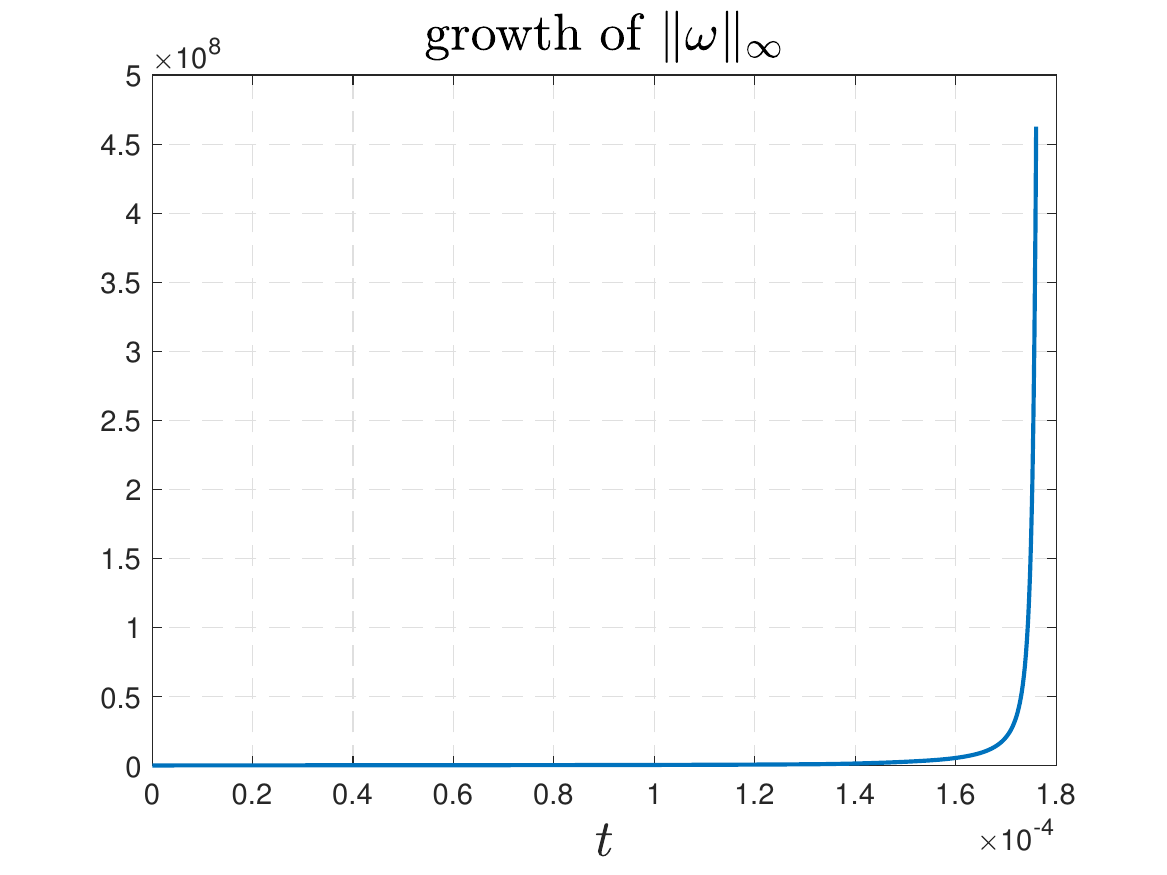}
    \includegraphics[width=0.32\textwidth]{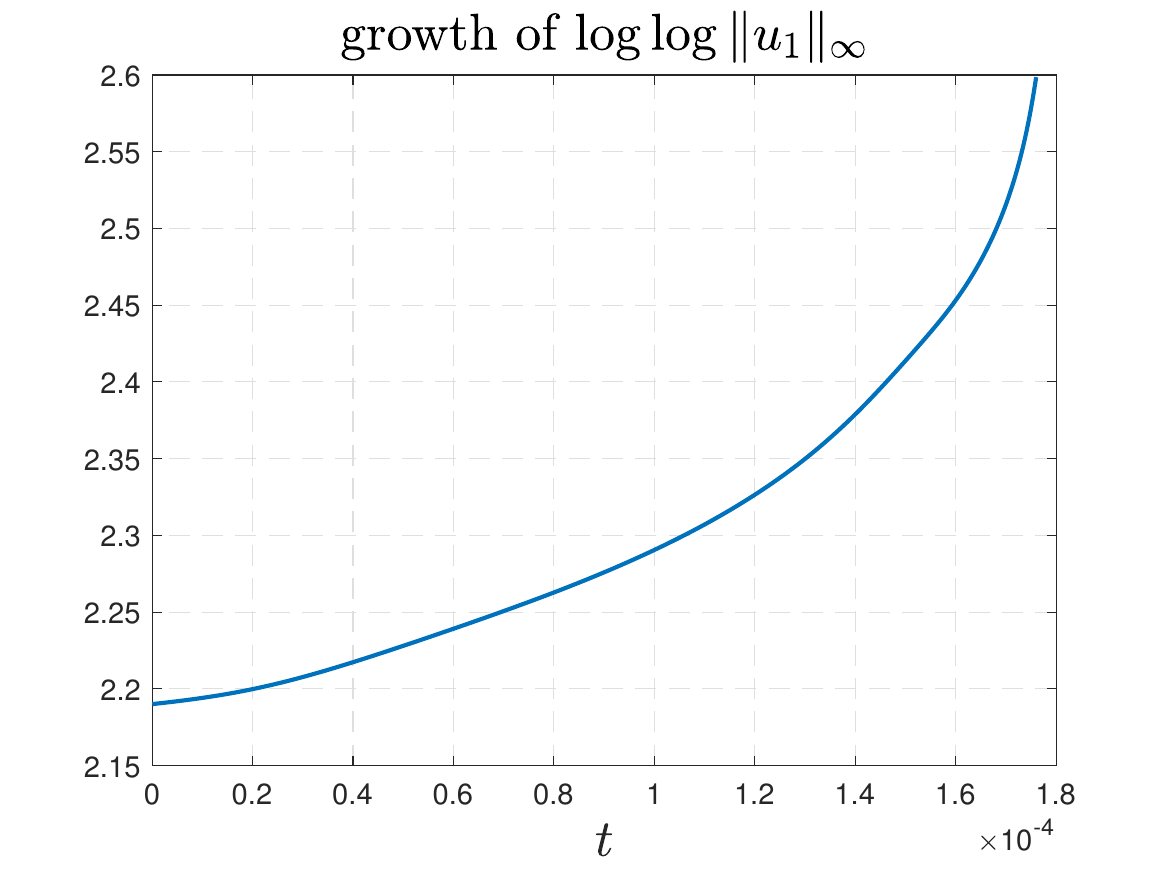}
    \includegraphics[width=0.32\textwidth]{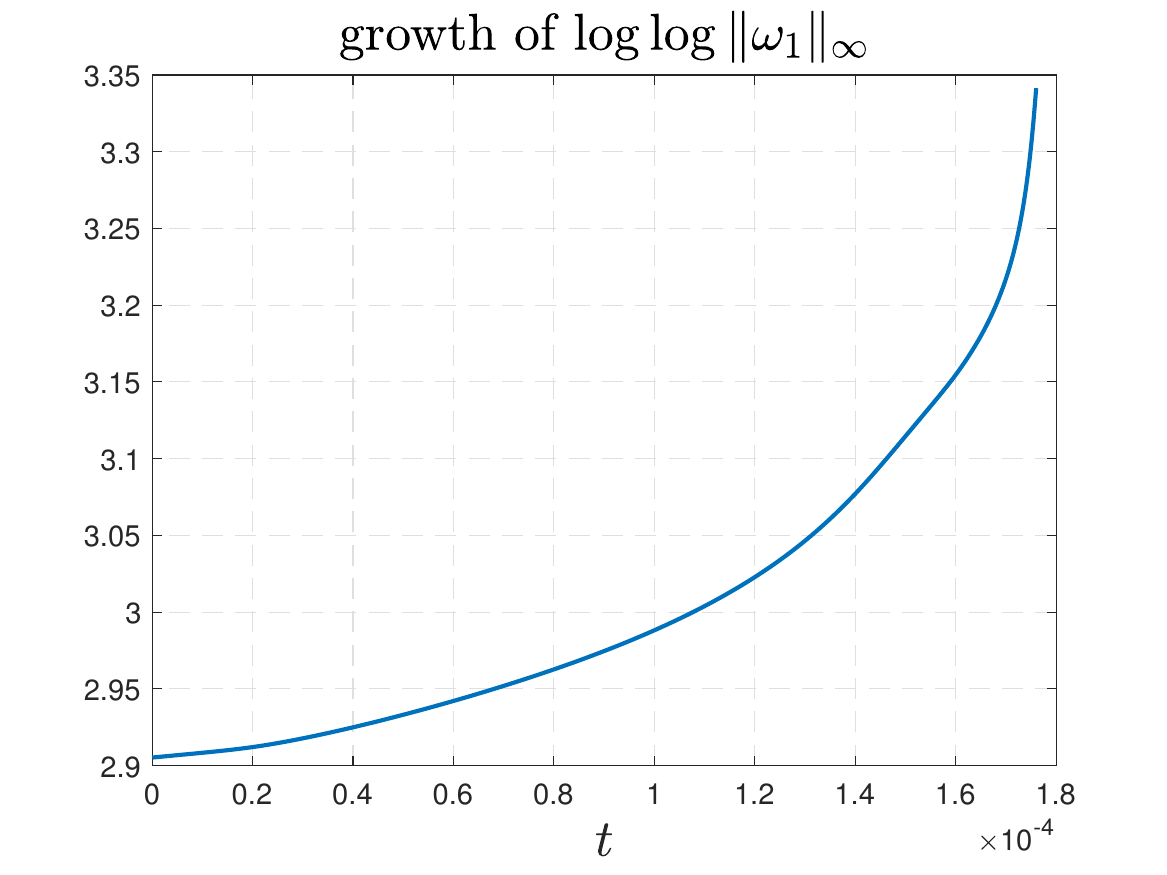}
    \includegraphics[width=0.32\textwidth]{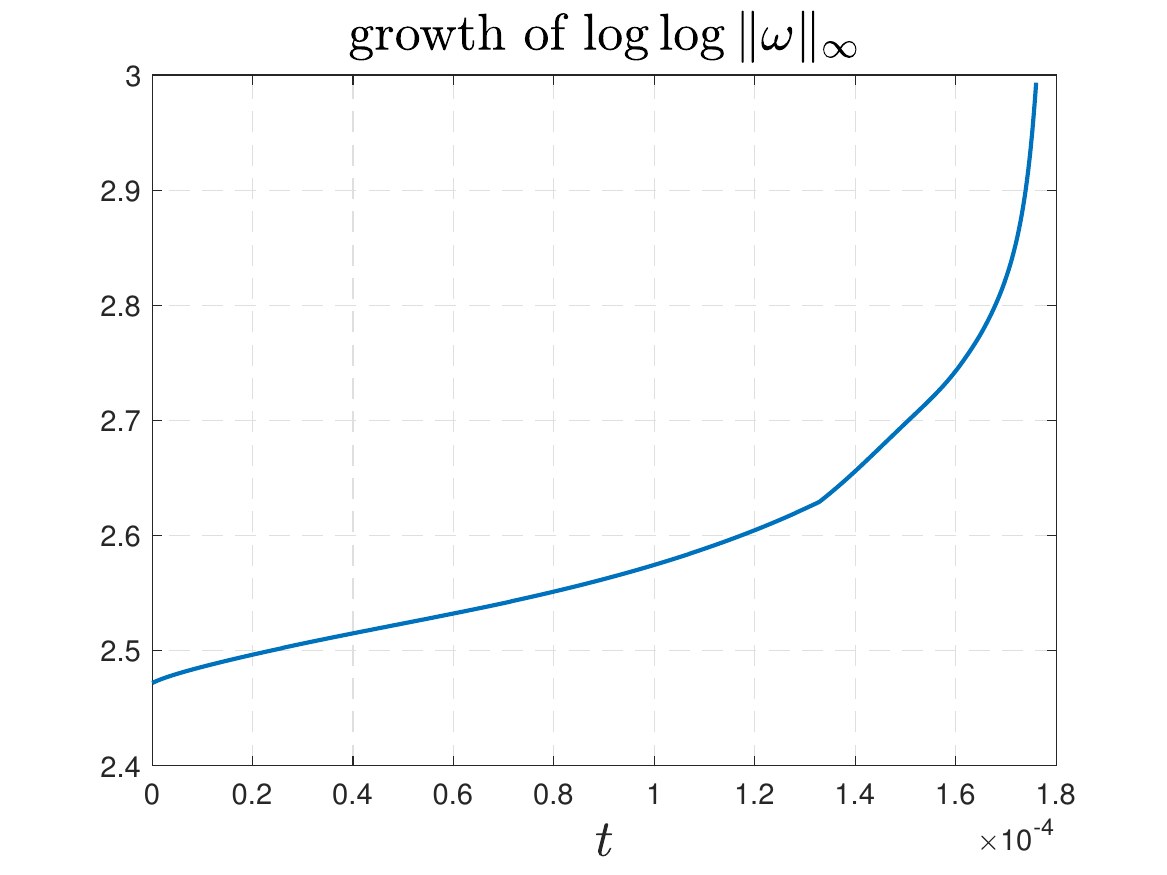}
    \caption[Rapid growth]{First row: the growth of $\|u_1\|_{L^\infty}$, $\|\om_1\|_{L^\infty}$ and $\|\vom\|_{L^\infty}$ as functions of time. Second row: $\log\log\|u_1\|_{L^\infty}$, $\log\log\|\om_1\|_{L^\infty}$ and $\log\log\|\vom\|_{L^\infty}$.} 
    \label{fig:rapid_growth}
\end{figure}

The most important observation in our computation is the rapid growth of the solution. The maximums of $|u_1|,|\om_1|$ and $|\vom|$ as functions of time are reported in Figure \ref{fig:rapid_growth}. Here 
\[\vom = (\om^\theta,\om^r,\om^z)^T = (r\om_1\,,\,-r u_{1,z}\,,\,2u_1+ru_{1,r})^T\]
is the vorticity vector, and 
\[|\vom| = \sqrt{(\om^\theta)^2+(\om^r)^2 + (\om^z)^2}.\]
We can see that these variables grow rapidly in time. In particular, they grow rapidly in the stable phase ($t\in[1.6\times 10^{-4},1.75\times 10^{-4}]$). Moreover, the second row in Figure \ref{fig:rapid_growth} shows that the solution grows much faster than a double-exponential rate. 

The rapid growth of the maximum vorticity $\|\vom\|_{L^\infty}$ is an important indicator of a finite-time singularity. In fact, the famous Beale--Kato--Majda criterion \cite{beale1984remarks} states that the solution to the standard Euler equations ceases to exist in some regularity class $H^s$ (for $s\geq 3$) at some finite time $T_*$ if and only if
\begin{equation}\label{eq:BKM}
\int_0^{T_*}\|\vom(t)\|_{L^\infty}\idiff t = +\infty. 
\end{equation}

Although the Beale--Kato--Majda criterion does not apply to the case of degenerate viscosity coefficients directly, we can still use an argument similar to the $H^s$ estimate of $\vu$ in \cite{beale1984remarks} to show that the solution to the Euler equations \eqref{eq:NSE_vc} with a degenerate viscosity coefficient $\nu$ ceases to exist in some regularity class $H^s$ ($s \geq 3$) if and only if 
\[\int_0^{T_*}\|\nabla\vu(t)\|_{L^\infty}\idiff t = +\infty.\]
Moreover, it is clear that $\|\vom\|_{L^\infty}\lesssim \|\nabla\vu\|_{L^\infty}$. Therefore, the rapid growth of maximum vorticity $\|\vom\|_{L^\infty}$ is still a good indicator for a finite-time singularity even in the case of a degenerate viscosity coefficient. We thus still view $\|\vom\|_{L^\infty}$ as a quantity of interest in our discussions. We will demonstrate in Section~\ref{sec:scaling_study} that the growth of $\|\vom\|_{L^\infty}$ has a nice fitting (with $R^2$ value greater that $0.9999$) to an inverse power law
\[\|\vom(t)\|_{L^\infty} \approx (T-t)^{-\gamma}\]
for some finite time $T$ and some power $\gamma>1$ (see Section \ref{sec:growth_fitting}). This then implies that the solution shall develop a potential singularity at some finite time $T$.

\subsection{Velocity field} In this subsection, we investigate the geometric structure of the velocity field. We first study the $3$D velocity field $\vu = u^r\vct{e}_r + u^\theta \vct{e}_\theta + u^z \vct{e}_z$ (denoted by $(u^r,u^\theta,u^z)$) by looking at the induced streamlines. An induced streamline $\{\Phi(s;X_0)\}_{s\geq0}\subset \mathbb{R}^3$ is completely determined by the background velocity $\vu$ and the initial point $X_0 = (x_0,y_0,z_0)^T$ through the initial value problem
\[\frac{\partial}{\partial s}\Phi(s;X_0) = \vu(\Phi(s;X_0)),\quad s\geq 0;\quad  \Phi(0;X_0) = X_0.\]
We remark that the induced streamlines do not give the particle trajectories in the real computation; they only characterize the geometric structure of the velocity field $\vu(t)$ for a fixed physical time $t$. The parameter $s$ dose not correspond to the physical time $t$. 

We will generate different streamlines with different initial points $X_0 = (r_0\cos(2\pi \theta),r_0\sin(2\pi\theta),z_0)^T$. Since the velocity field $\vu$ is now axisymmetric, the geometry of the streamline only depends on $(r_0,z_0)$. Varying the angular parameter $\theta$ only demonstrates the rotational symmetry of the streamlines.  

\subsubsection{A tornado singularity} 
\begin{figure}[!ht]
\centering
    \begin{subfigure}[b]{0.40\textwidth}
        \centering
        \includegraphics[width=1\textwidth]{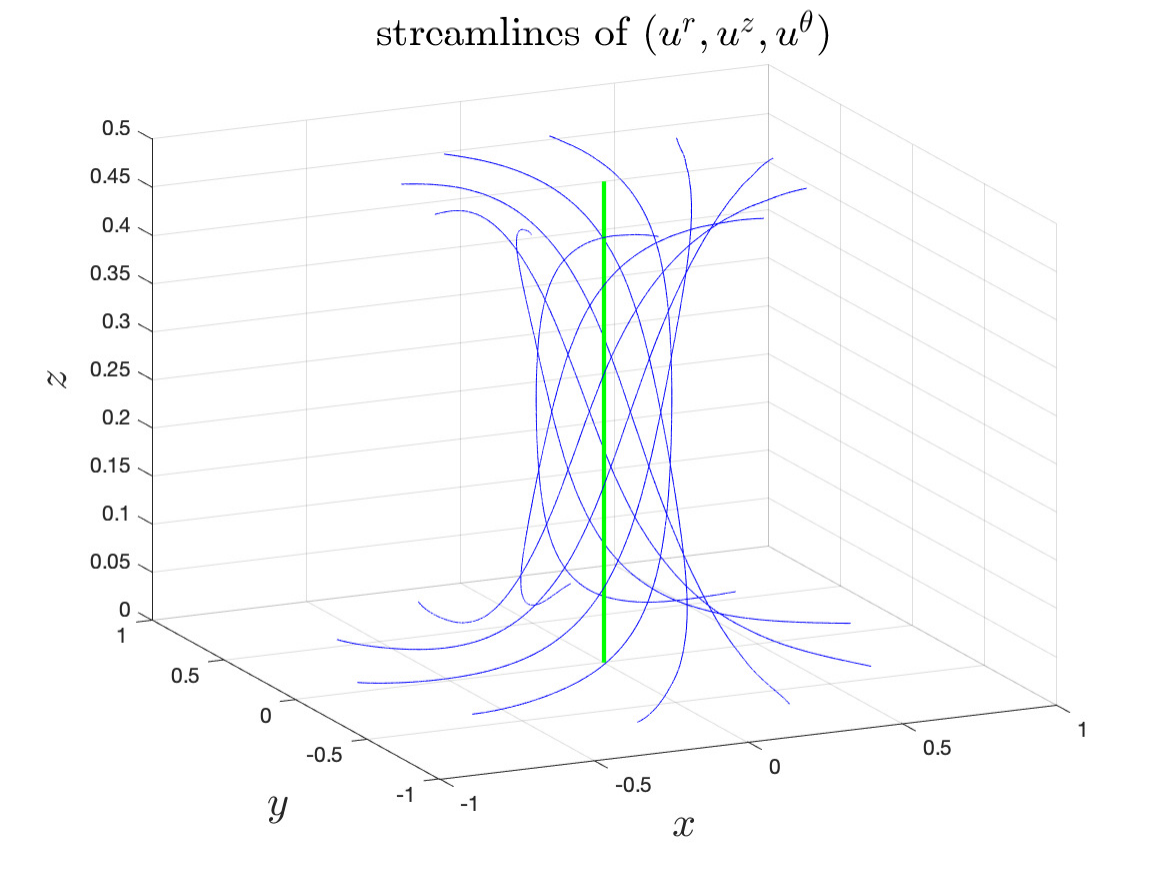}
        \caption{$r_0 = 0.8$, $z_0 = 0.01$}
    \end{subfigure}
    \begin{subfigure}[b]{0.40\textwidth}
        \centering
        \includegraphics[width=1\textwidth]{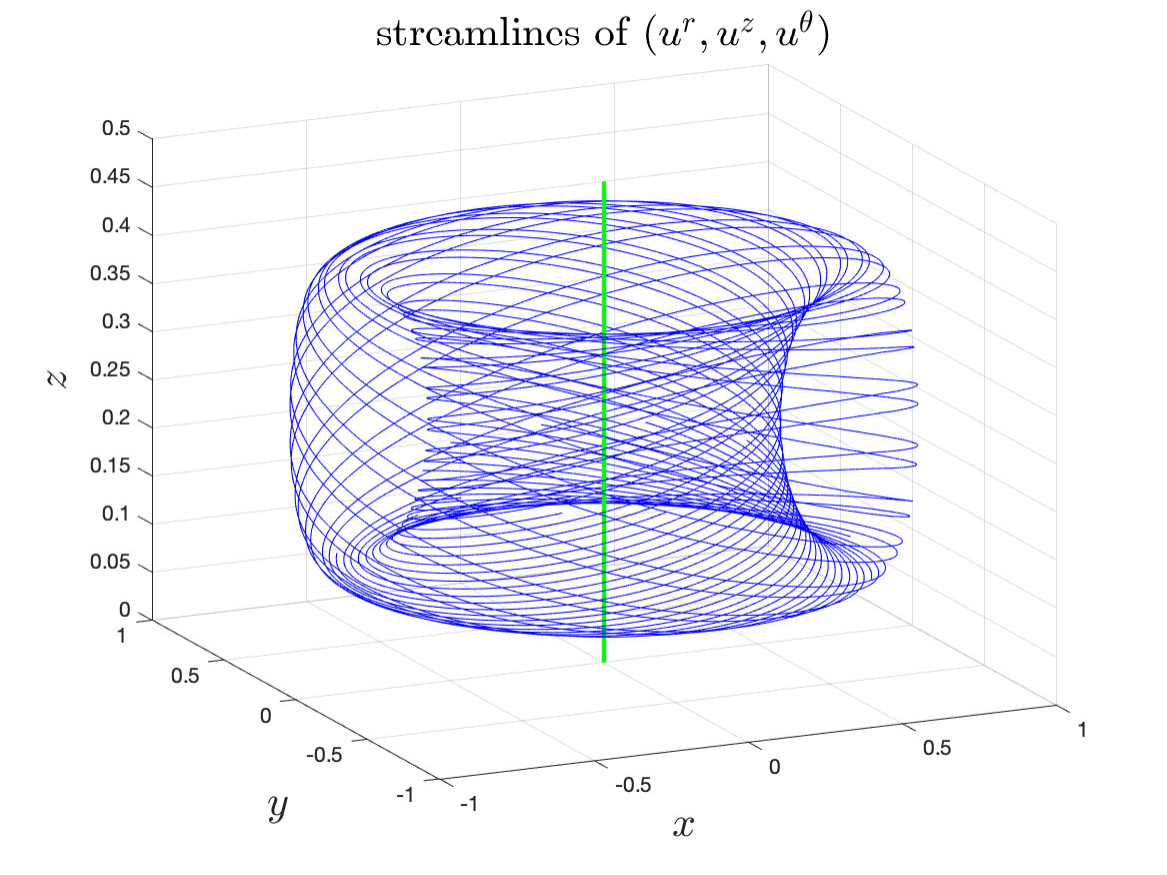}
        \caption{$r_0 = 0.8$, $z_0 = 0.1$}
    \end{subfigure}
    \caption[Global streamline]{The streamlines of $(u^r(t),u^\theta(t),u^z(t))$ at time $t=1.7\times 10^{-4}$ with initial points given by (a) $(r_0,z_0) = (0.8,0.01)$ and (b) $(r_0,z_0) = (0.8,0.1)$. The green pole is the symmetry axis $r=0$.}  
     \label{fig:streamline_3D_global}
        \vspace{-0.05in}
\end{figure}

Figure~\ref{fig:streamline_3D_global} shows the streamlines induced by the velocity field $\vu(t)$ at $t = 1.7\times 10^{-4}$ in a macroscopic scale (the whole cylinder domain $\mathcal{D}_1\times [0,2\pi]$) for different initial points with (a) $(r_0,z_0) = (0.8,0.01)$ and (b) $(r_0,z_0) = (0.8,0.1)$. The velocity field resembles that of a tornado spinning around the symmetry axis (the green pole). If the streamline starts near the ``ground'' ($z=0$) as in Figure~\ref{fig:streamline_3D_global}(a), it will first travel towards the symmetry axis, then move upward towards the ``ceiling'' ($z=1/2$), and at last turn outward away from the symmetry axis. In the meanwhile, it spins around the symmetry axis. On the other hand, if the initial point is higher (in the $z$ coordinate) as in Figure~\ref{fig:streamline_3D_global}(b), the streamline will not get very close the symmetry axis. Instead, it will travel in an ``inward-upward-outward-downward'' cycle in the $rz$-coordinates and, in the mean time, circle around the symmetry axis.

\begin{figure}[!ht]
\centering
    \begin{subfigure}[b]{0.32\textwidth}
        \centering
        \includegraphics[width=1\textwidth]{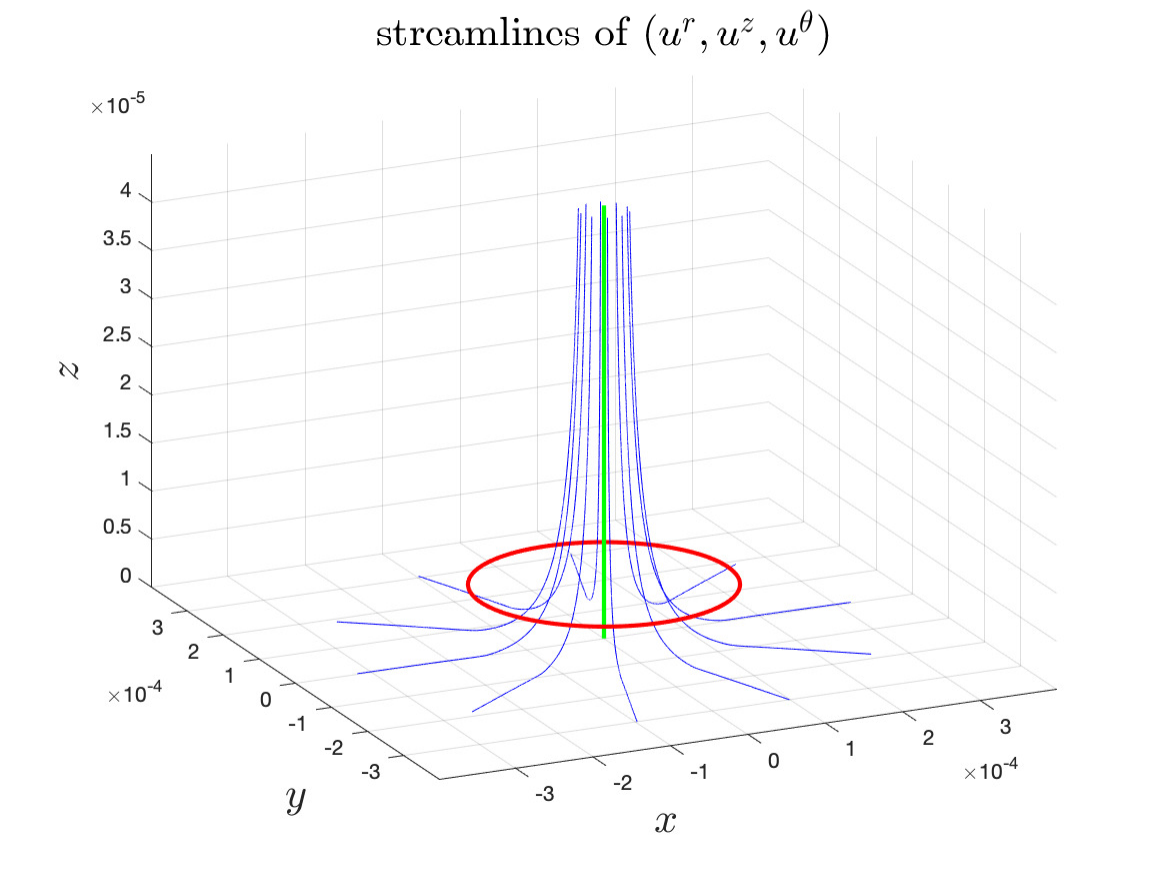}
        \caption{$r_0 = 2R(t)$, $z_0 = 0.01Z(t)$}
    \end{subfigure}
    \begin{subfigure}[b]{0.32\textwidth}
        \centering
        \includegraphics[width=1\textwidth]{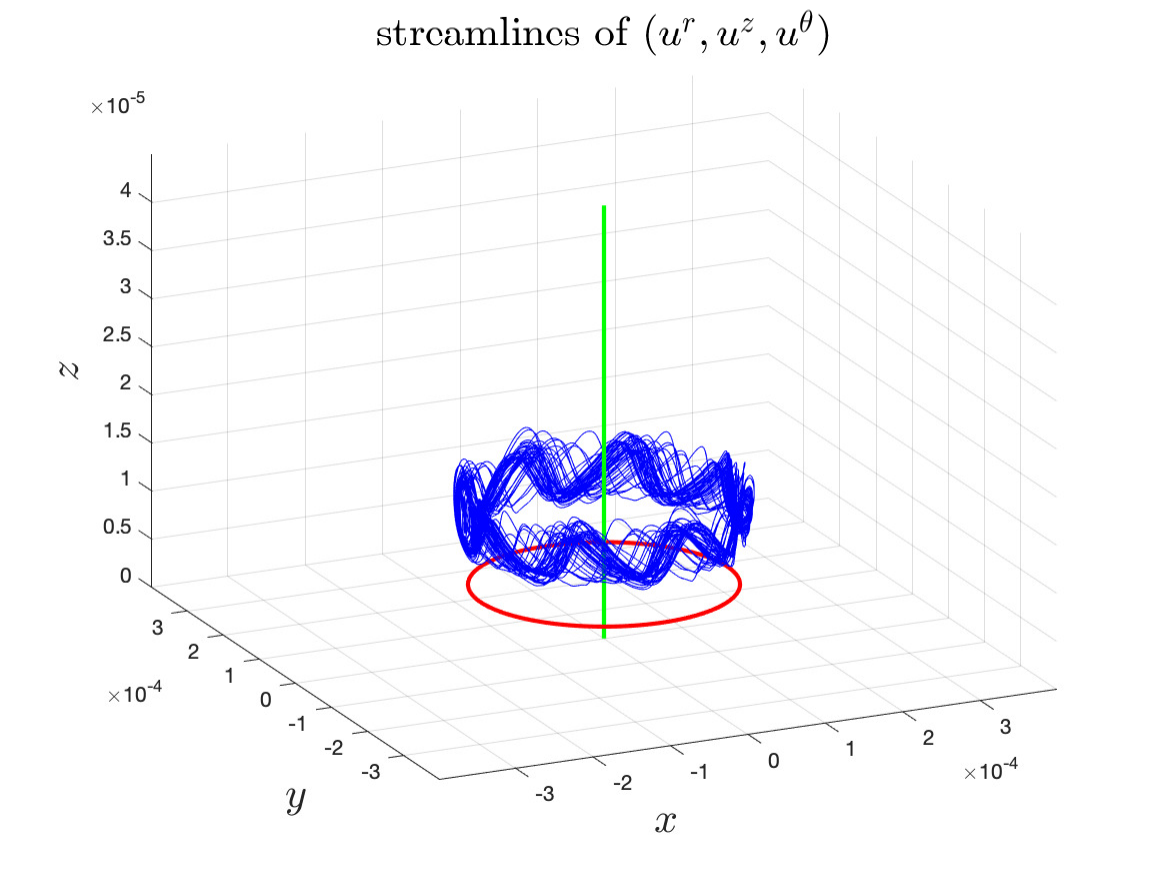}
        \caption{$r_0 = 1.05R(t)$, $z_0 = 2Z(t)$}
    \end{subfigure}
    \begin{subfigure}[b]{0.32\textwidth}
        \centering
        \includegraphics[width=1\textwidth]{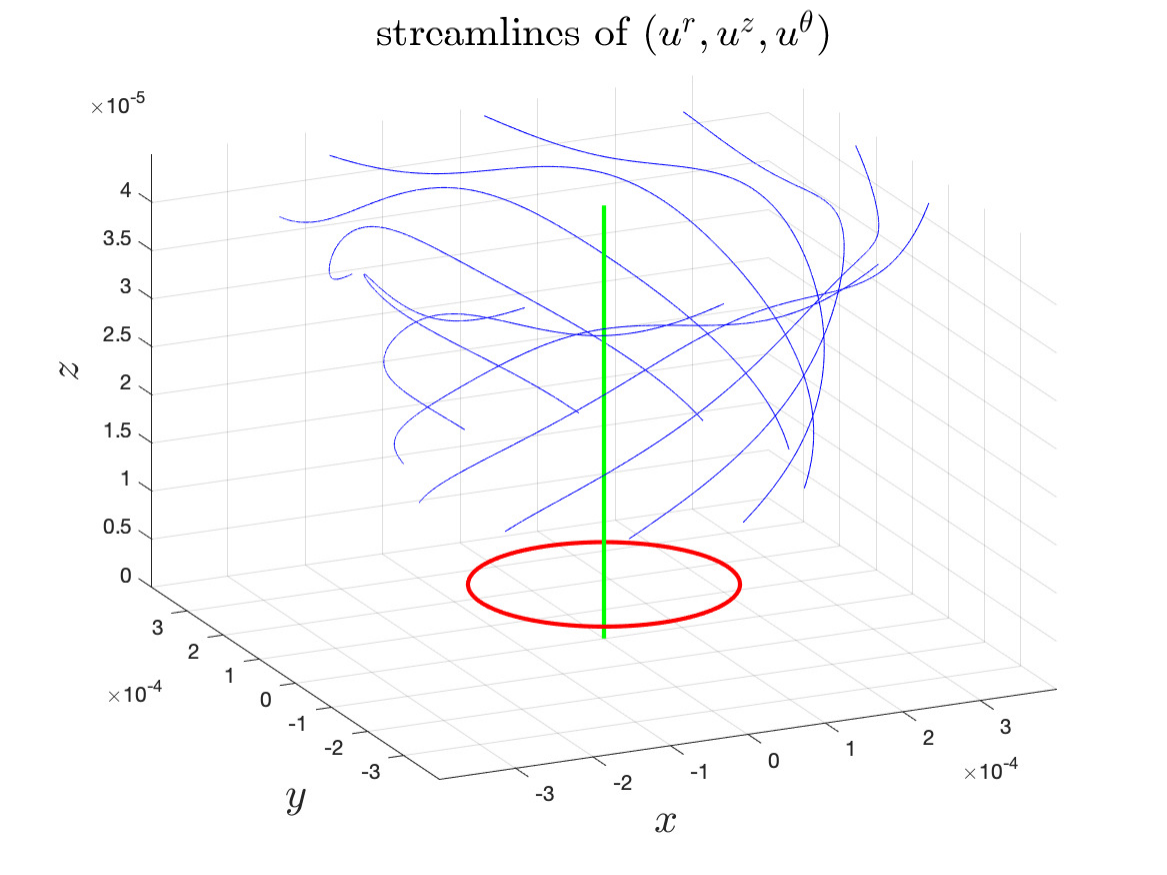}
        \caption{$r_0 = 1.5R(t)$, $z_0 = 3Z(t)$}
    \end{subfigure}
    \caption[Local streamline]{The streamlines of $(u^r(t),u^\theta(t),u^z(t))$ at time $t=1.7\times 10^{-4}$ with initial points given by (a) $(r_0,z_0) = (2R(t),0.01Z(t))$, (b) $(r_0,z_0) = (1.05R(t),2Z(t))$ and (c) $(r_0,z_0) = (1.5R(t),3Z(t))$. $(R(t),Z(t))$ is the maximum point of $u_1(t)$, indicated by the red ring. The green pole is the symmetry axis $r=0$.}  
     \label{fig:streamline_3D_zoomin}
        \vspace{-0.05in}
\end{figure}

Next, we take a closer look at the blowup region near the sharp front. Figure~\ref{fig:streamline_3D_zoomin} shows the streamlines at time $t=1.7\times 10^{-4}$ for different initial points near the maximum point $(R(t),Z(t))$ of $u_1(t)$. The red ring represents the location of $(R(t),Z(t))$, and the green pole is the symmetry axis $r=0$. The $3$ settings of $(r_0,z_0)$ are as follows.
\begin{itemize}
\item[(a)] $(r_0,z_0) = (2R(t),0.01Z(t))$. The streamline starts near the ``ground'' $z=0$ and below the red ring $(R(t),Z(t))$. It first travels towards the symmetry axis and then travels upward away from $z=0$. The spinning is weak since $u^\theta = ru_1$ is small in the corresponding region.
\item[(b)] $(r_0,z_0) = (1.05R(t),2Z(t))$. The streamline starts right above the ring $(R(t),Z(t))$. It gets trapped in a local region, oscillating and spinning around the symmetry axis periodically. The spinning is strong.
\item[(c)] $(r_0,z_0) = (1.5R(t),3Z(t))$.  The streamline starts even higher and away from the ring $(R(t),Z(t))$. It spins upward and outward, traveling away from the blowup region.
\end{itemize}

\subsubsection{The $2$D flow} 
\begin{figure}[!ht]
\centering
    \begin{subfigure}[b]{0.40\textwidth}
        \centering
        \includegraphics[width=1\textwidth]{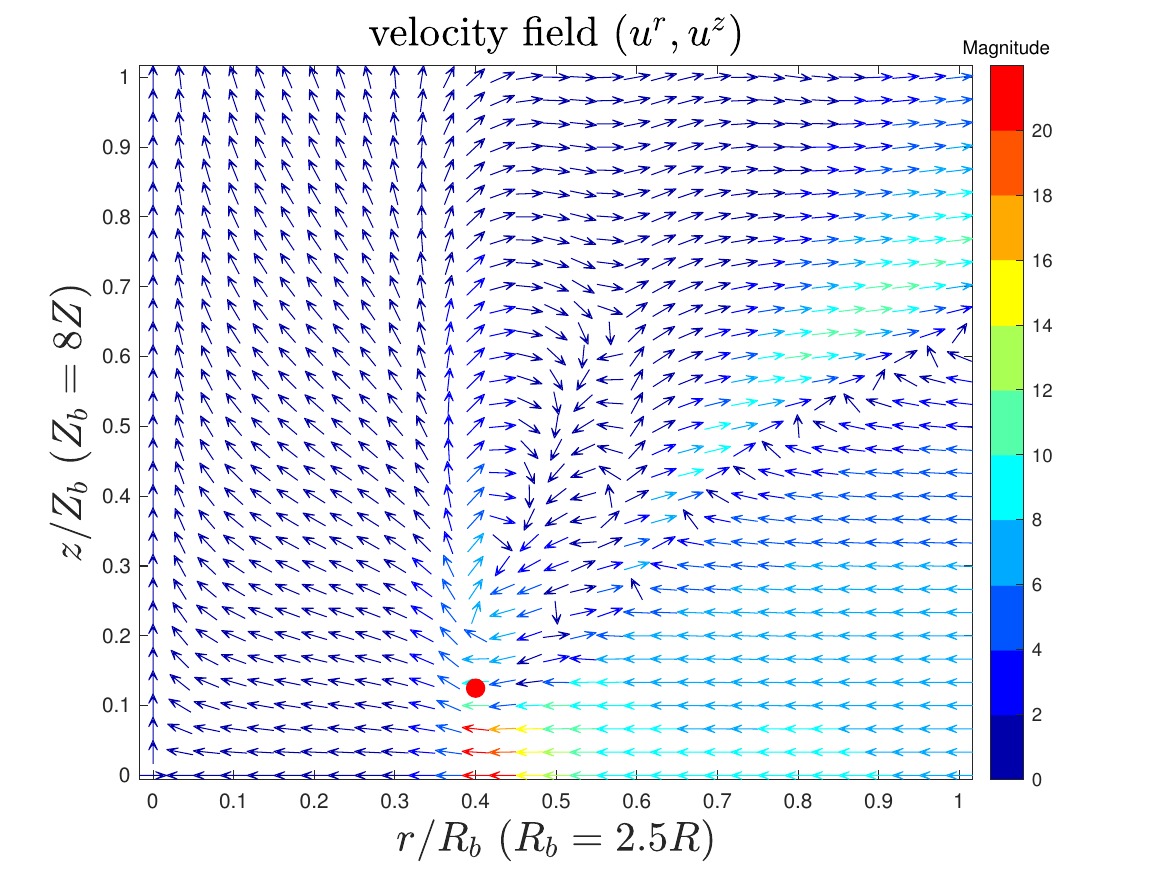}
        \caption{the velocity field $(u^r,u^z)$}
    \end{subfigure}
    \begin{subfigure}[b]{0.40\textwidth}
        \centering
        \includegraphics[width=1\textwidth]{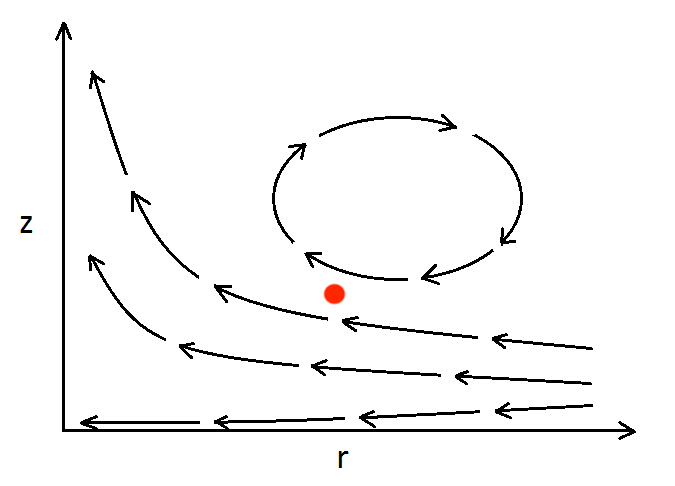} 
        \caption{a schematic}
    \end{subfigure}
    \caption[Global streamline]{(a) The velocity field $(u^r(t),u^z(t))$ near the maximum point $(R(t),Z(t))$ of $u_1(t)$ (the red point) at $t=1.7\times 10^{-4}$. The color corresponds to the magnitude of $\sqrt{(u^r)^2+(u^z)^2}$. The size of the domain has been rescaled. (b) A schematic of the vector field near the point $(R(t),Z(t))$.}  
     \label{fig:velocity_field}
        \vspace{-0.05in}
\end{figure}

To understand the phenomena in the blowup region as shown in Figure~\ref{fig:streamline_3D_zoomin}, we look at the $2$D velocity field $(u^r,u^z)$ in the computational domain $\mathcal{D}_1$. Figure~\ref{fig:velocity_field}(a) shows the vector field of $(u^r(t),u^z(t))$ in a local microscopic domain $[0,R_b]\times [0,Z_b]$, where $R_b = 2.5R(t) \approx 3.97\times 10^{-4}$ and $Z_b = 8Z(t) = 4.50\times 10^{-5}$. The domain has been rescaled in the figure for better visualization. Figure~\ref{fig:velocity_field}(b) is a schematic for the vector field in Figure~\ref{fig:velocity_field}(a). 

We can see that the streamline below $(R(t),Z(t))$ first travels towards $r=0$ and then move upward away from $z=0$, bypassing the sharp front near $(R(t),Z(t))$, which again demonstrates the ``two-phase'' feature of the flow. As the flow gets close to $r=0$, the strong axial velocity $u^z$ transports $u_1$ from near $z=0$ upward along the $z$ direction. Due to the odd symmetry of $u_1$, the angular velocity $u^\theta = ru^1$ is almost $0$ in the region near $z=0$. As a consequence, the upward stream dynamically generates a no-spinning region between the sharp front of $u_1$ and the symmetry axis $r=0$. This no-spinning region resembles the calm eye of a tornado, an area of relatively low wind speed near the center of the vortex. This explains why the streamlines almost do not spin around the symmetry axis in this region, as illustrated in Figure~\ref{fig:streamline_3D_zoomin}(a). 

Moreover, the velocity field $(u^r(t),u^z(t))$ forms a closed circle right above $(R(t),Z(t))$ as illustrated in Figure~\ref{fig:velocity_field}(b). The corresponding streamline is hence trapped in the circle region in the $rz$-plane. Since $u_1$ is large in this region (see Figure~\ref{fig:zoomin_profile}(a) and use the red point as a reference), the fluid spins fast around the symmetry axis $r=0$. Consequently, the corresponding streamline travels fast inside a $3$D torus surrounding the symmetry axis. This explains the oscillating and circling in Figure~\ref{fig:streamline_3D_zoomin}(b). This local circle structure of the $2$D velocity field is critical in stabilizing the blowup process, as it keeps the major profiles of $u_1,\om_1$ traveling towards the origin instead of being pushed upward.

\begin{figure}[!ht]
\centering
    \includegraphics[width=0.40\textwidth]{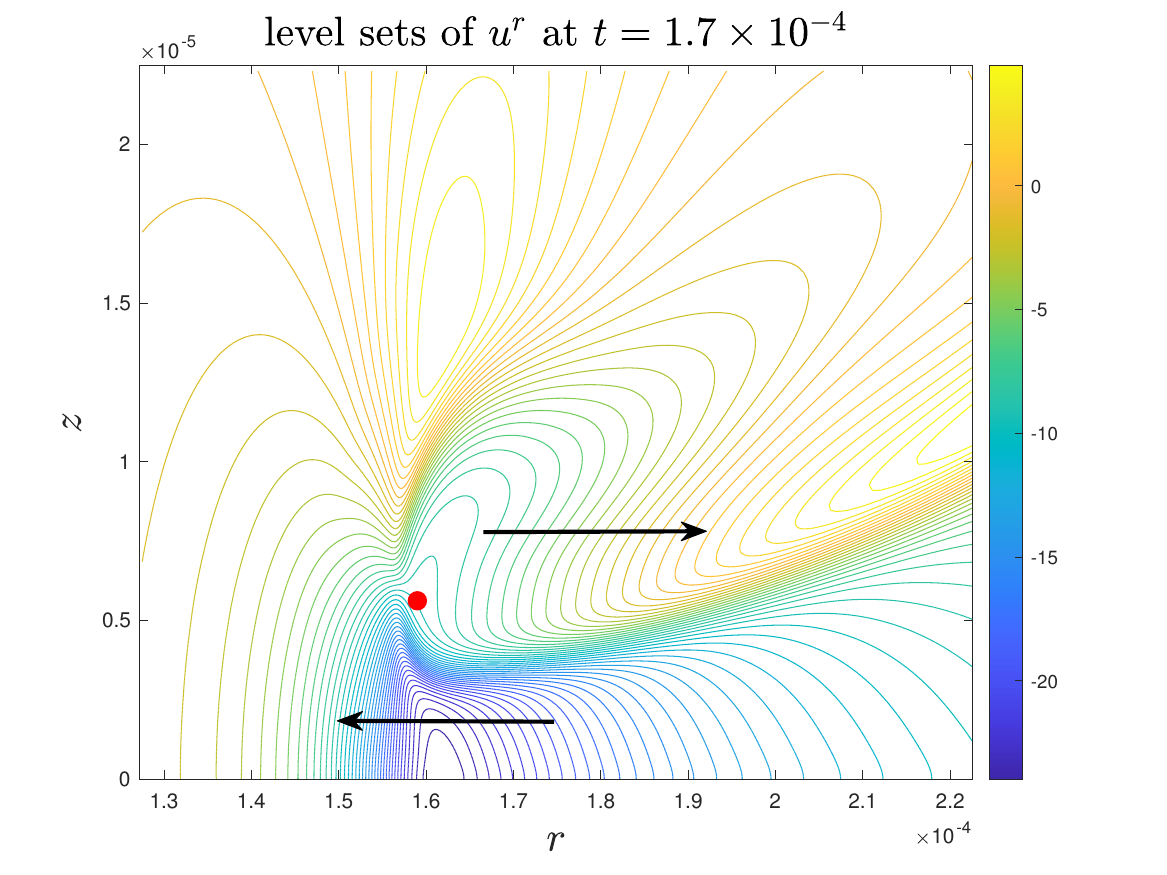}
    \includegraphics[width=0.40\textwidth]{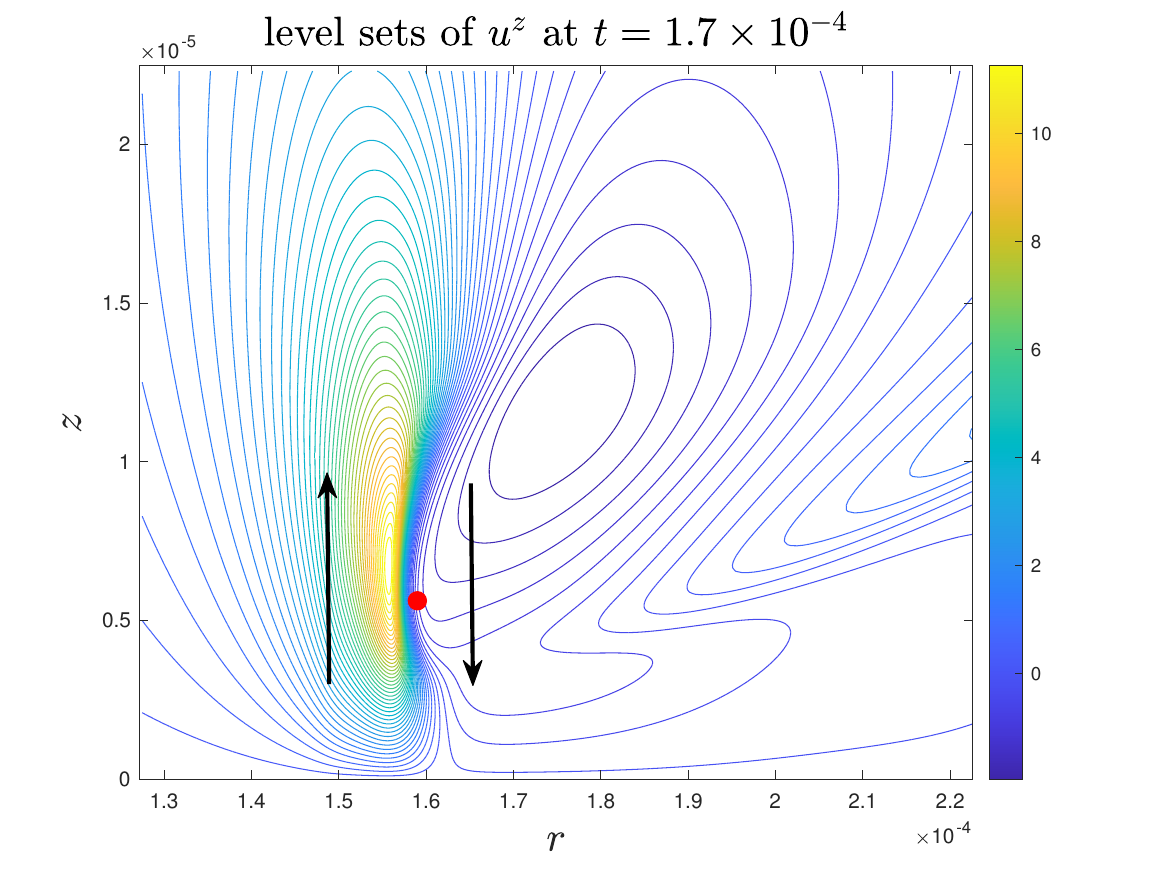} 
    \caption[Velocity level sets]{The level sets of $u^r$ (left) and $u^z$ (right) at $t = 1.7\times 10^{-4}$. The red point is the maximum point $(R(t),Z(t))$ of $u_1(t)$.}  
     \label{fig:velocity_levelset}
        \vspace{-0.05in}
\end{figure}

The velocity field $(u^r(t),u^z(t))$ can also explain the sharp structures of $u_1,\om_1$ in their local profiles (as shown in Figure~\ref{fig:zoomin_profile}(a),(b)). Figure~\ref{fig:velocity_levelset} shows the level sets of $u^r,u^z$ at $t = 1.7\times 10^{-4}$. One can see that the radial velocity $u^r$ has a strong shearing layer below $(R(t),Z(t))$ (the red point). This shearing contributes to the sharp gradient of $u_1$ in the $z$ direction (e.g., see Figure~\ref{fig:zoomin_profile}(a). Similarly, the axial velocity $u^z$ also has a strong shearing layer close to the point $(R(t),Z(t))$. This shearing explains the sharp front of $u_1$ in the $r$ direction. We will also explain in Section \ref{sec:asymptotic_analysis} the formation of a sharp front in the $r$ direction from a different perspective. 

\subsection{Understanding the blowup mechanism}\label{sec:mechanism} In this subsection, we elaborate our understanding of the potential blowup by examining several critical factors that lead to a sustainable blowup solution. 

\subsubsection{Vortex dipole and hyperbolic flow} Though we have only shown numerical results in the half-period domain $\mathcal{D}_1 = \{(r,z); 0\leq r\leq 1, 0\leq z\leq 1/2\}$, one should keep in mind that the real meaningful physics happens in the whole period $\{(r,z); 0\leq r\leq 1, -1/2\leq z\leq 1/2\}$. Moreover, the $2$D velocity field $(u^r, u^\theta)$ can be extended to the negative $r$ plane as an even function of $r$. The odd symmetry (in $z$) of the initial data of $\om_1$ leads to a dipole structure of the angular vorticity $\om^\theta$, which then induces a hyperbolic flow in the $rz$-plane with a pair of antisymmetric (with respect to $z$) local circulations. This pair of antisymmetric convective circulations is the cornerstone of our blowup scenario, as it has the desirable property of pushing the solution near $z=0$ towards $r=0$.

\begin{figure}[!ht]
\centering
    \includegraphics[width=0.40\textwidth]{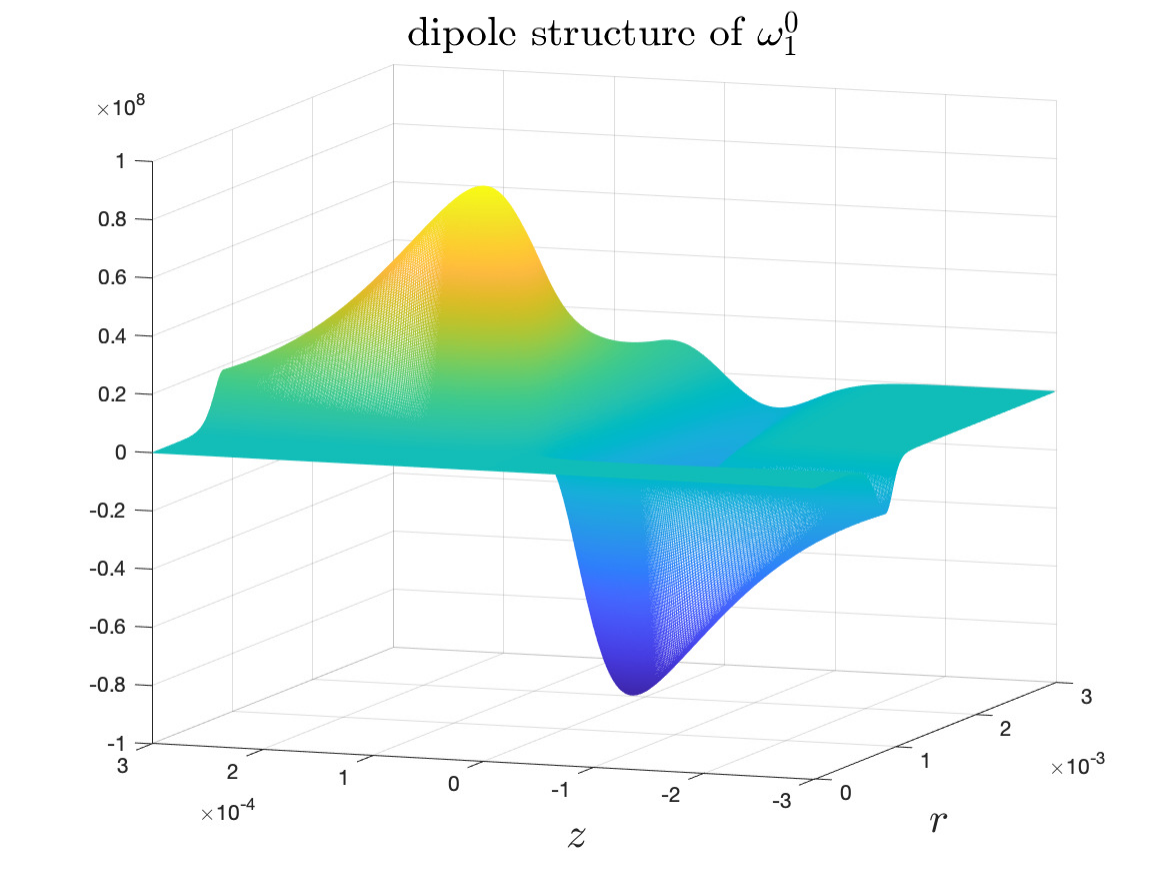}
    \includegraphics[width=0.40\textwidth]{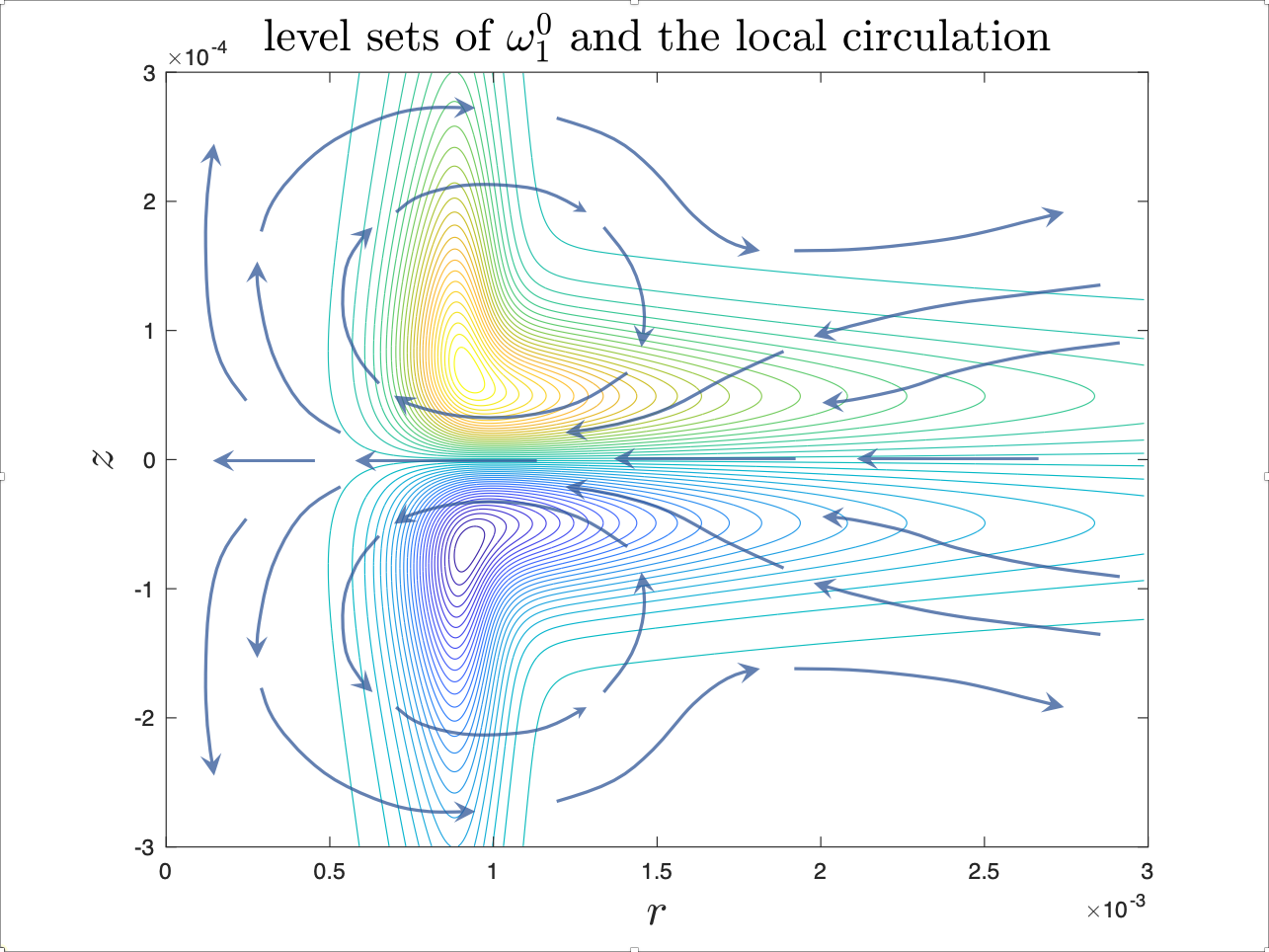} 
    \caption[Dipole]{The dipole structure of the initial data $\om_1^0$ and the induced local velocity field.}  
     \label{fig:dipole}
        \vspace{-0.05in}
\end{figure}

Figure~\ref{fig:dipole} presents the dipole structure of the initial data $\om_1^0$ in a local symmetric region $(r,z)\in [0,3\times 10^{-3}]\times [-3\times 10^{-4},3\times 10^{-4}]$ and the hyperbolic velocity field induced by it. The negative radial velocity near $z=0$ induced by the antisymmetric vortex dipole pushes the solution towards $r=0$, which is one of the driving mechanisms for a singularity on the symmetry axis. However, we also need another mechanism that squeezes the solution towards $z=0$, so that it can be driven by the inward velocity. Otherwise, the upward velocity away from $z=0$ may destroy the blowup trend. This critical squeezing mechanism comes from the odd symmetry of $u_1$.

\subsubsection{The odd symmetry and sharp gradient of $u_1$} 

The odd symmetry of $u_1$ is not required for $\om_1$ to be odd at $z=0$. The reason we construct $u_1$ to be an odd function of $z$ is that it ensures that $u_1^2$ has a large gradient in the $z$ direction near $z=0$.  

It is clear from the $\om_1$ equation \eqref{eq:as_NSE_1_b} that the driving force for $\om_1$ to blow up is the vortex stretching term $2(u_1^2)_z$. The odd symmetry of $u_1$ ensures that $u_1(r,0,t) = 0$ for all $t\geq0$. Therefore, $(u_1^2)_z$ is positive and large somewhere between $z = Z(t)$ and $z=0$, which drives $\om_1$ to grow fast near $z=0$. The growth of $\om_1$ then feeds to the growth of $u^r$ (in absolute value) around $z=0$, as a stronger dipole structure of the angular vorticity $\om^\theta$ induces a stronger inward flow in between the dipole (as demonstrated in Figure \ref{fig:dipole}). Note that $u^r$ being negative means $\psi_{1,z} = -u^r/r$ is positive, and the growth of $u^r$ around $z=0$ implies the growth of $\psi_{1,z}$ there, especially near $r = R(t)$. This in turn contributes to the rapid growth of $u_1$ in the critical region near $z=0$ through the vortex stretching term $2\psi_{1,z}u_1$ in the $u_1$ equation \eqref{eq:as_NSE_1_a}. 

Moreover, since $\psi_1=0$ along $z=0$ (by the odd condition), the Poisson equation \eqref{eq:as_NSE_1_c} can be approximated by $\psi_{1,zz} \approx -\om_1$ near $z=0$. This means that $\psi_{1,z}$, as a function of $z$, achieves its local maximum at $z=0$ in a neighborhood where $\text{sign}(\om_1) = \text{sign}(z)$. The rapid growth of $\psi_{1,z}$ and the nonlinear vortex stretching term $2\psi_{1,z}u_1$ in the $u_1$ equation induce a traveling wave for $u_1$ propagating towards $z=0$, which drags the maximum point of $u_1$ towards $z=0$. The traveling wave is so strong that it overcomes the stabilizing effect of advection along the $z$ direction, which pushes the flow upward away from $z=0$. The fact that the maximum point of $u_1$ traveling towards $z=0$ generates an even sharper gradient of $u_1^2$ in the $z$ direction. The whole coupling mechanism above forms a positive feedback loop, 
\begin{equation}\label{eq:mechanism}
(u_1^2)_z \uparrow \quad \Longrightarrow \quad\om_1 \uparrow \quad \Longrightarrow \quad \psi_{1,z} \uparrow \quad \Longrightarrow \quad u_1 \uparrow \quad \Longrightarrow \quad (u_1^2)_z\uparrow,
\end{equation}
that leads to a sustainable blowup solution shrinking towards $z=0$ and traveling towards $r=0$. 

To trigger this mechanism, it is important that the initial data have the proper symmetry and a strong alignment between $u_1$ and $\om_1$ as described in Appendix \ref{sec:initial_data}. The maximum point of $\om_1$ should align with the location where $u_{1,z}$ is positive and large, which is slightly lower (in $z$) than the maximum point of $u_1$. This is one of the guiding principles in the construction of our initial data.

\begin{figure}[!ht]
\centering
	\begin{subfigure}[b]{0.32\textwidth}
    \includegraphics[width=1\textwidth]{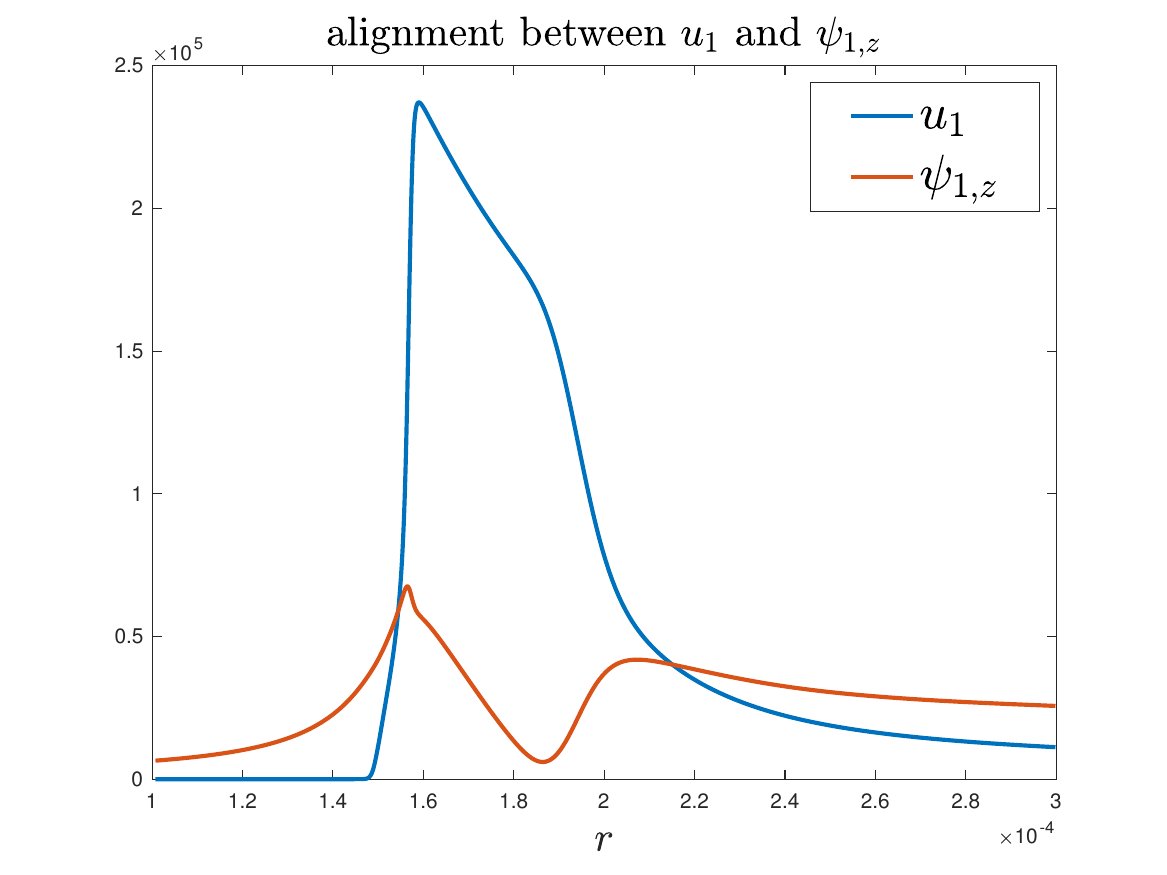}
    \caption{$r$ cross sections of $u_1,\psi_{1,z}$}
    \end{subfigure}
    \begin{subfigure}[b]{0.32\textwidth}
    \includegraphics[width=1\textwidth]{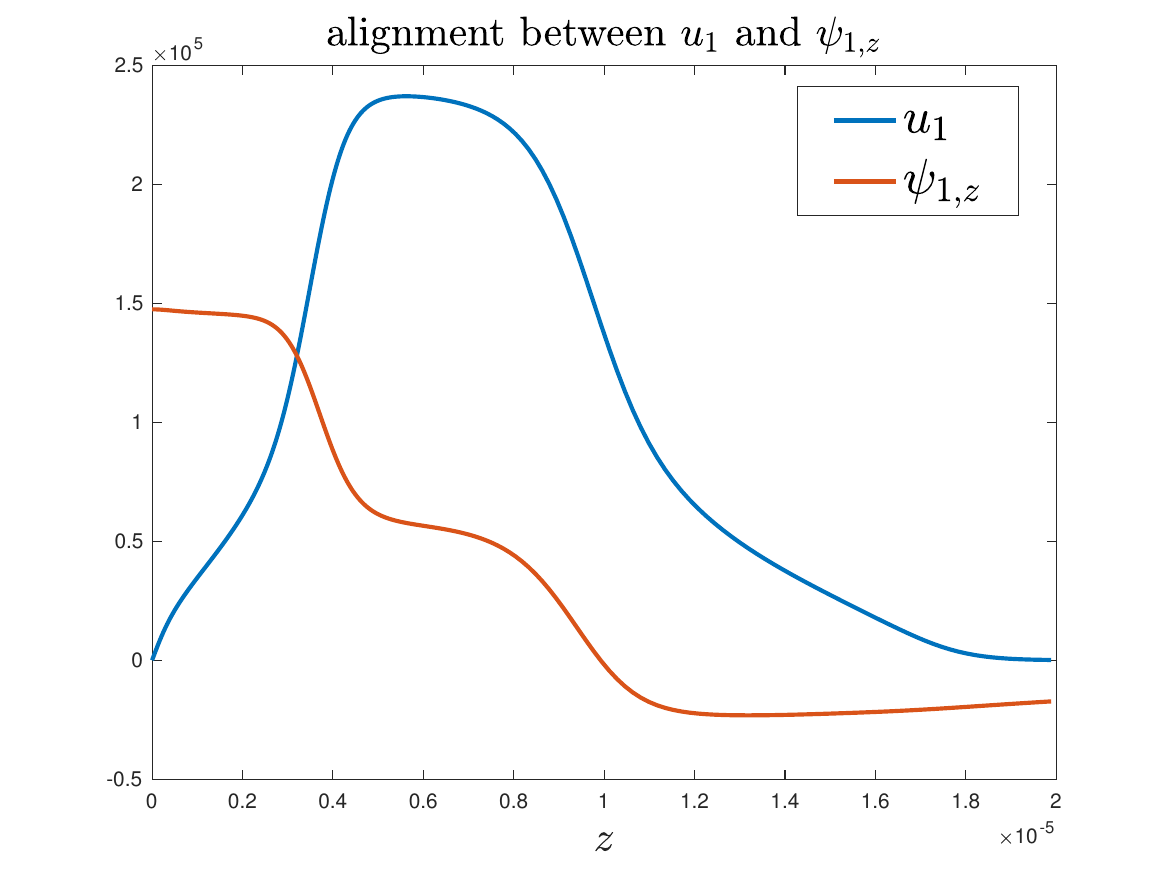}
    \caption{$z$ cross sections of $u_1,\psi_{1,z}$}
    \end{subfigure}
  	\begin{subfigure}[b]{0.32\textwidth}
    \includegraphics[width=1\textwidth]{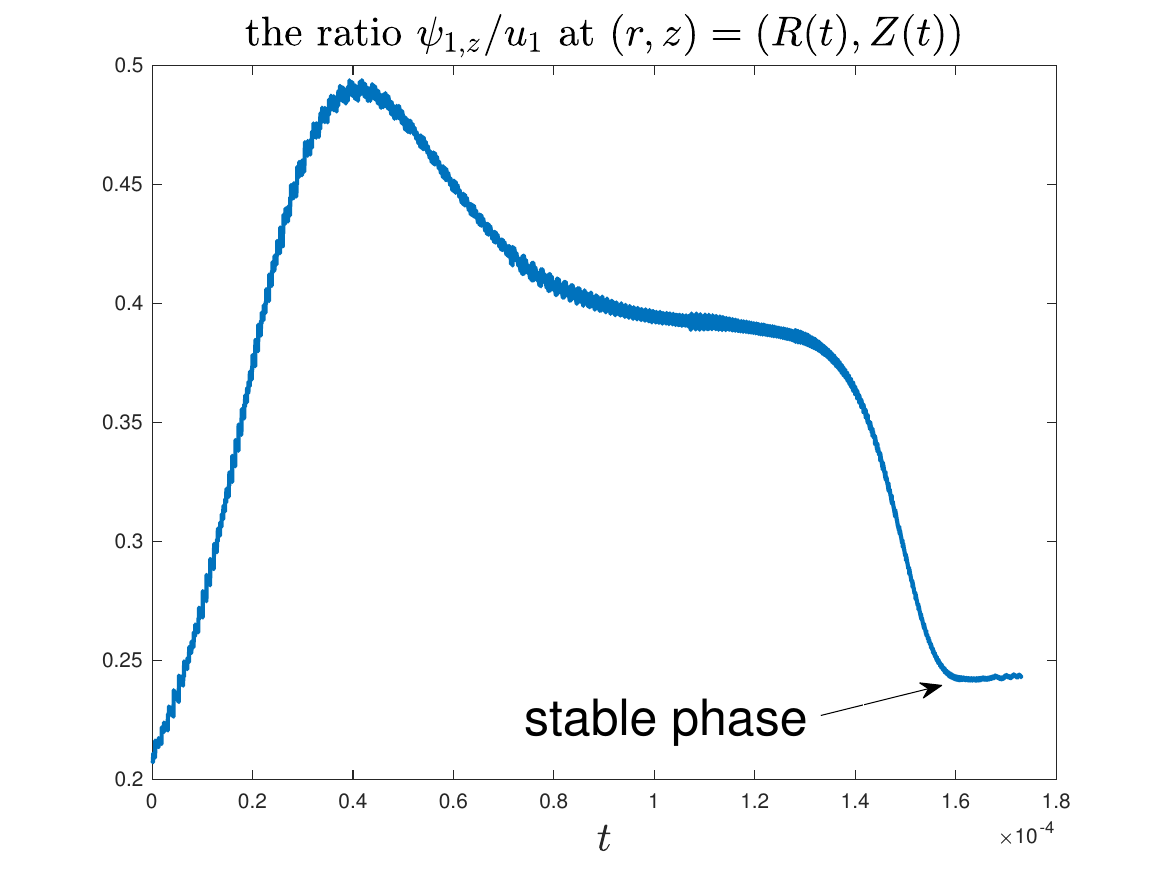}
    \caption{$\psi_{1,z}/u_1$ as a function of time}
    \end{subfigure}
    \caption[Alignment]{The alignment between $u_1$ and $\psi_{1,z}$. (a) and (b): cross sections of $u_1$ and $\psi_{1,z}$ through the point $(R(t),Z(t))$ at $t=1.7\times 10^{-4}$. (c): the ratio $\psi_{1,z}/u_1$ at the point $(R(t),Z(t))$ as a function of time up to $t=1.75\times 10^{-4}$.}  
     \label{fig:alignment}
        \vspace{-0.05in}
\end{figure}

Figure \ref{fig:alignment} demonstrates the alignment between $\psi_{1,z}$ and $u_1$. Figure \ref{fig:alignment}(b) shows the cross section of $u_1,\psi_{1,z}$ in the $z$ direction through $(R(t),Z(t))$ at $t = 1.7\times10^{-4}$. We can see that $\psi_{1,z}(R(t),z,t)$ is monotonically decreasing for $z \in [0, 2Z(t)]$ and relatively flat for $z \in [0, 0.5 Z(t)]$. Moreover, $\psi_{1,z}$ is large, positive, and comparable to $u_1$ in magnitude, which leads to the rapid growth of $u_1$ and pushes $Z(t)$ moving towards $z=0$. Figure \ref{fig:alignment}(c) shows the alignment ratio $\psi_{1,z}(R(t),Z(t),t)/u_1(R(t),Z(t),t)$, i.e. the alignment between $\psi_{1,z}$ and $u_1$ at the maximum point of $u_1$. One can see that the ratio $\psi_{1,z}/u_1$ settles down to a stable value at $(R(t),Z(t))$ in the stable phase which is characterized by the time interval $[1.6\times10^{-4},1.75\times10^{-4}]$; that is $\psi_{1,z}(R(t),Z(t),t) \sim u_1(R(t),Z(t),t)$ in the stable phase. Consequently, the vortex stretching term in the $u_1$ equation is formally quadratic at the maximum point of $u_1$ if we ignore the small viscosity: 
\[\frac{\diff\,}{\diff t} u_{1,\max} \approx u_{1,\max}^2\,,\]
which implies that maximum $u_1$ should blow up like $(T-t)^{-1}$ for some finite time $T$. We will see more clearly this observation in Section~\ref{sec:scaling_study}.

We remark that the above discussion on the potential blowup mechanism also applies to the $3$D axisymmetric Euler equations in the same scenario. We therefore expect that the solution to the Euler equations (in Case $3$) would develop a similar blowup if we were able to resolve the small scale features of the solution.

\subsection{Beyond the stable phase} One of the major features of the solution beyond the stable phase is the instabilities that appear in the tail part of the profile. These tail oscillations partially result from under-resolution, as they can be suppressed by increasing the number of grid points or applying stronger numerical regularization to the computation. However, according to our numerical experiments, these measures can only delay the occurrence of the instabilities until a slightly later time, and oscillations will eventually appear and be amplified by the strong shearing of the velocity along the tail. This implies that the potential blowup solution could be physically unstable. We may approach the potential singularity, but the unstable modes would prevent us from getting arbitrarily close to the potential singularity. In fact, as we can see in Section \ref{sec:Euler} of the Euler case, such unstable oscillations occur much earlier in time without the degenerate viscosity and the numerical regularization. A study of the unstable behavior of the solution beyond the stable phase is presented in Appendix \ref{sec:beyond}.

\section{Numerical Results: Resolution Study}\label{sec:performance_study}
In this section, we perform resolution study and investigate the convergence property of our numerical methods. In particular, we will study 
\begin{itemize}
\item[(i)] the effectiveness of the adaptive mesh (Section \ref{sec:mesh_effectiveness}), and 
\item[(ii)] the convergence of the solutions when $h_\rho,h_\eta\rightarrow 0$ (Section \ref{sec:resolution_study}).
\end{itemize}

\subsection{Effectiveness of the adaptive mesh}\label{sec:mesh_effectiveness} As discussed in Appendix~\ref{apdx:adaptive_mesh}, to effectively compute a potential singularity with finite resolution, it is important that the adaptive mesh well resolves the solution in the whole domain, especially in the most singular region. In this subsection, we study the effectiveness of our adaptive mesh. 

\begin{figure}[!ht]
\centering
	\begin{subfigure}[b]{0.40\textwidth}
    \includegraphics[width=1\textwidth]{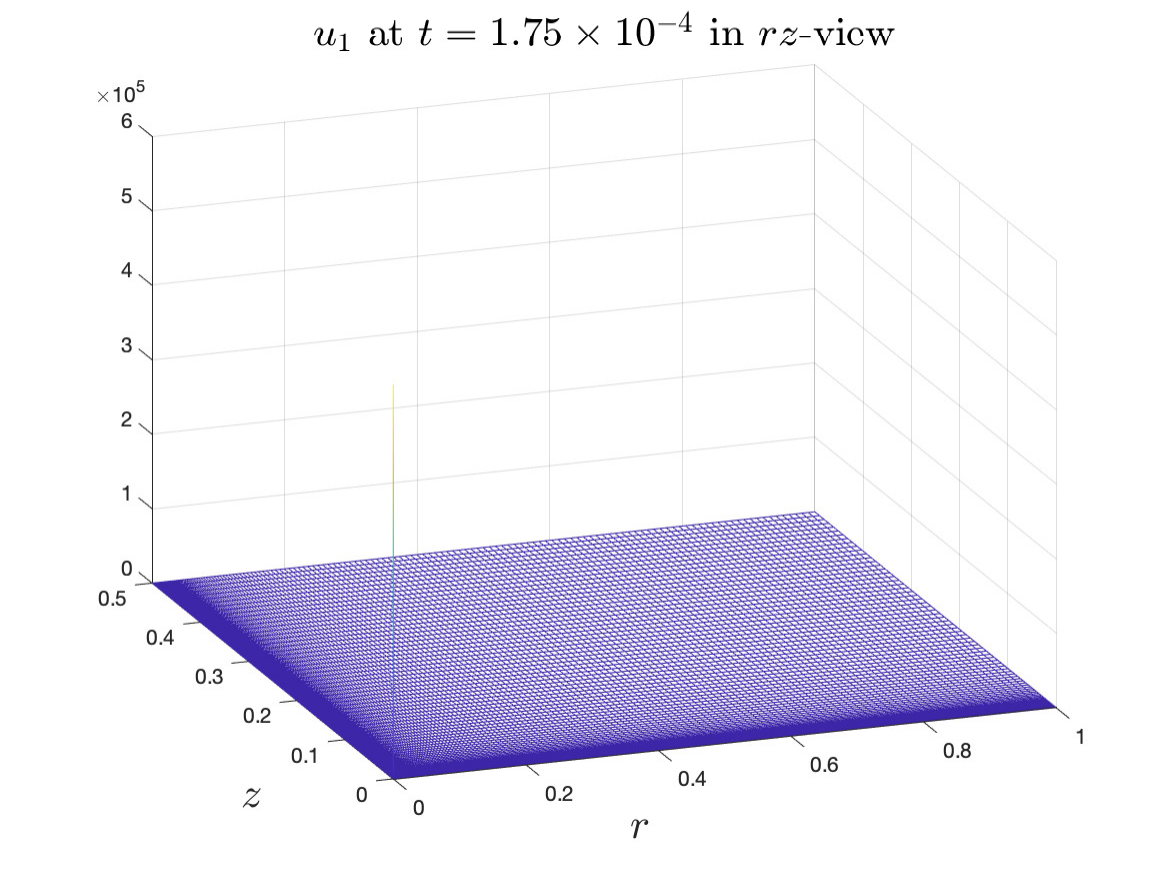}
    \caption{$u_1$ in the $rz$-plane}
    \end{subfigure}
  	\begin{subfigure}[b]{0.40\textwidth}
    \includegraphics[width=1\textwidth]{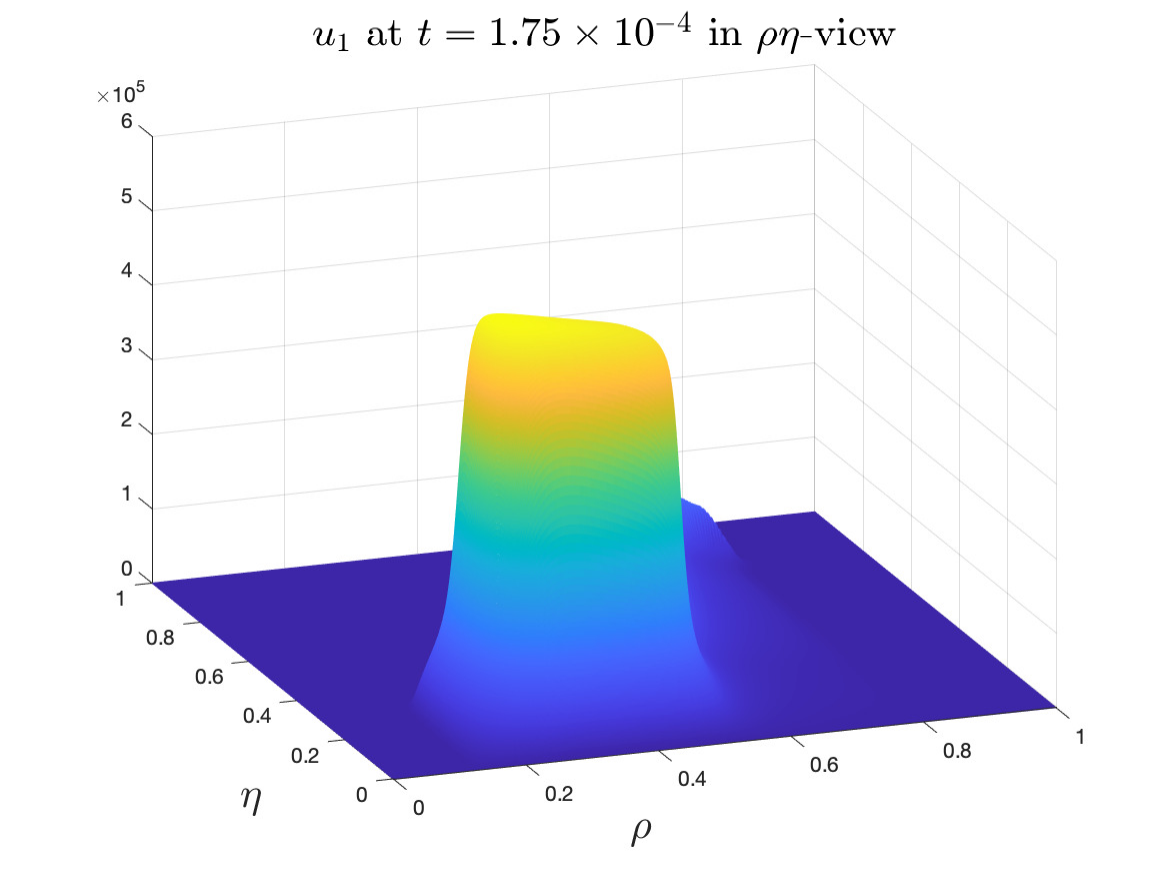}
    \caption{$u_1$ in the $\rho\eta$-plane}
    \end{subfigure}
    \begin{subfigure}[b]{0.40\textwidth}
    \includegraphics[width=1\textwidth]{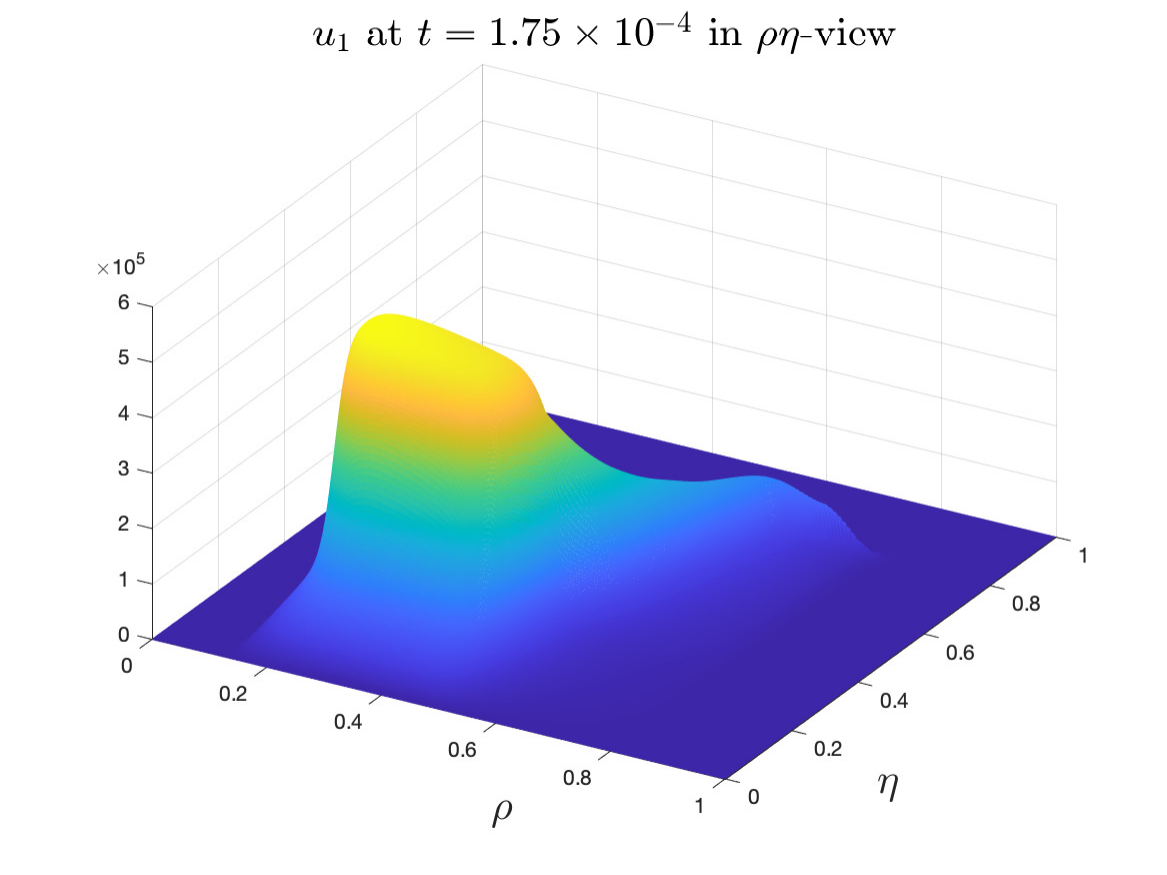}
    \caption{$u_1$ in the $\rho\eta$-plane (different angle)}
    \end{subfigure}
  	\begin{subfigure}[b]{0.40\textwidth}
    \includegraphics[width=1\textwidth]{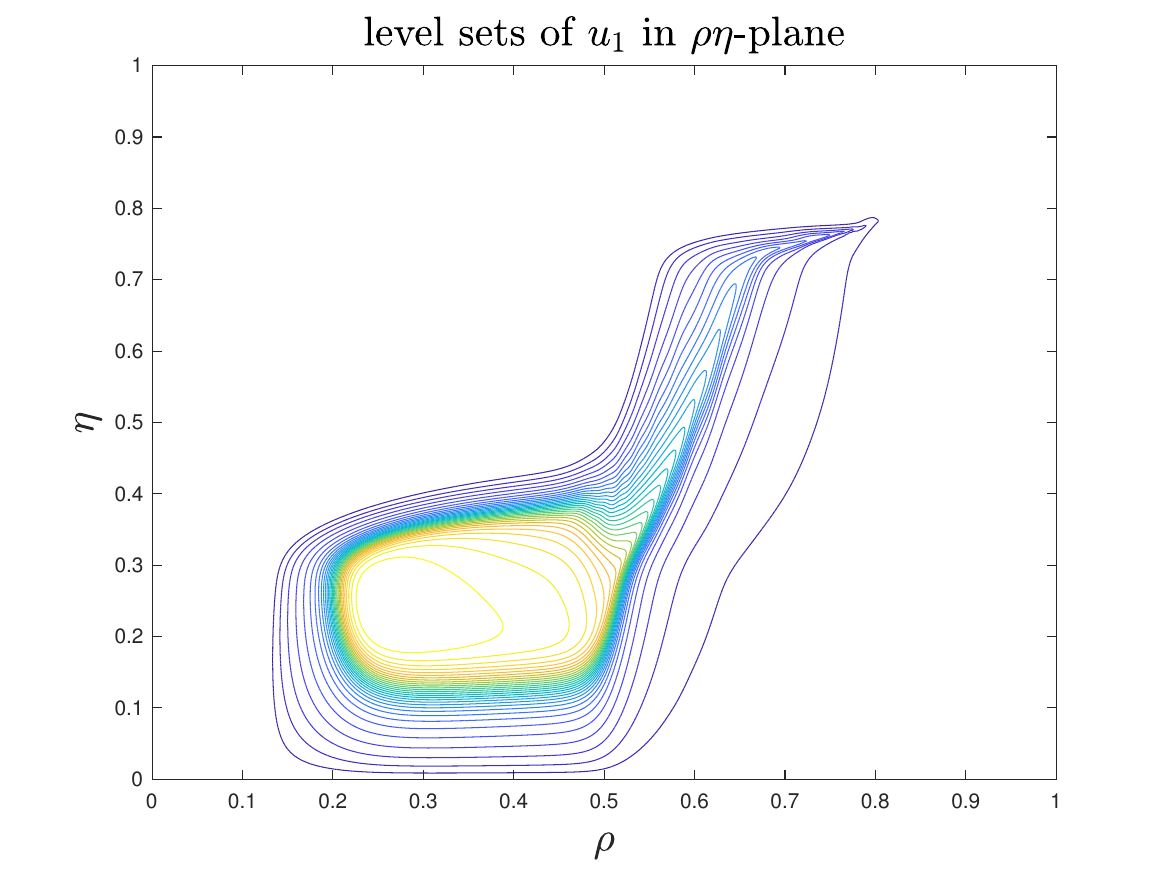}
    \caption{level sets of $u_1$ in the $\rho\eta$-plane}
    \end{subfigure}
    \caption[Mesh effectiveness]{The adaptive mesh resolves the solution in the $\rho\eta$-plane. (a) shows the focusing singularity structure of $u_1$ at $t=1.75\times10^{-4}$ in the $rz$-plane on the whole computational domain $\mathcal{D}_1$. (b) and (c) plot the profile of $u_1$ in the $\rho\eta$-plane from different angles. (d) shows the level sets of $u_1$ in the  $\rho\eta$-plane.}  
     \label{fig:mesh_effective}
        \vspace{-0.05in}
\end{figure}

To see how well the adaptive mesh resolves the solution, we first visualize how it transforms the solution from the $rz$-plane to the $\rho\eta$-plane. Figure~\ref{fig:mesh_effective}(a) shows the function $u_1$ at $t=1.75\times 10^{-4}$ in the original $rz$-plane. This plot suggests that the singular structure should be a focusing singularity at the original $(r,z) = (0,0)$. For comparison, Figure~\ref{fig:mesh_effective}(b) and (c) plot the profile of $u_1$ at the same time in the $\rho\eta$-plane from two different angles. It is apparent that our adaptive mesh resolves the potential point-singularity structure of the solution in the $(\rho,\eta)$ coordinates.

\begin{figure}[!ht]
\centering
	\begin{subfigure}[b]{0.40\textwidth}
    \includegraphics[width=1\textwidth]{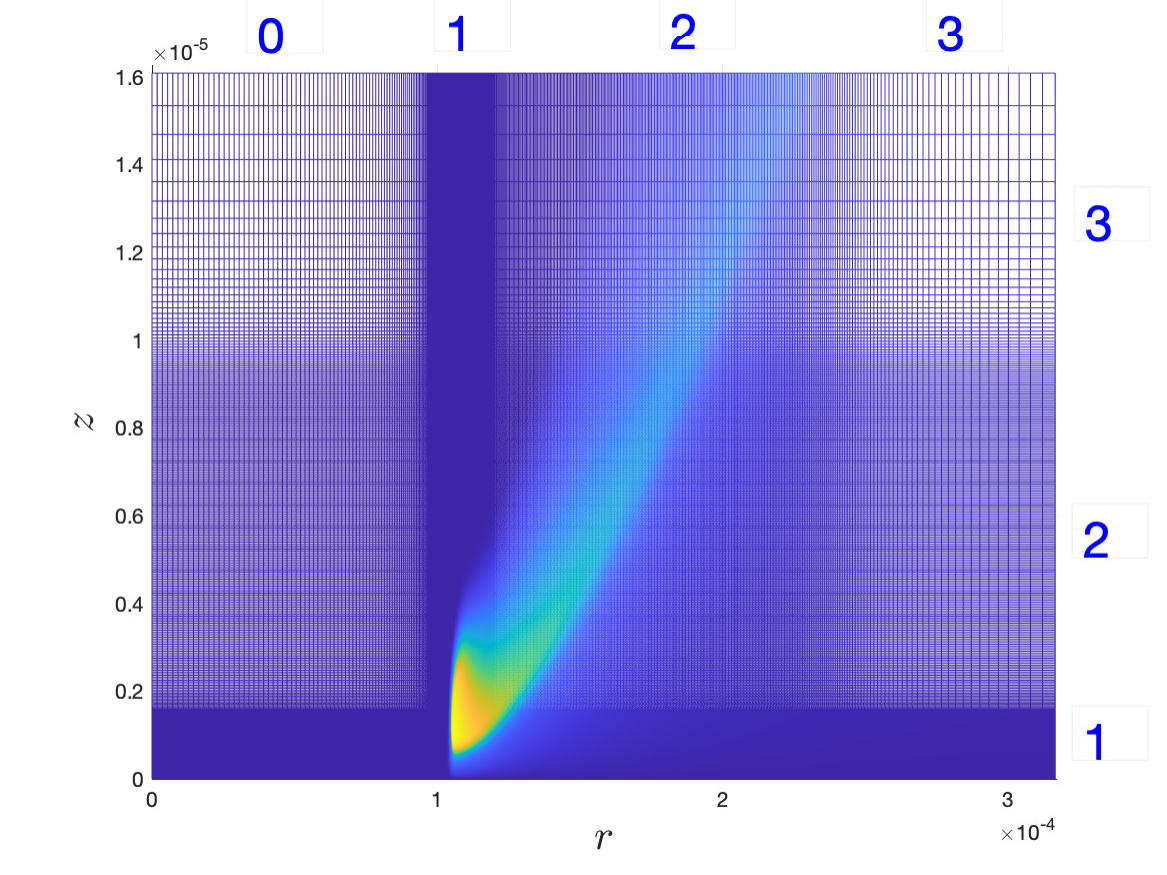}
    \caption{$u_1$ at $t=1.75\times 10^{-4}$}
    \end{subfigure}
  	\begin{subfigure}[b]{0.40\textwidth}
    \includegraphics[width=1\textwidth]{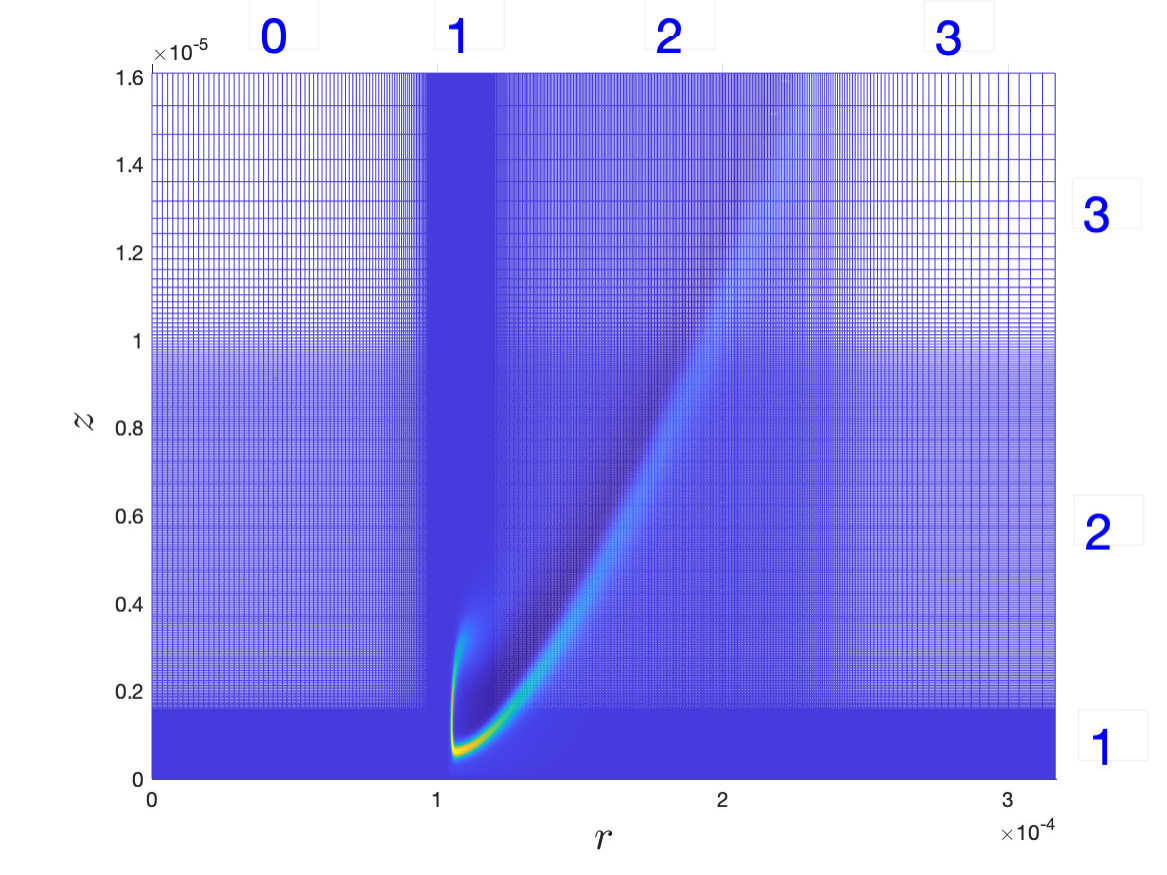}
    \caption{$\om_1$ at $t=1.75\times 10^{-4}$}
    \end{subfigure} 
    \caption[Mesh phases]{The adaptive mesh has different densities in different phases. The blue numbers are the labels of the phases.}  
     \label{fig:mesh_phases}
        \vspace{-0.05in}
\end{figure}

Figure~\ref{fig:mesh_phases} shows the top views of the profiles of $u_1,\om_1$ in a local domain at $t=1.75\times10^{-4}$. This figure demonstrates how the mesh points are distributed in different phases of the adaptive mesh. The blue numbers indicates the phase labels. By design, the adaptive meshes in phase $1$ in both directions, which have the most mesh points, capture the most singular part of the solution.

To obtain a quantitative measure of the maximum resolution power achieved by the adaptive mesh, we define the mesh effectiveness functions $ME_\rho(v),ME_\eta$ with respect to some solution variable $v$ as
\[ME_\rho(v) = \frac{h_\rho v_\rho}{\|v\|_{L^\infty}} = \frac{h_\rho r_\rho v_r}{\|v\|_{L^\infty}},\quad ME_\eta(v) = \frac{h_\eta v_\eta}{\|v\|_{L^\infty}} = \frac{h_\eta r_\eta v_z}{\|v\|_{L^\infty} },\]
and the corresponding mesh effectiveness measures (MEMs) as
\[ME_{\rho,\infty}(v) = \|ME_\rho(v)\|_{L^\infty},\quad ME_{\eta,\infty}(v) = \|ME_\eta(v)\|_{L^\infty}.\]
The MEMs quantify the largest relative growth of a function $v$ in one single mesh cell. The smaller the MEMs are, the better the adaptive mesh resolves the function $v$. Therefore, we can use these quantities to measure the effectiveness of our adaptive mesh.

\begin{figure}[!ht]
\centering
    \includegraphics[width=0.40\textwidth]{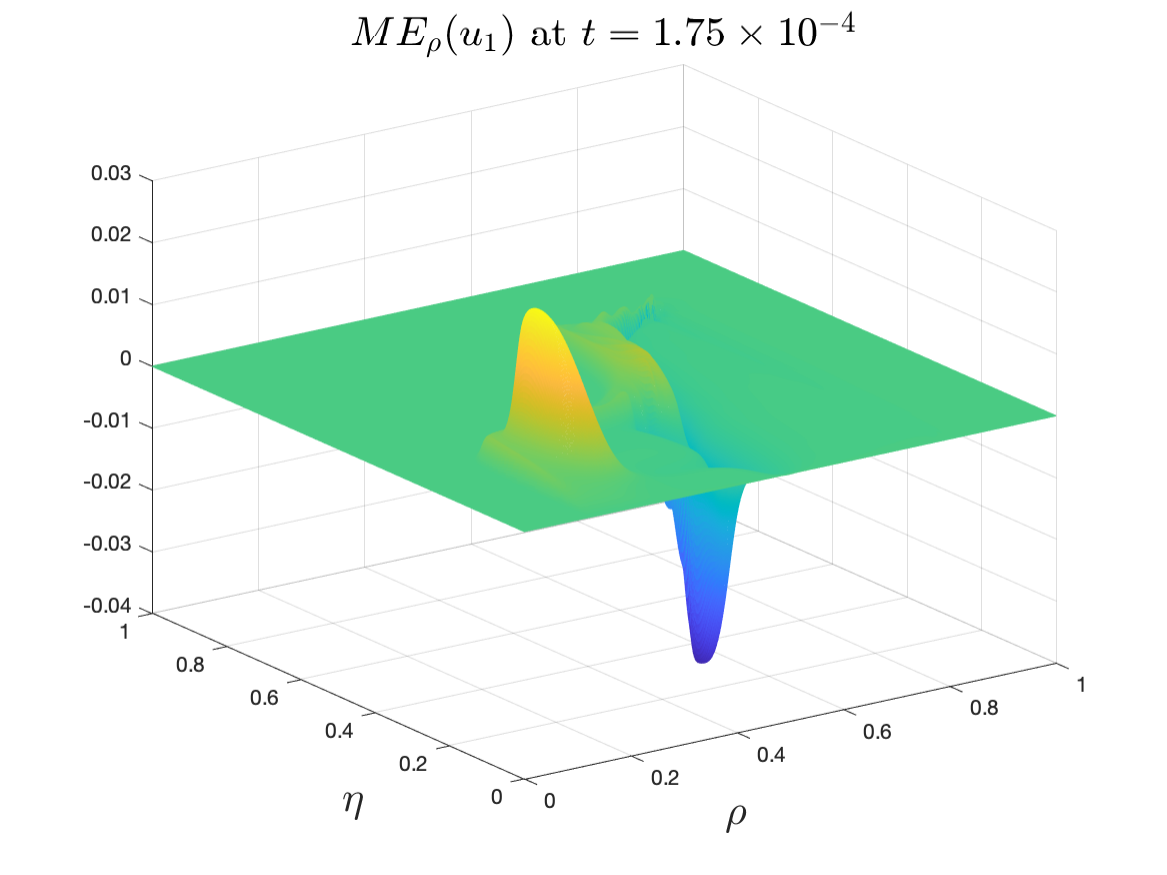}
    \includegraphics[width=0.40\textwidth]{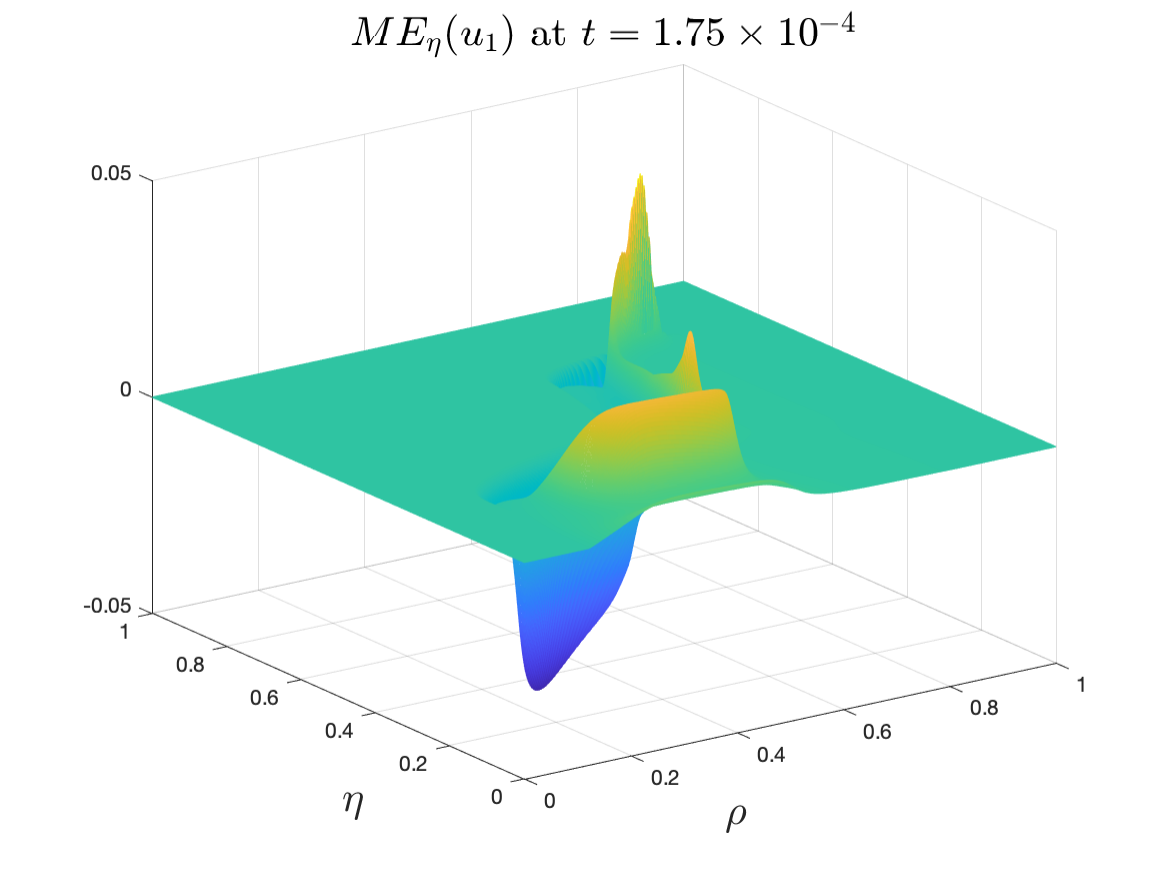} 
    \includegraphics[width=0.40\textwidth]{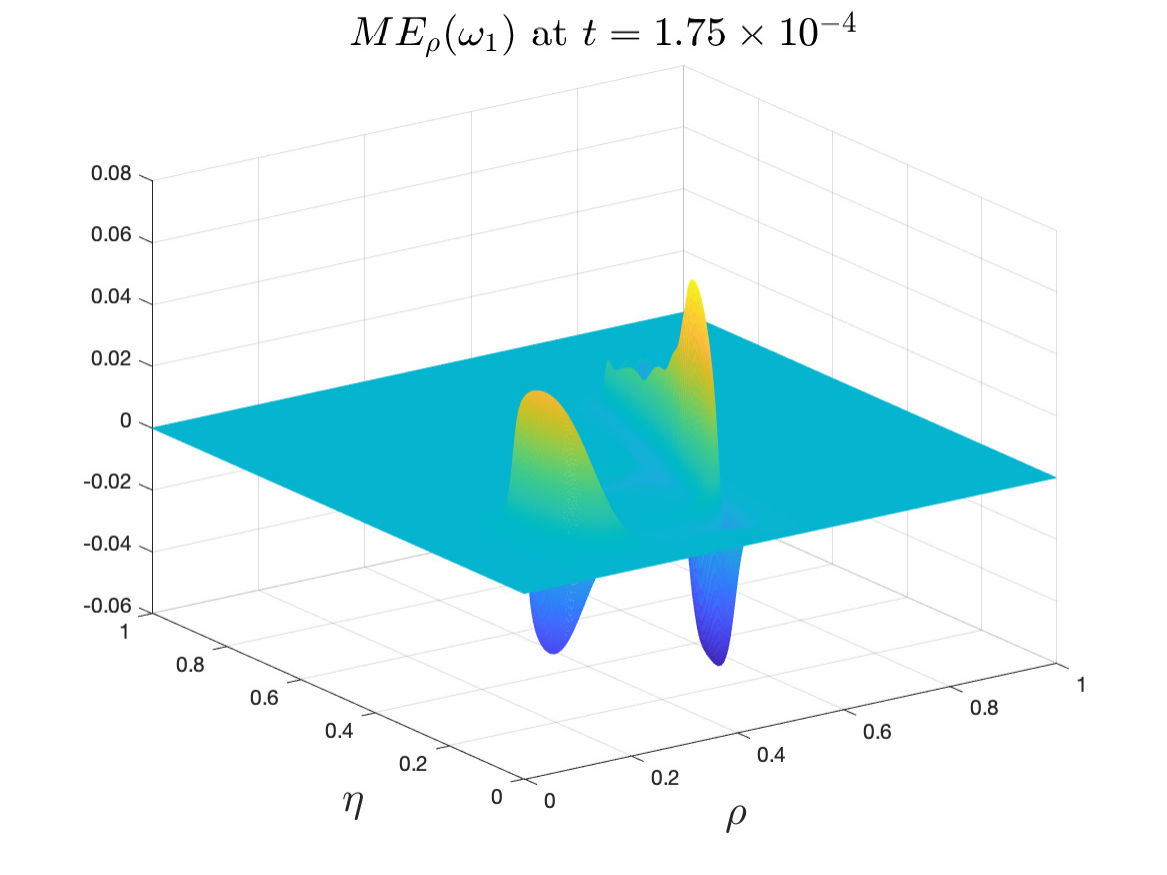}
    \includegraphics[width=0.40\textwidth]{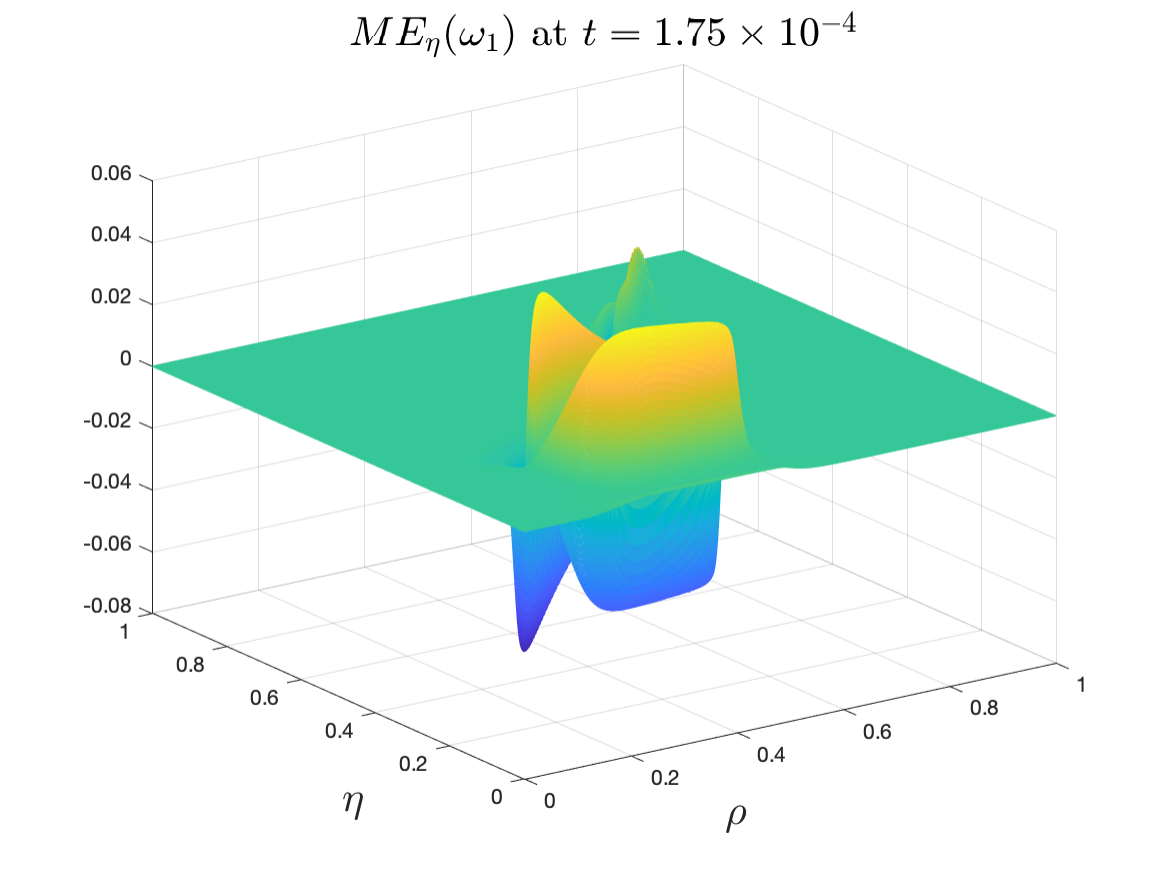}
    \caption[Mesh effectiveness functions]{First row: the mesh effectiveness functions of $u_1$ at $t=1.75\times10^{-4}$ with mesh dimension $(n,m) = (1024,512)$. Second row: the mesh effectiveness functions of $\om_1$ in the same setting.}  
     \label{fig:MEF}
        \vspace{-0.05in}
\end{figure}

Figure~\ref{fig:MEF} plots the mesh effectiveness functions of $u_1,\om_1$ at time $1.75\times 10^{-4}$ on the mesh of size $(n,m) = (1024,512)$. We can see that these functions are all uniformly bounded in absolute value by a small number (e.g., $0.1$). 

Table~\ref{tab:MEM_mesh} reports the MEMs of $u_1,\om_1$ at $t=1.75\times 10^{-4}$ on meshes of different sizes. The MEMs decrease as the grid sizes $h_\rho,h_\eta$ decrease, which is expected since the MEMs are proportional to $h_\rho,h_\eta$. Table~\ref{tab:MEM_time} reports the MEMs of $u_1,\om_1$ at different times with the same mesh size $(n,m) = (1024,512)$. We can see that the MEMs remain relatively small throughout this time interval, though with an increasing trend (especially $ME_{\rho,\infty}$) in time. From the above study, we can confirm that our adaptive mesh strategy is effective in resolving the potential singularity of the computed solution; it is well adaptive to the solution over the entire computational domain $\mathcal{D}_1$, especially in the local region near the sharp front where the solution is most singular.

We remark that though our mesh strategy can resolve the solution well before $t=1.76\times 10^{-4}$ in Case $1$, it may not work well all the way to the potential singularity time $T$. The limitation is mostly due to the unbounded growth of the ratio between the two separate scales in the solution as $t$ approaches $T$. To continue the computation in Case $1$ beyond $t=1.76\times 10^{-4}$, one may need to use a finer resolution or a more sophisticated adaptive mesh strategy. 

\begin{table}[!ht]
\centering
\footnotesize
\renewcommand{\arraystretch}{1.5}
    \begin{tabular}{|c|c|c|c|c|}
    \hline
    \multirow{2}{*}{Mesh size} & \multicolumn{4}{c|}{MEMs on mesh at $t=1.75\times 10^{-4}$} \\ \cline{2-5} 
    						   & $ME_{\rho,\infty}(u_1)$ & $ME_{\eta,\infty}(u_1)$ & $ME_{\rho,\infty}(\om_1)$ & $ME_{\eta,\infty}(\om_1)$  \\ \hline 
    $512\times 256$ & $0.086$ & $0.096$ & $0.136$ & $0.116$ \\ \hline 
    $1024\times 512$ & $0.038$ & $0.050$ & $0.067$ & $0.066$  \\ \hline 
    $1536\times 768$ & $0.026$ & $0.032$ & $0.046$ & $0.041$ \\ \hline 
    $2048\times 1024$ & $0.018$ & $0.026$ & $0.033$ & $0.030$  \\ \hline 
    \end{tabular}
    \caption{\small MEMs of $u_1,\om_1$ at $t=1.75\times 10^{-4}$ on the meshes of different sizes.}
    \label{tab:MEM_mesh}
    \vspace{-0.3in}
\end{table}

\begin{table}[!ht]
\centering
\footnotesize
\renewcommand{\arraystretch}{1.5}
    \begin{tabular}{|c|c|c|c|c|}
    \hline
    \multirow{2}{*}{Time} & \multicolumn{4}{c|}{MEMs on mesh $(n,m) = (1024,512)$} \\ \cline{2-5} 
    						   & $ME_{\rho,\infty}(u_1)$ & $ME_{\eta,\infty}(u_1)$ & $ME_{\rho,\infty}(\om_1)$ & $ME_{\eta,\infty}(\om_1)$  \\ \hline 
    $1.7\times 10^{-4}$ & $0.022$ & $0.042$ & $0.049$ & $0.070$ \\ \hline 
    $1.73\times 10^{-4}$ & $0.028$ & $0.047$ & $0.053$ & $0.068$  \\ \hline 
    $1.74\times 10^{-4}$ & $0.031$ & $0.041$ & $0.056$ & $0.059$ \\ \hline 
    $1.75\times 10^{-4}$ & $0.038$ & $0.050$ & $0.067$ & $0.066$  \\ \hline 
    $1.76\times 10^{-4}$ & $0.065$ & $0.059$ & $0.100$ & $0.069$  \\ \hline 
    \end{tabular}
    \caption{\small MEMs of $u_1,\om_1$ at different times on the mesh of size $(n,m) = (1024,512)$.}
    \label{tab:MEM_time}
    \vspace{-0.3in}
\end{table}

\subsection{Resolution study}\label{sec:resolution_study}
In this subsection, we perform resolution study on the numerical solutions of the initial-boundary value problem \eqref{eq:axisymmetric_NSE_1}--\eqref{eq:initial_data} at various time instants $t$ and in different cases of computation. We will estimate the relative error of some solution variable $f_p$ computed on the $256p\times 128p$ mesh by comparing it to a reference variable $\hat{f}$ that is computed on a finer mesh at the same time instant. If $f_p$ is a number, the relative error in absolute value is computed as 
$e_p = |f_p-\hat{f}|/|\hat{f}|$.
If $f_p$ is a spatial function, the reference variable $\hat{f}$ is first interpolated to the mesh on which $f$ is computed. Then the sup-norm relative error is computed as 
\begin{align*}
e_p &= \frac{\|f_p-\hat{f}\|_{\infty}}{\|\hat{f}\|_{\infty}}\quad \text{if $f$ is a scalar function,}\\
\text{and}\quad e_p &= \frac{\left\|\big|(f_p^\theta-\hat{f}_p^\theta,f_p^r-\hat{f}_p^r,f_p^z-\hat{f}_p^z)\big|\right\|_{\infty}}{\left\|\big|(\hat{f}_p^\theta,\hat{f}_p^r,\hat{f}_p^z)\big|\right\|_{\infty}}\quad \text{if $f$ is a vector function.}
\end{align*}
For all cases, the reference solution $\hat{f}$ is chosen to be the one computed at the same time instant on the finer mesh of size $256(p+1)\times 128(p+1)$; that is, $\hat{f} = f_{p+1}$.
The numerical order of the error is computed as 
\[\beta_p = \log_{\frac{p}{p-1}}\left(\frac{e_{p-1}}{e_p}\right) - 1.\]
Ideally, for a numerical method of order $\beta$, the error of a solution variable $f_p$ compared to the ground truth $f^*$ is proportional to $p^{-\beta}$. Suppose that $f_p$ converges to $f^*$ in a monotone fashion, then we should have $e_p\propto p^{-\beta} - (p+1)^{-\beta}$. Substituting this into the formula of $\beta_p$ yields
\[\beta_p = \log_{\frac{p}{p-1}}\left(\frac{(p-1)^{-\beta} - p^{-\beta}}{p^{-\beta} - (p+1)^{-\beta}}\right)-1.\]
One can then show that $\beta_p$ is monotone increasing in $p$ and will converge to the true order $\beta$ as $p\rightarrow \infty$. In particular, for our $2$nd-order method, $\beta_p$ should approach $2$ as $p$ increases.

\subsubsection{Case $1$} This case is the most fundamental one among all of our computations. We thus perform a more thorough resolution study for the numerical solutions corresponding to this case.  

We first study the sup-norm error of the solution, which is the most straightforward indication of the accuracy of our numerical method. Tables~\ref{tab:sup-norm_error_1-1}--\ref{tab:sup-norm_error_1-4} report the sup-norm relative errors and numerical orders of different solution variables at times $t = 1.65\times 10^{-4}$ and $t = 1.7\times 10^{-4}$, respectively. The results confirm that our method in Case $1$ is at least $2$nd-order accurate. We remark that the error in the solution mainly arises from the interpolation error near the sharp front, where the gradient of the solution is largest and becomes larger and larger in time.

\begin{table}[!ht]
\centering
\footnotesize
\renewcommand{\arraystretch}{1.5}
    \begin{tabular}{|c|c|c|c|c|c|c|}
    \hline
    \multirow{2}{*}{Mesh size} & \multicolumn{6}{c|}{Sup-norm relative error at $t=1.65\times10^{-4}$ in Case $1$} \\ \cline{2-7} 
    & Error of $u_1$ & Order & Error of $\omega_1$ & Order & Error of $\psi_1$ & Order \\ \hline 
    $512\times256$ & $1.2949\times10^{-1}$ & -- & $2.7872\times10^{-1}$ & -- & $2.4290\times10^{-2}$ & -- \\ \hline 
    $768\times384$ & $4.4825\times10^{-2}$ & $1.62$ & $9.8239\times10^{-2}$ & $1.57$ & $8.2729\times10^{-3}$ & $1.66$ \\ \hline 
    $1024\times512$ & $2.0467\times10^{-2}$ & $1.63$ & $4.4789\times10^{-2}$ & $1.73$ & $3.8283\times10^{-3}$ & $1.68$ \\ \hline 
    $1280\times640$ & $1.1264\times10^{-2}$ & $1.68$ & $2.5061\times10^{-2}$ & $1.60$ & $1.9905\times10^{-3}$ & $1.93$ \\ \hline 
    $1536\times768$ & $7.0304\times10^{-3}$ & $1.59$ & $1.5410\times10^{-2}$ & $1.67$ & $1.3228\times10^{-3}$ & $1.24$ \\ \hline 
    $1792\times896$ & $4.3618\times10^{-3}$ & $2.10$ & $9.5333\times10^{-3}$ & $2.12$ & $8.2235\times10^{-4}$ & $2.08$ \\
    \hline
    \end{tabular}
    \caption{\small Sup-norm relative errors and numerical orders of $u_1,\om_1,\psi_1$ at $t = 1.65\times 10^{-4}$ in Case $1$.}
    \label{tab:sup-norm_error_1-1}
    \vspace{-0.2in}
\end{table}

\begin{table}[!ht]
\centering
\footnotesize
\renewcommand{\arraystretch}{1.5}
    \begin{tabular}{|c|c|c|c|c|c|c|}
    \hline
    \multirow{2}{*}{Mesh size} & \multicolumn{6}{c|}{Sup-norm relative error at $t=1.65\times10^{-4}$ in Case $1$} \\ \cline{2-7} 
    & Error of $u^r$ & Order & Error of $u^z$ & Order & Error of $\vom$ & Order \\ \hline 
    $512\times256$ & $2.3116\times10^{-1}$ & -- & $6.5401\times10^{-2}$ & -- & $2.5887\times10^{-1}$ & -- \\ \hline 
    $768\times384$ & $8.0520\times10^{-2}$ & $1.60$ & $2.2632\times10^{-2}$ & $1.62$ & $9.1315\times10^{-2}$ & $1.57$ \\ \hline 
    $1024\times512$ & $3.6975\times10^{-2}$ & $1.71$ & $1.0477\times10^{-2}$ & $1.68$ & $4.1456\times10^{-2}$ & $1.74$ \\ \hline 
    $1280\times640$ & $2.0084\times10^{-2}$ & $1.74$ & $5.7374\times10^{-3}$ & $1.70$ & $2.3254\times10^{-2}$ & $1.59$ \\ \hline 
    $1536\times768$ & $1.2862\times10^{-2}$ & $1.44$ & $3.6477\times10^{-3}$ & $1.48$ & $1.4134\times10^{-3}$ & $1.73$ \\ \hline 
    $1792\times896$ & $7.9410\times10^{-3}$ & $2.13$ & $2.2479\times10^{-3}$ & $2.14$ & $8.7579\times10^{-3}$ & $2.11$ \\
    \hline
    \end{tabular}
    \caption{\small Sup-norm relative errors and numerical orders of $u^r,u^z,\vom$ at $t = 1.65\times 10^{-4}$ in Case $1$.}
    \label{tab:sup-norm_error_1-2}
    \vspace{-0.2in}
\end{table}

\begin{table}[!ht]
\centering
\footnotesize
\renewcommand{\arraystretch}{1.5}
    \begin{tabular}{|c|c|c|c|c|c|c|}
    \hline
    \multirow{2}{*}{Mesh size} & \multicolumn{6}{c|}{Sup-norm relative error at $t=1.7\times10^{-4}$ in Case $1$} \\ \cline{2-7} 
    & Error of $u_1$ & Order & Error of $\omega_1$ & Order & Error of $\psi_1$ & Order \\ \hline 
    $512\times256$ & $3.1543\times10^{-1}$ & -- & $5.7498\times10^{-1}$ & -- & $6.3176\times10^{-2}$ & -- \\ \hline 
    $768\times384$ & $1.1080\times10^{-1}$ & $1.58$ & $2.1421\times10^{-1}$ & $1.44$ & $2.1410\times10^{-2}$ & $1.67$ \\ \hline 
    $1024\times512$ & $5.3980\times10^{-2}$ & $1.50$ & $1.0487\times10^{-1}$ & $1.48$ & $1.0267\times10^{-2}$ & $1.55$ \\ \hline 
    $1280\times640$ & $2.8154\times10^{-2}$ & $1.92$ & $5.4662\times10^{-2}$ & $1.92$ & $5.3087\times10^{-3}$ & $1.96$ \\ \hline 
    $1536\times768$ & $1.8674\times10^{-2}$ & $1.25$ & $3.6419\times10^{-2}$ & $1.23$ & $3.5579\times10^{-3}$ & $1.19$ \\ \hline 
    $1792\times896$ & $1.1740\times10^{-2}$ & $2.01$ & $2.2894\times10^{-2}$ & $2.01$ & $2.2323\times10^{-3}$ & $2.02$ \\
    \hline
    \end{tabular}
    \caption{\small Sup-norm relative errors and numerical orders of $u_1,\om_1,\psi_1$ at $t = 1.7\times 10^{-4}$ in Case $1$.}
    \label{tab:sup-norm_error_1-3}
    \vspace{-0.2in}
\end{table}

\begin{table}[!ht]
\centering
\footnotesize
\renewcommand{\arraystretch}{1.5}
    \begin{tabular}{|c|c|c|c|c|c|c|}
    \hline
    \multirow{2}{*}{Mesh size} & \multicolumn{6}{c|}{Sup-norm relative error at $t=1.7\times10^{-4}$ in Case $1$} \\ \cline{2-7} 
    & Error of $u^r$ & Order & Error of $u^z$ & Order & Error of $\vom$ & Order \\ \hline 
    $512\times256$ & $5.1092\times10^{-1}$ & -- & $1.3200\times10^{-1}$ & -- & $5.0630\times10^{-1}$ & -- \\ \hline 
    $768\times384$ & $1.8141\times10^{-1}$ & $1.55$ & $4.6313\times10^{-2}$ & $1.58$ & $1.8885\times10^{-1}$ & $1.43$ \\ \hline 
    $1024\times512$ & $8.8443\times10^{-2}$ & $1.50$ & $2.2743\times10^{-2}$ & $1.47$ & $9.2927\times10^{-2}$ & $1.47$ \\ \hline 
    $1280\times640$ & $4.5596\times10^{-2}$ & $1.97$ & $1.2038\times10^{-2}$ & $1.85$ & $4.8958\times10^{-2}$ & $1.87$ \\ \hline 
    $1536\times768$ & $3.0860\times10^{-2}$ & $1.14$ & $7.8363\times10^{-3}$ & $1.35$ & $3.2255\times10^{-2}$ & $1.29$ \\ \hline 
    $1792\times896$ & $1.9350\times10^{-2}$ & $2.03$ & $4.9410\times10^{-3}$ & $1.99$ & $2.0271\times10^{-2}$ & $2.01$ \\
    \hline
    \end{tabular}
    \caption{\small Sup-norm relative errors and numerical orders of $u^r,u^z,\vom$ at $t = 1.7\times 10^{-4}$ in Case $1$.}
    \label{tab:sup-norm_error_1-4}
    \vspace{-0.2in}
\end{table}

\begin{figure}[!ht]
\centering
    \includegraphics[width=0.40\textwidth]{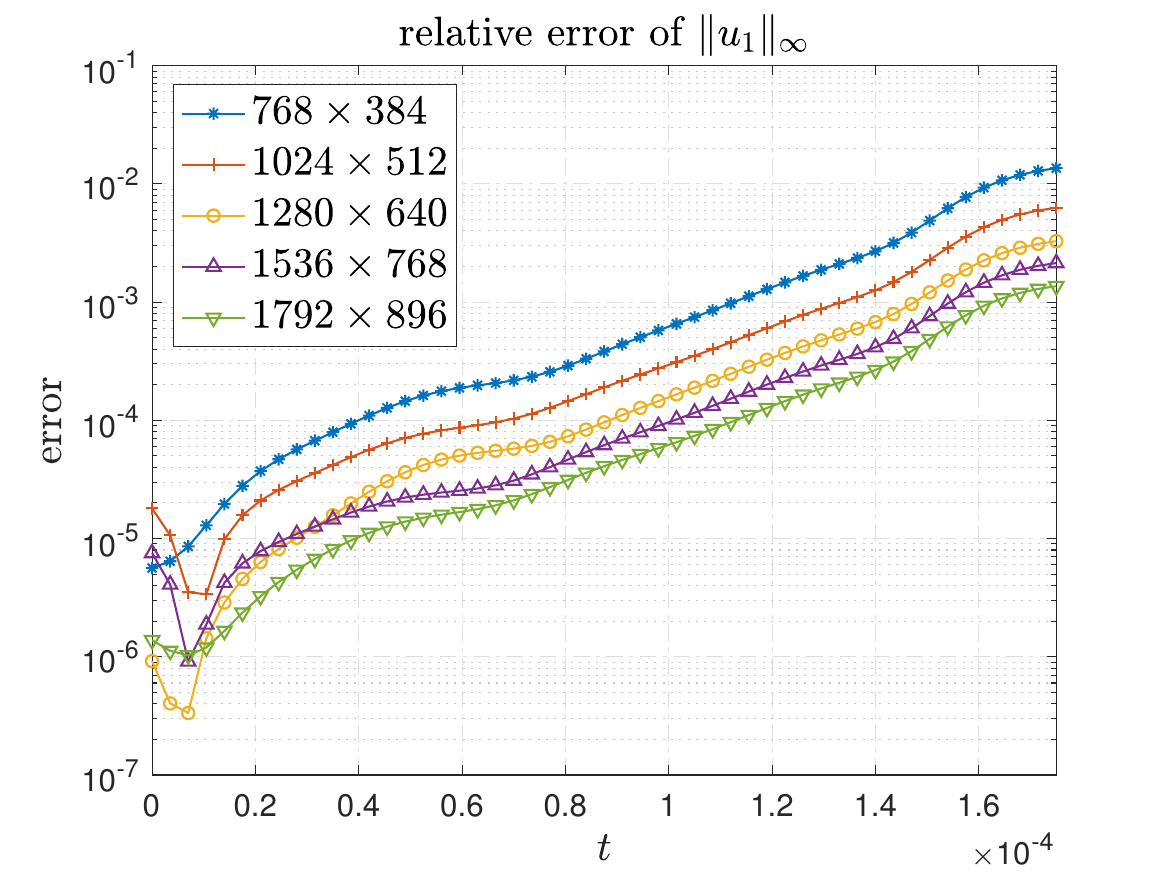}
    \includegraphics[width=0.40\textwidth]{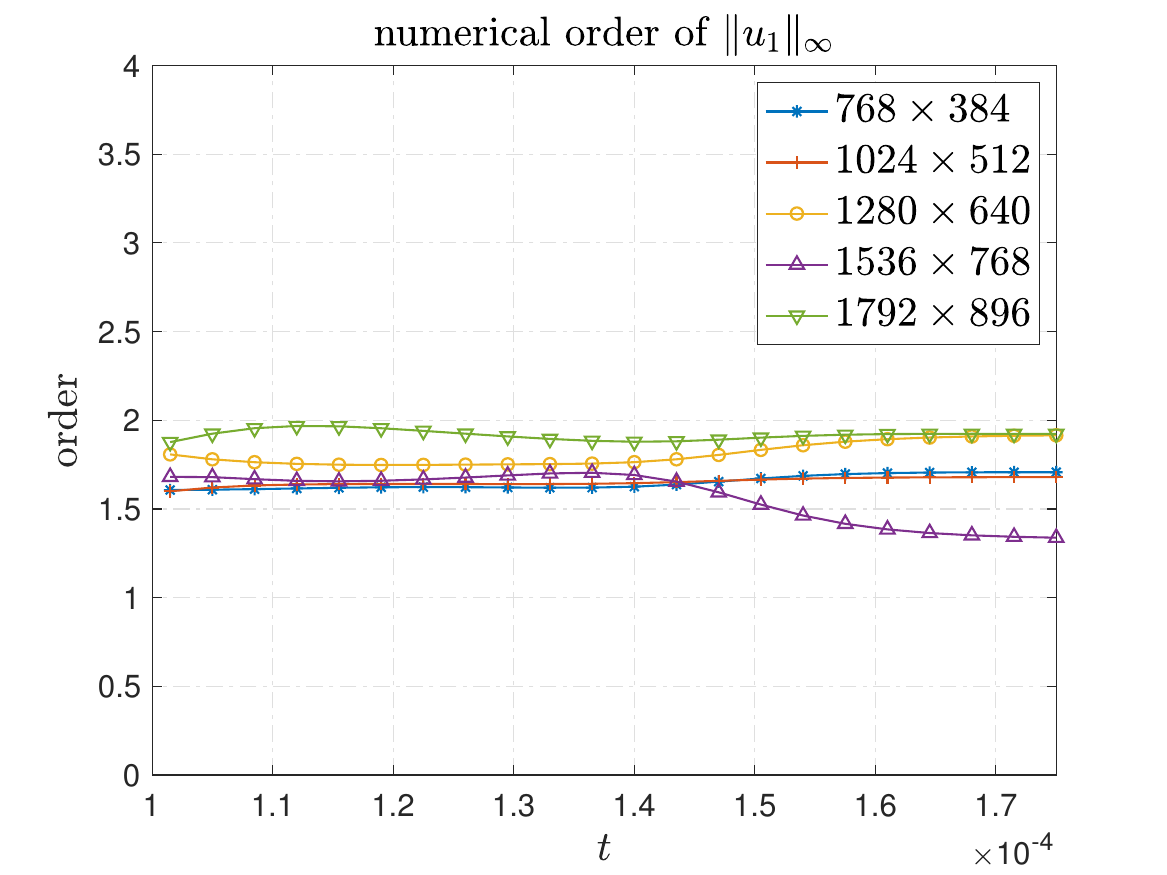} 
    \includegraphics[width=0.40\textwidth]{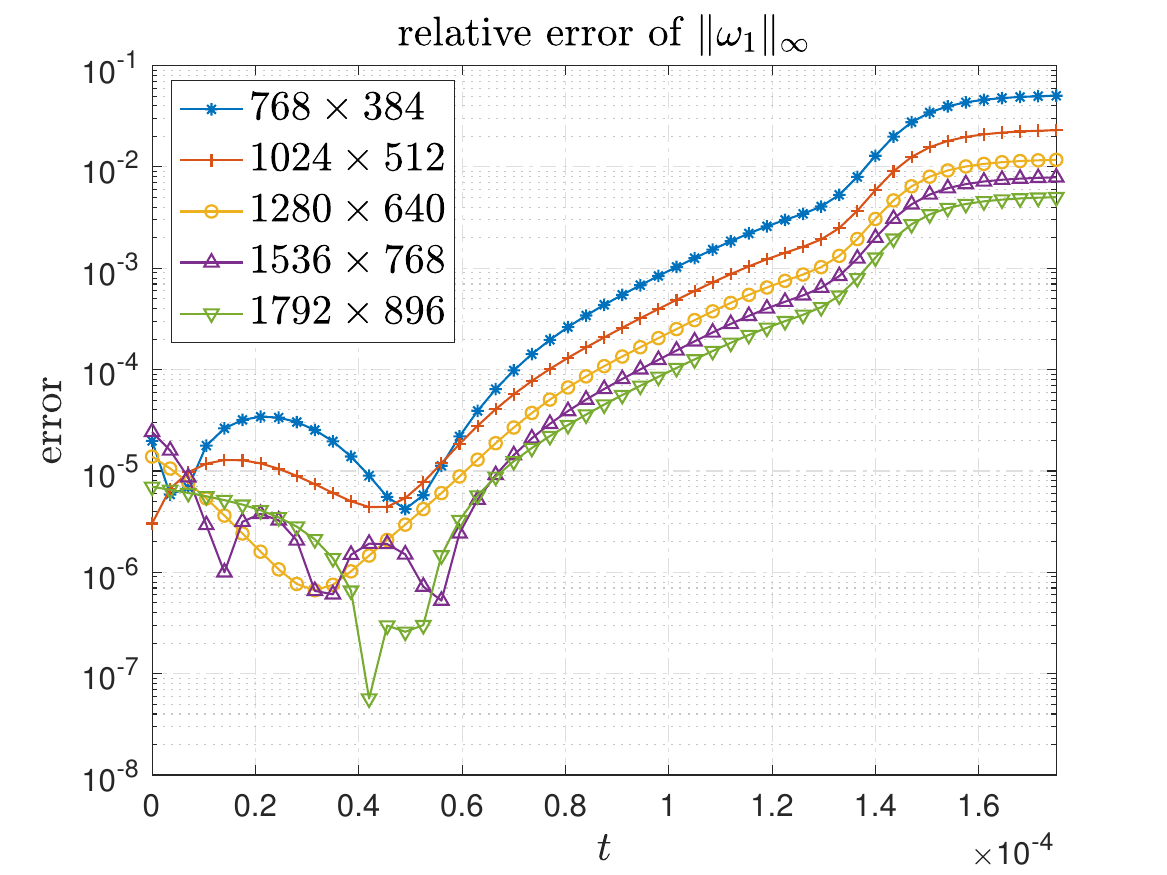}
    \includegraphics[width=0.40\textwidth]{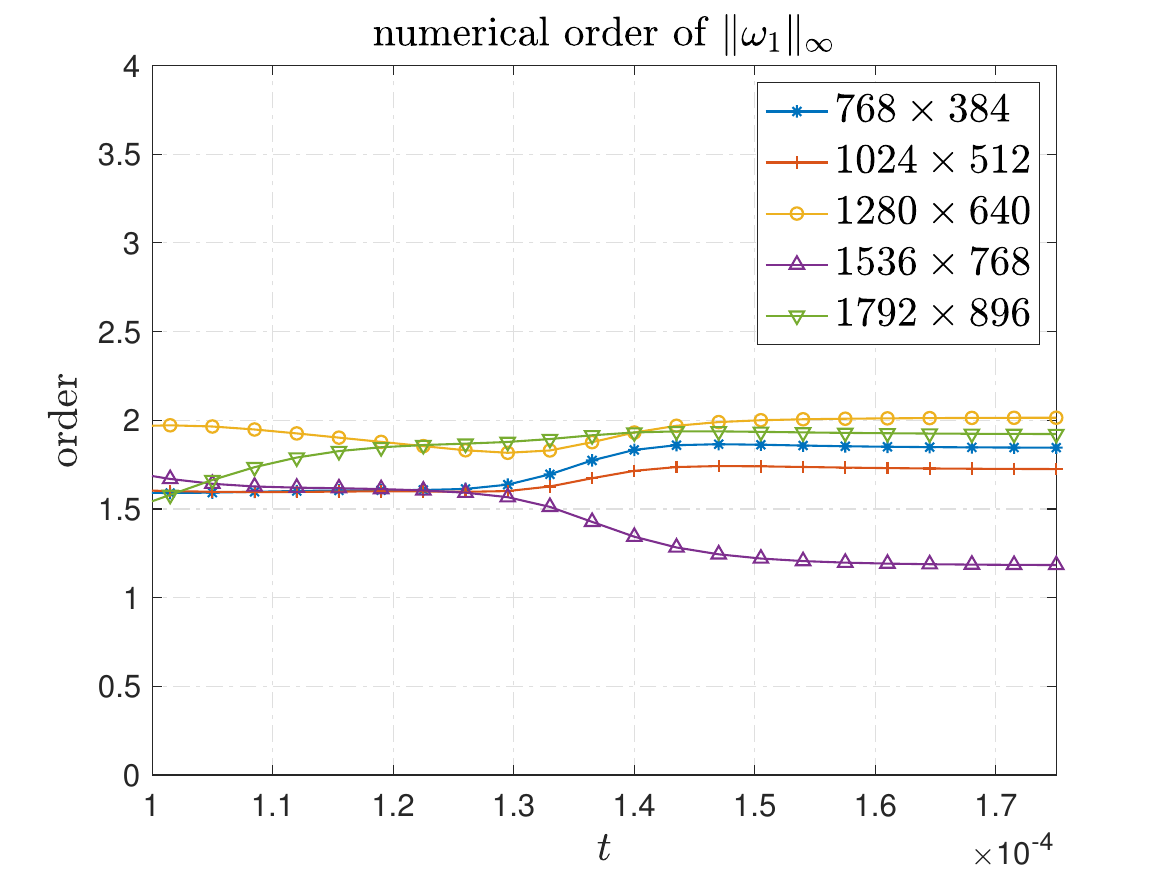}
    \caption[Relative error I]{First row: relative error and numerical order of $\|u_1(t)\|_{L^\infty}$. Second row: relative error and numerical order of $\|\om_1(t)\|_{L^\infty}$. The The last time instant shown in the figure is $t= 1.76\times10^{-4}$.}  
    \label{fig:relative_error_1}
\end{figure}

\begin{figure}[!ht]
\centering  
    \includegraphics[width=0.40\textwidth]{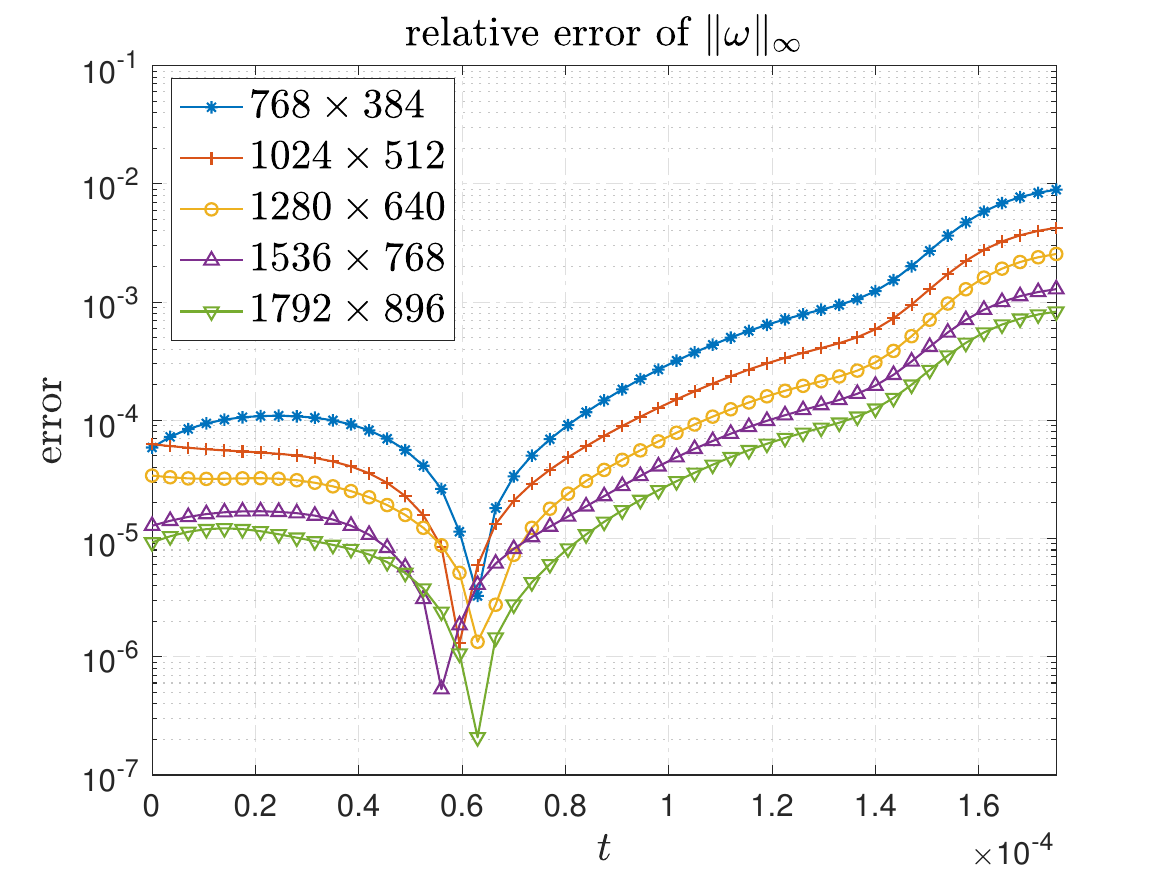}
    \includegraphics[width=0.40\textwidth]{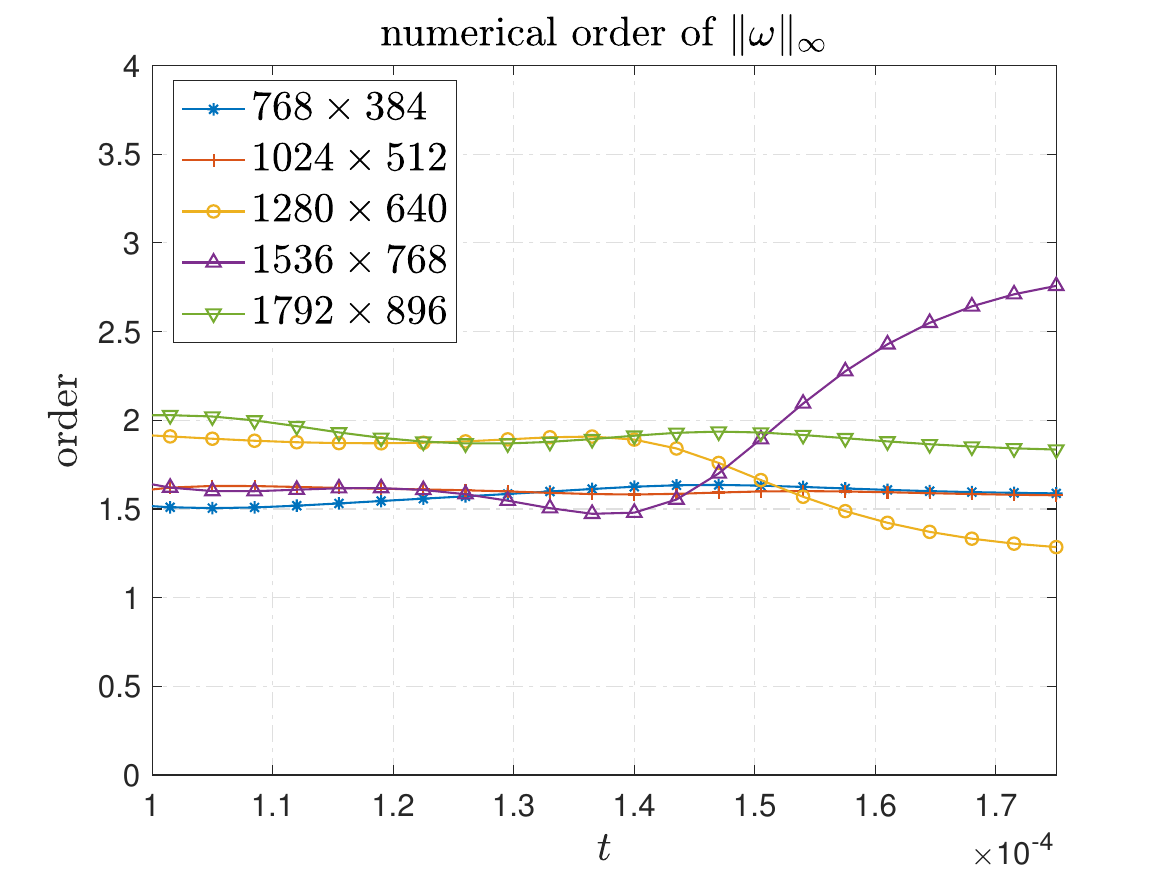} 
    \includegraphics[width=0.40\textwidth]{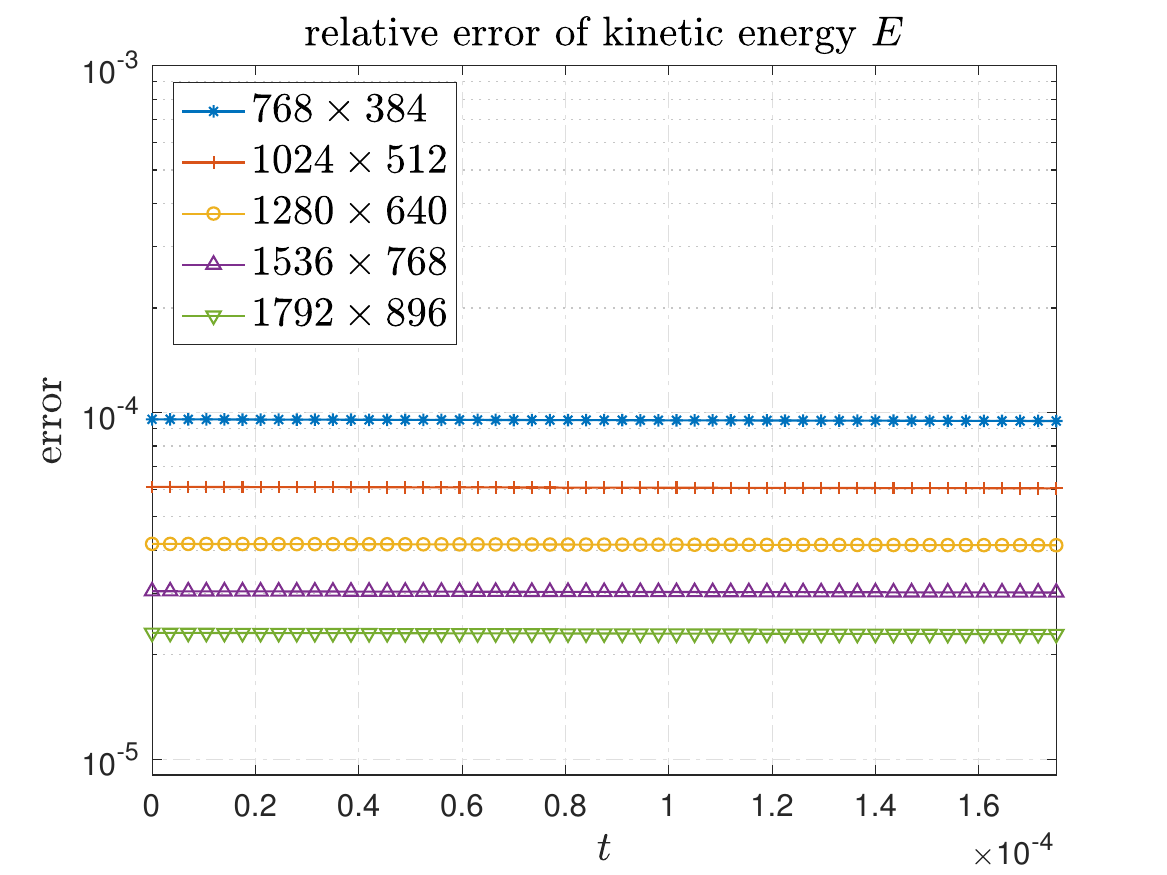}
    \includegraphics[width=0.40\textwidth]{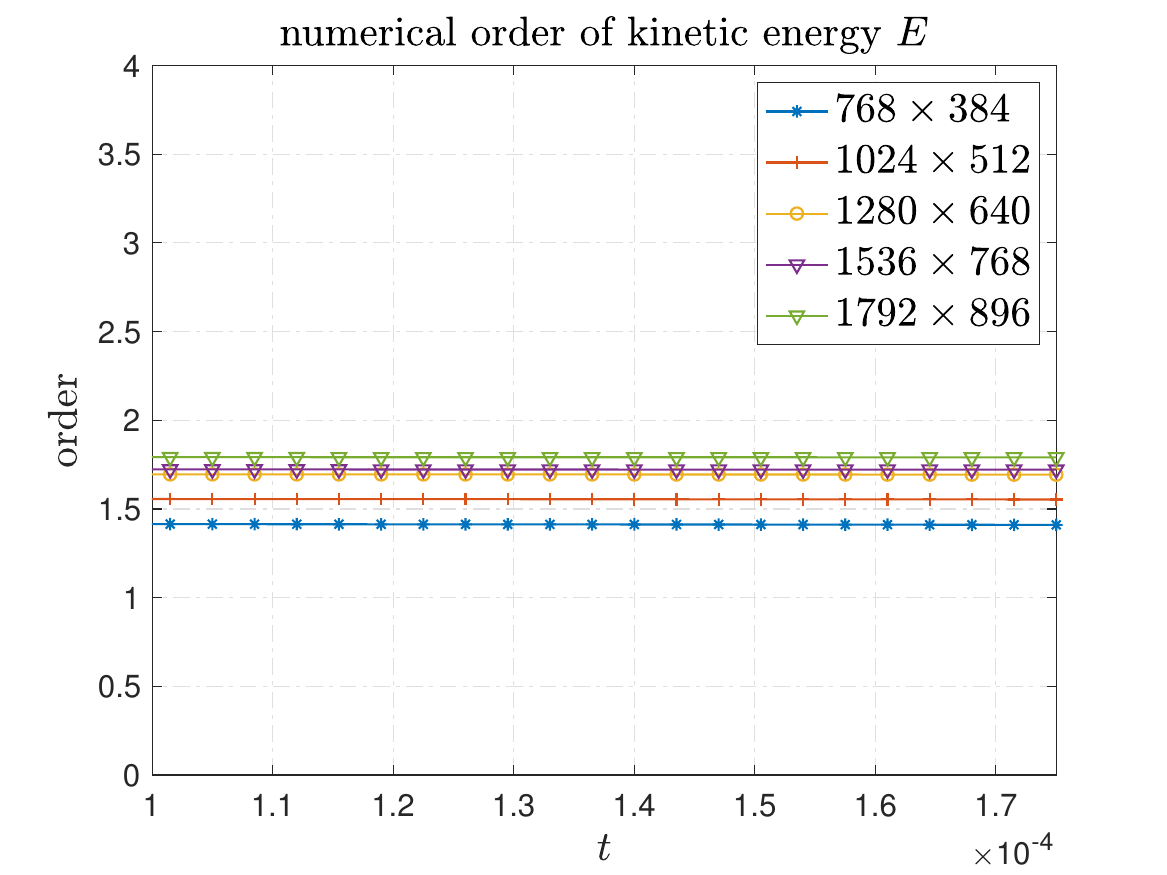}
    \caption[Relative error II]{First row: relative error and numerical order of $\|\vom(t)\|_{L^\infty}$. Second row: relative error and numerical order of $E(t)$. The last time instant shown in the figure is $t= 1.76\times10^{-4}$.}  
    \label{fig:relative_error_2}
       \vspace{-0.05in}
\end{figure}

We can also study the convergence of some variables as functions of time. In particular, we report the convergence of the quantities $\|u_1\|_{L^\infty}$, $\|\om_1\|_{L^\infty}$, $\|\vom\|_{L^\infty}$, and the kinetic energy $E$. Here the kinetic energy $E$ is given by 
\[E := \frac{1}{2}\int_{\mathcal{D}_1}|\vu|^2 \idiff x = \frac{1}{2}\int_0^1\int_0^{1/2}\left(|u^r|^2+ |u^\theta|^2 + |u^z|^2\right)r\idiff r \idiff z. \]
Since the viscosity term with variable coefficients in \eqref{eq:NSE_vc} is given in a conservative form, the kinetic energy is a non-increasing function of time: $E(t_1)\leq E(t_2)$ for $t_2\geq t_1\geq0$. Figures~\ref{fig:relative_error_1} and \ref{fig:relative_error_2} plot the relative errors and numerical orders of these quantities as functions of time. The results further confirm that our method is $2$nd-order in $h_\rho,h_\eta$. 

\begin{remark}
There are two main sources of errors in our computation: one is the global error of the numerical scheme, and the other is the interpolation error from one mesh to another when we change our adaptive mesh at certain time instants. Since the solution is quite smooth in the early stage of the computation, the scheme error has not accumulated much, and the total error may be dominated by the interpolation error. The change of mesh may happen more frequently for the computation with a finer mesh, and thus the total error can even be smaller on a coarser mesh, as we can see in the first row of Figure~\ref{fig:relative_error_1}. We can only see the expected order of accuracy when the discretization error accumulates to some level such that it dominates the total error.
\end{remark}

On the other hand, we also observe an increasing trend in the relative errors of $\|u_1\|_{L^\infty}$, $\|\om_1\|_{L^\infty}$, and $\|\vom\|_{L^\infty}$, which implies that our numerical method with a fixed mesh size will not work for all time up to the anticipated singularity. As we mentioned in Section~\ref{sec:mesh_effectiveness}, our adaptive mesh strategy may lose its power to resolve the solution as the two scales in the solution becomes more and more separated. Indeed, the sharp front in the $r$ direction becomes thinner and thinner as $t$ approaches the potential singularity time, which makes it more and more difficult to construct an adaptive mesh with a fixed number of grid points that provides a small approximation error in the entire domain. Therefore, to obtain a well-resolved solution sufficiently close to the potential singularity time, one must use an extremely large number of grid points, which is, unfortunately, beyond the capacity of our current computational resources.

\subsubsection{Case $2$} Recall that in the Case $2$ computation, the only change we make is to replace the variable viscosity coefficients \eqref{eq:viscosity_coefficient} by a constant viscosity coefficient $\mu$, so we actually solve the original Navier--Stokes equations. Since the numerical methods in Case $1$ and Case $2$ are the same, we expect to see similar convergence behaviors of the solutions. Nevertheless, we still perform a resolution study in Case $2$ to confirm the $2$nd-order accuracy of our method, and we only present the results from the case $\nu^r=\nu^z = 10^{-5}$. We remark that the solution in Case $2$ evolves in a completely different way; in particular, it does not lead to a finite-time blowup. We can thus compute the solution to a much later time. The sup-norm relative errors and numerical orders of $u_1,\om_1,\vom$ at $t = 2\times 10^{-4}$ and $t = 2.5\times 10^{-4}$ are reported in Tables~\ref{tab:sup-norm_error_2-1} and \ref{tab:sup-norm_error_2-2}, respectively. 

\begin{table}[!ht]
\centering
\footnotesize
\renewcommand{\arraystretch}{1.5}
    \begin{tabular}{|c|c|c|c|c|c|c|}
    \hline
    \multirow{2}{*}{Mesh size} & \multicolumn{6}{c|}{Sup-norm relative error at $t=2\times10^{-4}$ in Case $2$} \\ \cline{2-7} 
    & Error of $u_1$ & Order & Error of $\omega_1$ & Order & Error of $\vom$ & Order \\ \hline 
    $512\times256$ & $1.6283\times10^{-3}$ & -- & $2.0004\times10^{-3}$ & -- & $2.4642\times10^{-3}$ & -- \\ \hline 
    $768\times384$ & $5.7153\times10^{-4}$ & $1.58$ & $7.0144\times10^{-4}$ & $1.58$ & $8.6683\times10^{-4}$ & $1.57$ \\ \hline 
    $1024\times512$ & $2.6353\times10^{-4}$ & $1.69$ & $3.2570\times10^{-4}$ & $1.67$ & $3.9648\times10^{-4}$ & $1.72$ \\ \hline 
    $1280\times640$ & $1.4402\times10^{-4}$ & $1.70$ & $1.7634\times10^{-4}$ & $1.75$ & $2.1780\times10^{-4}$ & $1.68$ \\ \hline 
    $1536\times768$ & $8.6198\times10^{-5}$ & $1.82$ & $1.0735\times10^{-4}$ & $1.72$ & $1.2873\times10^{-4}$ & $1.88$ \\ \hline 
    $1792\times896$ & $5.5432\times10^{-5}$ & $1.86$ & $6.8494\times10^{-5}$ & $1.91$ & $8.3230\times10^{-5}$ & $1.83$ \\
    \hline
    \end{tabular}
    \caption{\small Sup-norm relative errors and numerical orders of $u_1,\om_1,\vom$ at $t = 2\times 10^{-4}$ in Case $2$.}
    \label{tab:sup-norm_error_2-1}
    \vspace{-0.2in}
\end{table}

\begin{table}[!ht]
\centering
\footnotesize
\renewcommand{\arraystretch}{1.5}
    \begin{tabular}{|c|c|c|c|c|c|c|}
    \hline
    \multirow{2}{*}{Mesh size} & \multicolumn{6}{c|}{Sup-norm relative error at $t=2.5\times10^{-4}$ in Case $2$} \\ \cline{2-7} 
    & Error of $u_1$ & Order & Error of $\omega_1$ & Order & Error of $\vom$ & Order \\ \hline 
    $512\times256$ & $1.3726\times10^{-3}$ & -- & $1.6971\times10^{-3}$ & -- & $2.4229\times10^{-3}$ & -- \\ \hline 
    $768\times384$ & $4.7140\times10^{-4}$ & $1.63$ & $5.8979\times10^{-4}$ & $1.61$ & $8.3826\times10^{-4}$ & $1.58$ \\ \hline 
    $1024\times512$ & $2.1351\times10^{-4}$ & $1.75$ & $2.7200\times10^{-4}$ & $1.69$ & $3.7361\times10^{-4}$ & $1.85$ \\ \hline 
    $1280\times640$ & $1.2303\times10^{-4}$ & $1.47$ & $1.5144\times10^{-4}$ & $1.62$ & $2.0920\times10^{-4}$ & $1.60$ \\ \hline 
    $1536\times768$ & $7.1537\times10^{-5}$ & $1.97$ & $9.0012\times10^{-5}$ & $1.85$ & $1.2082\times10^{-4}$ & $2.01$ \\ \hline 
    $1792\times896$ & $4.6082\times10^{-5}$ & $1.85$ & $5.7634\times10^{-5}$ & $1.89$ & $7.8138\times10^{-5}$ & $1.83$ \\
    \hline
    \end{tabular}
    \caption{\small Sup-norm relative errors and numerical orders of $u_1,\om_1,\vom$ at $t = 2.5\times 10^{-4}$ in Case $2$.}
    \label{tab:sup-norm_error_2-2}
     \vspace{-0.2in}
\end{table}

\subsubsection{Case $3$} This can be seen as a special case of Case $2$ with $\nu^r =\nu^z = 0$. Therefore, it is expected that our computation in Case $3$ also enjoys $2$nd-order convergence as in Case $2$. Nevertheless, we still present a resolution study for Case $3$ at an early time instant $t=1.5\times 10^{-4}$ when the solution is still reasonably resolved. The sup-norm relative errors and numerical orders of $u_1,\om_1,\vom$ at $t = 1.5\times 10^{-4}$ are reported in Tables~\ref{tab:sup-norm_error_3}. We can see that the numerical orders are somewhat off expectation, which may be due to the more singular behavior of the Euler solution.

\begin{table}[!ht]
\centering
\footnotesize
\renewcommand{\arraystretch}{1.5}
    \begin{tabular}{|c|c|c|c|c|c|c|}
    \hline    
    \multirow{2}{*}{Mesh size} & \multicolumn{6}{c|}{Sup-norm relative error at $t=1.5\times10^{-4}$ in Case $3$} \\ \cline{2-7} 
    & Error of $u_1$ & Order & Error of $\omega_1$ & Order & Error of $\vom$ & Order \\ \hline 
    $512\times256$ & $5.2272\times10^{-2}$ & -- & $2.3259\times10^{-1}$ & -- & $2.4051\times10^{-1}$ & -- \\ \hline 
    $768\times384$ & $2.2214\times10^{-2}$ & $1.11$ & $1.1263\times10^{-1}$ & $0.79$ & $1.1639\times10^{-1}$ & $0.79$ \\ \hline 
    $1024\times512$ & $1.1326\times10^{-2}$ & $1.33$ & $6.3290\times10^{-2}$ & $1.00$ & $6.5100\times10^{-2}$ & $1.02$ \\ \hline 
    $1280\times640$ & $6.1266\times10^{-3}$ & $1.75$ & $3.7231\times10^{-2}$ & $1.38$ & $3.8318\times10^{-2}$ & $1.38$ \\ \hline 
    $1536\times768$ & $2.5521\times10^{-3}$ & $3.80$ & $1.1869\times10^{-2}$ & $5.27$ & $1.4159\times10^{-2}$ & $4.46$ \\ \hline 
    $1792\times896$ & $1.1643\times10^{-3}$ & $4.09$ & $6.9587\times10^{-3}$ & $2.47$ & $7.1981\times10^{-3}$ & $3.39$ \\
    \hline
    \end{tabular}
    \caption{\small Sup-norm relative errors and numerical orders of $u_1,\om_1,\vom$ at $t = 1.5\times 10^{-4}$ in Case $3$.}
    \label{tab:sup-norm_error_3}
    \vspace{-0.2in}
\end{table}

\section{Comparison with the Original Navier--Stokes Equations}\label{sec:original_NSE}
In this section, we compare the solution to the equations \eqref{eq:axisymmetric_NSE_1} with the variable viscosity coefficients \eqref{eq:viscosity_coefficient} (Case $1$) and the solution to the original Navier--Stokes Equations (Case $2$) using the same initial-boundary conditions \eqref{eq:BC}, \eqref{eq:initial_data}. This comparison will explain why the degeneracy of the variable viscosity coefficients is crucial for the solution to develop a potential finite-time singularity. In fact, we observe that the Navier--Stokes equations with a constant viscosity coefficient will destroy the critical two-scale feature of the solution and eventually prevent the finite-time blowup. 

\subsection{Profile evolution in Case $2$} In Section \ref{sec:first_sign}, we studied the evolution of the solution in Case $1$, and observed a stable blowup with a two-scale feature. Here, we investigate how the solution evolves differently when the degenerate viscosity coefficients $\nu^r,\nu^z$ are replaced by a constant $\mu$. As an illustration, we will focus our study on the case where $\mu = 10^{-5}$. In what follows, Case $2$ refers to the computation of the Navier--Stokes equations with constant viscosity coefficient $\mu = 10^{-5}$ without further clarification. Similar phenomena have been observed in Case $2$ when $\mu$ takes different values. 

\begin{figure}[!ht]
\centering
    \includegraphics[width=0.32\textwidth]{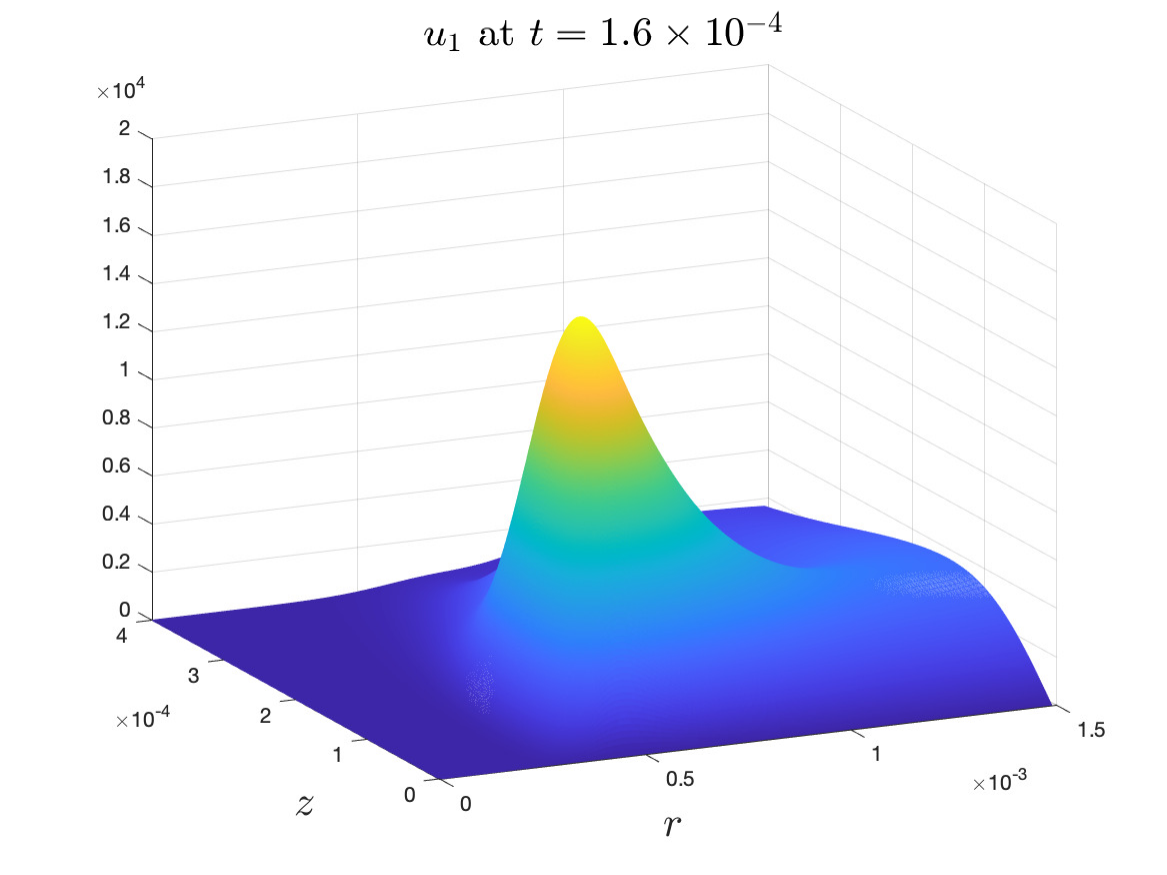}
    \includegraphics[width=0.32\textwidth]{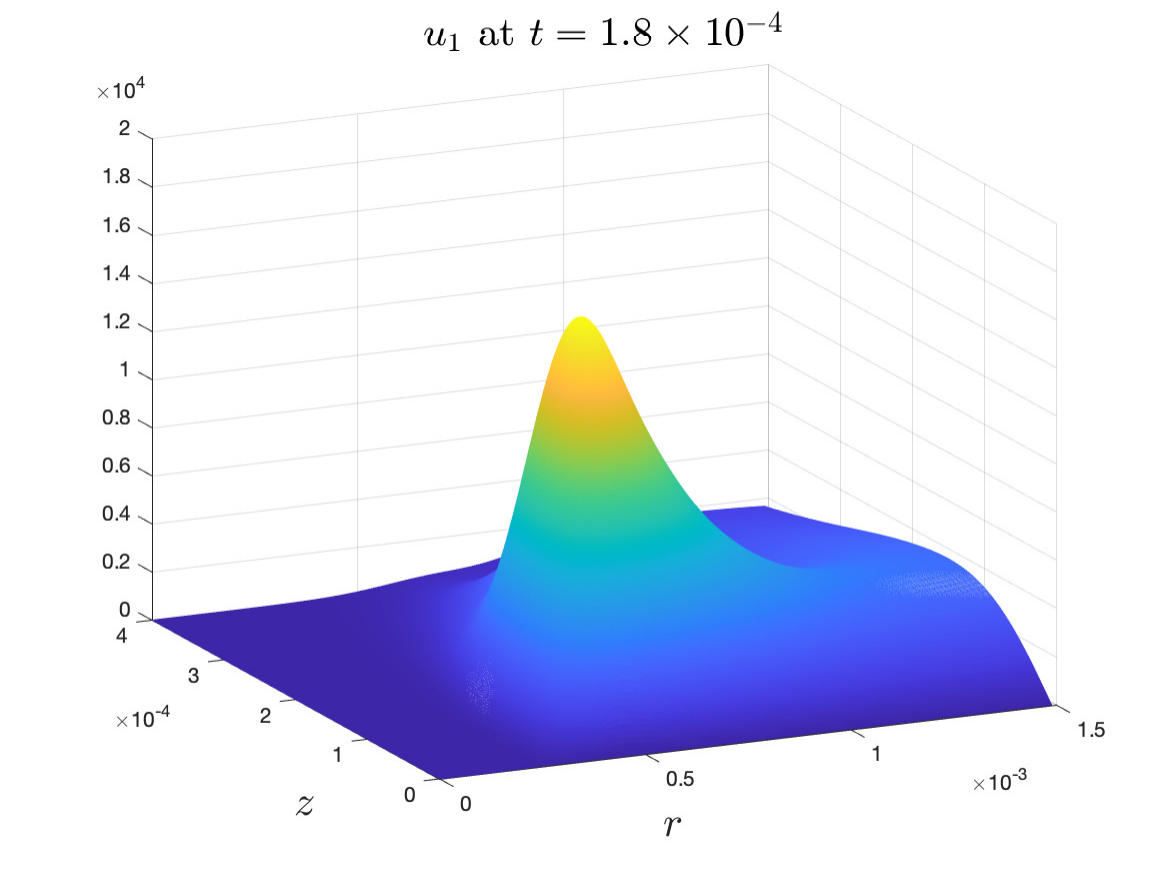}
    \includegraphics[width=0.32\textwidth]{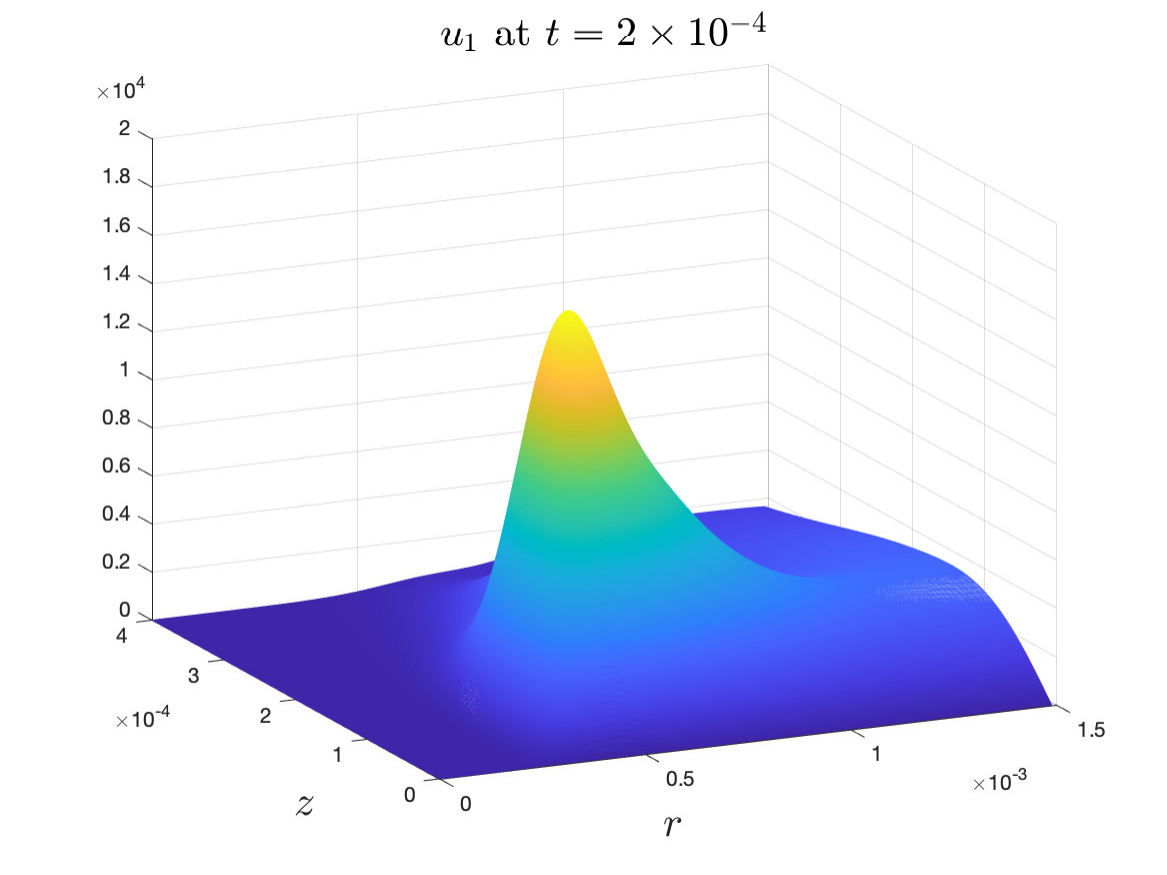}
    \includegraphics[width=0.32\textwidth]{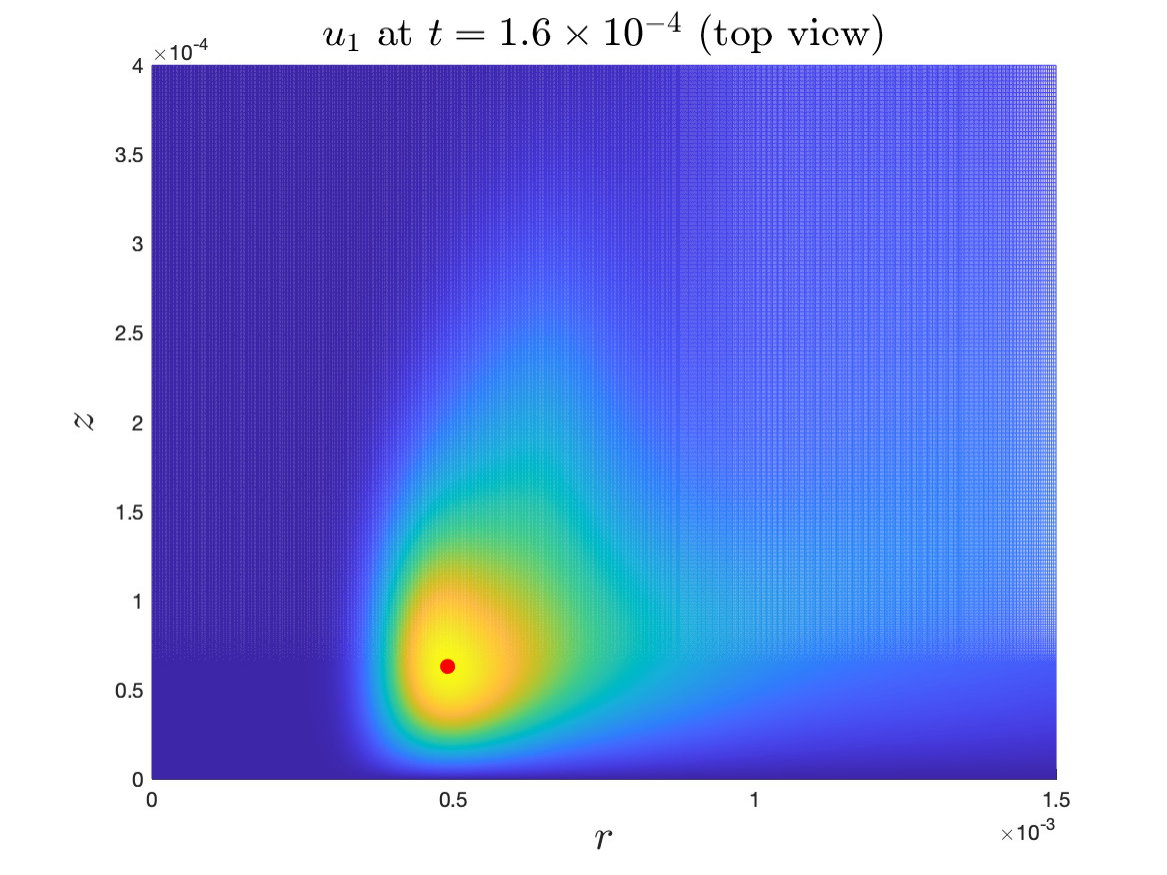}
    \includegraphics[width=0.32\textwidth]{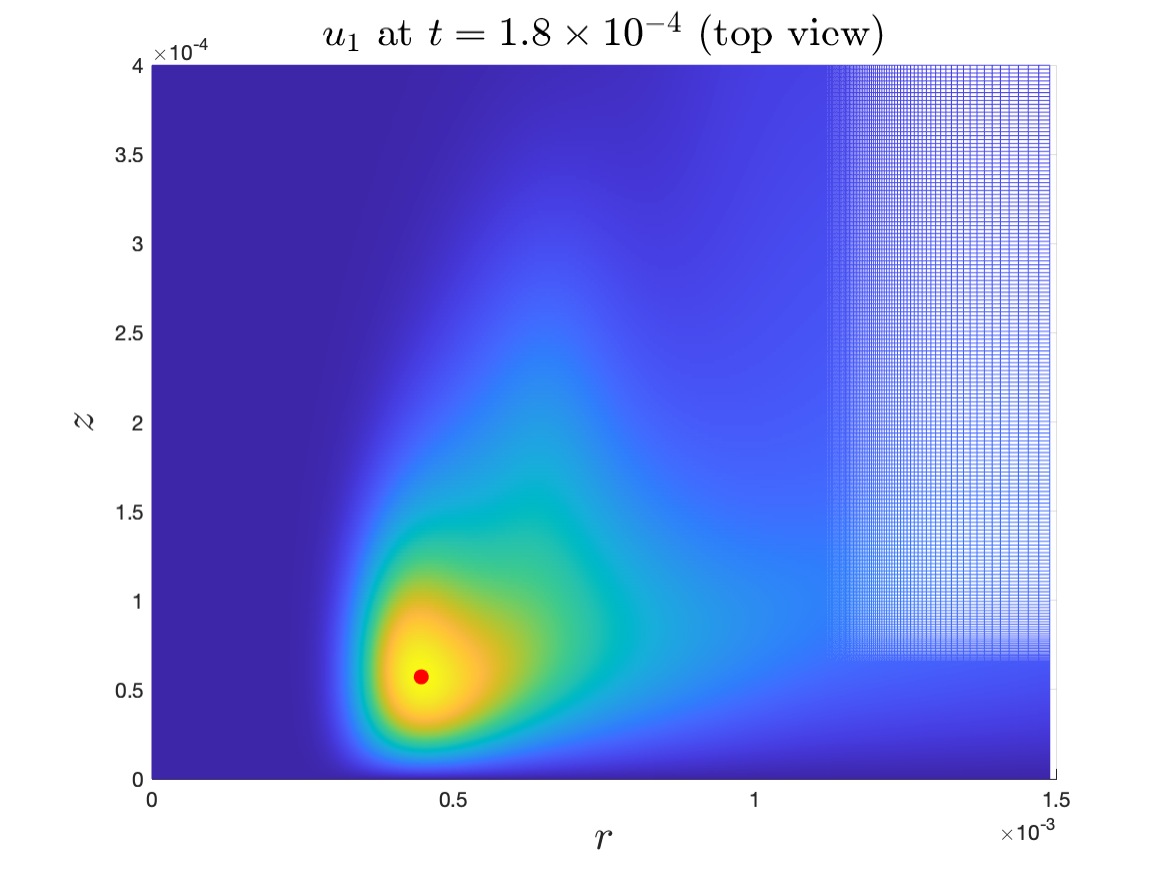}
    \includegraphics[width=0.32\textwidth]{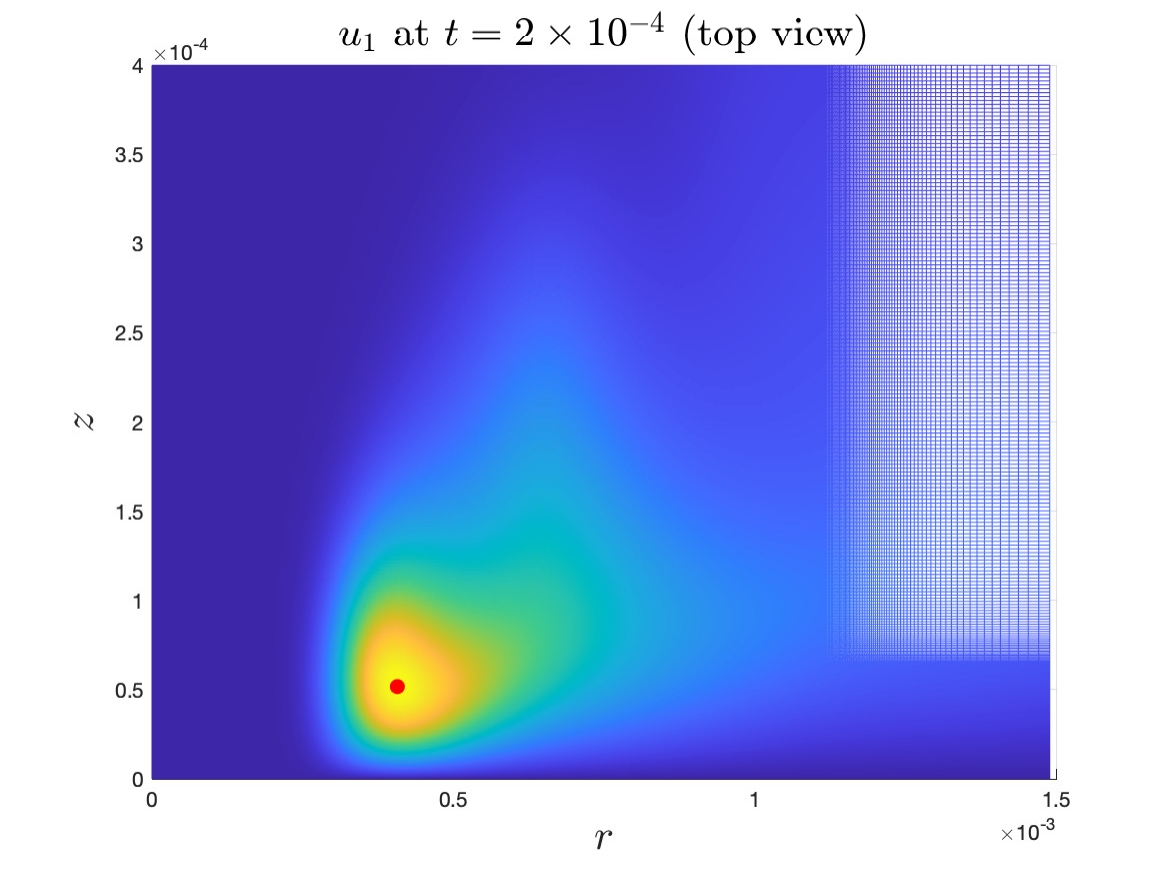}
    \includegraphics[width=0.32\textwidth]{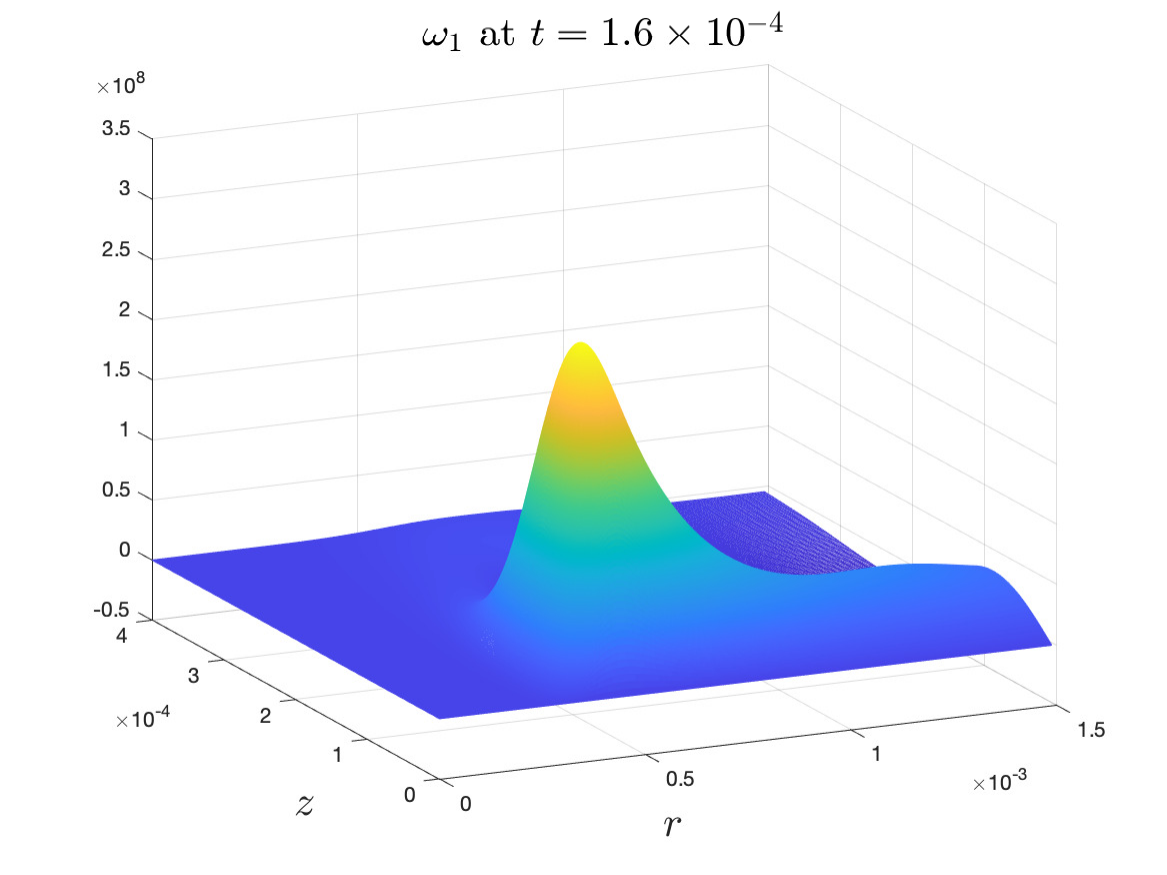}
    \includegraphics[width=0.32\textwidth]{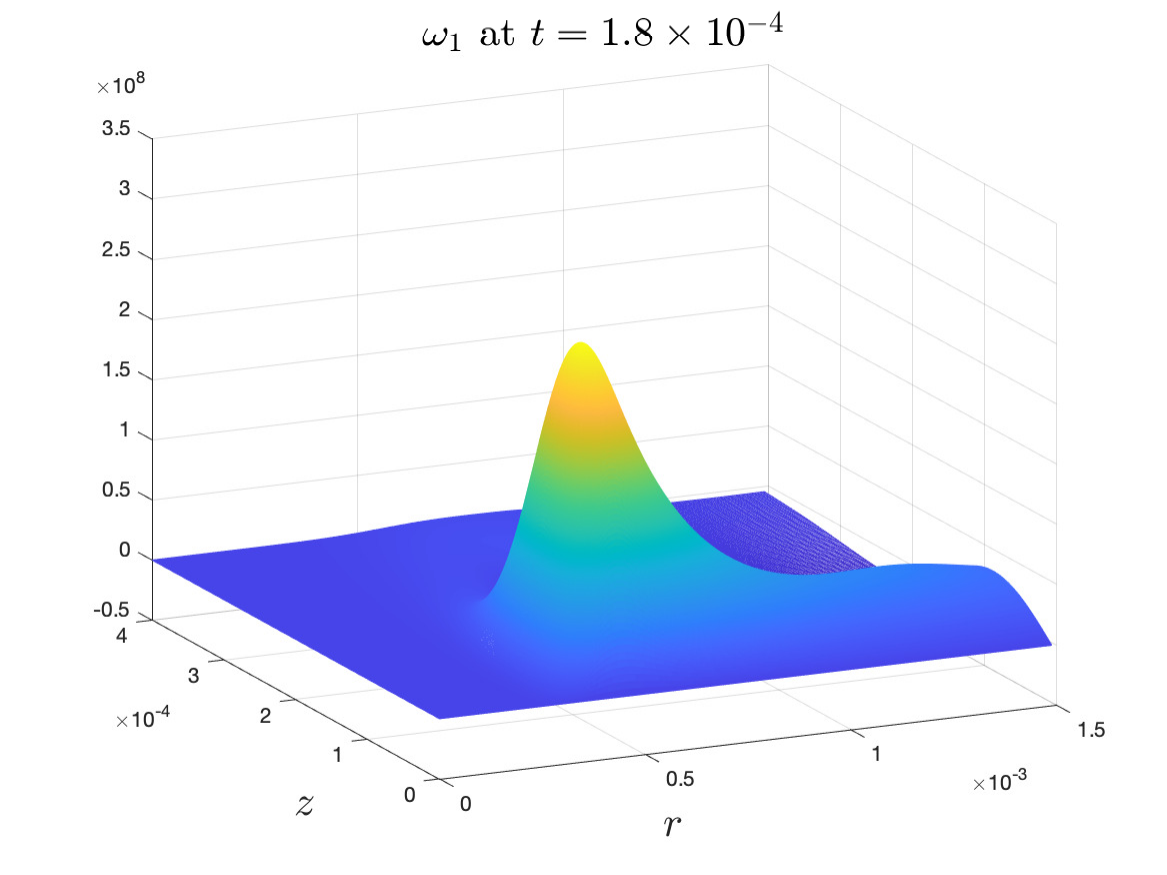}
    \includegraphics[width=0.32\textwidth]{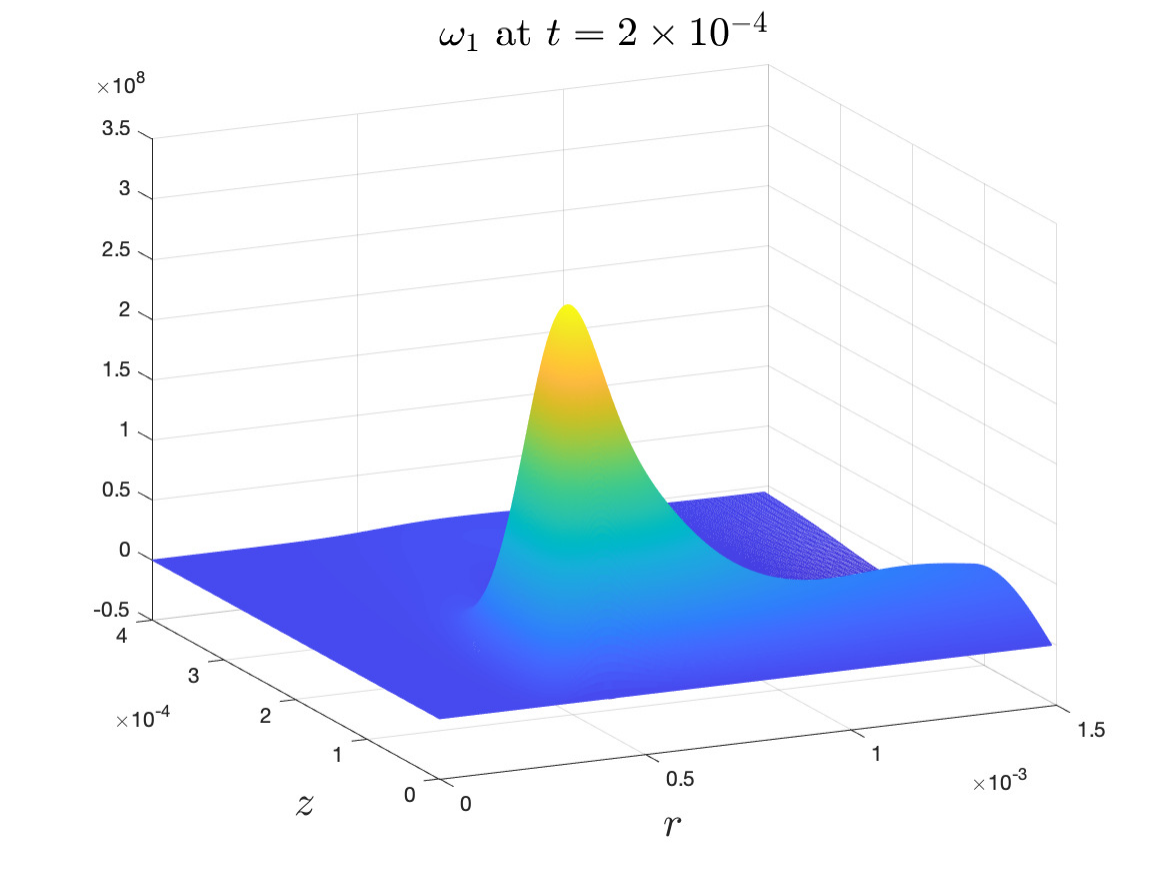}
    \includegraphics[width=0.32\textwidth]{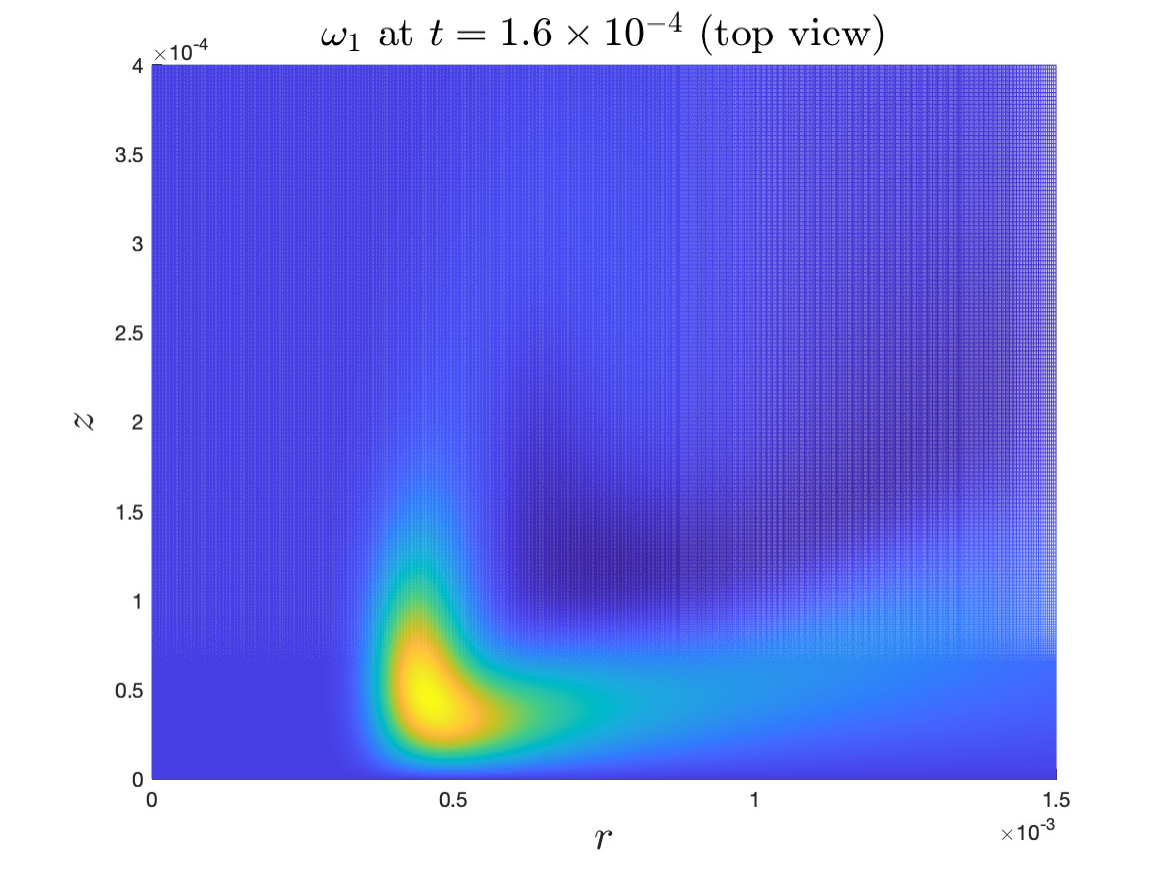}
    \includegraphics[width=0.32\textwidth]{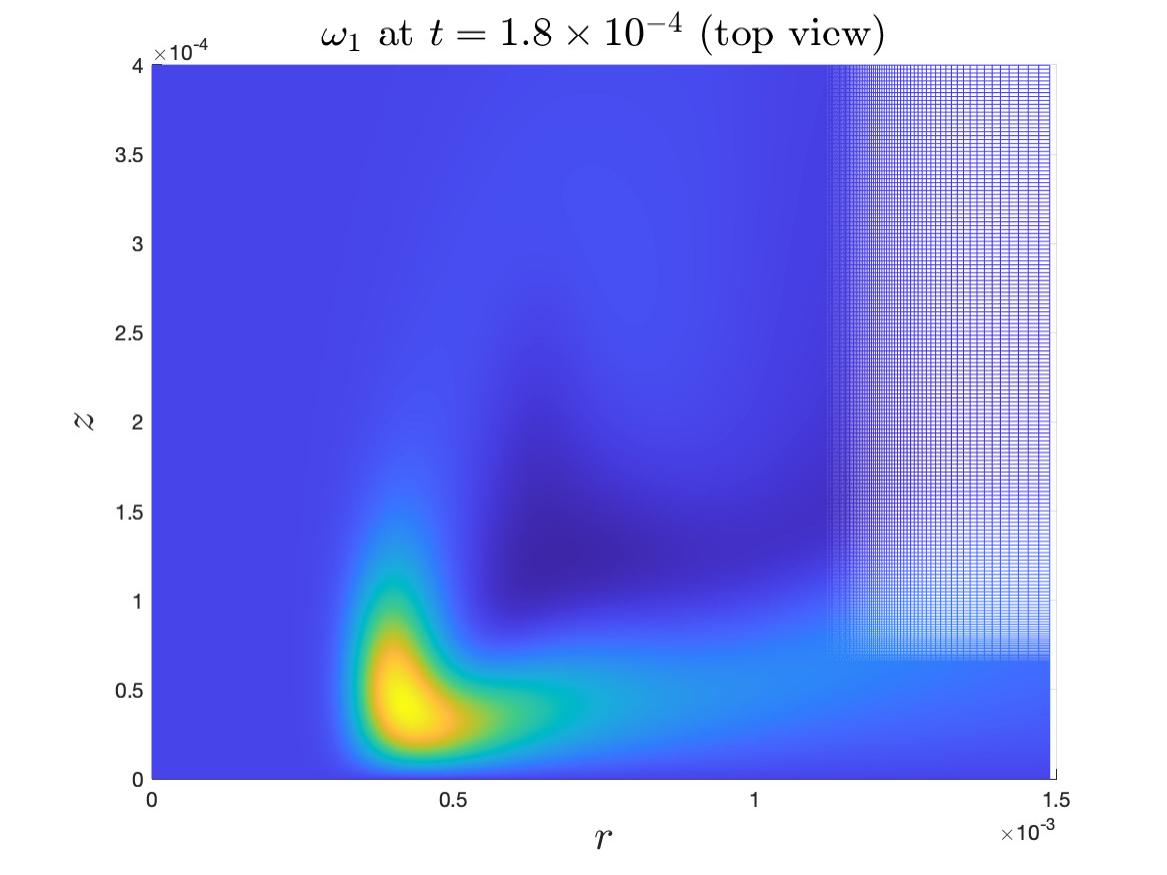} 
    \includegraphics[width=0.32\textwidth]{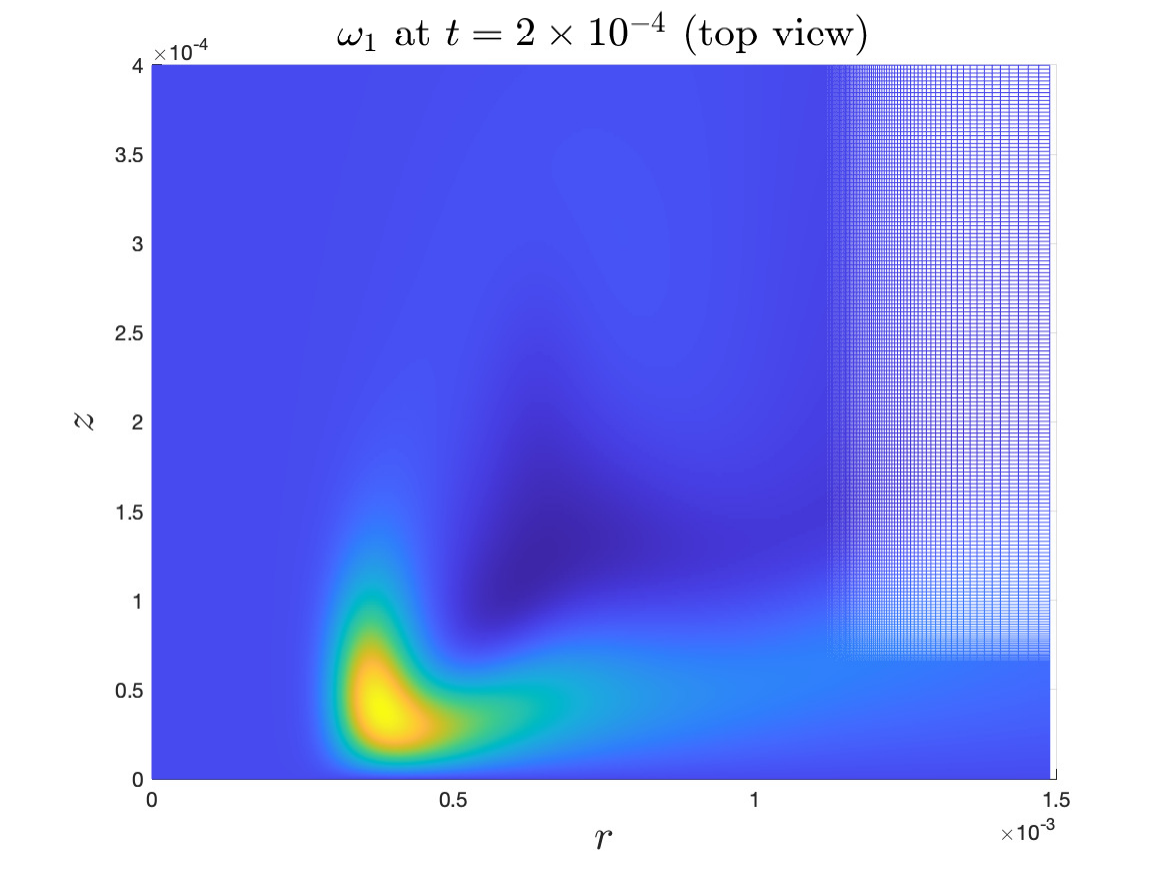} 
    \caption[Profile evolution in Case $2$]{The evolution of the profiles of $u_1$ (row $1$ and $2$) and $\om_1$ (row $3$ and $4$) in Case $2$ with $\nu^r=\nu^z = 10^{-5}$. Line $1$ and $3$ are the profiles of $u_1,\om_1$ at three different times;  Line $2$ and $4$ are the corresponding top-views. The red dot is the location of maximum point of $u_1$.}  
    \label{fig:profile_evolution_case3}
\end{figure}

Figure \eqref{fig:profile_evolution_case3} demonstrates the evolution of the solution in Case $2$ from $t=1.6\times 10^{-4}$ to $t=2\times 10^{-4}$. One should notice the obvious difference in behavior between the solution in Case $1$ and that in Case $2$ when comparing Figure \eqref{fig:profile_evolution_case3} with Figure \eqref{fig:profile_evolution}. Below we list some of our observations. 
\begin{itemize}
\item Unlike in Case $1$, the computation in Case $2$ can be continued to a much later time, and the solution still remains quite smooth. 
\item The solution does not change much from $t=1.6\times 10^{-4}$ to $t=2\times 10^{-4}$. In particular, it does not develop a two-scale spatial structure. Instead, it maintains a profile with a single scale that is comparable to $R(t)$, the distance between the maximum point of $u_1$ and the symmetry axis $r=0$. Moreover, the profile of $u_1$ does not form a sharp gradient in the $z$ direction or a sharp front in the $r$ direction, and the profile of $\om_1$ does not develop a thin structure.
\item The maximums of the solution $u_1$ and $\omega_1$ only grow modestly in the early stage and eventually decrease in the late stage. From $t=0$ to $t=2\times 10^{-4}$, $\|u_1\|_{L^\infty}$ increases only  by a factor of $2.34$, and $\|\om_1\|_{L^\infty}$ increases only by a factor of $3.67$.
\end{itemize}
These observations suggest that the solution in Case $2$ does not develop a finite-time blowup, at least not in the same way as in Case $1$. The main reason for such difference in behavior is that the viscosity term with a constant viscosity coefficient is so strong that it regularizes the smaller scale $Z(t)$ in the two-scale solution profile that we observed in Section \ref{sec:mechanism}, thus destroys the critical blowup mechanism. We will explain in Section \ref{sec:asymptotic_analysis} why the degenerate viscosity coefficient is crucial for a two-scale blowup to appear and persist. 

\begin{figure}[!ht]
\centering
    \includegraphics[width=0.40\textwidth]{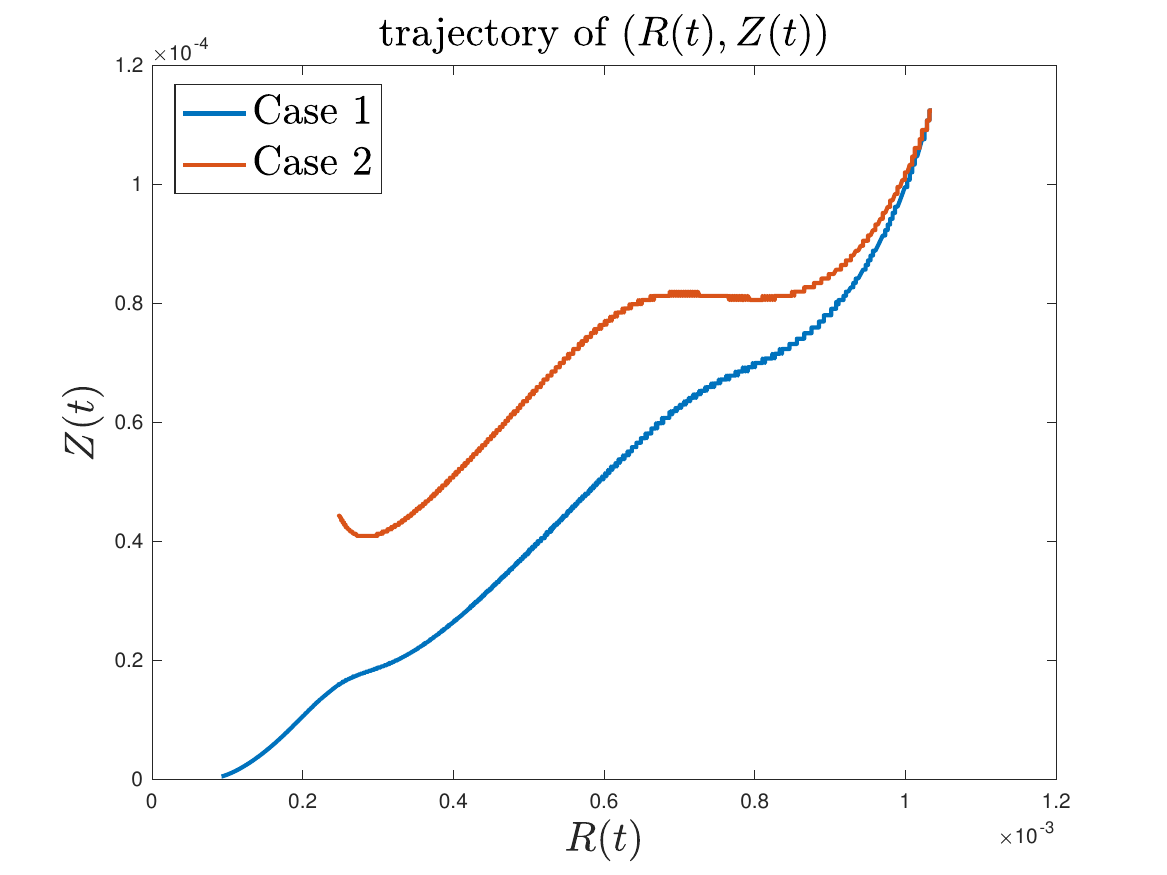}
    \includegraphics[width=0.40\textwidth]{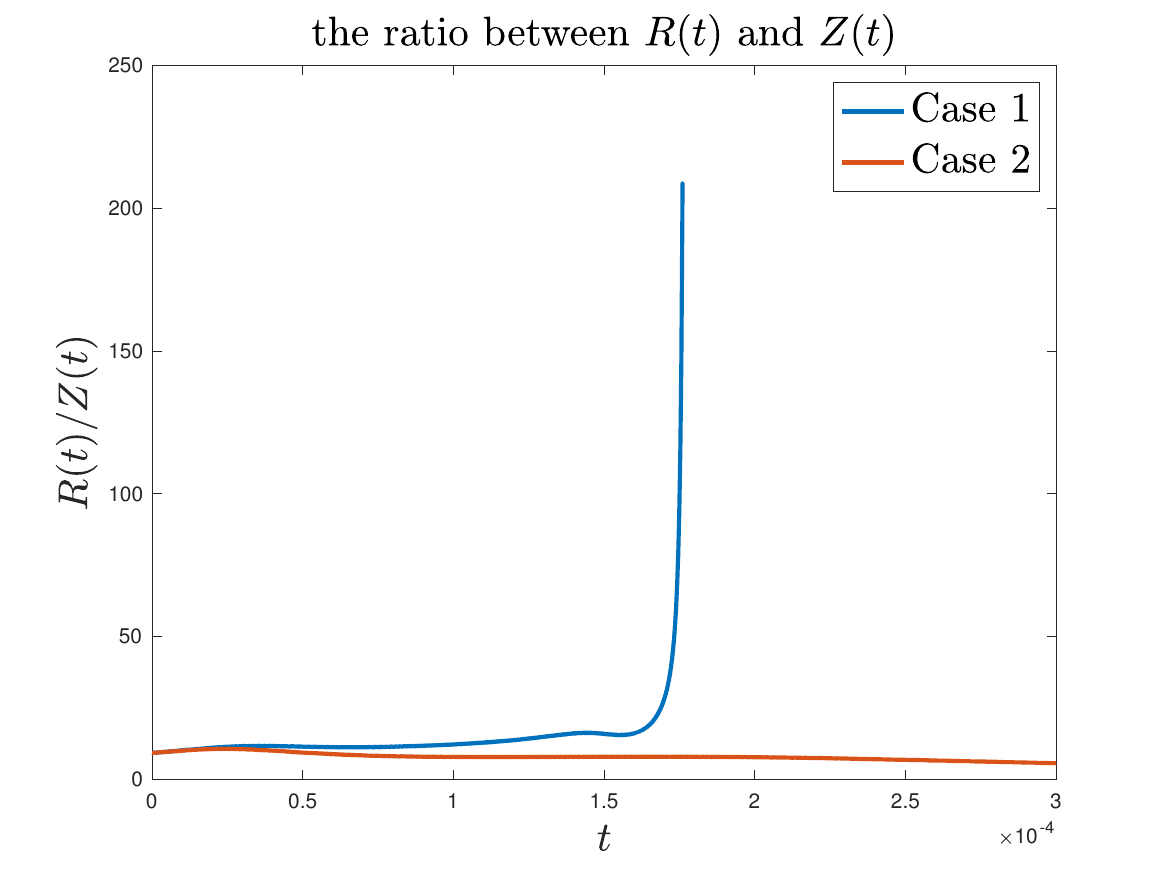}
    \caption[Trajectory compare]{Left: trajectories of $(R(t),Z(t))$ in Case $1$ and in Case $2$. Right: ratios between $R(t),Z(t)$ in Case $1$ and in Case $2$. Blue curves: Case $1$ for $t\in[0,1.76\times 10^{-4}]$. Red curves: Case $2$ for $t\in[0,3\times 10^{-4}]$.}  
    \label{fig:trajectory_compare}
       \vspace{-0.05in}
\end{figure}

Figure \eqref{fig:trajectory_compare} compares the trajectories of $(R(t),Z(t))$ and the ratios $R(t)/Z(t)$ in Case $1$ for $t\in[0,1.76\times 10^{-4}]$ and in Case $2$ for $t\in[0,3\times 10^{-4}]$. We can see that, due to the effect of the stronger viscosity, the point $(R(t),Z(t))$ in Case $2$ does not travel towards the symmetry axis $r=0$ or towards the symmetry plane $z=0$ as fast as in Case $1$. The ratio $R(t)/Z(t)$ in Case $2$ does not blow up rapidly; instead, the two coordinates remain comparable to each other. This again confirms that the solution does not develop the critical two-scale feature in Case $2$. More interestingly, the red trajectory turns upward after some time, suggesting that there will be no blowup focusing at the origin $(r,z)=(0,0)$ in Case $2$. This is consistent with our discussion in Section \ref{sec:mechanism} that when there is no two-scale feature, the maim profile of the solution will eventually be pushed away from the ``ground'' $z=0$ by the upward flow. As a result, the critical blowup mechanism in our scenario will be destroyed.  

\subsection{Growth of some key quantities in Case $2$} 
\begin{figure}[!ht]
\centering
    \includegraphics[width=0.32\textwidth]{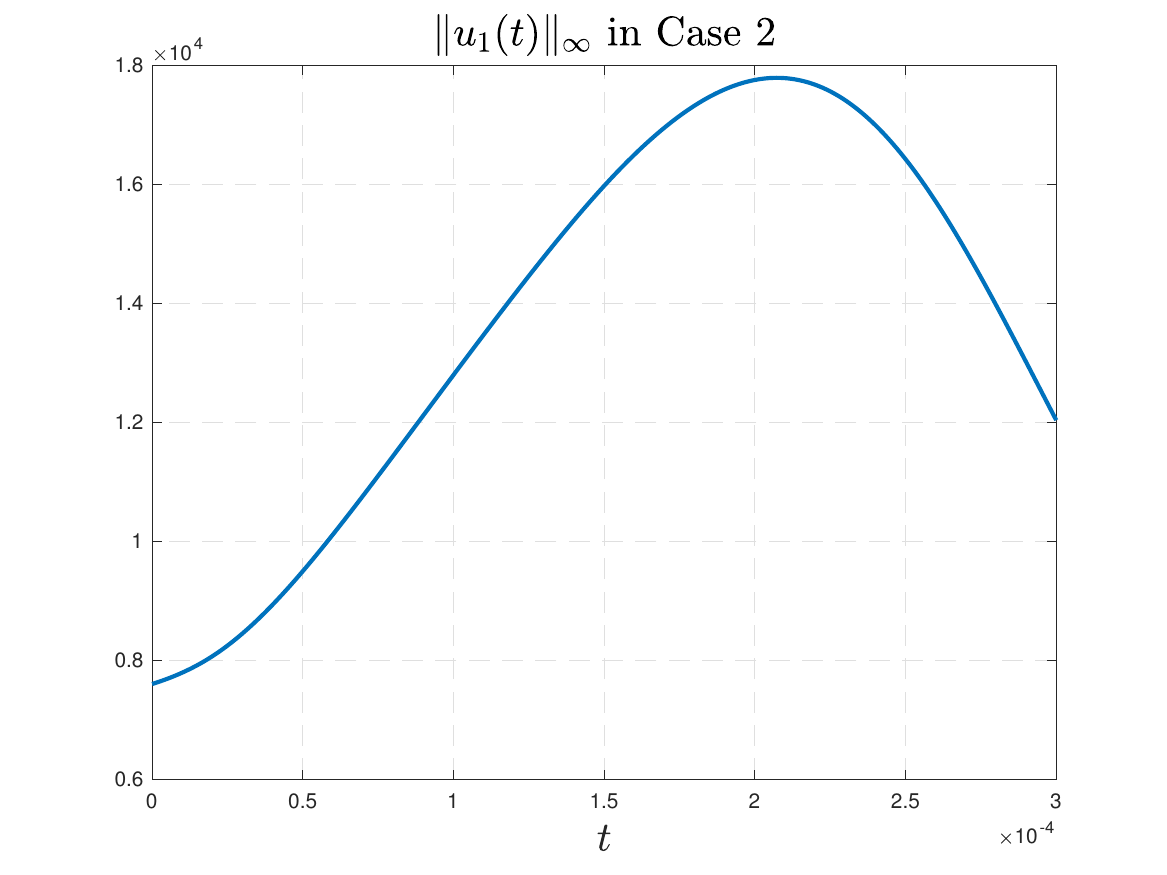}
    \includegraphics[width=0.32\textwidth]{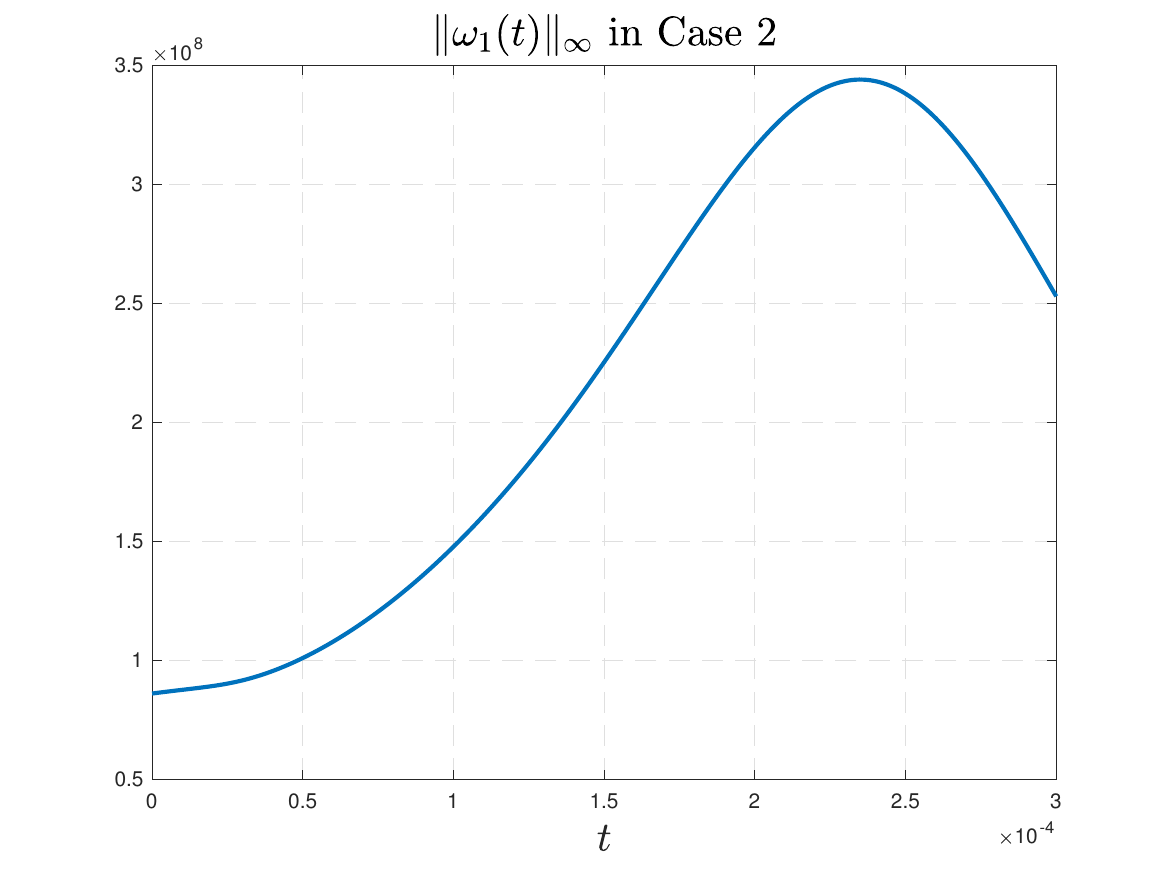} 
    \includegraphics[width=0.32\textwidth]{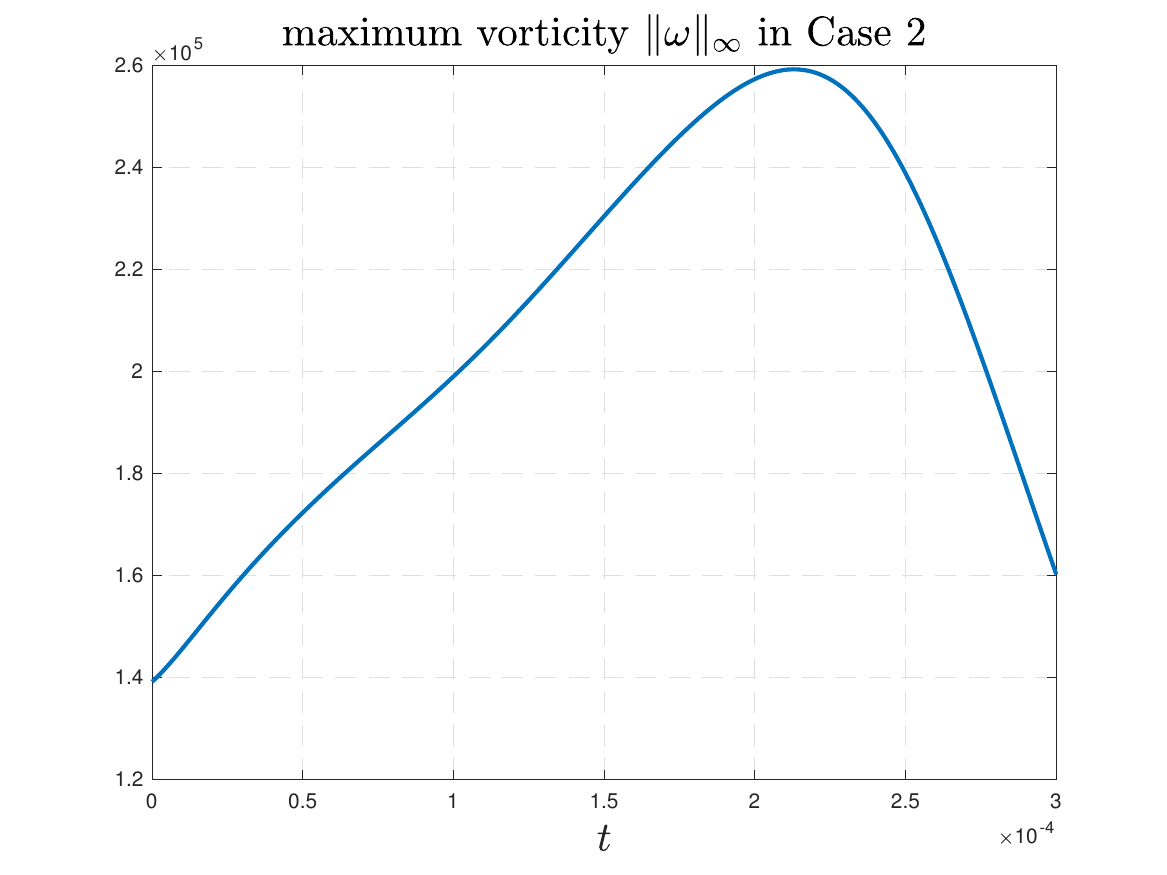}
    \caption[No growth in Case $2$]{The values of $\|u_1\|_{L^\infty}$, $\|\om_1\|_{L^\infty}$ and $\|\vom\|_{L^\infty}$ as functions of time in Case $2$.} 
    \label{fig:no_growth}
\end{figure}

To further illustrate that the solution will remain regular in Case $2$, we directly study the growth of different solution variables. Figure \ref{fig:no_growth} plots the maximums of $u_1,\om_1,|\vom|$ as functions of time. We can see that these quantities do not increase rapidly as in Case $1$ (compared to the dramatic growth shown in Figure \ref{fig:rapid_growth}); moreover, they all start to decrease after some time. Note that the Beale--Kato--Majda criterion (see Section \ref{sec:rapid_growth}) also applies to the Navier--Stokes equations: the solution develops a singularity at some finite time $T$ if and only if $\int_0^T\|\vom\|_{L^\infty}\idiff t = +\infty$. From Figure \ref{fig:no_growth} we can see that the maximum vorticity $\|\vom\|_{L^\infty}$ tends to remain bounded, at least for the duration of our computation. This observation strongly suggests that the solution to the equations \eqref{eq:axisymmetric_NSE_1} with a constant viscosity coefficient (namely the Navier--Stokes equations) will not blow up under the initial-boundary conditions \eqref{eq:BC},\eqref{eq:initial_data}. 

To understand why the maximum of $u_1$ does not rise rapidly and eventually drops in the later stage in Case $2$, we study the competition between the vortex stretching and the viscosity. Since $2u_1$ is the leading order part of the axial vorticity $\om^z = 2u_1 + ru^r$ for $r$ near 0, the forcing term $2\psi_{1,z}u_1$ in the $u_1$ equation \eqref{eq:axisymmetric_NSE_1} can be considered as a vortex stretching term. This term is the driving force for the growth of $u_1$. On the contrary, the viscosity term $f_{u_1}$ (given by \eqref{eq:viscosity_u1}), which is always negative at $(R(t),Z(t))$, damps the maximum of $u_1$. If the vortex stretching dominates the viscosity near $(R(t),Z(t))$, $\|u_1\|_{L^\infty}$ should grow; otherwise, $\|u_1\|_{L^\infty}$ will drop. 

\begin{figure}[!ht]
\centering
    \includegraphics[width=0.40\textwidth]{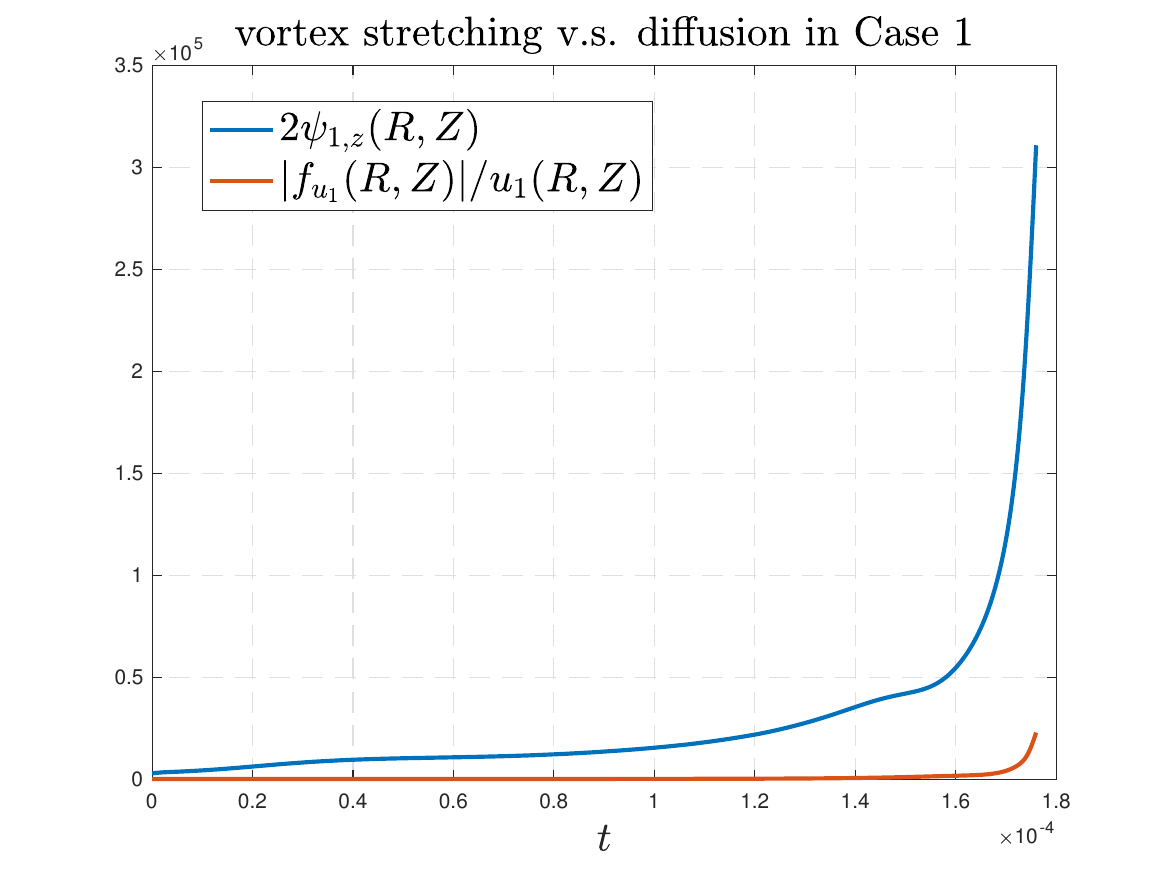}
    \includegraphics[width=0.40\textwidth]{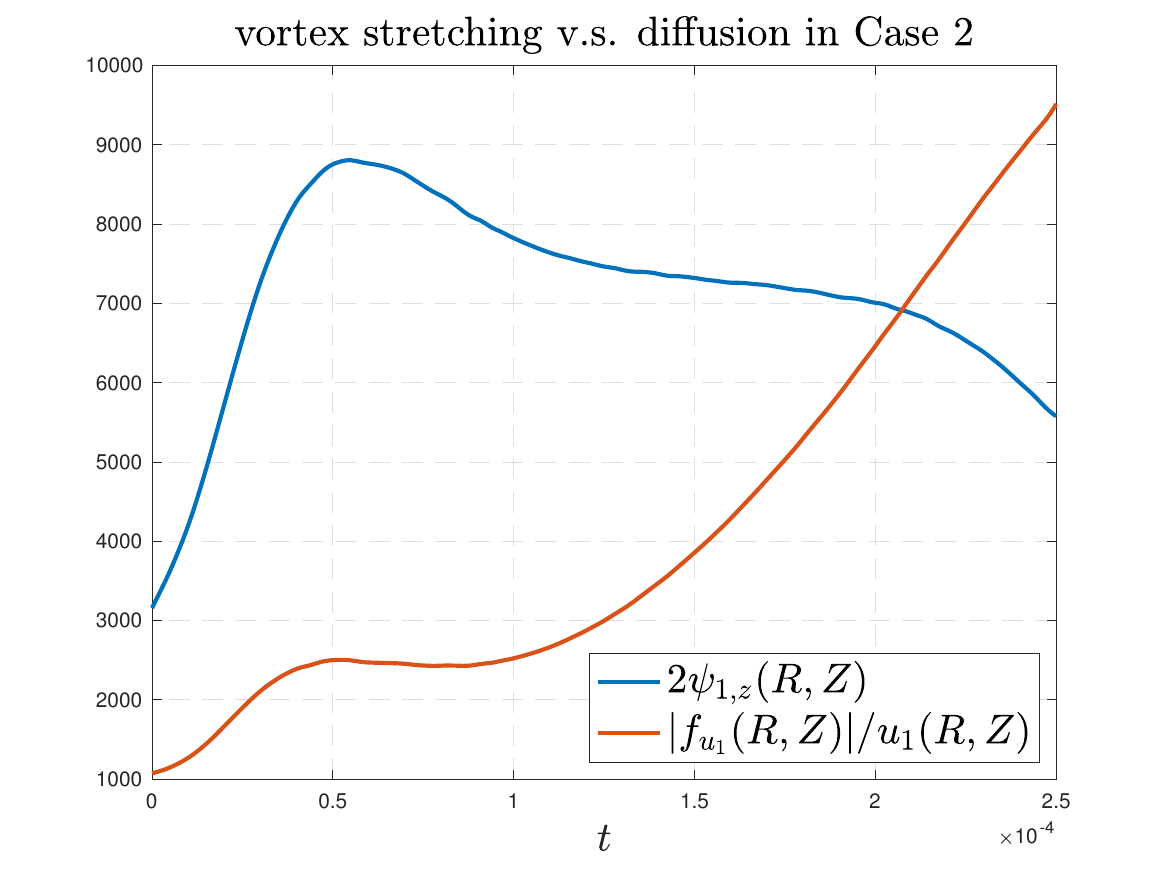} 
    \caption[Competition between vortex stretching and viscosity]{The relative magnitudes of the vortex stretching $2\psi_{1,z}$ and the viscosity $|f_{u_1}|/u_1$ at the point $(R(t),Z(t))$ in Case $1$ (left) and in Case $2$ (right).}
    \label{fig:competition}
       \vspace{-0.05in}
\end{figure}

In Figure \ref{fig:competition} we plot the relative magnitudes of the vortex stretching $2\psi_{1,z}$ and the viscosity $|f_{u_1}|/u_1$ at $(R(t),Z(t))$ in Case $1$ (left) and in Case $2$ (right). It is clear that the vortex stretching keeps growing and always dominates the viscosity term in the $u_1$ equations at $(R(t),Z(t))$ in Case $1$; thus we observe a rapid growth of $\|u_1\|_{L^\infty}$ in time. This is the consequence of (i) the good alignment between $\psi_{1,z}$ and $u_1$ that relies on the thin structure (the smaller scale) of the solution in the $z$ direction as described in Section \ref{sec:mechanism} and (ii) the fact that the viscosity coefficients are degenerate at the origin. On the contrary, we observe in Case 2 that the relative strength of the vortex stretching starts to decrease after some time and is dominated by the viscosity term in later time, which leads to the decrease of $\|u_1\|_{L^\infty}$. This is caused by the strong viscosity from two aspects. On the one hand, if the blowup mechanism in Section \ref{sec:mechanism} tries to generate a thinner scale in the solution, then the viscosity with a constant/non-degenerate coefficient will become too strong for the smaller scale to survive, and thus the alignment between $\psi_{1,z}$ and $u_1$ is not sustainable. On the other hand, if the solution does not develop a two-scale structure, $\psi_{1,z}$ and $u_1$ cannot cannot develop a strong alignment for the coupling mechanism \eqref{eq:mechanism} to last. This dilemma prevents a sustainable blowup to occur in Case $2$.

We remark that we have carried out computations in Case $2$ with different values of $\mu$ in the range $[10^{-7},10^{-4}]$, and we have qualitatively similar observations in all trials: there is no sign of finite-time blowup for all tested values of $\mu$. For a smaller $\mu$, the solution in Case $2$ in the early stage of the computation is very similar the solution in Case $1$, and a two-scale feature seems to develop. However, the viscosity with a constant coefficient will eventually take dominance and eliminate the potential two-scale blowup, and the solution starts to drop afterwards. If $\mu$ is even smaller ($\mu< 10^{-7}$), the viscosity term will be too weak to regularize the sharp fronts in the early stage of the computation and cannot effectively control the mild instability in the intermediate field and far field where the mesh is not as dense as in the near field. The solution quickly becomes under-resolved. It is still unclear whether the solution to the original $3$D axisymmetric Navier--Stokes equations can develop a focusing blowup at the symmetry axis in a different manner when $\mu$ is sufficiently small. Yet we conjecture that this cannot happen in the two-scale manner described in Case 1. 

\subsection{On the viscosity effect} 
We remark that the driving mechanism for the potential finite-time singularity is from the incompressible Euler equations. But it is much harder to obtain a convincing Euler singularity due to the development of extremely sharp fronts in the early stage of the computation and the shearing instability in the tail region.
The Euler equations with degenerate viscosity coefficients capture the main effect of this potential Euler singularity and it is easier to resolve computationally due to the viscous regularization. As we will show in our scaling analysis, the degenerate viscosity coefficients also select a length scale for the smaller scale $Z(t)$, which plays an important role in generating a stable two-scale locally self-similar blowup. Without any viscous regularization, it is not clear what is the mechanism for selecting the length scale of $Z(t)$.

We need to choose the degenerate viscosity coefficients carefully to fulfill two purposes. On the one hand, the viscosity coefficients should be strong enough to control the shearing instability associated with the Euler equations so that the potential blowup can be stable and robust. On the other hand, the viscosity coefficients must not be too strong to suppress the intrinsic blowup mechanism of the Euler equations. The strength of the viscosity coefficients needs to maintain a delicate balance, especially in the crucial blowup region near the sharp front of the solution.    

Out of such consideration, we propose to use the variable viscosity coefficients of the kind \eqref{eq:viscosity_coefficient}, which consists of a space-dependent part (SDP) and a time-dependent part (TDP). The SDP is globally small and degenerate at the origin with order $O(r^2) + O(z^2)$. Recall that the most important part of the blowup solution, i.e. the local profile near the sharp front, travels towards the origin $(r,z)=(0,0)$. Correspondingly, the effective value of the SDP that affects the shrinking region of interest is decreasing to $0$. Thus the blowup mechanism will not be hindered by this degenerate part of the viscosity. As we will see in Section \ref{sec:asymptotic_analysis}, the order of degeneracy of the SDP is compatible with the smaller scale of the solution, so the SDP actually helps stabilize the two-scale blowup. In the meanwhile, the outer part of the solution (the far field) is under the influence of the non-degenerate part of the SDP, which suppresses some mild instability in the far field.

\begin{figure}[!ht]
\centering
    \includegraphics[width=0.40\textwidth]{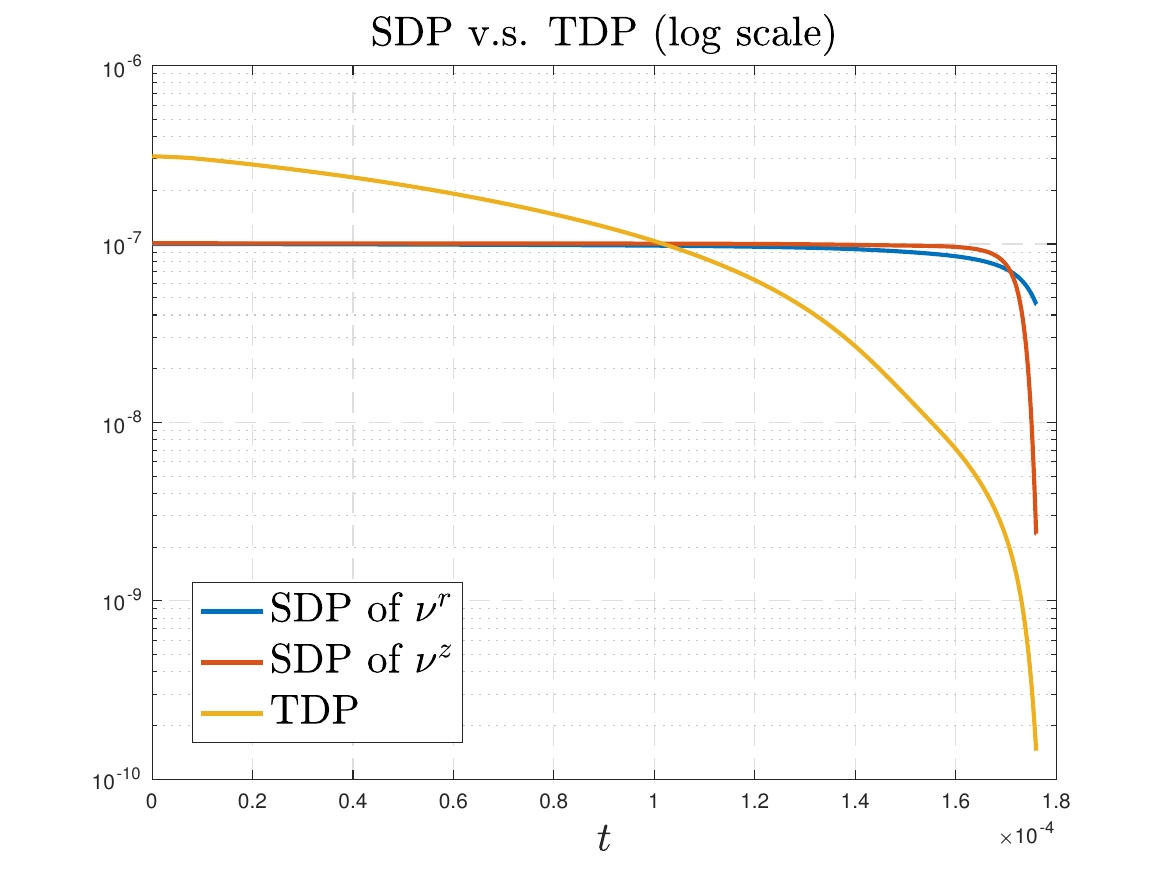}
    \caption[SDP v.s. TDP]{The SDP of $\nu^r$ and the SDP of $\nu^z$ evaluated at the point $(R(t),Z(t))$ versus the TDP.}
    \label{fig:SDPvsTDP}
\end{figure}

The TDP, which is equal to $0.025/\|\om^\theta\|_{L^\infty}$, is relatively large in the early stage of the computation, so that it can help regularize the thin structure of the solution in the warm-up phase ($t\in [0,1.6\times10^{-4}]$). However, as $\|\om^\theta\|_{L^\infty}$ grows rapidly in time, the TDP will drops quickly and be dominated by the SDP in the critical blowup region near $(R(t),Z(t))$, and thus it will not harm the development of the focusing singularity in our blowup scenario. Figure \ref{fig:SDPvsTDP} plots the SDPs of $\nu^r$ and $\nu^z$ evaluated at the point $(R(t),Z(t))$ and the TDP. We can see that the TDP drops blow the SDP at an early time around $t=1\times 10^{-4}$ and is much smaller than the SDP in the stable phase ($t\in [1.6\times10^{-4},1.76\times 10^{-4}]$). In fact, the TDP has no essential contribution to the development of the two-scale blowup. Removing the TDP in the late stage of the computation has almost no influence on the behavior of the solution. 

In summary, the effective values of the variable viscosity coefficients in the focusing blowup region near the sharp front are decreasing in time. As a result, the viscosity in this region of interest is not strong enough to prevent the blowup of the solution. In particular, this weak viscosity cannot prevent the occurrence and development of the two-scale feature in our scenario, which is crucial to the sustainability and stability of the blowup.

On the contrary, the viscosity with a constant coefficient will be too strong for a two-scale blowup to maintain in the late stage when the ratio between the large scale and the smaller scale becomes huge. If there were a two-scale structure, then the smaller scale would have generated a strong viscosity effect (with a constant viscosity coefficient) due to large second derivatives, which, paradoxically, would have smoothed out the smaller scale. This explains why we have observed that the two-scale feature is not sustainable in Case $2$ computation and that the solution does not blowup in the end. In fact, the viscosity with a constant coefficient in the Navier--Stokes equations imposes a constraint on the spatial scaling of the blowup solution, if there is a stable blowup scaling. Such constraint denies the existence of a two-scale structure for the solution, which, however, appears to be crucial for the potential blowup in our computation in Case $1$. We shall come back to this point in the next Section.

\begin{remark}
There were some unpublished computational results obtained by Dr. Guo Luo and the first author for $3$D axisymmetric Euler and Navier--Stokes equations that showed some early sign of a tornado like singularity. However, the solution became unstable dynamically. In the case of the Euler equations, the shear layer in $u_1$ rolled up into several small vortices near the upper edge of the sharp front. The Navier--Stokes equations with a constant viscosity coefficient $\mu=10^{-5}$ regularized these small vortices, but the solution developed into a ``turbulent flow'' in the late stage, which is very different from our computational results for the original Navier--Stokes equations with $\mu=10^{-5}$. Unfortunately, these unpublished computational results are not reproducible due to the loss of Dr. Luo's computer hard disk that contained all the computer codes and the initial data for these results.
\end{remark}

\section{Comparison with the Euler equations}\label{sec:Euler} 
In this section, we will discuss our potential blowup scenario in Case $3$ of the Euler setting. That is, we study the evolution of the solution to the initial-boundary value problem \eqref{eq:axisymmetric_NSE_1}--\eqref{eq:initial_data} with $\nu^r = \nu^z = 0$. As an overview, the Euler solution behaves very similarly to the solution in our main Case $1$ in the warm-up phase. This is not surprising as the critical blowup mechanism discussed in Section \ref{sec:mechanism} relies only on the Euler part of the equations. In particular, the Euler solution grows faster than the solution in Case $1$ (with degenerate viscosity coefficients) during the warm-up phase. However, the Euler solution also quickly develops unfavorable oscillations in the critical blowup region, which is likely due to under-resolution of the extremely sharp structure in the profile. 

\subsection{Profile evolution} 
\begin{figure}[!ht]
\centering
    \includegraphics[width=0.40\textwidth]{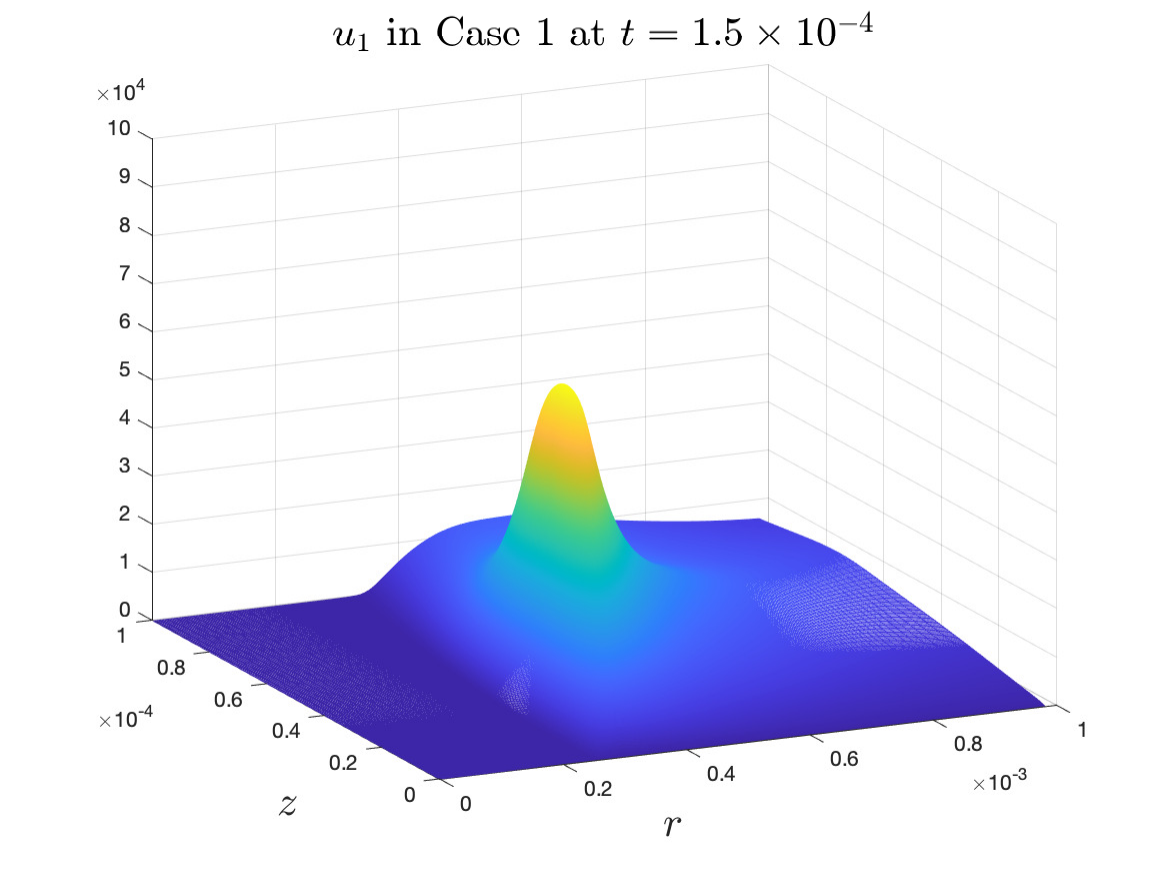}
    \includegraphics[width=0.40\textwidth]{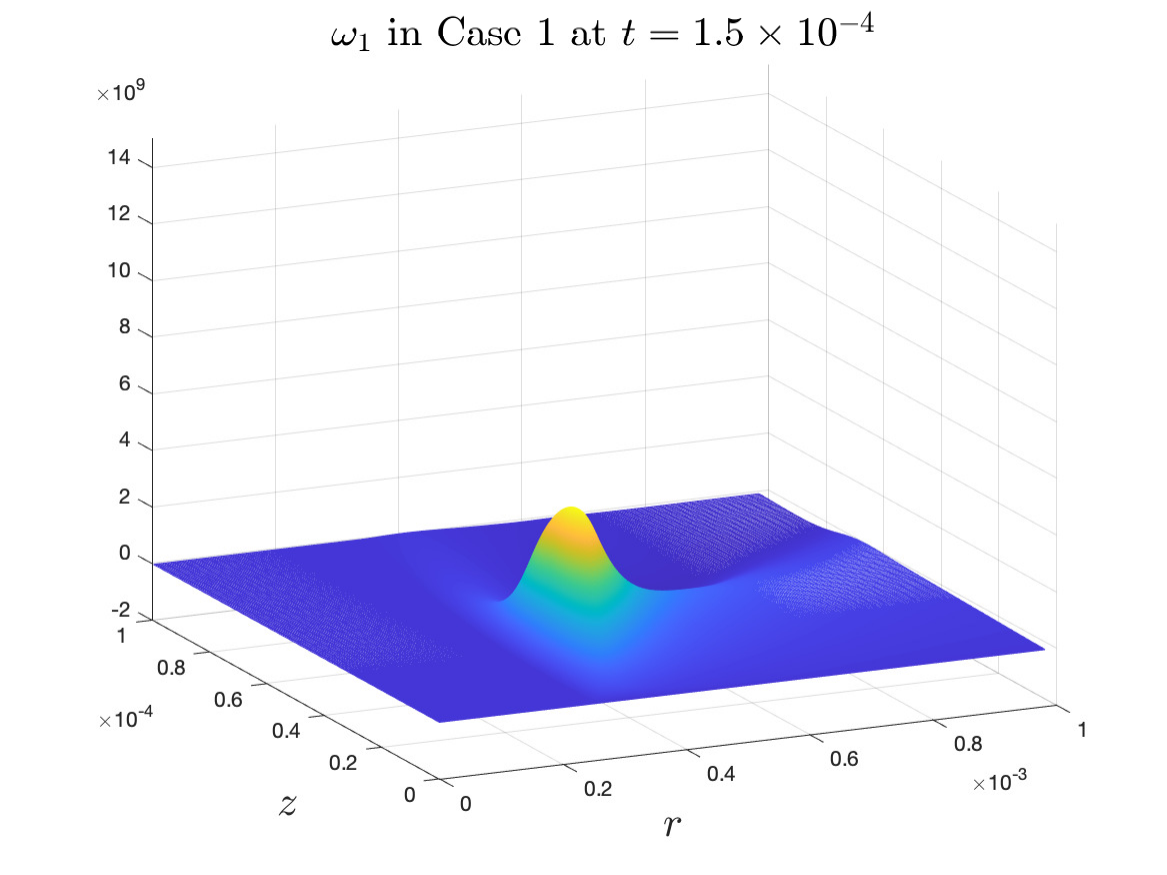} 
    \includegraphics[width=0.40\textwidth]{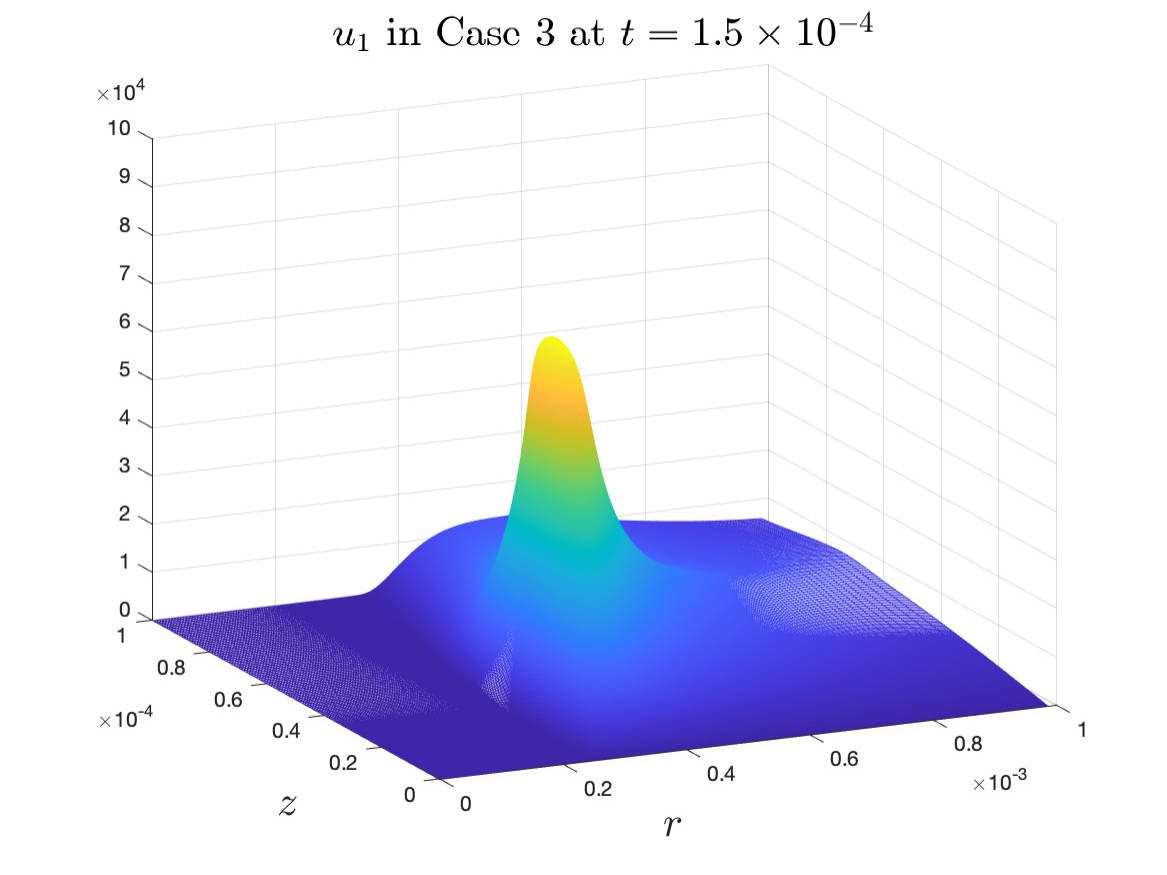}
    \includegraphics[width=0.40\textwidth]{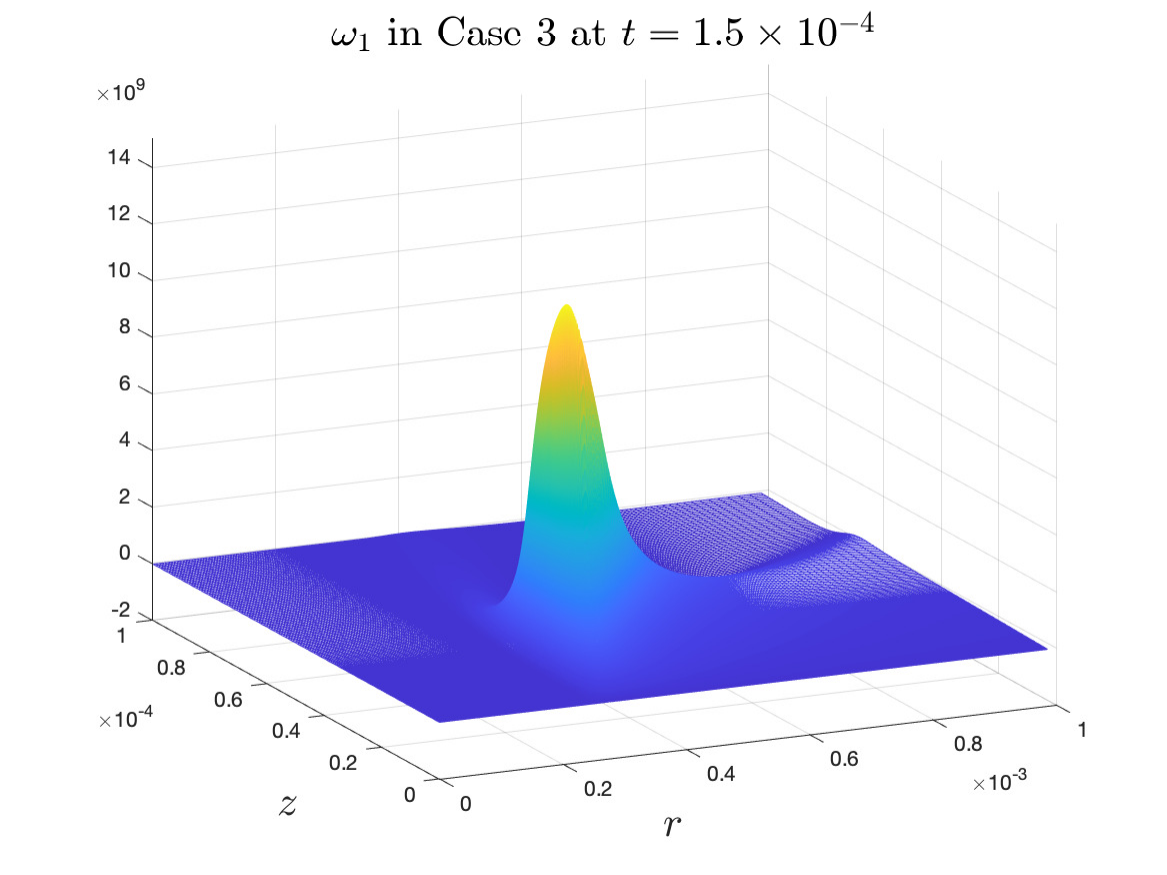}
    \caption[Compare with Euler: Profile]{The profiles of $u_1$ and $\om_1$ at $t=1.5\times 10^{-4}$ in Case $1$ (first row) and in Case $3$ (second row).}  
    \label{fig:Euler_profile_compare}
\end{figure}

To demonstrate that the solutions in Case $1$ and Case $3$ behave similarly in the warm-up phase, we compare their profiles at the same time instant. In Figure \ref{fig:Euler_profile_compare}, we plot the profiles of $u_1$ and $\om_1$ at $t=1.5\times 10^{-4}$ in Case $1$ (first row) and Case $3$ (second row), respectively. We can see that the solution profiles in both cases are qualitatively similar, except that the solution in Case $3$ grows faster. The solution also develops a sharp front in the $r$ direction and a no-spinning region between the front and the symmetry axis $r=0$.

\begin{figure}[!ht]
\centering
    \includegraphics[width=0.40\textwidth]{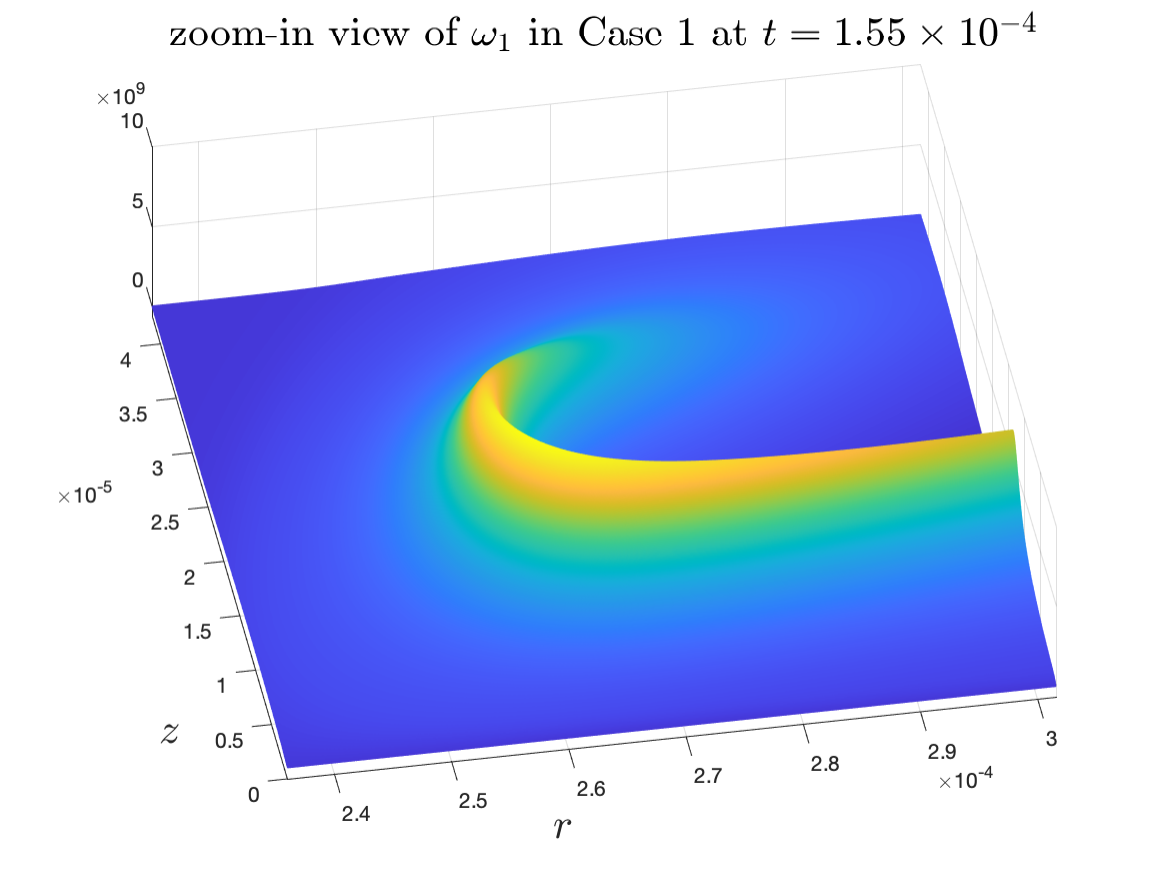}
    \includegraphics[width=0.40\textwidth]{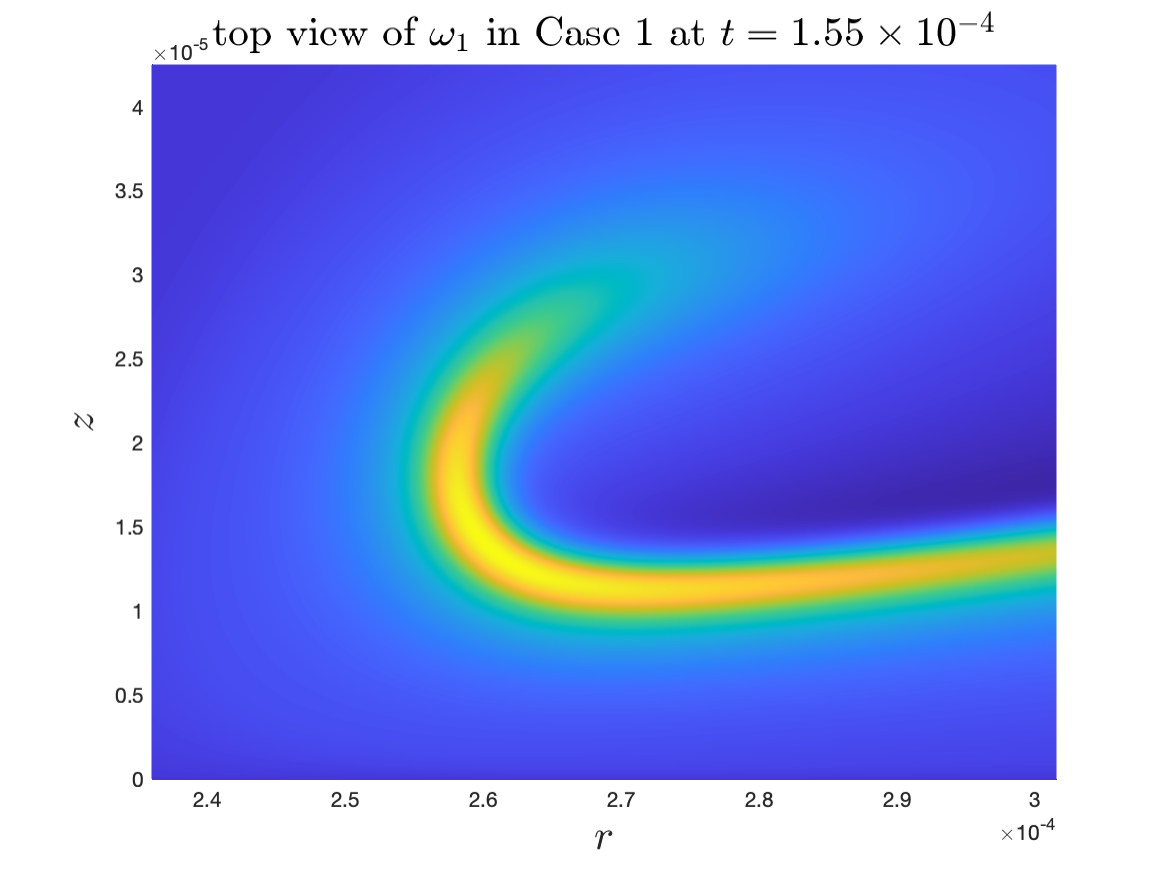} 
    \includegraphics[width=0.40\textwidth]{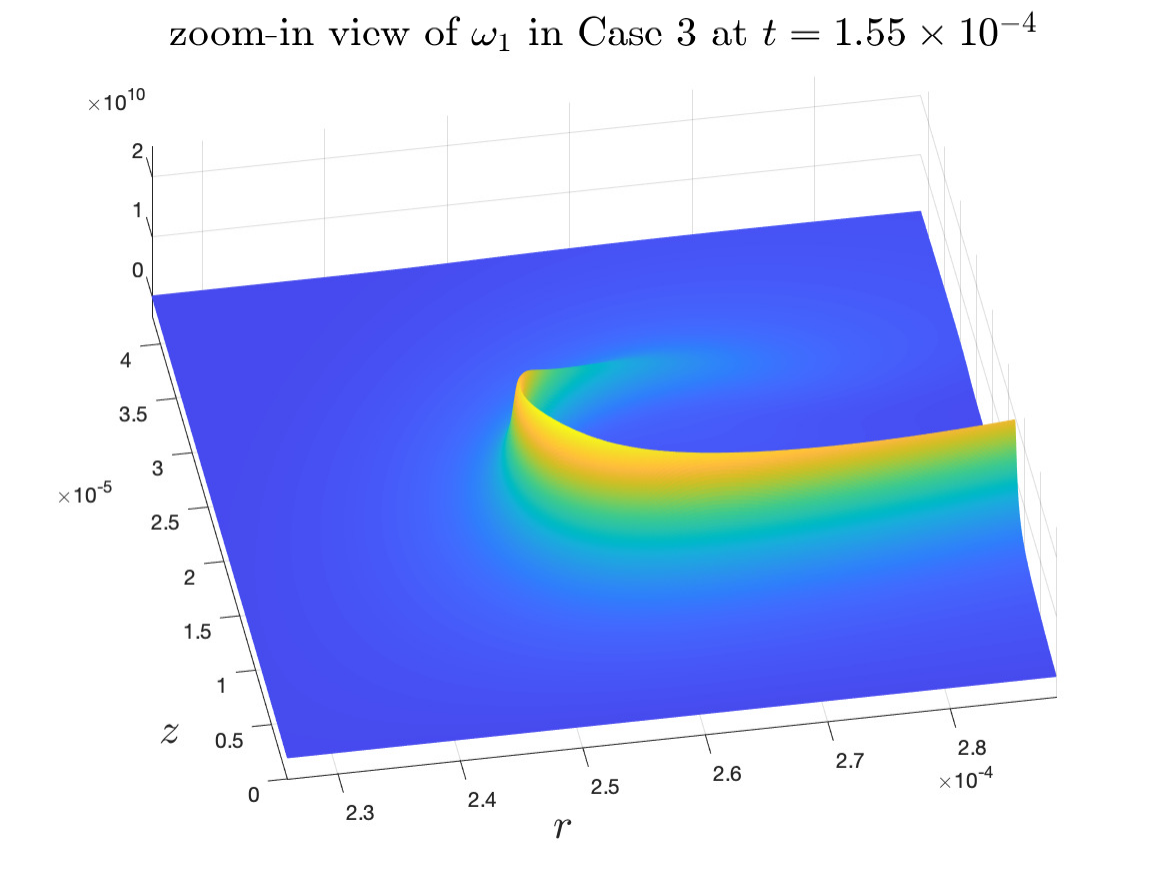}
    \includegraphics[width=0.40\textwidth]{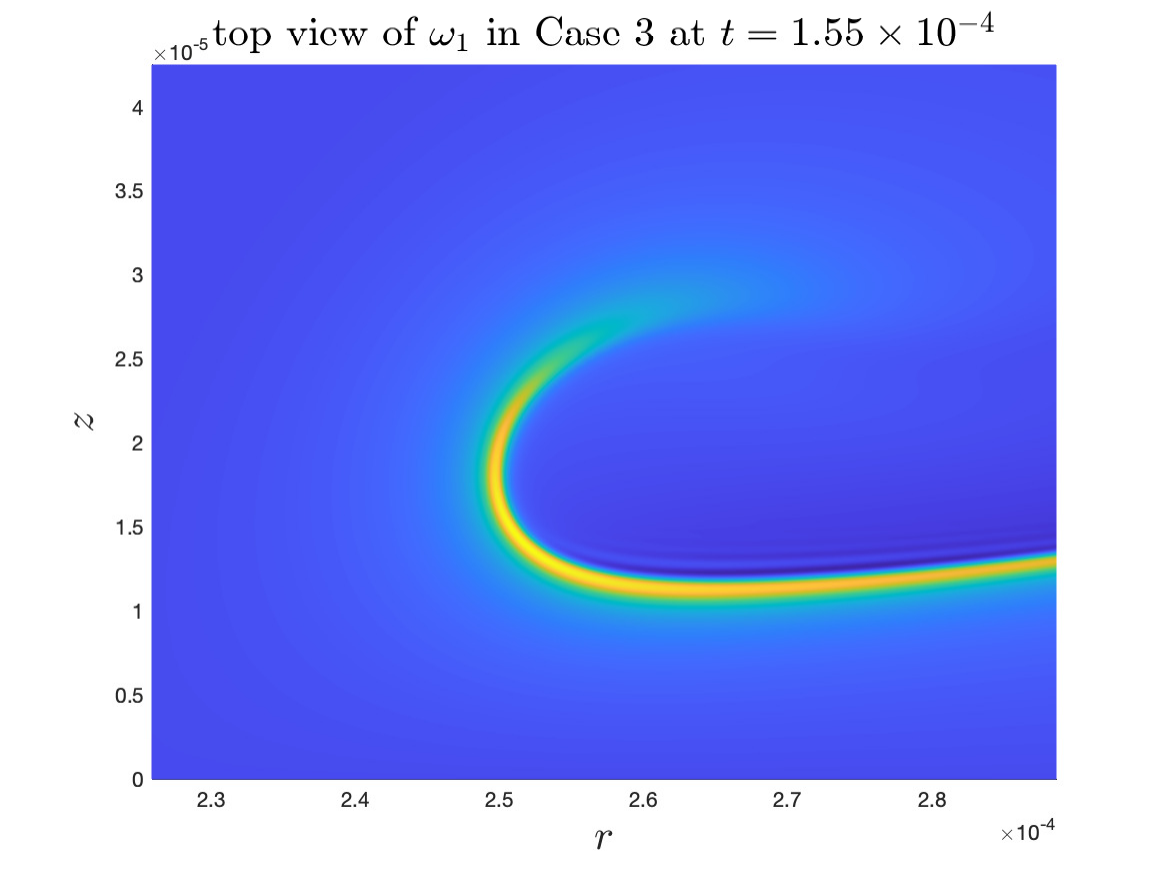}
    \caption[Compare with Euler: Zoomin]{The zoom-in profile and the top view of $\om_1$ at $t=1.55\times 10^{-4}$ in Case $1$ (first row) and in Case $3$ (second row).}  
    \label{fig:Euler_zoomin_compare}
\end{figure}

As in Case $1$, the solution in Case $3$ also demonstrates two-scale features: a long tail in the $r$ direction and a thin structure in the $z$ direction. If we zoom into the front part of the solution, we can also see local isotropic profiles that are similar to those in Case $1$. Figure \ref{fig:Euler_zoomin_compare} compares the local profiles of $\om_1$ near the front part in Case $1$ with those in Case $3$ at $t=1.55\times 10^{-4}$. Again, these profiles are qualitatively similar. However, one can see that the $\om_1$ profile in Case $3$ is much thinner at this early time, due to the absence of the regularization of the degenerate viscosity. Recall that, in Case $1$, the curved structure of $\om_1$ only becomes very thin at a much later time (see Figure \ref{fig:zoomin_profile}).

\subsection{Even faster growth} 
\begin{figure}[!ht]
\centering
    \includegraphics[width=0.32\textwidth]{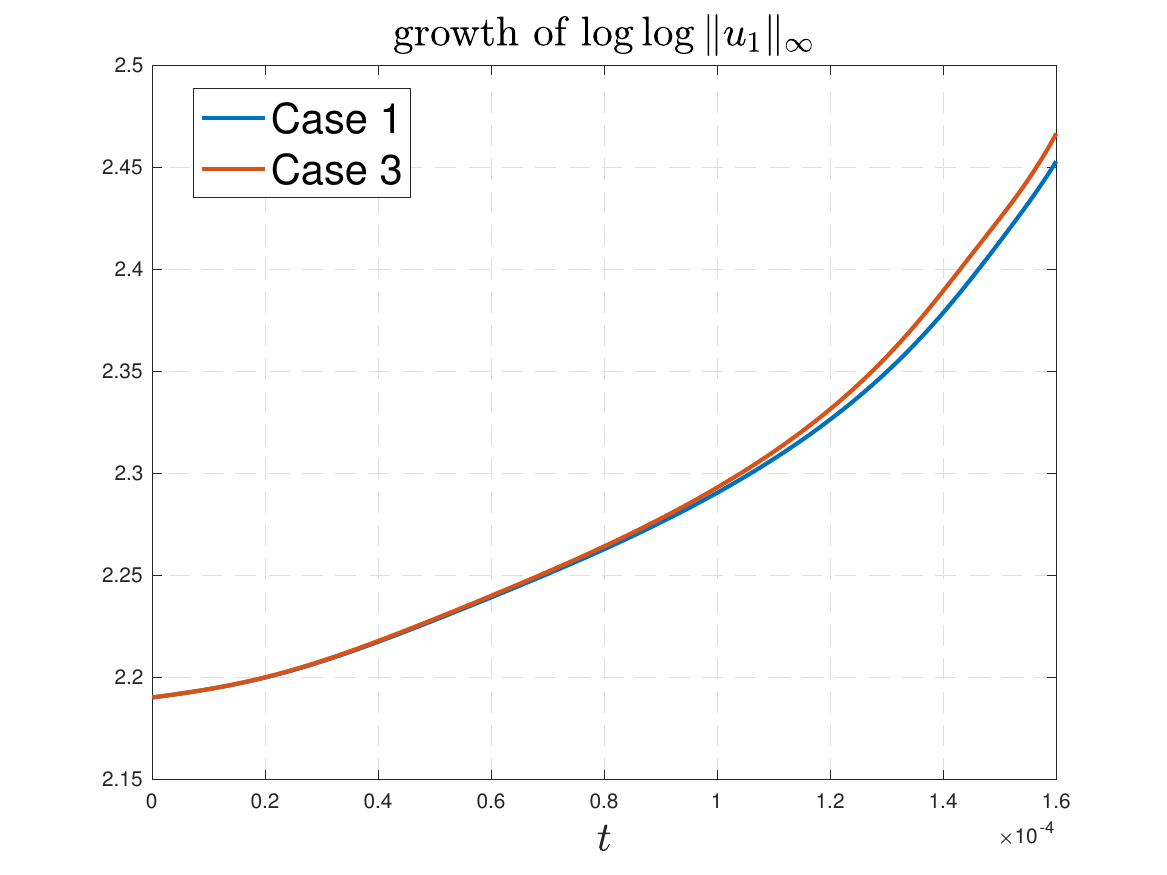}
    \includegraphics[width=0.32\textwidth]{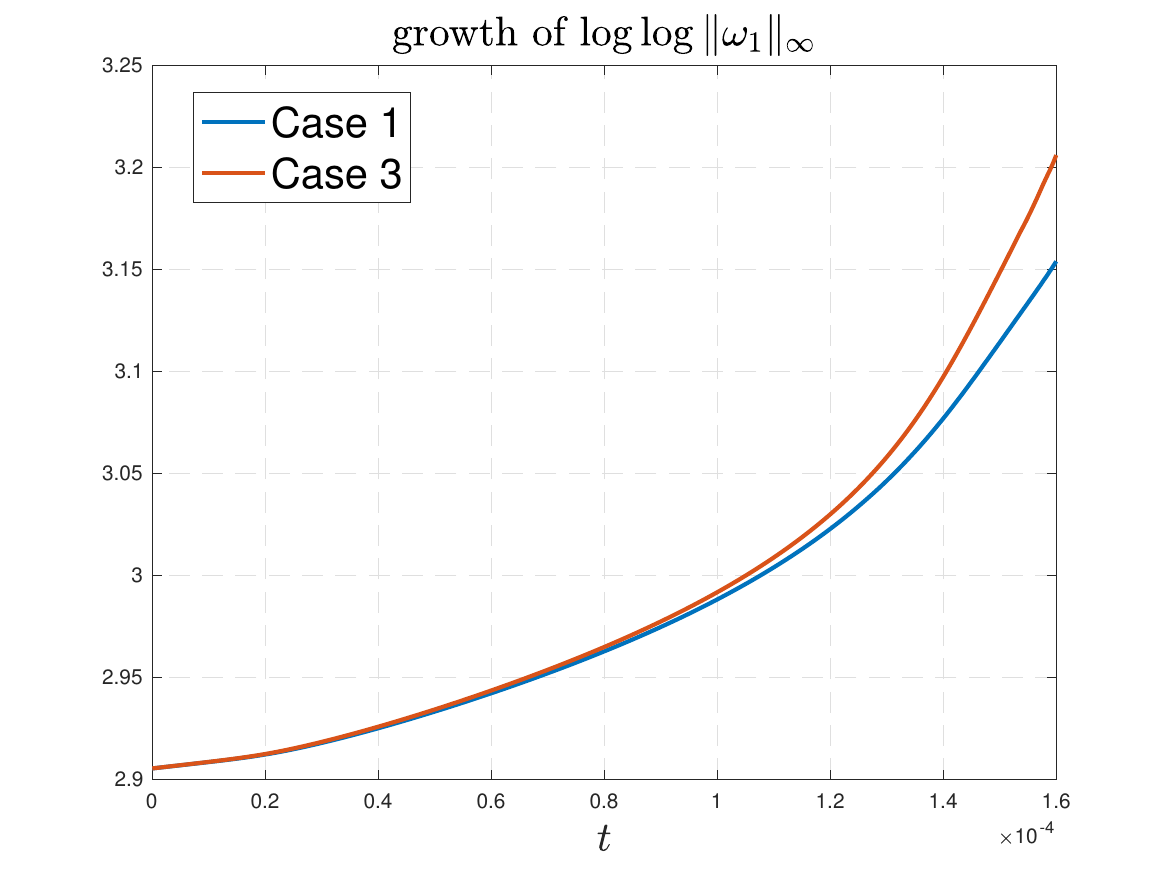}
    \includegraphics[width=0.32\textwidth]{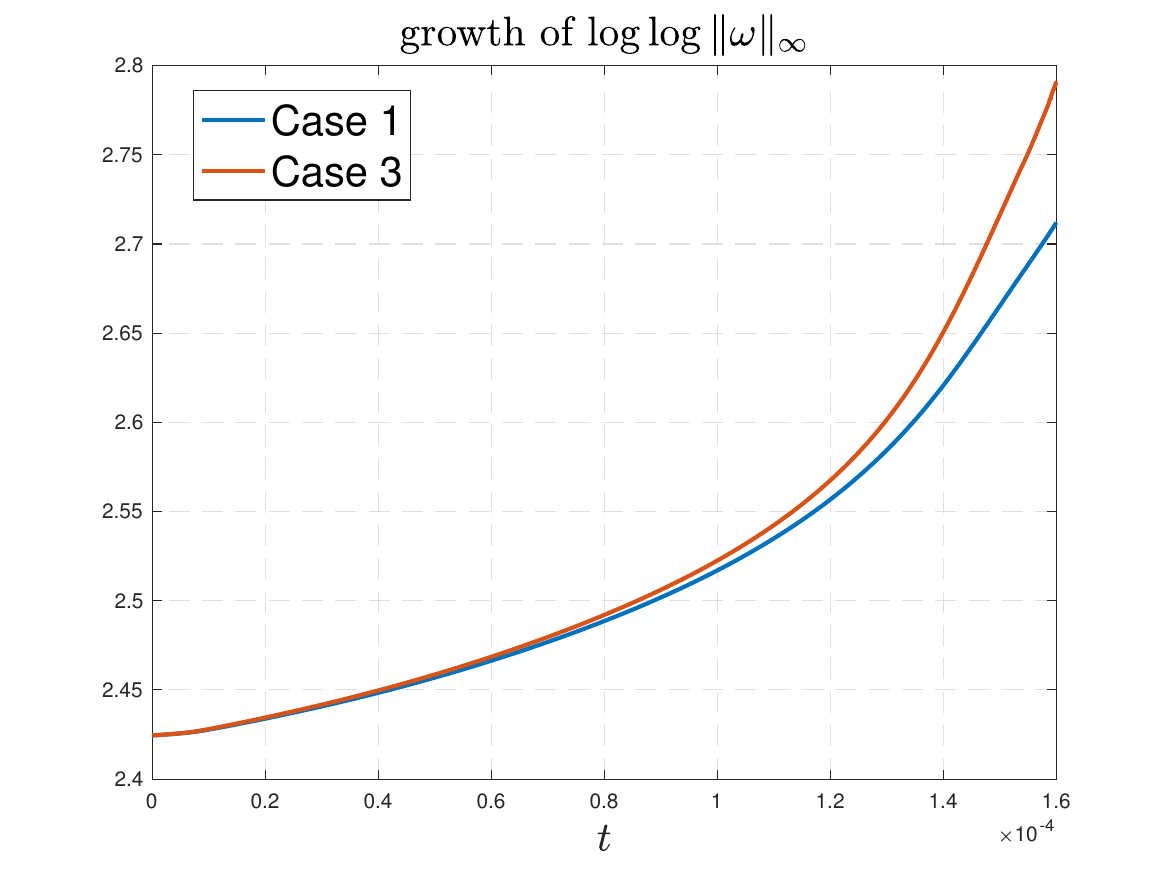}
    \caption[Even faster growth]{Comparison between the growth of $\|u_1\|_{L^\infty}$, $\|\om_1\|_{L^\infty}$ and $\|\vom\|_{L^\infty}$ in Case $1$ and Case $3$.} 
    \label{fig:rapid_growth_compare}
\end{figure}

Without the viscosity regularization, the solution in Case $3$ grows even faster than that in Case $1$. In Figure \ref{fig:rapid_growth_compare}, we compare the growth of $\|u_1\|_{L^\infty}$, $\|\om_1\|_{L^\infty}$ and $\|\vom\|_{L^\infty}$ (in double-log scale) in Case $1$ and Case $3$. The plots stop at $1.6\times 10^{-4}$ when the solution in Case $3$ is still resolved. We can see that these variables computed in Case $3$ grow faster than a double-exponential rate, even above the corresponding growth curve in Case $3$. This result implies that the solution in Case $1$ may also develop a similar blowup at a finite time in a fashion similar to that of Case $3$.

In fact, the Euler solution in Case $3$ also enjoys the critical blowup mechanism discussed in Section \ref{sec:mechanism}, which does not rely on the viscosity terms. Intuitively, the viscosity terms should slow down the blowup instead of promoting it. From this point of view, the Euler solution is more likely to blow up at a finite time.

\subsection{Under-resolution at early time}
\begin{figure}[!ht]
\centering
    \includegraphics[width=0.40\textwidth]{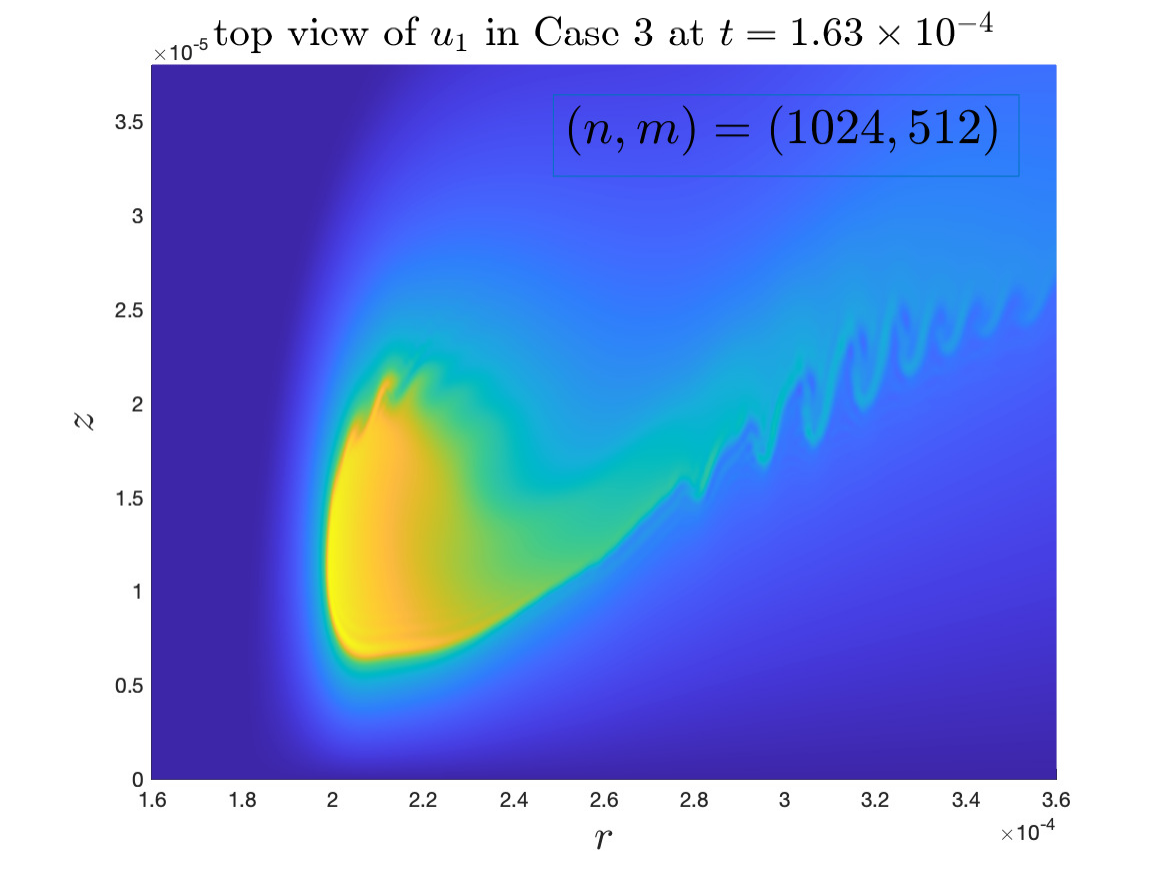} 
    \includegraphics[width=0.40\textwidth]{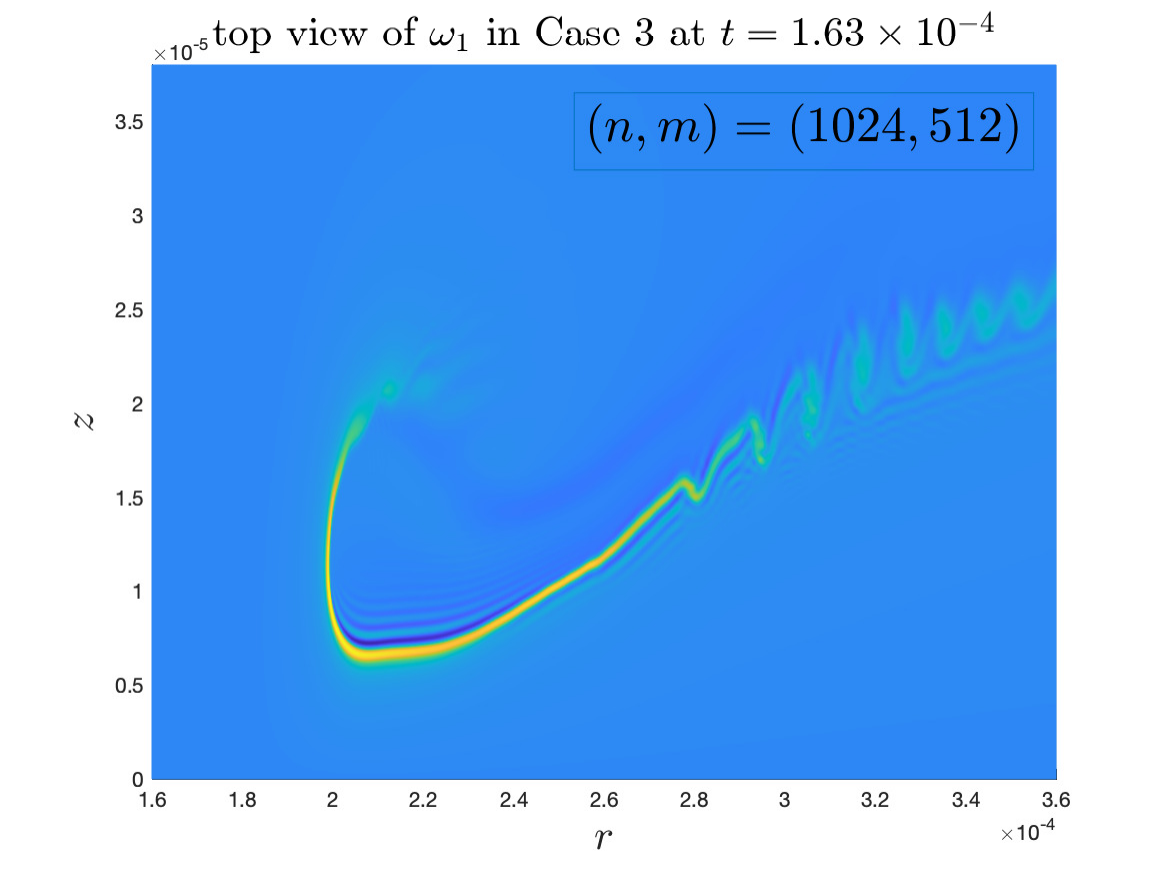}
    \includegraphics[width=0.40\textwidth]{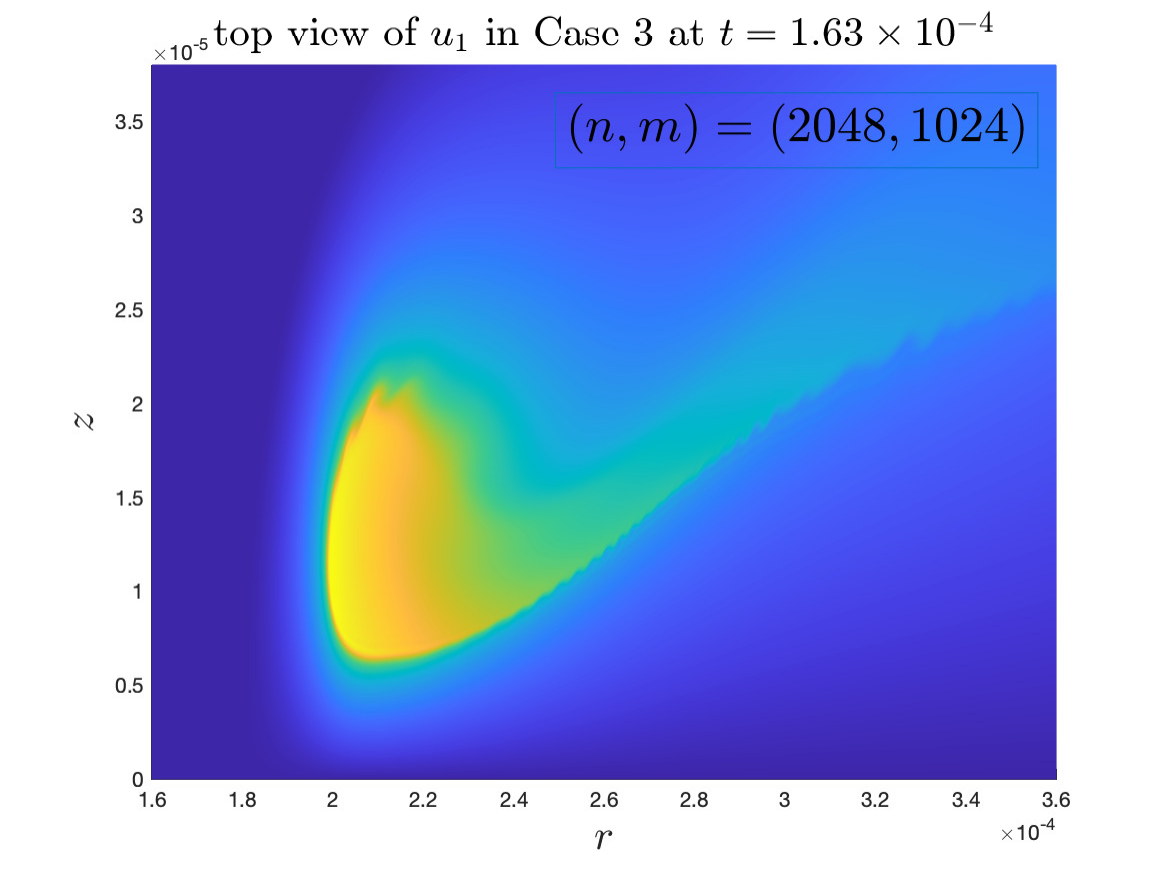} 
    \includegraphics[width=0.40\textwidth]{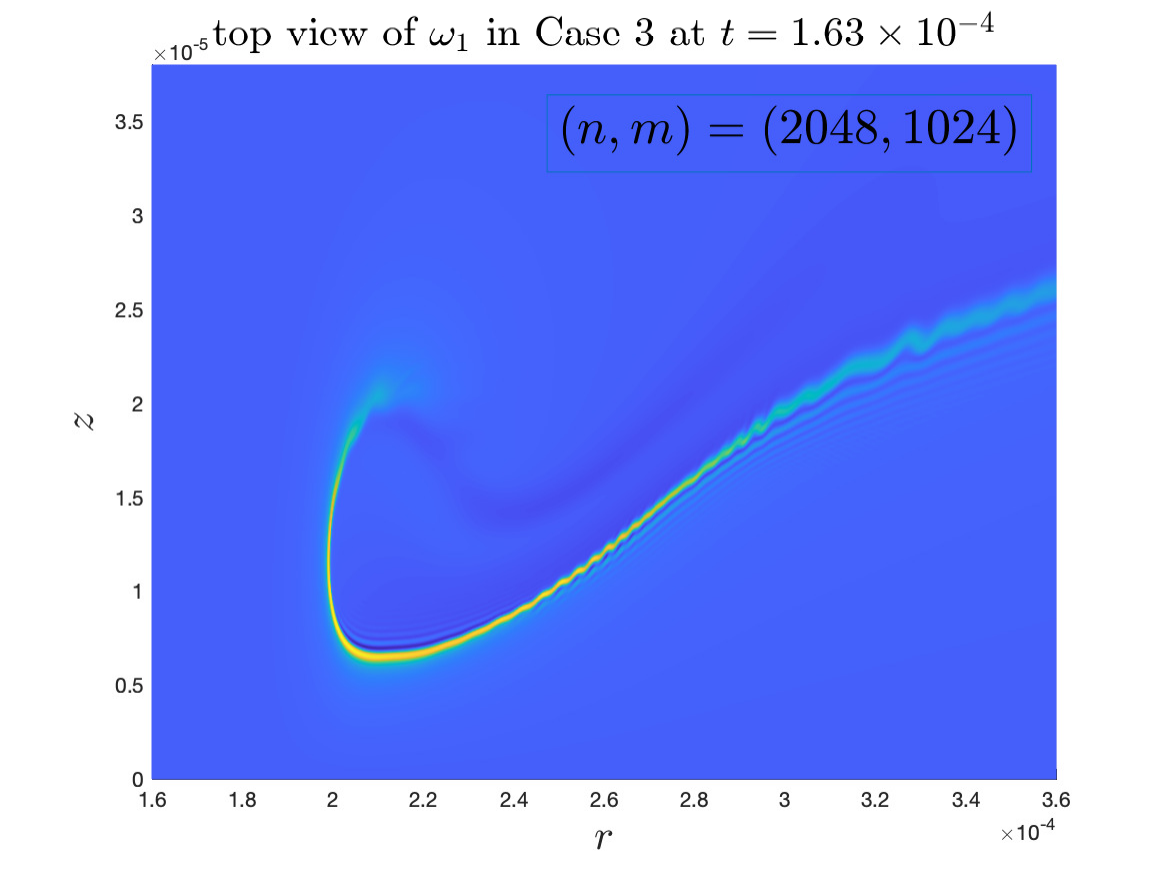}
    \caption[Euler: under-resolution]{Top views of the profiles of $u_1$ and $\om_1$ in Case $3$ at $1.63\times 10^{-4}$, computed with different resolutions. First row: $(n,m) = (1024, 512)$; second row: $(n,m) = (2048, 1024)$.  Oscillations appear all over the profiles before the stable phase can last for a long time.}  
    \label{fig:Euler_underresolution}
\end{figure}

Currently, we are not able to study thoroughly the potential blowup of the Euler solution in Case $3$ for a longer time because the solution quickly develops visibly oscillations in the critical region when it enters the stable phase (beyond $t=1.6\times 10^{-4}$). Figure \ref{fig:Euler_underresolution} shows the top views of the profiles of $u_1$ and $\om_1$ in Case $3$ at $t=1.63\times 10^{-4}$, computed with $(n,m) = (1024, 512)$ (first row) and $(n,m) = (2048, 1024)$ (second row), respectively. One can see that the oscillations appear not only in the tail part but also in the front part of the solution, which may disturb the crucial alignment between $u_1$ and $\psi_{1,z}$ near the maximum point of $u_1$. In fact, the oscillations already start to occur at an earlier time $t=1.61\times 10^{-4}$. Increasing the resolution can help suppress the oscillations (the plots in the second row of Figure \ref{fig:Euler_underresolution} are less oscillatory than those in the first row), which implies that this phenomenon is a consequence of under-resolution of the Euler solution. However, even if we use a finer mesh, the oscillations still appear quickly before we can obtain convincing numerical evidences of locally self-similar blowup at a finite time. 

\begin{figure}[!ht]
\centering
	\begin{subfigure}[b]{0.32\textwidth}
    \includegraphics[width=1\textwidth]{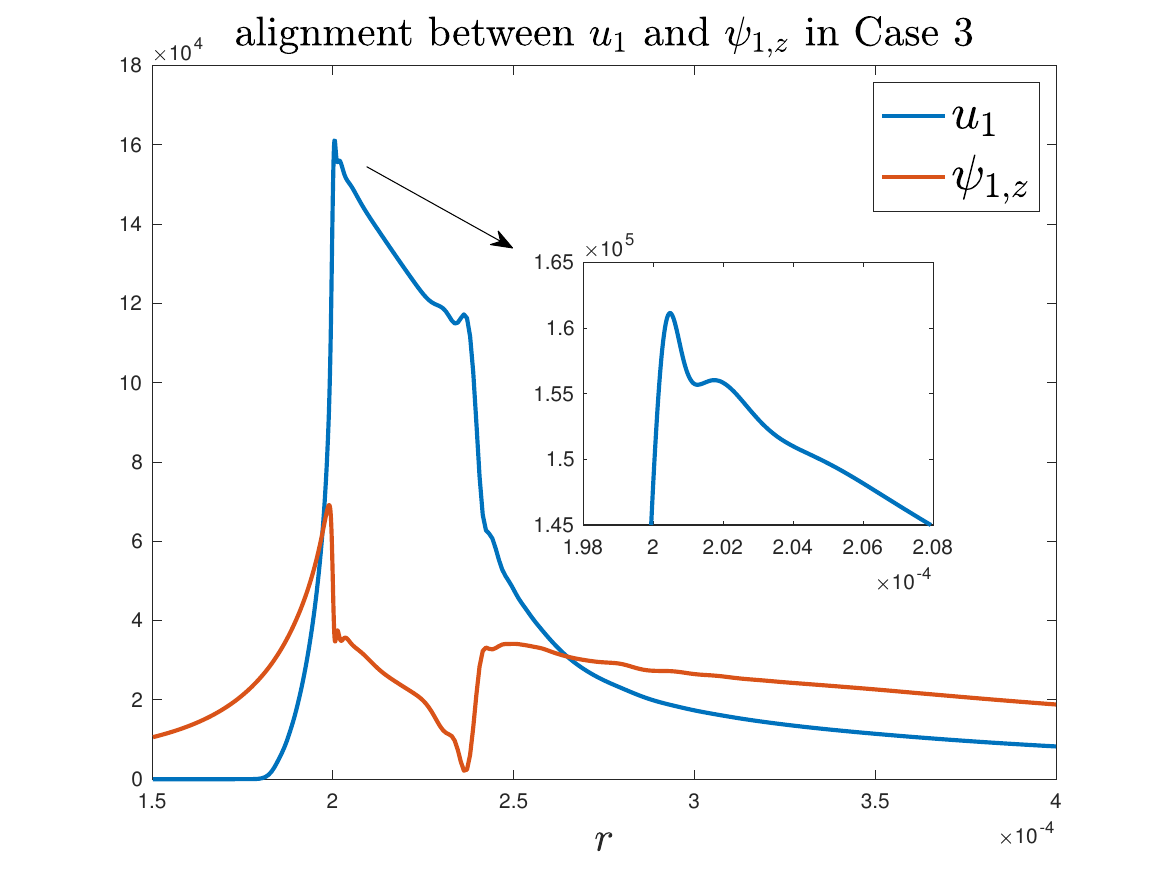}
    \caption{$r$ cross sections of $u_1,\psi_{1,z}$}
    \end{subfigure}
    \begin{subfigure}[b]{0.32\textwidth}
    \includegraphics[width=1\textwidth]{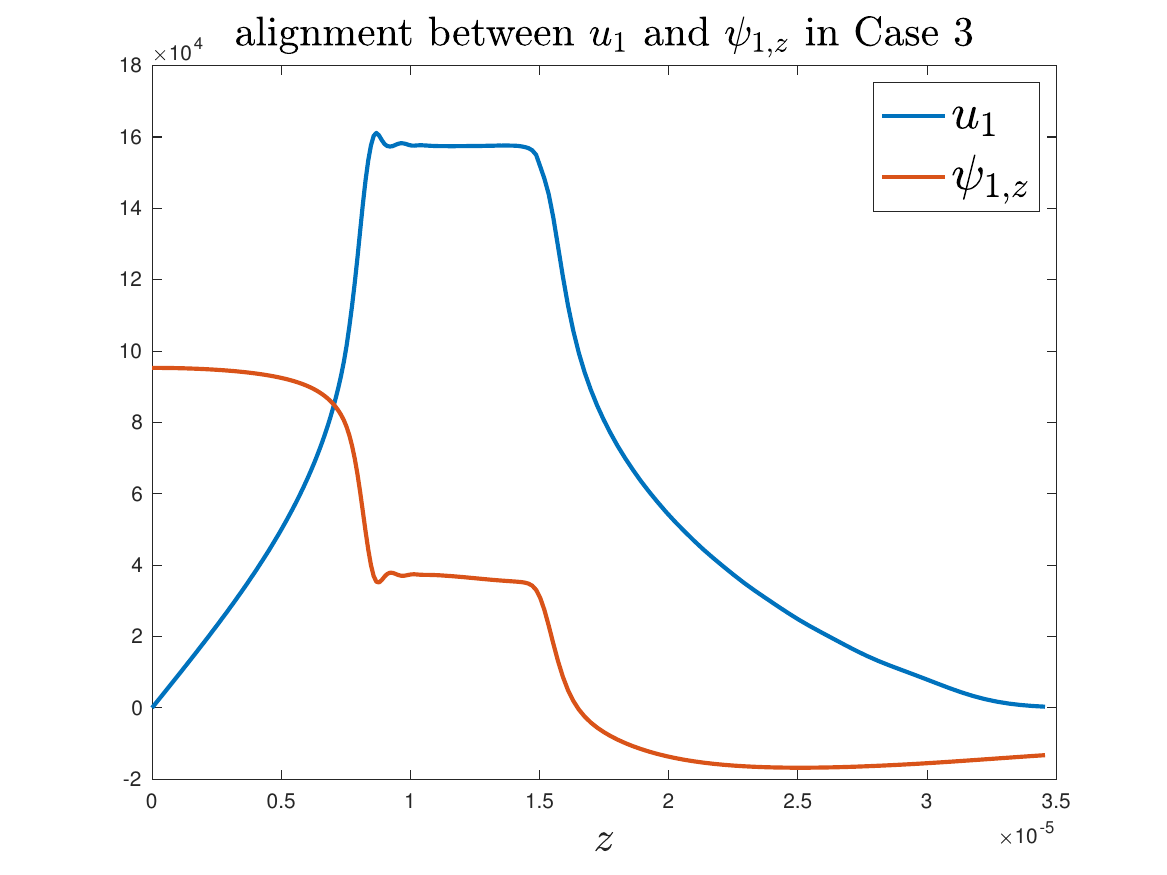}
    \caption{$z$ cross sections of $u_1,\psi_{1,z}$}
    \end{subfigure}
  	\begin{subfigure}[b]{0.32\textwidth}
    \includegraphics[width=1\textwidth]{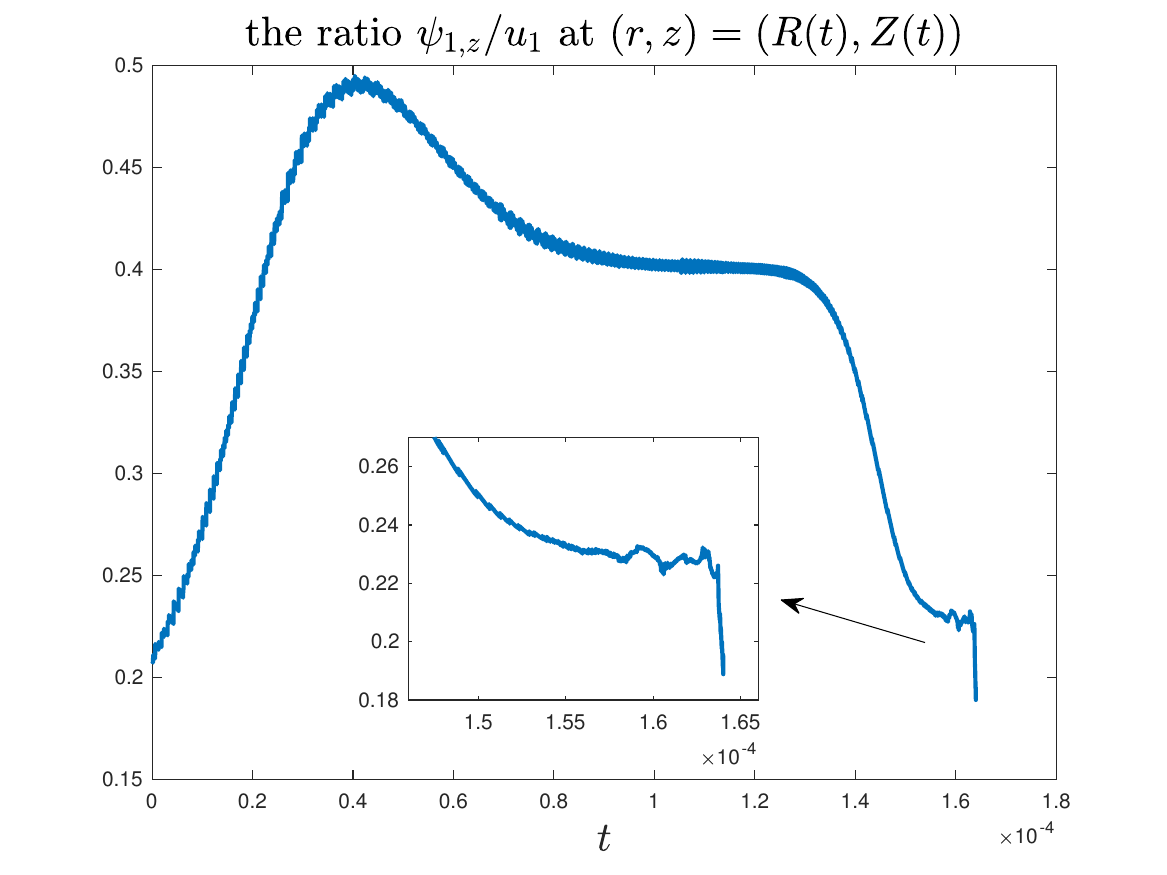}
    \caption{$\psi_{1,z}/u_1$ as a function of time}
    \end{subfigure}
    \caption[Euler: Alignment]{The alignment between $u_1$ and $\psi_{1,z}$. (a) and (b): cross sections of $u_1$ and $\psi_{1,z}$ through the point $(R(t),Z(t))$ at $t=1.63\times 10^{-4}$. (c): the ratio $\psi_{1,z}/u_1$ at the point $(R(t),Z(t))$ as a function of time up to $t=1.64\times 10^{-4}$.}  
     \label{fig:Euler_alignment}
        \vspace{-0.05in}
\end{figure}

FIn igure \ref{fig:Euler_alignment}(a) and (b), we plot the cross sections of $u_1$ and $\psi_{1,z}$ through the point $(R(t),Z(t))$ at $t=1.63\times 10^{-4}$ in Case $3$. We can see that the Euler solution also enjoys the favorable nonlinear alignment between $u_1$ and $\psi{1,z}$ near the maximum point of $u_1$ as described in Section \ref{sec:mechanism}. One should compare these plots with those in Figure \ref{fig:alignment}. However, the under-resolution of the Euler solution leads to oscillations in the front part of the $u_1$ profile, which may compromise the critical blowup mechanism. In Figure \ref{fig:Euler_alignment}(c), we plot the ratio $\psi_{1,z}/u_1$ at $(R(t),Z(t))$ against time. The alignment between $u_1$ and $\psi_{1,z}$ begins to decrease before $1.64\times 10^{-4}$ due to under-resolution.

A possible reason for the Euler solution to become under-resolved at an early time is that the local geometric structure of the solution becomes too singular to be resolved by our current adaptive mesh strategy. The front of $u_1$ is much sharper and the structure of $\om_1$ is much thinner than that in Case $3$ at the same time instant. If we treat the thickness of the thin structure of $\om_1$, denoted by $d(t)$, as an additional spatial scale, then this scale is even smaller than the scale of $Z(t)$. That is, the Euler solution demonstrates three separate spatial scales $d(t), Z(t), R(t)$ (from small to large), each converges to $0$ at a different rate. However, our adaptive mesh strategy is only powerful enough to handle the high contrast between two separate scales in the critical region of the solution over the stable phase. The three-scale feature of the Euler solution in Case $3$ is beyond our current computational capacity. Moreover, the thin $1$D-like structure of $\om_1$ induces strong shearing instabilities that will amplify the errors from under-resolution and lead to the visible rolling oscillations. 

In summary, the Euler solution in Case $3$ quickly develops an even more singular structure that is extremely difficult to resolve with our current computational capacity. This is why we adopt the degenerate viscosity coefficients in our main Case $1$: the degenerate viscosity is strong enough to prevent the occurrence of a third scale but also not too strong to suppress the two-scale blowup. We believe that the Euler solution may develop a locally self-similar blowup as in Case $1$. To obtain a convincing numerical evidence of a potential $3$D Euler blowup, we need to develop a more effective adaptive mesh strategy and have access to larger computational resources. We will leave this to our future work.

\section{Potential Blowup Scaling Analysis}\label{sec:scaling_study}
In this section, we will quantitatively examine the features of the potential blowup in our computation. We will first provide adequate numerical evidences that the growth and the spatial scaling of the solution obey some (inverse) power laws, which suggest that a finite-time singularity exists in a locally self-similar form. In particular, we employ a linear fitting procedure to estimate the blowup rates and scalings. Then we will perform an asymptotic analysis of the potential blowup based on a two-scale self-similar ansatz. We show that the results of the asymptotic analysis are highly consistent with our numerical results, supporting the existence of a locally self-similar blowup.

\subsection{Linear fitting procedure}\label{sec:fitting_procedure} The most straightforward way to numerically identify a finite-time blowup is to study the growth rate of the solution. For a solution quantity $v(t)$ that is expected to blow up at some finite time $T$, a typical asymptotic model is the inverse power law:
\begin{equation}\label{eq:inverse_power_law}
v(t) \sim \alpha (T-t)^{-c_v}\quad \text{as}\quad t\rightarrow T,
\end{equation}
where $c_v>0$ is the blowup rate and $\alpha>0$ is some constant. To verify that $v(t)$ satisfies an inverse power law and to learn the power $c_v$, we follow the idea of Luo--Hou \cite{luo2014toward} and study the time derivative of the logarithm:
\[\frac{\diff\,}{\diff t}\log v(t) = \frac{v'(t)}{v(t)} \sim \frac{c_v}{T-t}.\]
This naturally leads to the linear regression model
\begin{equation}\label{eq:regression_1}
y(t;v) := \frac{v(t)}{v'(t)} \sim -\frac{1}{c_v}(t-T) =: \tilde{a}t + \tilde{b},
\end{equation}
with response variable $y$, explanatory variable $t$, and model parameters $\tilde{a} = -1/c_v,\tilde{b}=T/c_v$. Then the blowup rate $c_v$ can be estimated via a standard least-squares procedure. The quality of the fitting using this model can be measured by the coefficient of determination (the $R^2$):
\[R^2 = 1 - \frac{SS_{\text{err}}}{SS_{\text{tot}}},\]
with a value close to $1$ indicating a high quality fitting. Here 
$SS_{\text{tot}} = \sum_{i} (y_i - \bar{y})^2$
is the total sum of squares, 
$SS_{\text{err}} = \sum_{i} (y_i - \hat{y}_i)^2$
is the residual sum of squares, $y_i,\hat{y}_i$ denote the observed and predicted values of the response
variable $y$, respectively, and $\bar{y}$ denotes the mean of the observed data $y_i$.

To have a convincing estimate of the blowup rate $c_v$, it is important that the fitting procedure is performed in a proper time interval $[t_1,t_2]$. First of all, this time interval must lie in the asymptotic regime of the inverse power law \eqref{eq:inverse_power_law} if such scaling does exist. Secondly, the solution must be well resolved in this time interval $[t_1,t_2]$. As we have observed, the blowup settles down to a stable phase at around $t = 1.6\times 10^{-4}$, after which the evolution of the solution begins to have a stable pattern. It is likely that the solution enters the asymptotic regime of the blowup after this time instant. In addition, we have learned in Section \ref{sec:resolution_study} that the numerical solution is well resolved before $t = 1.76\times 10^{-4}$. Therefore, according to the two criteria, we should place the fitting interval $[t_1,t_2]$ within the time interval $[1.6\times10^{-4}, 1.76\time10^{-4}]$. Moreover, the interval should not to be too short; otherwise, any curve may look like a straight line. In particular, we choose $[t_1,t_2]=[1.60\times10^{-4}, 1.75\times10^{-4}]$. We will denote by $\tilde{c}_v = -1/\tilde{a},\tilde{T}_v = -\tilde{b}/\tilde{a}$ the approximate blowup rate and blowup time obtained from this fitting procedure.

Since the quantities for which we would like to obtain the potential blowup rates are mostly the $L^\infty$ norms of some solution functions, their values are sensitive to the discretization methods, the choice of the adaptive mesh, and the interpolation operations, especially when the maximum points are traveling as in our scenario. Therefore, the model \eqref{eq:regression_1} may not yield an ideal fitting even if the inverse power law \eqref{eq:inverse_power_law} does exist, and the resulting $\tilde{c}_v$ may not reflect the true blowup rate $c_v$, though it should still be a good approximation. To obtain a better approximation of $c_v$, we will perform a local search near the crude estimate $\tilde{c}_v$ and find a value $\bar{c}_v$ such that the model 
\begin{equation}\label{eq:regression_2}
\gamma(t;v) := v(t)^{-1/\bar{c}_v} \sim \alpha^{-1/\bar{c}_v}(T - t)^{c_v/\bar{c}_v} \sim \alpha^{-1/c_v}(T - t) =: \bar{a} t + \bar{b}
\end{equation}
has the best linear regression fitness (the $R^2$) against a least-square test. More precisely, we put a uniform mesh (with $100$ points) over the interval $[\tilde{c}_v-0.1,\tilde{c}_v+0.1]$, compute the $R^2$ value of the model \eqref{eq:regression_2} for each candidate blowup rate on the mesh, and find the blowup rate $\bar{c}_v$ such that the model \eqref{eq:regression_2} has the optimal fitness over all candidates. If $\bar{c}_v$ falls into one of the endpoints of the interval, i.e. $\tilde{c}_v\pm0.1$, we will perform a local search again around $\bar{c}_v$ and update the value of $\bar{c}_v$. After this procedure, the resulting $\bar{c}_v$ should be a better approximation of $c_v$. Corresponding, $\bar{T} := -\bar{b}/\bar{a}$ is an approximation of the true blowup time $T$. Note that for the fitting of the model \eqref{eq:regression_2}, we use the original recorded values of $v(t)$ rather than the time-interpolated values. We remark that the fitting of model \eqref{eq:regression_2} is more reliable than the fitting of model \eqref{eq:regression_1} in reflecting a potential inverse power law, as it is directly applied to the quantity of interest without taking the time derivative of the logarithm of this quantity. It is much harder for the logarithm of a blowup quantity to fall into the asymptotic regime in comparison with the blowup quantity itself. In other words, the fitting based on \eqref{eq:regression_2} is a refinement of the fitting based on \eqref{eq:regression_1}.

\subsection{Fitting of the growth rate}\label{sec:growth_fitting}
\begin{figure}[!ht]
\centering
    \includegraphics[width=0.4\textwidth]{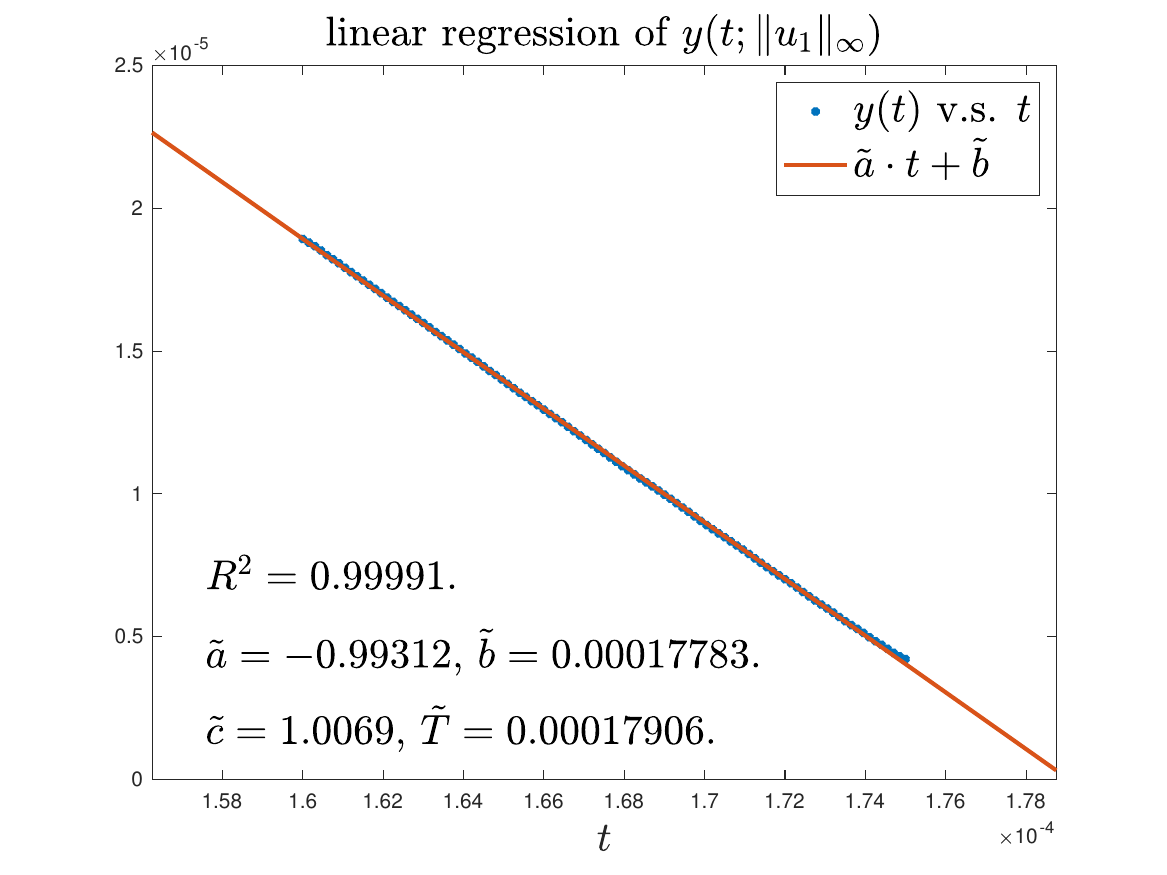}
    \includegraphics[width=0.4\textwidth]{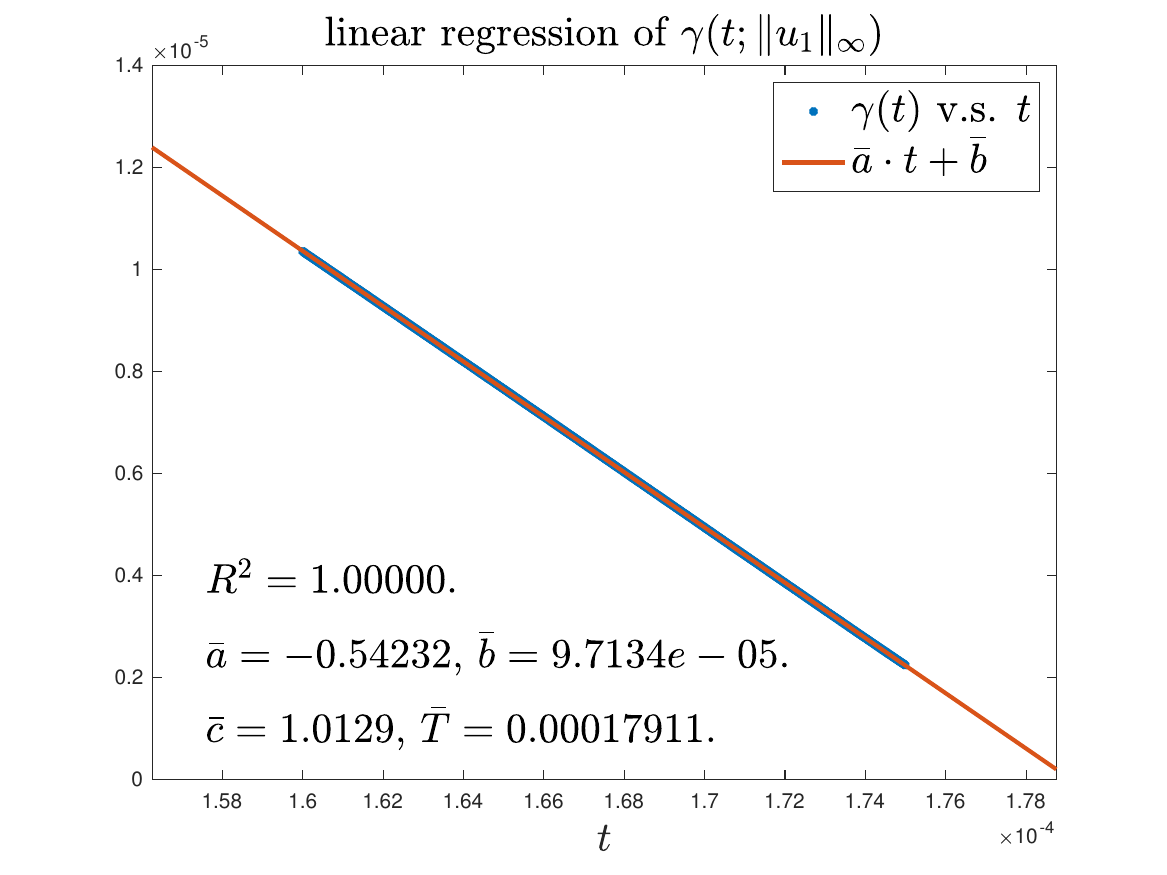}
    \caption[Linear regression $u_1$]{The linear regression of $y(t;\|u_1\|_{L^\infty})$ (left) and $\gamma(t;\|u_1(t)\|_{L^\infty})$ (right) on the time interval $[t_1,t_2] = [1.65\times 10^{-4},1.75\times 10^{-4}]$. The blue points are the data points obtained from our computation, and the red lines are the linear models. We plot the linear models on a larger interval.}   
    \label{fig:linear_regression_u1}
       \vspace{-0.05in}
\end{figure}

We are now ready to apply the above fitting procedures to the numerical solutions obtained in our computation. Figure\eqref{fig:linear_regression_u1} shows the fitting results for the quantity $\|u_1(t)\|_{L^\infty}$ (in Case $1$) on the time interval $[t_1,t_2] = [1.6\times 10^{-4},1.75\times10^{-4}]$. We can see that both models, $y(t;\|u_1\|_{L^\infty})$ and $\gamma(t;\|u_1\|_{L^\infty})$, have excellent linear fitness with $R^2$ values very close to $1$ (the $R^2$ value for the model $\gamma(t;\|u_1\|_{L^\infty})$ is greater than $1-10^{-6}$). Note that the blowup rates (or the blowup times) learned from the two models are close to each other, cross-verifying the validity of both models. This strongly implies that $\|u_1\|_{L^\infty}$ has a finite-time singularity of an inverse power law with a blowup rate very close to $1$, which is consistent with our observation and analysis in Section \ref{sec:mechanism}. Recall that we have observed a strong positive alignment between $\psi_{1,z}$ and $u_1$ around the maximum point $R(t),Z(t)$ of $u_1$. In particular, $\psi_{1,z}(R(t),Z(t),t)\sim u_1(R(t),Z(t),t) $ in the stable phase $[1.6\times 10^{-4},1.75\times10^{-4}]$. If we ignore the degenerate viscosity, then the equation of $\|u_1(t)\|_{L^\infty}$ can be approximated by 
\[\frac{\diff \,}{\diff t} \|u_1(t)\|_{L^\infty} = 2\psi_{1,z}(R(t),Z(t),t)\cdot u_1(R(t),Z(t),t) \sim c_0 \|u_1(t)\|_{L^\infty}^2,\]
which then implies that 
$\|u_1(t)\|_{L^\infty} \sim  (T-t)^{-1} $
for some finite time $T$. This asymptotic analysis is now supported by our linear fitting results.

\begin{figure}[!ht]
\centering
    \includegraphics[width=0.4\textwidth]{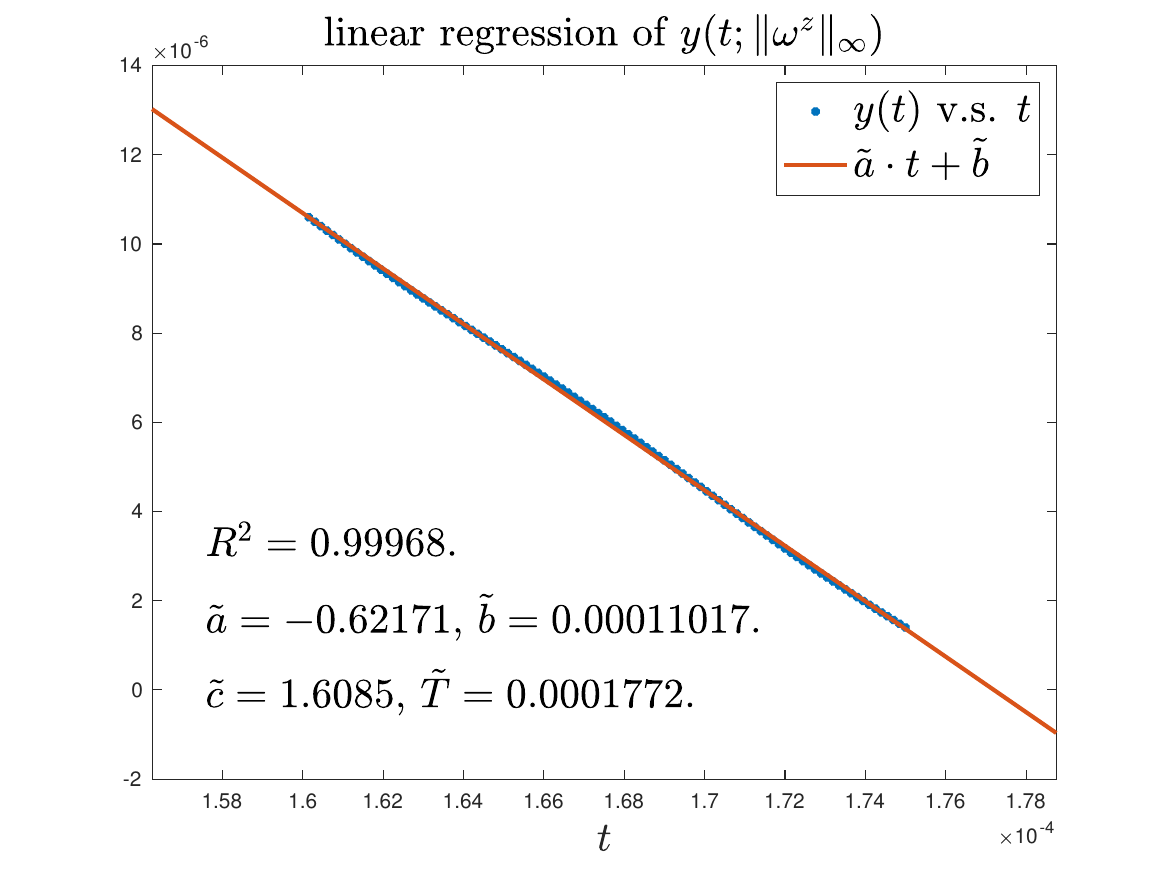}
    \includegraphics[width=0.4\textwidth]{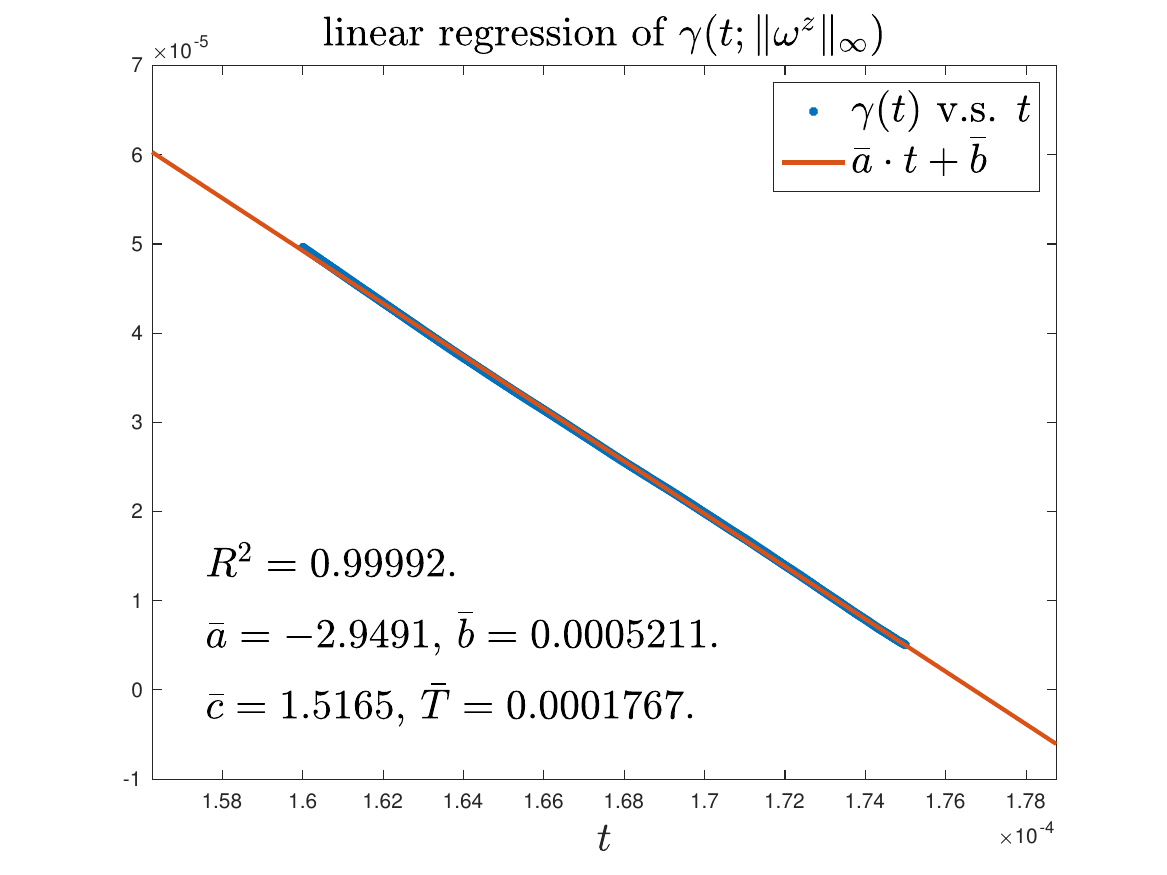}
    \caption[Linear regression $\vom$]{The linear regression of $y(t;\|\om^z\|_{L^\infty})$ (left) and $\gamma(t;\|\om^z\|_{L^\infty})$ (right) on the time interval $[t_1,t_2] = [1.65\times 10^{-4},1.75\times 10^{-4}]$. The blue points are the data points obtained from our computation, and the red lines are the linear models. We plot the linear models on a larger interval.}   
    \label{fig:linear_regression_omegaz}
       \vspace{-0.05in}
\end{figure}

Next, we study the growth of the maximum vorticity $\|\vom\|_{L^\infty}$, which is an important indicator of a finite-time singularity. However, the growth of $\|\vom\|_{L^\infty}$ may not reflect a clean inverse power law, since the maximum point of the vector magnitude $|\vom|$ and the maximum points of the components $\om^\theta,\om^r,\om^z$ are distinct in general. Therefore, we directly apply the fitting procedure to the maximums of the vorticity components instead of to the maximum vorticity. As an illustration, we present the study of the inverse power law of the axial vorticity component $\om^z$. Figure \ref{fig:linear_regression_omegaz} shows the linear fitting for the associated models of $\|\om^z(t)\|_{L^\infty}$ (in Case $1$) on the time interval $[t_1,t_2] = [1.6\times 10^{-4},1.75\times10^{-4}]$. Again, we see that both models $y(t;\|\om^z\|_{L^\infty})$ and $\gamma(t;\|\om^z\|_{L^\infty})$ have good linear fitness, which provides evidences of the finite-time blow of $\|\om^z(t)\|_{L^\infty}$ in the form of an inverse power law
\[\|\om^z(t)\|_{L^\infty} \sim (T-t)^{-\bar{c}_{\om^z}}.\] 
More importantly, the estimated blowup rate is approximately equal to $1.5$, which implies that 
\[\int_0^T\|\vom(t)\|_{L^\infty} \geq \int_0^T\|\om^z(t)\|_{L^\infty}  \idiff t = +\infty.\]
According to our discussion in Section \ref{sec:rapid_growth}, this strongly supports the existence of a finite-time singularity of the solution to the initial-boundary value problem \eqref{eq:axisymmetric_NSE_1}--\eqref{eq:initial_data}. 

To further illustrate the existence of a potential finite-time blowup, we perform the linear fitting procedure on more blowup quantities computed with different mesh resolutions. We only present the fitting results of model \eqref{eq:regression_2}. Table \ref{tab:fitting_resolution} reports the linear fitting results of different solution quantities computed in Case $1$ with different mesh sizes. Again, the fitting time interval is $[t_1,t_2] = [1.6\times 10^{-4},1.75\times 10^{-4}]$. We can see that all the quantities reported in the table have excellent fitting to some inverse power law, and the fitting results are consistent over different resolutions.

\begin{table}[!ht]
\centering
\footnotesize
\renewcommand{\arraystretch}{1.5}
    \begin{tabular}{|c|c|c|c|c|c|c|c|c|c|}
    \hline
    \multirow{2}{*}{Mesh size} & \multicolumn{3}{c|}{$\|u_1\|_{L^\infty}$} & \multicolumn{3}{c|}{$\|\om_1\|_{L^\infty}$}  & \multicolumn{3}{c|}{$\|\psi_{1,r}\|_{L^\infty}$}\\ \cline{2-10} 
    & $\bar{c}$ & $\bar{T}/10^{-4}$ & $R^2$ & $\bar{c}$ & $\bar{T}/10^{-4}$ & $R^2$ & $\bar{c}$ & $\bar{T}/10^{-4}$ & $R^2$ \\ \hline 
    $1024\times 512$ & $1.0129$ & $1.7911$ & $1.00000$ & $1.9773$ & $1.7908$ & $1.00000$ & $1.1154$ & $1.7998$ & $0.99999$\\ \hline
    $1536\times 768$ & $1.0126$ & $1.7914$ & $1.00000$ & $2.0366$ & $1.7966$ & $1.00000$ & $1.1151$ & $1.7975$ & $0.99999$\\ \hline
    $2048\times 1024$ & $1.0125$ & $1.7916$ & $1.00000$ & $2.0217$ & $1.7956$ & $1.00000$ & $1.1129$ & $1.7966$ & $1.00000$ \\ \hline
    \hline
    \multirow{2}{*}{Mesh size} & \multicolumn{3}{c|}{$\|\psi_{1,z}\|_{L^\infty}$} & \multicolumn{3}{c|}{$\|u_{1,r}\|_{L^\infty}$}  & \multicolumn{3}{c|}{$\|u_{1,z}\|_{L^\infty}$}\\ \cline{2-10} 
    & $\bar{c}$ & $\bar{T}/10^{-4}$ & $R^2$ & $\bar{c}$ & $\bar{T}/10^{-4}$ & $R^2$ & $\bar{c}$ & $\bar{T}/10^{-4}$ & $R^2$ \\ \hline 
    $1024\times 512$ & $0.9744$ & $1.7952$ & $0.99990$ & $2.0427$ & $1.7712$ & $0.99993$ & $1.9752$ & $1.7687$ & $0.99998$\\ \hline
    $1536\times 768$ & $0.9730$ & $1.7954$ & $0.99990$ & $2.0445$ & $1.7718$ & $0.99993$ & $1.9690$ & $1.7690$ & $0.99998$\\ \hline
    $2048\times 1024$ & $0.9724$ & $1.7954$ & $0.99991$ & $2.0438$ & $1.7719$ & $0.99993$ & $1.9666$ & $1.7692$ & $0.99998$ \\ \hline
    \hline
    \multirow{2}{*}{Mesh size} & \multicolumn{3}{c|}{$\|\om^\theta\|_{L^\infty}$} & \multicolumn{3}{c|}{$\|\om^r\|_{L^\infty}$}  & \multicolumn{3}{c|}{$\|\om^z\|_{L^\infty}$}\\ \cline{2-10} 
    & $\bar{c}$ & $\bar{T}/10^{-4}$ & $R^2$ & $\bar{c}$ & $\bar{T}/10^{-4}$ & $R^2$ & $\bar{c}$ & $\bar{T}/10^{-4}$ & $R^2$ \\ \hline 
    $1024\times 512$ & $1.5045$ & $1.7950$ & $1.00000$ & $1.4878$ & $1.7652$ & $0.99997$ & $1.5165$ & $1.7670$ & $0.99992$\\ \hline
    $1536\times 768$ & $1.5407$ & $1.7937$ & $1.00000$ & $1.4912$ & $1.7658$ & $0.99997$ & $1.5108$ & $1.7648$ & $0.99992$\\ \hline
    $2048\times 1024$ & $1.5270$ & $1.7911$ & $1.00000$ & $1.4924$ & $1.7661$ & $0.99997$ & $1.5097$ & $1.7623$ & $0.99993$ \\ \hline
    \end{tabular}
    \caption{\small Linear fitting results of model \eqref{eq:regression_2} on time interval $[t_1,t_2] = [1.6\times 10^{-4},1.75\times 10^{-4}]$ for different solution quantities computed with different mesh sizes.}
    \label{tab:fitting_resolution}
    \vspace{-0.2in}
\end{table}

However, the estimated blowup times obtained from the fitting of different blowup quantities agree only up to the third digit. This may be due to the fact that the estimated blowup time is an extrapolated information from the fitting model, which can be very sensitive to the recorded observations and the model parameters. In fact, if we force the blowup rate of $\|u_1\|_{L^\infty}$ to be $\bar{c} = 0.9$ (instead of $\bar{c} = 1.0129$) in model \eqref{eq:regression_2}, then it still yields a fairly good linear regression with fitness $R^2 = 0.99945$, but the estimated blowup time drops to $\bar{T} = 1.7802$ (compared to $1.7911$). Therefore, the estimated blowup time may not be a robust approximation of the true blowup time even if the finite-time singularity does exist. It is the estimated blowup rate that is more interesting in our analysis. 

For the blowup rates reported in Table \ref{tab:fitting_resolution}, we have the following observations and discussions. First, the fittings of $\|u_1\|_{L^\infty},\|\psi_{1,z}\|_{L^\infty}$ again confirm the strong alignment between $u_1$ and $\psi_{1,z}$. As we can see, the blowup rates of $\psi_{1,z}$ and $u_1$ are both close to $1$, which is consistent with our observation that $\psi_{1,z}\sim u_1$ in the critical region around $(R(t),Z(t))$. This result supports our discussions on the blowup mechanism in Section \ref{sec:mechanism}.

Moreover, there seems to be a pattern in the blowup rates of the quantities reported in Table \ref{tab:fitting_resolution}. More precisely, we have the following approximation that agrees with the fitting results up to the first two digits:
\begin{equation}\label{eq:pattern}
\bar{c}_{u_1} \approx \bar{c}_{\psi_{1,r}} \approx \bar{c}_{\psi_{1,z}} \approx 1,\quad \bar{c}_{\om^\theta} \approx \bar{c}_{\om^r} \approx \bar{c}_{\om^z} \approx 1.5,\quad \bar{c}_{\om_1} \approx \bar{c}_{u_{1,r}} \approx \bar{c}_{u_{1,z}} \approx 2, 
\end{equation}
where we us $\bar{c}_f$ to denote the estimated blowup rate of $\|f(t)\|_{L^\infty}$ for a function $f(r,z,t)$. In the next section, we will argue that this pattern implies the existence of power laws for the two spatial scales of the solution, which is consistent to the direct linear fitting of these scales. Furthermore, this pattern reflects the possibility that the solution blows up in a locally self-similar way.

\subsection{Fitting of spatial scalings}\label{sec:scale_fitting}
Recall that we have observed a two-scale property of the solution in our scenario: the smaller scale (featured by $Z(t)$) measures the length scale of the local solution profile near the point $(R(t),Z(t))$ (or the sharp front), and the larger scale (featured by $R(t)$) measures the distance between the sharp front and the symmetry axis $r=0$. Our numerical observations suggest that these two scales are separated and both converge to $0$, which characterizes the focusing nature of the blowup. The next step is to quantitatively investigate how fast these two scales go down to $0$. Just like how we characterize the growth of the solution, we assume that these two spatial scales of the solution also admit some power laws:
\begin{equation}\label{eq:power_law}
C_s(t) = (T-t)^{c_s},\quad C_l(t) = (T-t)^{c_l}, 
\end{equation}
where $C_s$ denotes the larger scale (with a smaller power $c_s>0$) and $C_l$ the smaller scale (with a larger power $c_l>0$). We will use two methods to check the validity of the power laws \eqref{eq:power_law} and to learn the powers $c_s,c_l$, and we will compare the results from both methods to see if they are consistent.   

The first method to learn the scalings of $C_s,C_l$ is to extract the spatial information of the solution from the growth of different quantities. In fact, if the blowup solution develops a local profile of a isotropic spatial scale $C_l(t)$, it should be reflected by the spatial derivatives. More precisely, for a blowup function $f(r,z,t)$ that is smooth with respect to the scale $C_l(t)$, we should have
\[C_l(t) \sim \frac{\|f(t)\|_{L^\infty}}{\|f_r(t)\|_{L^\infty}} \sim \frac{\|f(t)\|_{L^\infty}}{\|f_z(t)\|_{L^\infty}},\]
which is equivalent to the relation 
\[c_l = c_{f_r} - c_f = c_{f_z} - c_f,\]
if $f(t)$ also admits an inverse power law \eqref{eq:inverse_power_law}. We have seen in Section \ref{sec:two_scale} that the profile of $u_1$ is smooth when observed in a local neighborhood around $(R(t),Z(t))$, so we may use $u_1$ to extract the scaling of $C_l$. Our data in the previous subsection, particularly the pattern in \eqref{eq:pattern}, show that 
\[\bar{c}_{u_{1,r}} - \bar{c}_{u_1} \approx \bar{c}_{u_{1,z}} - \bar{c}_{u_1} \approx 1,\]
which suggests that the solution has a local spatial scale 
\begin{equation}\label{eq:Cl_power}
C_l(t) \sim Z(t) \sim T-t.
\end{equation}

Similarly, we can learn the power of the larger scale $C_s$ from our fitting of the growth of different blowup quantities. By the definition of $u_1,\om_1$ and $\vom$, we have 
\begin{equation}\label{eq:vorticity_formula}
\om^\theta = r \om_1,\quad \om^r =  -ru_{1,z},\quad \om^z = 2u_1 + ru_{1,r}.
\end{equation}
We observe that the maximums of $\om_1$, $u_{1,z}$, $u_{1,r}$, $\om^\theta$, $\om^r$ and $\om^z$ are all achieved inside a local region around $(R(t),Z(t))$, whose length scale is featured by $Z(t)\ll R(t)$. It thus follows from \eqref{eq:vorticity_formula} that 
\[R(t) \sim \frac{\|\om^\theta\|_{L^\infty}}{\|\om_1\|_{L^\infty}} \sim \frac{\|\om^r\|_{L^\infty}}{\|u_{1,z}\|_{L^\infty}} \sim \frac{\|\om^z\|_{L^\infty}}{\|u_{1,r}\|_{L^\infty}},\]
which, ideally, is equivalent to the relations 
\[c_s = c_{\om^\theta} - c_{\om_1} = c_{\om^r} - c_{u_{1,z}} = c_{\om^z} - c_{u_{1,r}}.\]
The data in Table \ref{tab:fitting_resolution} or the approximations in \eqref{eq:pattern} yield that 
\[\bar{c}_{\om^\theta} - \bar{c}_{\om_1} \approx \bar{c}_{\om^r} - \bar{c}_{u_{1,z}} \approx \bar{c}_{\om^z} - \bar{c}_{u_{1,r}} \approx 0.5,\]
which then suggests that 
\begin{equation}\label{eq:Cs_power}
C_s(t) \sim R(t) \sim (T-t)^{0.5}.
\end{equation}

\begin{figure}[!ht]
\centering
    \includegraphics[width=0.4\textwidth]{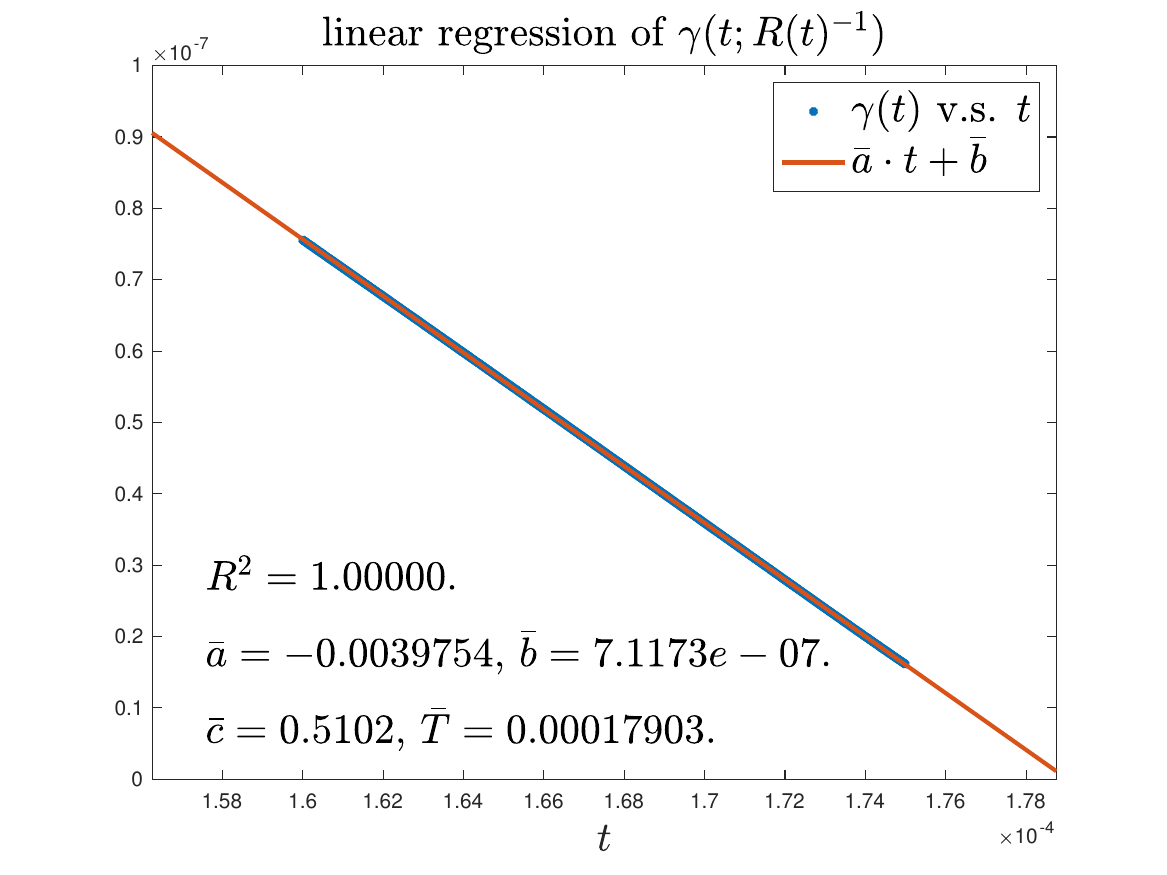}
    \includegraphics[width=0.4\textwidth]{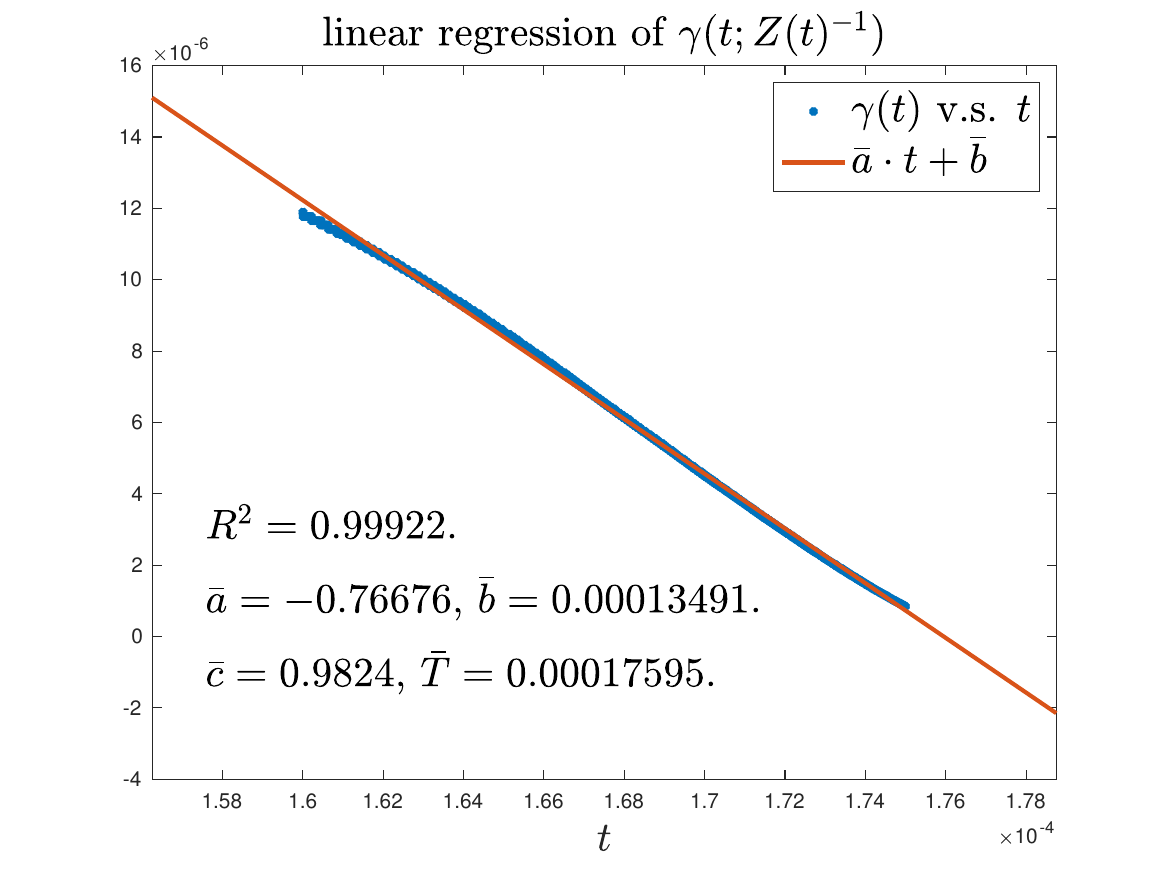}
    \caption[Linear regression $R^{-1},Z^{-1}$]{The linear regression of $\gamma(t;R(t)^{-1})$ (left) and $\gamma(t;Z(t)^{-1})$ (right) on the time interval $[t_1,t_2] = [1.65\times 10^{-4},1.75\times 10^{-4}]$. The blue points are the data points obtained from our computation, and the red lines are the linear models. We plot the linear models on a larger interval.}   
    \label{fig:linear_regression_RZ}
       \vspace{-0.05in}
\end{figure}

The second method to verify the power laws \eqref{eq:power_law} and to learn the powers $c_s,c_l$ is by applying the fitting procedure in Section \ref{sec:fitting_procedure} directly to the blowup quantities $R(t)^{-1},Z(t)^{-1}$. Figure \ref{fig:linear_regression_RZ} presents the linear regression of model \eqref{eq:regression_2} for $R(t)^{-1},Z(t)^{-1}$. We can see in Figure \ref{fig:linear_regression_RZ}(a) that $R(t)$ has an excellent fitting to a power law with a rate $\bar{c}_s \approx 0.5$, which is very close to the conjectured power law \eqref{eq:Cs_power} obtained from the first method.

Figure \ref{fig:linear_regression_RZ}(b) shows that $Z(t)$ also approximately satisfies the power law with a rate $\bar{c}_l \approx 1$, though the fitness is not as good as that of $R(t)$. This lower fitness may be due to the issue that the numerical recording of $Z(t)$ is sensitive to the construction of the adaptive mesh and the interpolation operation between different meshes. Since $Z(t)\ll R(t)$, the relative error in $Z(t)$ due to discretization and interpolation is expected to be much larger than the relative error in $R(t)$. Nevertheless, the estimated power of $Z(t)$, which is close to $1$, is consistent with the conjectured power law \eqref{eq:Cl_power}. 

The consistency between the numerical fitting results in this subsection and in the last subsection further confirms the validity of the (inverse) power laws \eqref{eq:inverse_power_law},\eqref{eq:power_law} of the solution, which provides additional supporting evidence for the existence of a finite-time singularity.

\subsection{Locally self-similar structure} It is well know that the $3$D Euler equations have the scaling invariance property that, if $\vu(x,t)$ is a solution to the equations, then 
\[\vu_{\lambda,\tau}(x,t) := \frac{\lambda}{\tau}\vu\left(\frac{x}{\lambda},\frac{t}{\tau}\right),\quad \forall \lambda,\tau \in \mathbb{R},\]
is still a solution. Similarly, the $3$D Navier--Stokes equations satisfy a more restricted scaling invariance property that, if $\vu(x,t)$ is a solution to the equations, then 
\[\vu_{\lambda}(x,t) := \frac{1}{\lambda}\vu\left(\frac{x}{\lambda},\frac{t}{\lambda^2}\right),\quad \forall \lambda \in \mathbb{R},\]
is still a solution. Regarding these scaling properties, a fundamental question is whether the Euler equations or the Navier--Stokes equations have self-similar solutions of the form 
\begin{equation}\label{eq:self_similar_old}
\vu(x,t) = \frac{1}{(T-t)^\gamma}\vct{U}\left(\frac{x-x_0}{(T-t)^\beta}\right),
\end{equation}
where $\vct{U}$ is a self-similar vector profile and $\beta,\gamma>0$ are scaling powers. If such a solution exists, it will imply that the Euler equations or the Navier--Stokes equations can develop a focusing self-similar singularity at the point $x_0$ at a finite time $T$.

The scaling properties of the Euler equations or the Navier--Stokes equations cannot hold in our scenario due to the existence of the cylinder boundary at $r = 1$ and the variable viscosity coefficients. Nevertheless, a focusing self-similar blowup can still exist asymptotically and locally near the symmetry axis $r = 0$ for two reasons. First of all, since the potential singularity is a focusing one around the origin $(r,z) = 0$, the solid boundary at $r=1$ has no essential contribution to the blowup and can be neglected as a far field. Secondly, the variable viscosity coefficients in our scenario are degenerate near the origin and have an asymptotic behavior in the critical region that may be compatible with a self-similar blowup. More importantly, as we have seen in the previous subsections, the (inverse) power law fitting for the growth and the spatial scales of the solution and the consistency among the fitting results strongly suggest that the solution develops a finite-time self-similar singularity of the form \eqref{eq:self_similar_old}. 

However, the conventional self-similar ansatz \eqref{eq:self_similar_old} with a single spatial scaling is not suitable to  characterize the potential two-scale blowup in our computation, since we have observed two separate spatial scales in the solution. To describe the locally self-similar two-scale blowup scenario, we propose the following self-similar ansatz with two spatial scales in the axisymmetric setting: 
\begin{subequations}\label{eq:self-similar_ansatz}
\begin{align}
u_1(r,z,t) &\sim (T-t)^{-c_{u}}\bar{U}\left(\frac{r-R(t)}{(T-t)^{c_l}}\,,\,\frac{z}{(T-t)^{c_l}}\right), \label{eq:self-similar_u1}\\
\om_1(r,z,t) &\sim (T-t)^{-c_{\om}}\bar{\Omega}\left(\frac{r-R(t)}{(T-t)^{c_l}}\,,\,\frac{z}{(T-t)^{c_l}}\right), \label{eq:self-similar_w1}\\
\psi_1(r,z,t) &\sim (T-t)^{-c_{\psi}}\bar{\Psi}\left(\frac{r-R(t)}{(T-t)^{c_l}}\,,\,\frac{z}{(T-t)^{c_l}}\right), \label{eq:self-similar_psi1}\\
R(t) &\sim (T-t)^{c_s}R_0. \label{eq:self-similar_R}
\end{align}
\end{subequations}
Here $\bar{U},\bar{\Omega},\bar{\Psi}$ denote the self-similar profiles of $u_1,\om_1,\psi_1$ respectively. For notational simplicity, we use $c_u,c_\om,c_\psi$ for $c_{u_1},c_{\om_1},c_{\psi_1}$. As in our previous setting, the reference point $R(t)$ is chosen to be $r$-coordinate of the maximum point of $u_1$, and $R_0>0$ is some normalization constant. This ansatz depicts that, the solution develops an asymptotically self-similar blowup focusing at the point $(R(t),0)$ with a local spatial scaling $(T-t)^{c_l}$, and in the mean time, the point $(R(t),0)$ travels towards the origin with a different length scaling $(T-t)^{c_s}$. 

In what follows, we will carry out a numerical study of the solution profile in a local region around the dynamic location $(R(t),Z(t))$ (the maximum point of $u_1$) to further demonstrate the existence of a locally self-similar blowup of the form \eqref{eq:self-similar_ansatz}. After that, we will analyze the existence conditions and the potential properties of the self-similar profiles $U,\Omega,\Psi$ using a dynamic rescaling formulation of the equations \eqref{eq:axisymmetric_NSE_1}.  

\vspace{-0.05in}
\subsection{Numerical evidence of locally self-similar profiles}\label{sec:evidence_of_self-similar}
An straightforward but useful way to identify a self-similar solution is by comparing the properly normalized profiles of the solution at different time instants. As we have mentioned, the self-similar ansatz \eqref{eq:self-similar_ansatz} is not supposed to hold globally in the entire computational domain $\mathcal{D}_1$; it should only characterize the asymptotic blowup behavior of the solution in a local critical region. Therefore, it is important that we focus our study on the solution profile in a small-scale neighborhood of a reference point. In particular, the reference point is chosen as usual to be $(R(t),Z(t))$, the maximum point of $u_1$.

\begin{figure}[!ht]
\centering
    \includegraphics[width=1\textwidth]{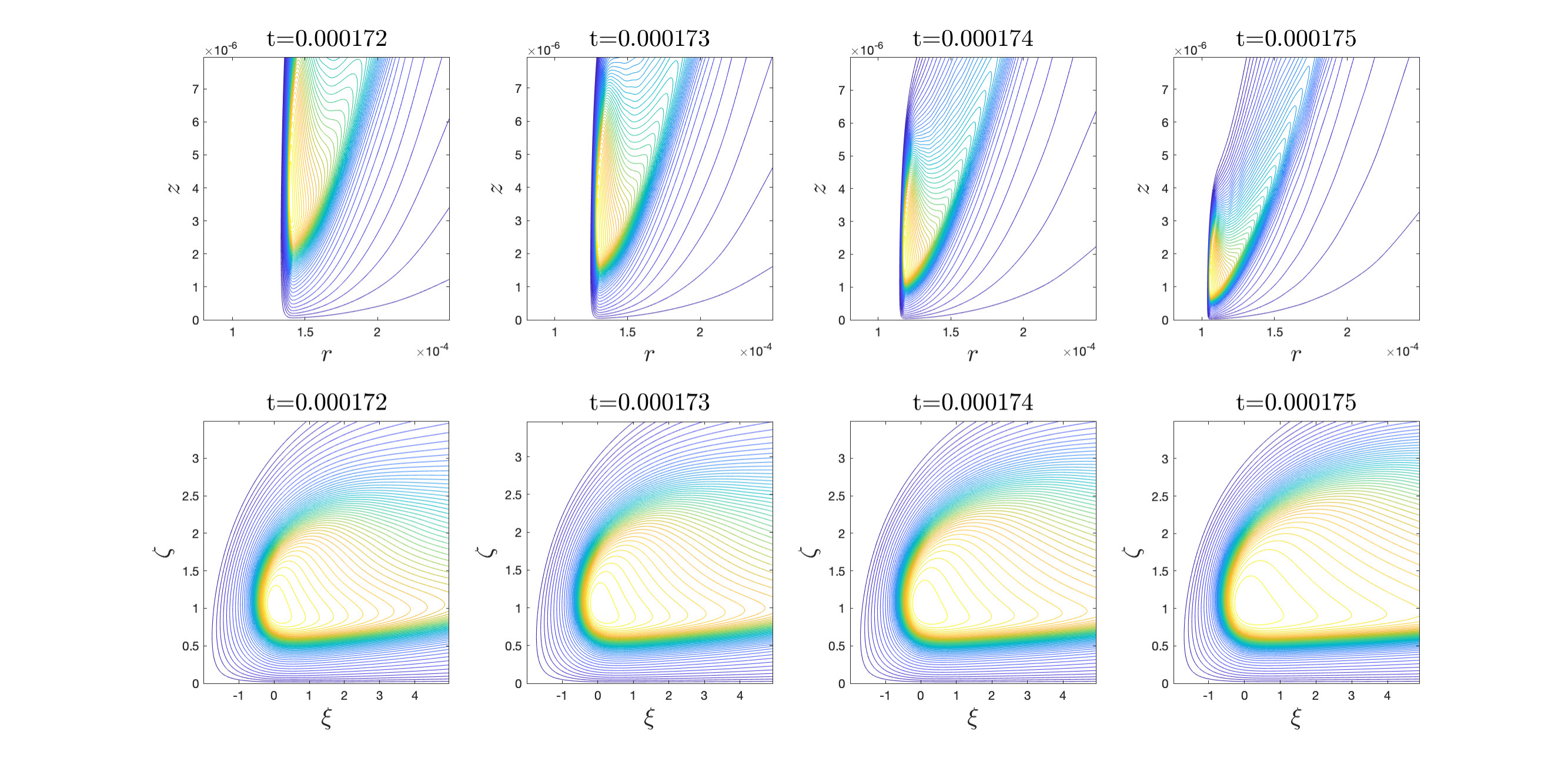}
    \caption{Comparison of the level sets of $u_1$ at different time instants. First row: original level sets of $u_1$ in the domain $(r,z)\in[0.8\times 10^{-4},2.5\times10^{-4}]\times[0,8\times 10^{-6}]$. Second row: rescaled level sets of $u_1$ as a function of $(\xi,\zeta)$ in the domain $(\xi,\zeta)\in[-2,5]\times [0,3.5]$.} 
    \label{fig:levelset_compare_u1} 
\end{figure}

Figure \ref{fig:levelset_compare_u1} compares the level sets of $u_1$ at different time instants. The fisrt row of Figure \ref{fig:levelset_compare_u1} plots the level sets of $u_1$ in a local domain $(r,z)\in[0.8\times 10^{-4},2.5\times10^{-4}]\times[0,8\times 10^{-6}]$. We can see that in a short time interval, from $t = 1.72\times 10^{-4}$ to $t = 1.75\times 10^{-4}$, the profile of $u_1$ changes remarkably. The main part of the profile shrinks in space and travels towards $z=0$ in the $z$ direction and towards the symmetry axis $r=0$ in the $r$ direction. The sharp front also becomes thinner and thinner. However, if we plot the level sets of the spatially rescaled function 
\begin{equation}\label{eq:stretch_u1}
\tilde{u}_1(\xi,\zeta,t) = u_1(Z(t)\xi + R(t),Z(t)\zeta,t)
\end{equation}
as in the second row of Figure \ref{fig:levelset_compare_u1}, we can see that the landscape of $\tilde{u}_1$ (in the $\xi\zeta$-plane) is almost static in the presented time interval. Here 
\[\xi = \frac{r - R(t)}{Z(T)} \sim \frac{r - R(t)}{(T-t)^{c_l}},\quad \zeta = \frac{z}{Z(T)} \sim \frac{z}{(T-t)^{c_l}}\]
are the shifted and rescaled coordinates motivated by the self-similar ansatz \eqref{eq:self-similar_ansatz}. This observation suggests that the asymptotic behavior \eqref{eq:self-similar_u1} of $u_1$ is valid and a self-similar profile $U(\xi,\zeta)$ exists.

\begin{figure}[!ht] 
\centering
    \includegraphics[width=1\textwidth]{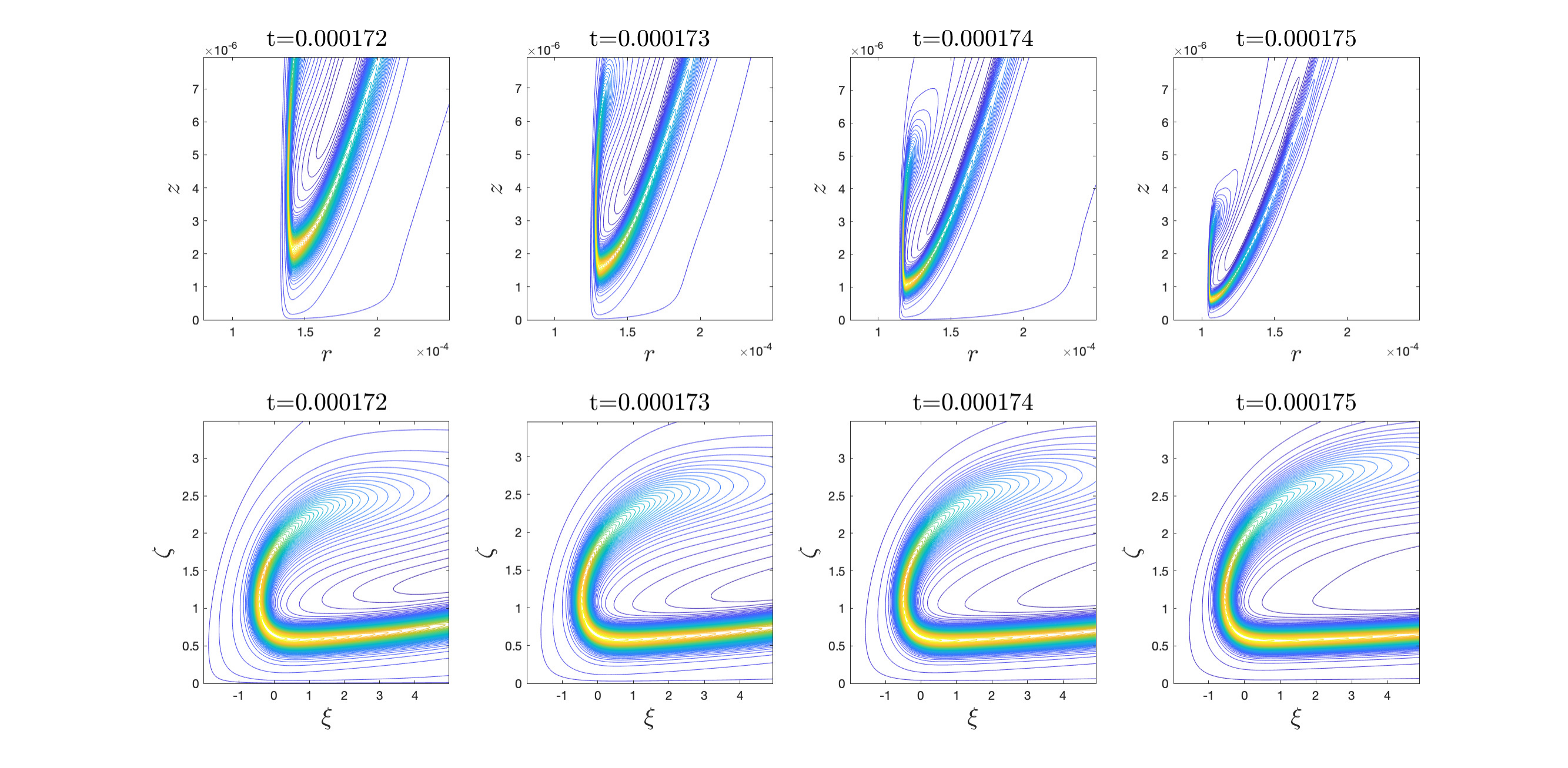}
    \caption{Comparison of the level sets of $\om_1$ at different time instants. First row: original level sets of $u_1$ in the domain $(r,z)\in[0.8\times 10^{-4},2.5\times10^{-4}]\times[0,8\times 10^{-6}]$. Second row: rescaled level sets of $\om_1$ as a function of $(\xi,\zeta)$ in the domain $(\xi,\zeta)\in[-2,5]\times [0,3.5]$.} 
    \label{fig:levelset_compare_w1} 
\end{figure}

Figure \ref{fig:levelset_compare_w1} compares the level sets of $\om_1$ and the level sets of the spatially rescaled function
\begin{equation}\label{eq:stretch_w1}
\tilde{\om}_1(\xi,\zeta,t) = \om_1(Z(t)\xi + R(t),Z(t)\zeta,t)
\end{equation}
in a similar manner. Again, we can see that this profile of $\om_1$ changes a lot in the presented time interval $t\in[1.72\times 10^{-4},1.75\times 10^{-4}]$, while the landscape of $\tilde{\om}_1$ seems to converge. This further suggests the validity of the self-similar conjecture \eqref{eq:self-similar_ansatz} for the solution.

\begin{figure}[!ht]
\centering
    \includegraphics[width=0.35\textwidth]{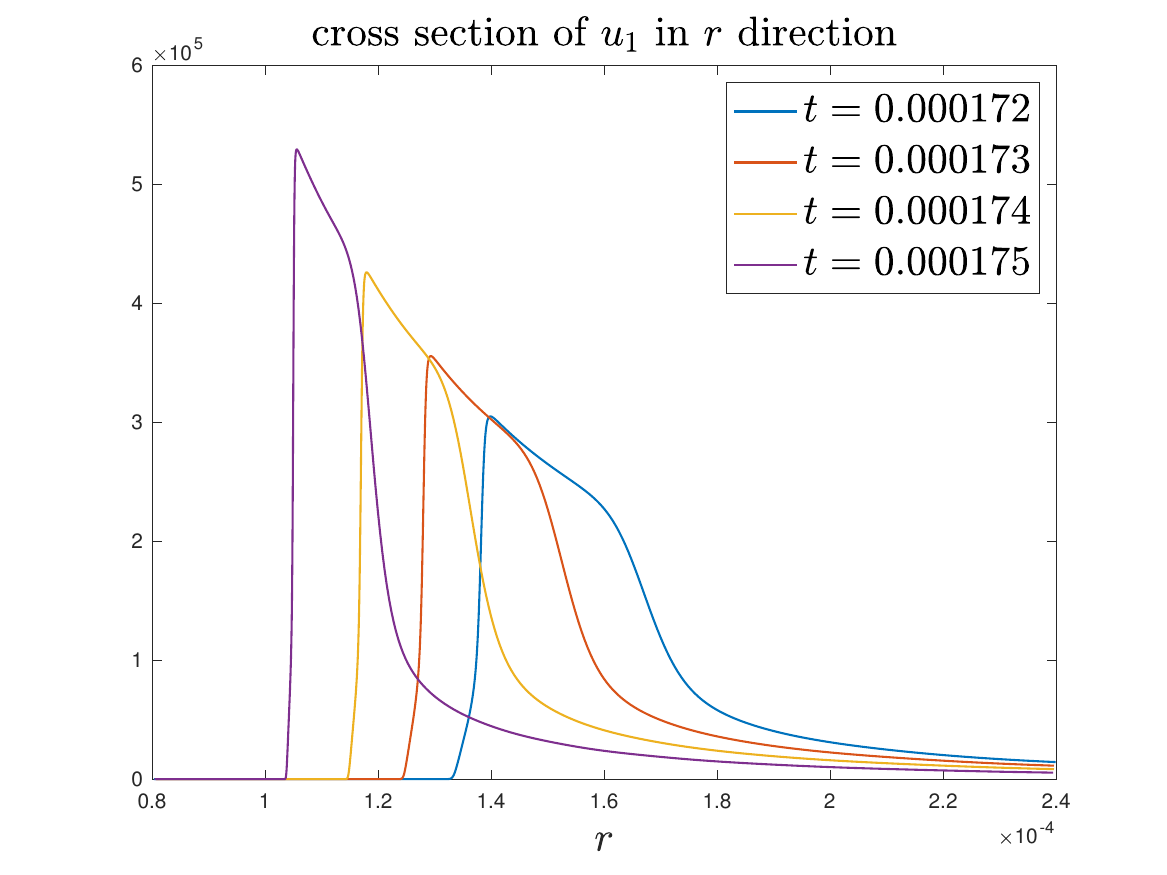}
    \includegraphics[width=0.35\textwidth]{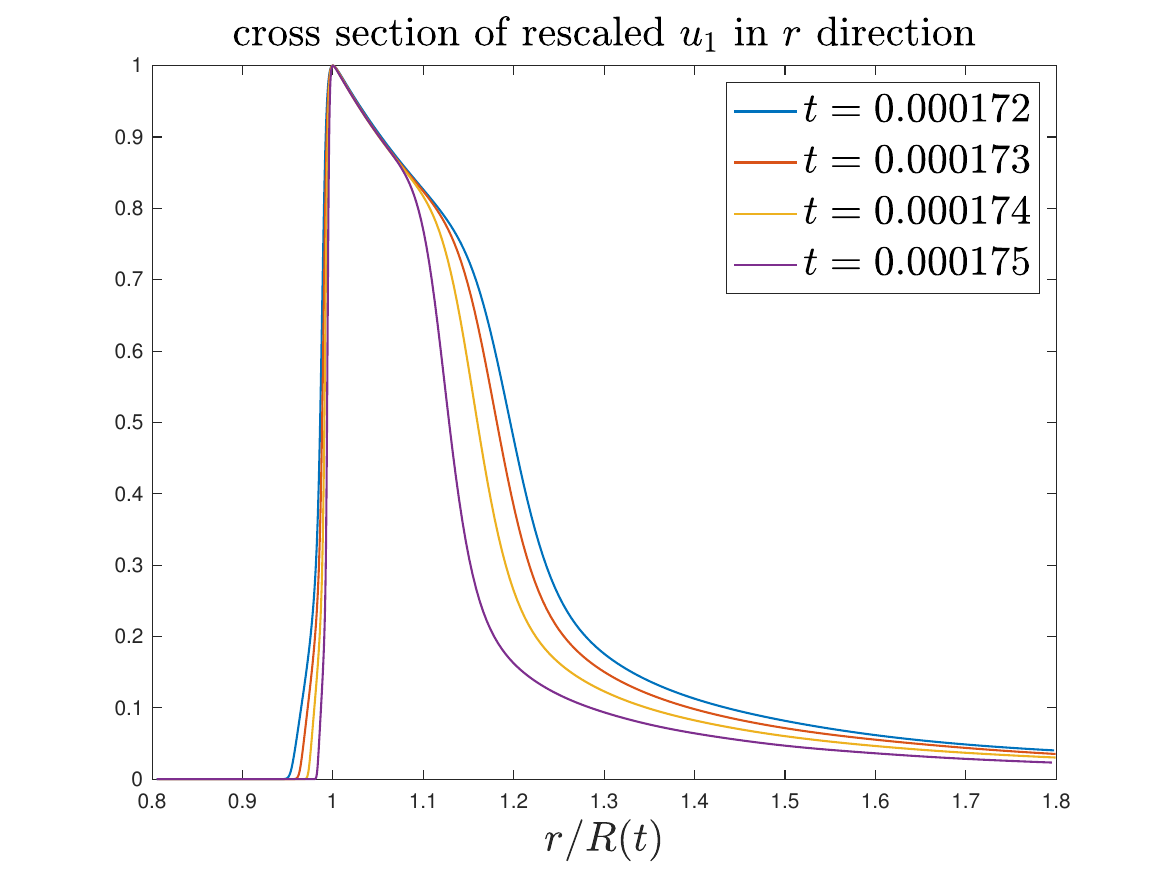}
    \includegraphics[width=0.35\textwidth]{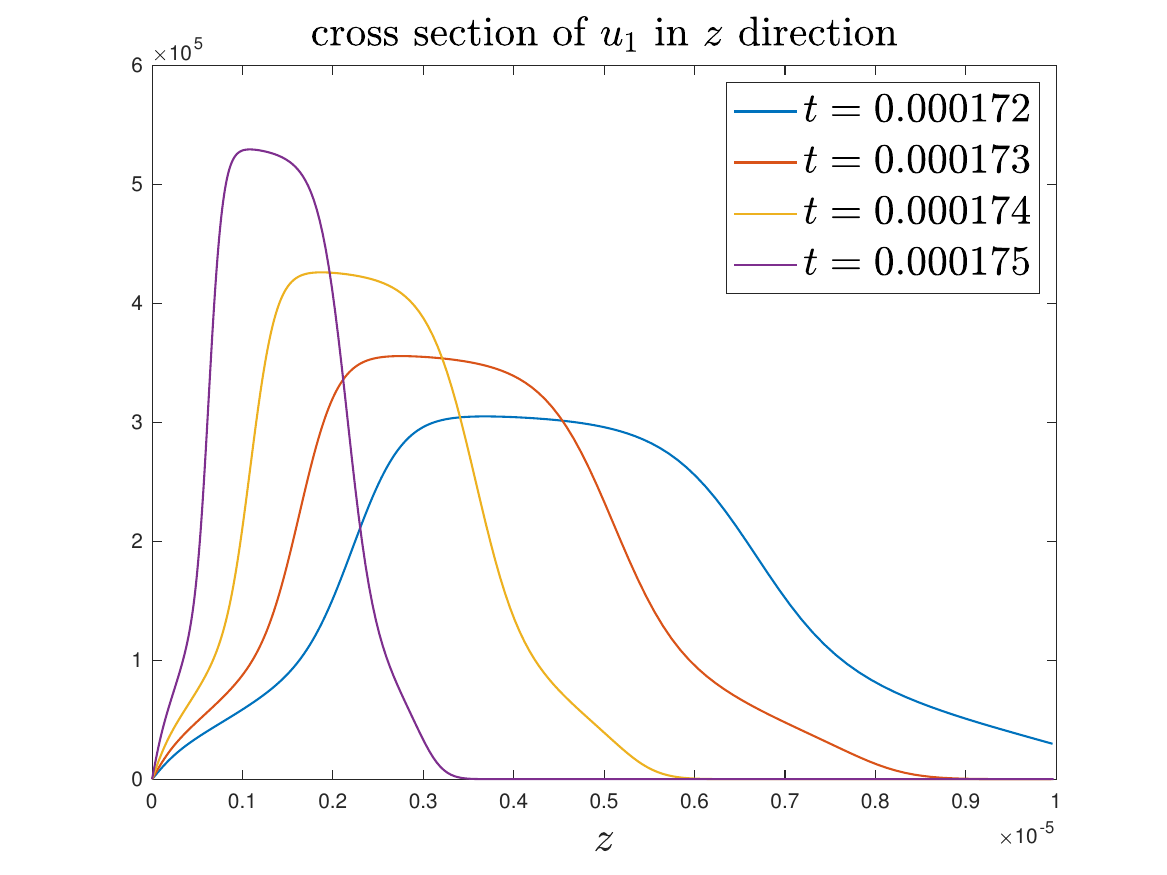}
    \includegraphics[width=0.35\textwidth]{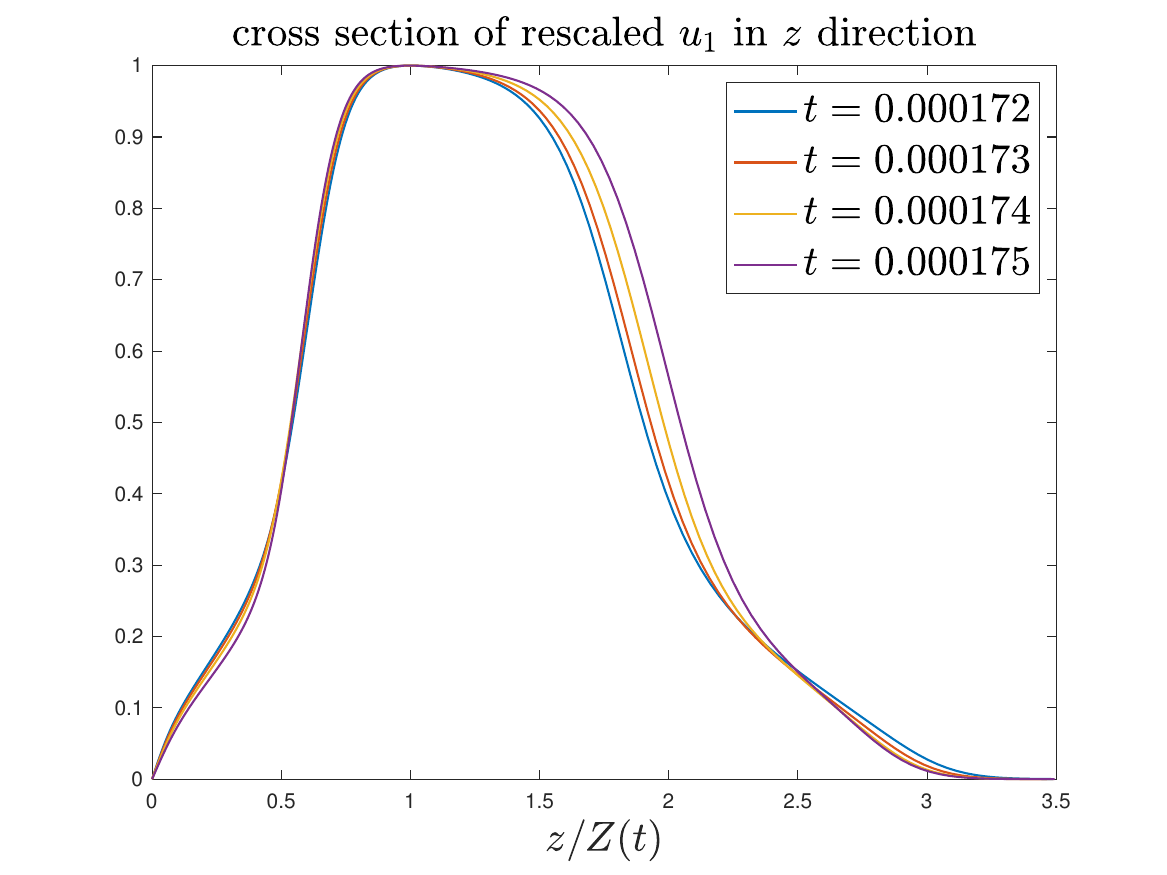}
    \caption[Cross section compare]{Cross sections and rescaled cross sections of $u_1$ through the point $R(t),Z(t)$ in both directions at different time instants instants. (a) Cross sections in the $r$ direction. (b) Rescaled cross sections in the $r$ directions. (c) Cross sections in the $z$ direction. (d) Rescaled cross sections in the $z$ directions.}   
    \label{fig:cross_section_compare}
       \vspace{-0.05in}
\end{figure}

We can also compare the cross sections of the solution at different time instants to study the potential self-similar blowup. As an example, Figure \ref{fig:cross_section_compare}(a) and (c) present the evolution of the cross sections of $u_1$ through the point $R(t),Z(t)$ in both directions. The length scale of the profile shrinks in both directions, and the sharp front travels towards $r=0$. For comparison, Figure \ref{fig:cross_section_compare}(b) and (d) plot the corresponding cross sections of the normalized function $u_1/\|u_1\|_{L^\infty}$ in the rescaled coordinates $(r/R(t),z/Z(t))$, illustrating the potential convergence of the properly rescaled profile of $u_1$. Note that the rescaled cross sections in the $r$ direction seem to converge only in a small neighborhood of $R(t)$, i.e. within $|r-R(t)| = O(Z(t))$, implying that the asymptotic self-similar behavior \eqref{eq:self-similar_u1} is only valid locally. These results again support the existence of a locally self-similar profile of the solution near the reference point $(R(t),Z(t))$. 

\subsection{Asymptotic analysis of self-similar blowup}\label{sec:asymptotic_analysis}
In this subsection, we carry out an asymptotic analysis based on the self-similar ansztz \eqref{eq:self-similar_ansatz} to provide a possible understanding of the two-scales features, blowup rates and self-similar behaviors that we observed numerically. 

A standard method to study a self-similar blowup is by substituting the self-similar ansztz \eqref{eq:self-similar_ansatz} into the physical equations \eqref{eq:axisymmetric_NSE_1} and deriving equations for the potential self-similar profiles $\bar{U},\bar{\Omega},\bar{\Psi}$, based on the fundamental assumption that these profiles exist and are smooth functions. More appropriately, we will introduce time-dependent profile solutions $U,\Omega,\Psi$ and treat the potential self-similar profiles $\bar{U},\bar{\Omega},\bar{\Psi}$ as the steady state of $U,\Omega,\Psi$. Thus, we can relate ($u_1, \omega_1,\psi_1)$ to $(U,\Omega,\Psi)$ by a dynamic change of variables given below:
\begin{subequations}\label{eq:dynamic_profile}
\begin{align}
u_1(r,z,t) &= (T-t)^{-c_{u}}U\left(\xi,\zeta,\tau(t)\right), \label{eq:dynamic_profile_u1}\\
\om_1(r,z,t) &= (T-t)^{-c_{\om}}\Omega\left(\xi,\zeta,\tau(t)\right), \label{eq:dynamic_profile_w1}\\
\psi_1(r,z,t) &= (T-t)^{-c_{\psi}}\Psi\left(\xi,\zeta,\tau(t)\right), \label{eq:dynamic_profile_psi1}
\end{align}
where 
\begin{equation}\label{eq:dynamic_coordinate}
\xi := \frac{r - R(t)}{C_l(t)},\quad \zeta := \frac{z}{C_l(t)},
\end{equation}
and $\tau(t)$ is a rescaled time satisfying 
\begin{equation}\label{eq:dynamic_time}
\tau'(t) = (T-t)^{-1}.
\end{equation}
Now the self-similar ansatz \eqref{eq:self-similar_ansatz} asserts that the profile solutions $U(\xi,\zeta,\tau),\Omega(\xi,\zeta,\tau),\Psi(\xi,\zeta,\tau)$ in the $\xi\zeta$-coordinates should converge to some time-independent profiles $\bar{U}(\xi,\zeta),\bar{\Omega}(\xi,\zeta),\bar{\Psi}(\xi,\zeta)$ as $\tau\rightarrow \infty$ (i.e. $t\rightarrow T$). In particular, $\bar{U},\bar{\Omega},\bar{\Psi}$ should be smooth functions of $\xi, \zeta$.
\end{subequations}

Before we derive the equations for the profile solutions $U,\Omega,\Psi$, we first make some preparations to simplify our argument, so that we can focusing on delivering the main idea. For simplicity, we only keep the viscosity terms \eqref{eq:viscosity} to their leading-order terms:
\[f_{u_1} \approx \nu^r\left(u_{1,rr} + \frac{3}{r}u_{1,r}\right) + \nu^zu_{1,zz},\quad f_{\om_1} \approx \nu^r\left(\om_{1,rr} + \frac{3}{r}\om_{1,r}\right) + \nu^z\om_{1,zz}.\]
We remark that this approximation will not change the resulting equations for the self-similar profiles in the asymptotic analysis, since the terms that we have dropped are asymptotically small under our ansatz. Moreover, we assume that the ansatz \eqref{eq:self-similar_R} is actually an identity:
\[R(t) = (T-s)^{c_s}R_0.\]
Guided by our numerical observations, we make the two-scale assumption:
\begin{equation}\label{eq:two-scale_assumption}
c_s < c_l,\quad \text{or equivalently,} \quad R(t)/C_l(t) \rightarrow +\infty\quad \text{as}\quad t\rightarrow T.
\end{equation}
Recall that $C_l := (T-t)^{c_l}$ (see the definition \eqref{eq:power_law}). We will also use the following notations
\[C_1(t) := (T-t)^{-1},\quad C_u(t) := (T-t)^{-c_u}, \quad C_\om(t) := (T-t)^{-c_\om},\quad C_\psi(t) := (T-t)^{-c_\psi}.\]

\vspace{-0.1in}
\subsubsection{Substituting the self-similar ansztz}
Now we substitute the change of variables \eqref{eq:dynamic_profile} into the equations \eqref{eq:axisymmetric_NSE_1} (with the simplified viscosity terms). For clarity, we do this term by term. For the $u_1$ equation \eqref{eq:as_NSE_1_a}, we have
\begin{subequations}\label{eq:substitution}
\begin{equation}\label{eq:substitution_u1}
\begin{split}
u_{1,t} &= C_1C_uU_\tau + c_uC_1C_uU + c_l C_1 C_u(\xi U_\xi + \zeta U_\zeta) + \underline{c_sC_1C_uC_l^{-1}RU_\xi} ,\\
u^ru_{1,r} + u^zu_{1,z} &= C_\psi C_u C_l^{-1}\big( - \xi\Psi_\zeta U_\xi +(2\Psi + \xi\Psi_\xi)U_\zeta\big) + \underline{C_\psi C_u C_l^{-2}R\big( - \Psi_\zeta U_\xi  + \Psi_\xi U_\zeta\big)},\\
2\psi_{1,z} u_1 &= 2 C_\psi C_u C_l^{-1}\Psi_\zeta U, \\ 
f_{u_1} &= C_u C_l^{-2} \big( \nu^r U_{\xi\xi} + 3\nu^r (\xi + R C_l^{-1})^{-1} U_\xi  + \nu^zU_{\zeta\zeta}\big).
\end{split}
\end{equation}
Note that we have used the expressions of $u^r,u^z$ in \eqref{eq:as_NSE_1_d}. We have also used the relation \eqref{eq:dynamic_time}: $\tau'(t)=(T-t)^{-1} = C_1$. Similarly, for the $\om_1$ equation \eqref{eq:as_NSE_1_b}, we have
\begin{equation}\label{eq:substitution_w1}
\begin{split}
\om_{1,t} &= C_1C_\om \Omega_\tau + c_\om C_1C_\om \Omega + c_l C_1 C_\om(\xi \Omega_\xi + \zeta \Omega_\zeta) + \underline{c_sC_1C_\om C_l^{-1}R\Omega_\xi} ,\\
u^r\om_{1,r} + u^z\om_{1,z} &= C_\psi C_\om C_l^{-1}\big( - \xi\Psi_\zeta \Omega_\xi +(2\Psi + \xi\Psi_\xi)\Omega_\zeta\big) + \underline{C_\psi C_\om C_l^{-2}R\big( - \Psi_\zeta \Omega_\xi  + \Psi_\xi \Omega_\zeta\big)},\\
2u_{1,z} u_1 &= 2 C_u^2 C_l^{-1} U_\zeta U, \\ 
f_{\om_1} &= C_\om C_l^{-2} \big( \nu^r \Omega_{\xi\xi} + 3\nu^r (\xi + R C_l^{-1})^{-1} \Omega_\xi  + \nu^z \Omega_{\zeta\zeta}\big).
\end{split}
\end{equation}
Finally, for the Poisson equation \eqref{eq:as_NSE_1_c}, we have
\begin{equation}\label{eq:substitution_Poisson}
-\left(\partial_r^2+\frac{3}{r}\partial_r+\partial_z^2\right)\psi_1 = \om_1\quad \Longrightarrow\quad - C_\psi C_l^{-2}\left (\partial_\xi^2+\frac{3}{\xi + R C_l^{-1}} \partial_\xi+ \partial_\zeta^2\right)\Psi = C_\om \Omega.
\end{equation}
\end{subequations}

\subsubsection{Balancing the equations} The next step is to determine the relations between the quantities $C_u$, $C_\om$, $C_\psi$, $C_l$ and $R$ by balancing the terms in each equation of \eqref{eq:substitution} in the asymptotic regime $t\rightarrow T$, based on the assumption that the limit profiles $\bar{U},\bar{\Omega},\bar{\Psi}$ are smooth regular functions of $\xi,\zeta$ and are independent of time $t$. We also assume that the viscosity term are of the same order as the vortex stretching term. This balance is crucial in determining the length scale for 
$C_l$ or $Z(t)$.

We have underlined some terms in \eqref{eq:substitution_u1} and \eqref{eq:substitution_w1} for some reason to be clarified later. For those terms that are not underlined in \eqref{eq:substitution_u1}, the balance among various terms as $t\rightarrow T$ requires
\[C_1C_u = C_\psi C_u C_l^{-1} \sim \nu C_u C_l^{-2}.\]
Similarly, for those terms that are not underlined in \eqref{eq:substitution_w1}, the balance among various terms as $t\rightarrow T$ enforces
\[C_1C_\om = C_\psi C_\om C_l^{-1} = C_u^2C_l^{-1} \sim \nu C_\om C_l^{-2}.\]
Finally, for the Poisson equation \eqref{eq:substitution_Poisson} to balance each other as $t\rightarrow T$, we must have
\[C_\psi C_l^{-2} = C_\om.\]
Summarizing up these relations, we obtain 
\begin{equation}\label{eq:relation_1}
\left\{\begin{array}{l}
C_u = C_1,\\
C_\om = C_1C_l^{-1},\\
C_\psi = C_1C_l,\\
\end{array}\right.
\quad \Longleftrightarrow \quad 
\left\{\begin{array}{l}
c_u = 1,\\
c_\om = 1 + c_l,\\
c_\psi = 1 - c_l,\\
\end{array}\right.
\end{equation}
and 
\begin{equation}\label{eq:relation_2}
\nu^r\sim\nu^z \sim C_1C_l^2 = (T-t)^{2c_l - 1}.
\end{equation}
Note that the relations \eqref{eq:relation_1} also imply that the underlined terms can balance with each other in \eqref{eq:substitution_u1} and in \eqref{eq:substitution_w1}.

So far, we have already obtained some meaningful information of the blowup rates. If the self-similar ansatz \eqref{eq:self-similar_ansatz} is true, then no matter what the spatial scalings $c_s,c_l$ are, the asymptotic blowup rates of $u_1,\psi_{1,r},\psi_{1,z}$ are always $1$:
\[\|u_{1}\|_{L^\infty} \sim C_u  = C_1 = (T-t)^{-1},\quad \|\psi_{1,r}\|_{L^\infty}\sim\|\psi_{1,z}\|_{L^\infty}\sim C_\psi C_l^{-1} = C_1 = (T-t)^{-1}. \]
This result of the asymptotic analysis is consistent with our observations and fitting results in Section \ref{sec:growth_fitting}, which confirms the validity of the inverse power law for $u_1$. To obtain the blowup rate of $\om_1$, we still need to determine the value of $c_l$.

\vspace{-0.15in}
\subsubsection{Conservation of circulation} To determine the values of $c_s,c_l$, we need to make use of the conservation of the total circulation, an important physical property of the axisymmetric Euler or Navier--Stokes equations. 

The total circulation is defined as 
\[\Gamma(r,z,t) := ru^\theta(r,z,t) = r^2 u_1(r,z,t).\] 
It is easy to derive the equation of $\Gamma$ from the $u^\theta$ equation \eqref{eq:as_NSE_a} of (or the $u_1$ equation \eqref{eq:as_NSE_1_a}):
\begin{equation}\label{eq:circulation}
\Gamma_t + u^r \Gamma_r + u^z \Gamma_z = rf_u.
\end{equation}
Recall that $f_u$ is the viscosity term in the $u^\theta$ equation \eqref{eq:as_NSE_a}. In our scenario, $\Gamma$ is a non-negative variable in the computational domain $\mathcal{D}_1$, since $u_1\geq0$ in $\mathcal{D}_1$. It is well known that the circulation function $\Gamma$ satisfies a maximum principle for the Euler equations or the Navier--Stokes equations. In fact, it also satisfies a maximum principle for our equations \eqref{eq:axisymmetric_NSE_1}. Let $(\tilde{R}(t),\tilde{Z}(t))$ be a local maximum point of $\Gamma$, then we have from \eqref{eq:circulation} that
\[\frac{\diff\,}{\diff t} \Gamma(\tilde{R},\tilde{Z},t) = \nu^r\Gamma_{rr} + \nu^z \Gamma_{zz} - \frac{\nu^r_r}{\tilde{R}}\Gamma\leq 0,\]
where we have used that $\Gamma_{rr},\Gamma_{zz}\leq 0$ at the local maximum point $(\tilde{R}(t),\tilde{Z}(t))$, and that $\nu^r_r\geq 0$ for the variable viscosity coefficient we use in our computation. This means that $\Gamma(\tilde{R}(t),\tilde{Z}(t),t)$ is always non-increasing in time. 

In fact, we have observed that the maximum point $(R(t),Z(t))$ of $u_1$ is also a local maximum point of $\Gamma$. Therefore, 
\[\frac{\diff\,}{\diff t} \Gamma(R(t),Z(t),t) = \nu^r\Gamma_{rr} + \nu^z \Gamma_{zz} - \frac{\nu^r_r}{R}\Gamma\leq 0.\]
The viscosity term $\nu^r\Gamma_{rr} + \nu^z \Gamma_{zz}$ cannot damp $\Gamma(R,Z,t)$ to $0$ in a finite time as long as the viscosity coefficients are bounded, which is the case in our computation. From the expression \eqref{eq:nu_r} of $\nu_r$, we have $\nu^r_r = O(r)$ in the critical blowup region and thus $\nu^r_r(R,Z)/R = O(1)$. Therefore, the term $-\nu^r_r\Gamma/R$ is a linear damping term with an $O(1)$ coefficient, which can only drive $\Gamma$ to $0$ as $t\rightarrow +\infty$. In summary, $\Gamma(R,Z,t)$ will not blow up or vanish to $0$ in any finite time $T$. That is, $\Gamma(R,Z,t) = R^2u_1(R,Z,t) = O(1)$ as $t\rightarrow T$. It follows from the blowup scaling of $u_1$ that 
\[R(t) \sim \|u_1\|_{L^\infty}^{-1/2} \sim (T-t)^{c_u/2} = (T-t)^{1/2},\]
which implies 
\[c_s = 1/2.\]
We remark that this property only relies on the conservation of the maximum circulation and the fact that $c_u = 1$, which are intrinsic to the equations \eqref{eq:axisymmetric_NSE_1}.  

\subsubsection{Determining the smaller scale} To determine $c_l$, we need to use again the asymptotic values of the variable viscosity coefficients, 
\[\nu^r\sim \nu^z = O(r^2) + O(z^2) + 0.025 \;\|\om^\theta\|_{L^\infty}^{-1},\]
for $r,z$ close to $0$, which follows from the expressions \eqref{eq:viscosity_coefficient}. In particular, in the critical region around the reference point $(R(t),Z(t))$ where the self-similar ansatz \eqref{eq:self-similar_ansatz} is assumed to be valid, we have
\begin{align*}
\nu^r\sim \nu^z &= O(R^2) + O(Z^2) + O(R^{-1}\|\om_1\|_{L^\infty}^{-1}) \\
&= O((T-t)^{2c_s}) + O((T-t)^{2c_l}) + O((T-t)^{1+c_l - c_s}) \sim (T-t)^1,
\end{align*}
where we have used $c_s = 1/2$ and the two-scale assumption \eqref{eq:two-scale_assumption} that $c_l>c_s$, so that $O((T-t)^{2c_l})$ and $O((T-t)^{1+c_l - c_s})$ are dominated by $O((T-t)^{2c_s})$. Comparing this with the relation \eqref{eq:relation_2}, we conclude that $2c_l - 1  = 1$, that is, $c_l = 1$.

\begin{subequations}\label{eq:asymptotic_powers}
We have now obtained all the blowup rates and the spatial scalings in the self-similar ansztz \eqref{eq:self-similar_ansatz}: 
\begin{equation}
c_s = 1/2,\quad c_l = 1,\quad c_u = 1,\quad c_\om = 1 + c_l = 2,\quad c_\psi = 1 - c_l = 0.
\end{equation}
Moreover, the derivative relations and product relations yield that
\begin{equation}
\begin{split}
c_{\psi_{1,r}} = c_{\psi_{1,z}} &= c_\psi + c_l = 1,\quad c_{u_{1,r}} = c_{u_{1,z}} = c_u + c_l = 2,\\
c_{\om^\theta} = c_{\om} - c_s &= 1.5,\quad c_{\om^r} = c_{u_{1,z}} - c_s = 1.5, \quad c_{\om^z} = c_{u_{1,r}} - c_s = 1.5.
\end{split}
\end{equation}
These results are surprisingly consistent with the fitting data in Sections \ref{sec:growth_fitting} and \ref{sec:scale_fitting}, especially the numerically observed pattern \eqref{eq:pattern}. The consistency between the blowup rates obtained by our numerical fitting procedures and those obtained by the asymptotic scaling analysis provides further support for the existence of a finite-time locally self-similar blowup of the form \eqref{eq:self-similar_ansatz}. 
\end{subequations} 

We remark that if we use a constant viscosity coefficient $\nu^r=\nu^z=\mu$, then through a similar balancing procedure we will obtain a different scaling result with $c_l=1/2$. This implies that there is no two-scale feature in the potential blowup solution, which is consistent with our numerical observations. This also explains why the two-scale blowup cannot survive the viscosity with a constant coefficient in Case $2$ computation, as we have seen in Section \ref{sec:original_NSE}.

\vspace{-0.1in}
\subsubsection{Equations for the self-similar profiles} Our previous analysis on the blowup rates is based on the fundamental assumption that the asymptotic self-similar profiles $\bar{U},\bar{\Omega},\bar{\Psi}$ exist, for which we present strong numerical evidences in the previous sections. To gain more insights into the locally self-similar blowup, we will derive some potential time-independent equations for the self-similar profiles, which may help us understand the properties and the existence conditions of these profiles.

Collecting the terms in \eqref{eq:substitution_u1}, \eqref{eq:substitution_w1} and \eqref{eq:substitution_Poisson} and using the relations in \eqref{eq:relation_1} and \eqref{eq:relation_2}, we first obtain the following time-dependent equations for the dynamic profiles:
\begin{subequations}\label{eq:profile_equation_1}
\begin{align}
&U_\tau + (c_l\xi-\xi\Psi_\zeta)U_\xi + (c_l\zeta+ 2\Psi + \xi\Psi_\xi) U_\zeta + \underline{ RC_l^{-1}\big((c_s-\Psi_\zeta)U_\xi + \Psi_\xi U_\zeta\big)} \nonumber \\
&\qquad\qquad = 2\Psi_\zeta U - c_uU + \tilde{\nu}^r\big(U_{\xi\xi} + 3(\xi+RC_l^{-1})U_\xi\big) + \tilde{\nu}^z U_{\zeta\zeta},\label{eq:profile_equation_u1_1}\\
&\Omega_\tau + (c_l\xi-\xi\Psi_\zeta)\Omega_\xi + (c_l\zeta+ 2\Psi + \xi\Psi_\xi) \Omega_\zeta + \underline{ RC_l^{-1}\big((c_s-\Psi_\zeta)\Omega_\xi + \Psi_\xi \Omega_\zeta\big)} \nonumber \\
&\qquad\qquad = 2U_\zeta U - c_\om\Omega + \tilde{\nu}^r\big(\Omega_{\xi\xi} + 3(\xi+RC_l^{-1})\Omega_\xi\big) + \tilde{\nu}^z \Omega_{\zeta\zeta}\label{eq:profile_equation_w1_1},\\
&-\big(\partial_{\xi\xi} + 3(\xi + RC_l^{-1})^{-1}\partial_\xi + \partial_{\zeta\zeta}\big)\Psi  = \Omega, \label{eq:profile_equation_Poisson_1}
\end{align}
\end{subequations}
where 
\[\tilde{\nu}^r(\xi,\zeta) := (T-t)^{-1}\nu^r(C_l\xi + R(t),C_l\zeta),\quad \tilde{\nu}^z(\xi,\zeta) := (T-t)^{-1}\nu^z(C_l\xi + R(t),C_l\zeta).\]
For the self-similar profiles to exist, it requires that the solution to the equations \eqref{eq:profile_equation_1} converge to some non-trivial steady state $\bar{U},\bar{\Omega},\bar{\Psi}$ as $\tau\rightarrow \infty$. We thus expect that each equation of \eqref{eq:profile_equation_1} is balanced in the limit $\tau\rightarrow \infty$. Note that from the relation \eqref{eq:dynamic_time} we have $\tau = -\log(T-t) + c$ for some constant $c$, so $\tau\rightarrow \infty$ means $t\rightarrow T$. Since we have argued that $R(t)/C_l(t)\sim (T-t)^{-1/2}\rightarrow \infty$ as $t\rightarrow T$, the underlined terms in equations \eqref{eq:profile_equation_u1_1} and \eqref{eq:profile_equation_w1_1} need to satisfy some extra conditions so that they can balance with the other terms. In particular, we should have
\[RC_l^{-1}\big((c_s-\Psi_\zeta)U_\xi + \Psi_\xi U_\zeta\big) \rightarrow G_u(\xi,\zeta), \quad RC_l^{-1}\big((c_s-\Psi_\zeta)\Omega_\xi + \Psi_\xi \Omega_\zeta\big) \rightarrow G_\om(\xi,\zeta) \]
for some smooth functions $G_u,G_\om = O(1)$ as $t\rightarrow T$. This further implies that 
\[(c_s-\Psi_\zeta)U_\xi + \Psi_\xi U_\zeta\rightarrow 0\quad\text{and}\quad (c_s-\Psi_\zeta)\Omega_\xi + \Psi_\xi \Omega_\zeta\rightarrow 0 \quad \text{as $t\rightarrow T$}.\]
Moreover, since $R(t)/C_l(t)\rightarrow \infty$ as $t \rightarrow T$, the lower order viscosity term $3(\xi + RC_l^{-1})^{-1}\partial_\xi$ in equations \eqref{eq:profile_equation_u1_1},\eqref{eq:profile_equation_w1_1} and in the Poisson equation \eqref{eq:profile_equation_Poisson_1} should vanish as $t\rightarrow T$.  

Based on the preceding discussions, we conjecture the following time-independent equations for the self-similar profiles $\bar{U},\bar{\Omega},\bar{\Psi}$:
\begin{subequations}\label{eq:profile_equation_2}
\begin{align}
(c_l\xi-\xi\bar{\Psi}_\zeta)\bar{U}_\xi + (c_l\zeta+ 2\bar{\Psi} + \xi\bar{\Psi}_\xi) \bar{U}_\zeta + G_u&= 2\bar{\Psi}_\zeta \bar{U} - c_u\bar{U} + \bar{\nu}^r \bar{U}_{\xi\xi}  + \bar{\nu}^z \bar{U}_{\zeta\zeta},\label{eq:profile_equation_u1_2}\\
(c_l\xi-\xi\bar{\Psi}_\zeta)\bar{\Omega}_\xi + (c_l\zeta+ 2\bar{\Psi} + \xi\bar{\Psi}_\xi) \bar{\Omega}_\zeta + G_\om &= 2\bar{U}_\zeta \bar{U} - c_\om\bar{\Omega} + \bar{\nu}^r\bar{\Omega}_{\xi\xi} + \bar{\nu}^z \bar{\Omega}_{\zeta\zeta}\label{eq:profile_equation_w1_2},\\
-(\partial_{\xi\xi} + \partial_{\zeta\zeta})\bar{\Psi}  &= \bar{\Omega}, \label{eq:profile_equation_Poisson_2}\\
(c_s-\bar{\Psi}_\zeta)\bar{U}_\xi + \bar{\Psi}_\xi \bar{U}_\zeta &= 0,\label{eq:profile_constraint_u1}\\
(c_s-\bar{\Psi}_\zeta)\bar{\Omega}_\xi + \bar{\Psi}_\xi \bar{\Omega}_\zeta &= 0,\label{eq:profile_constraint_w1}
\end{align}
\end{subequations}
where
\[\bar{\nu}^r(\xi,\zeta) := \lim_{t\rightarrow T}(T-t)^{-1}\nu^r(R(t),0),\quad \bar{\nu}^z(\xi,\zeta) := \lim_{t\rightarrow T}(T-t)^{-1}\nu^z(R(t),0).\]
The existence of the self-similar profiles of the above equations is beyond the scope of this paper. 
We remark that our preceding asymptotic scaling analysis is valid when we properly rescale the solution and zoom into an $O(C_l(t))$ neighborhood of the point $(R(t),0)$ in the $rz$-plane. From a $3$D macroscopic perspective, one may view the blowup region of the solution as a tubular ring surronding the symmetry axis with a decreasing radius $R(t)$ and a shrinking thickness $C_l(t)$. Correspondingly, the equations \eqref{eq:profile_equation_2} only characterize the asymptotic self-similar behavior of the solution on the scale of $C_l(t)$ in the limit $t\rightarrow T$ around the $1$D ring $\{(r,z,\theta)| r= R(t), z=0, \theta\in [0,2\pi)\}$. If we zoom out to an $O(R(t))$ region around the origin, we can only see the blowup region shrinking into a $1$D ring, and hence we cannot be see the effect of these equations. Therefore, we say that the potential two-scale blowup in our scenario is only locally self-similar with respect to the smaller scale $C_l(t)$.  

\vspace{-0.1in}
\subsubsection{A level set condition} Though we cannot use the equations \eqref{eq:profile_equation_2} to determine the self-similar profiles $\bar{U},\bar{\Omega},\bar{\Psi}$, we can still learn some properties of the profiles from them. Note that the equations \eqref{eq:profile_constraint_u1} and \eqref{eq:profile_constraint_w1} are independent of the unknown functions $G_u,G_\om$, it thus makes sense to study their implications. One can see these two equations as necessary conditions for the self-similar profiles to exist. The physical solution $u_1,\om_1,\psi_1$ can only develop the two-scale self-similar blowup when their profiles satisfy these two conditions locally (or after rescaling). If these two conditions are not satisfied locally, the underlined terms in equations \eqref{eq:profile_equation_1} will not be compatible with our preceding scaling analysis based on the balance of scales among various terms. 

Let $\Phi(\xi,\zeta) = \Psi(\xi,\zeta) - c_s\zeta$. Then the equations \eqref{eq:profile_constraint_u1} and \eqref{eq:profile_constraint_w1} can be written as 
\begin{equation}\label{eq:necessary_condition}
\Phi_\xi U_\zeta -\Phi_\zeta U_\xi = 0,\quad \Phi_\xi \Omega_\zeta - \Phi_\zeta \Omega_\xi = 0,
\end{equation}
which implies that the gradient of $U,\Omega,\Phi$ are parallel to each other:
\[(U_\xi,U_\zeta) \parallel (\Omega_\xi,\Omega_\zeta) \parallel (\Phi_\xi,\Phi_\zeta),\]
or that the level sets of $U,\Omega,\Phi$ have the same geometric contours. In other words, the profiles $U,\Omega,\Phi$ can be viewed as functions of each other:
\[U = U(\Phi),\quad \Omega = \Omega(\Phi).\]
The above relationship also implies that the velocity field induced by the modified stream function $\Phi$ is parallel to the level set of $U$ and $\Omega$. In other words, the large underlined advection terms in \eqref{eq:profile_equation_u1_1} and \eqref{eq:profile_equation_w1_1} enforce a condition that the local velocity field near the sharp front should be tangent to the sharp front.
This also provides us a way to numerically verify the condition \eqref{eq:necessary_condition}. Note that under the asymptotic ansatz \eqref{eq:self-similar_ansatz}, we have
\[(u_{1,r},u_{1,z}) \parallel (U_\xi,U_\zeta)\quad \text{and} \quad (\om_{1,r},\om_{1,z}) \parallel (\Omega_\xi,\Omega_\zeta), \]
which means that, if $U$ is a function of $\Omega$, then $u_1$ is also a function of $\omega_1$:
\begin{equation}\label{eq:necessary_condition_physical}
(u_{1,r},u_{1,z}) \parallel (\om_{1,r},\om_{1,z}).
\end{equation} 
We can thus examine the validity of the condition \eqref{eq:necessary_condition} by comparing the level sets of $u_1$ and $\omega_1$.

\begin{figure}[!ht]
\centering
    \includegraphics[width=0.32\textwidth]{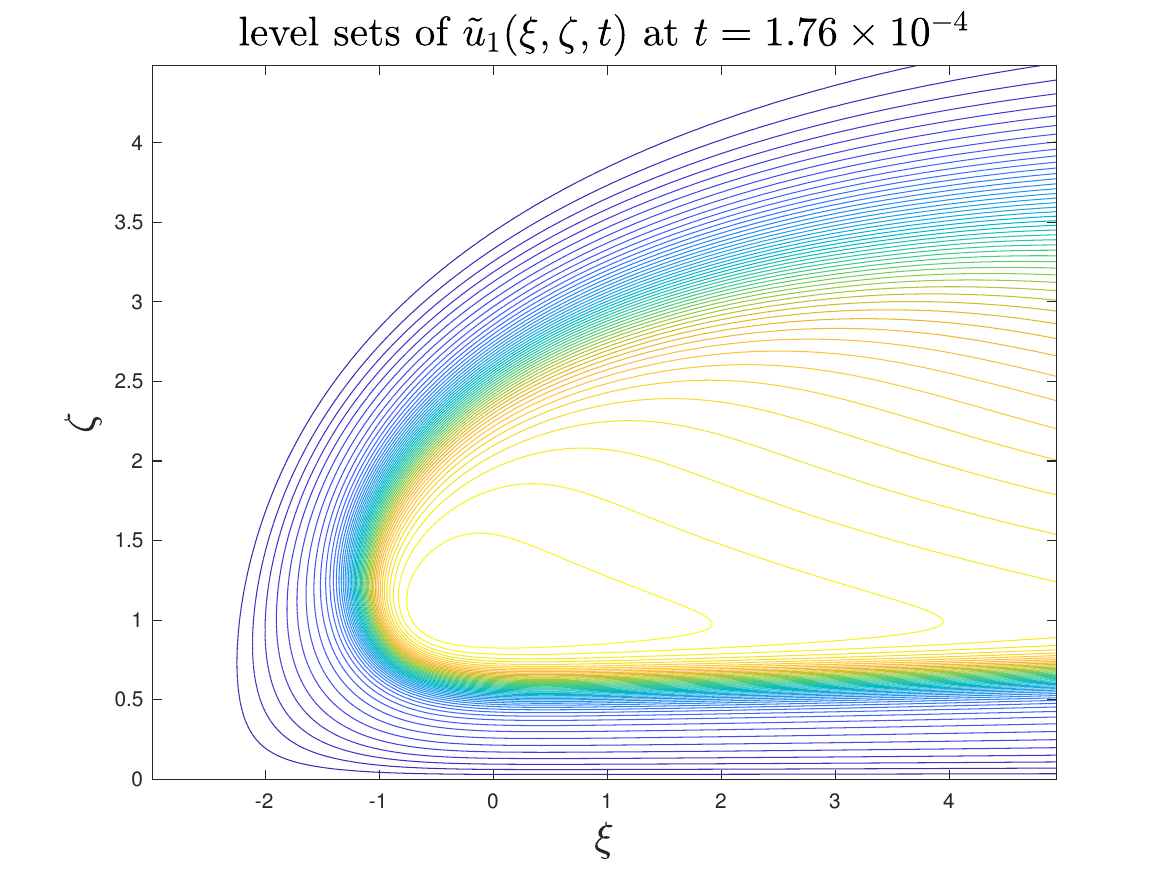}
    \includegraphics[width=0.32\textwidth]{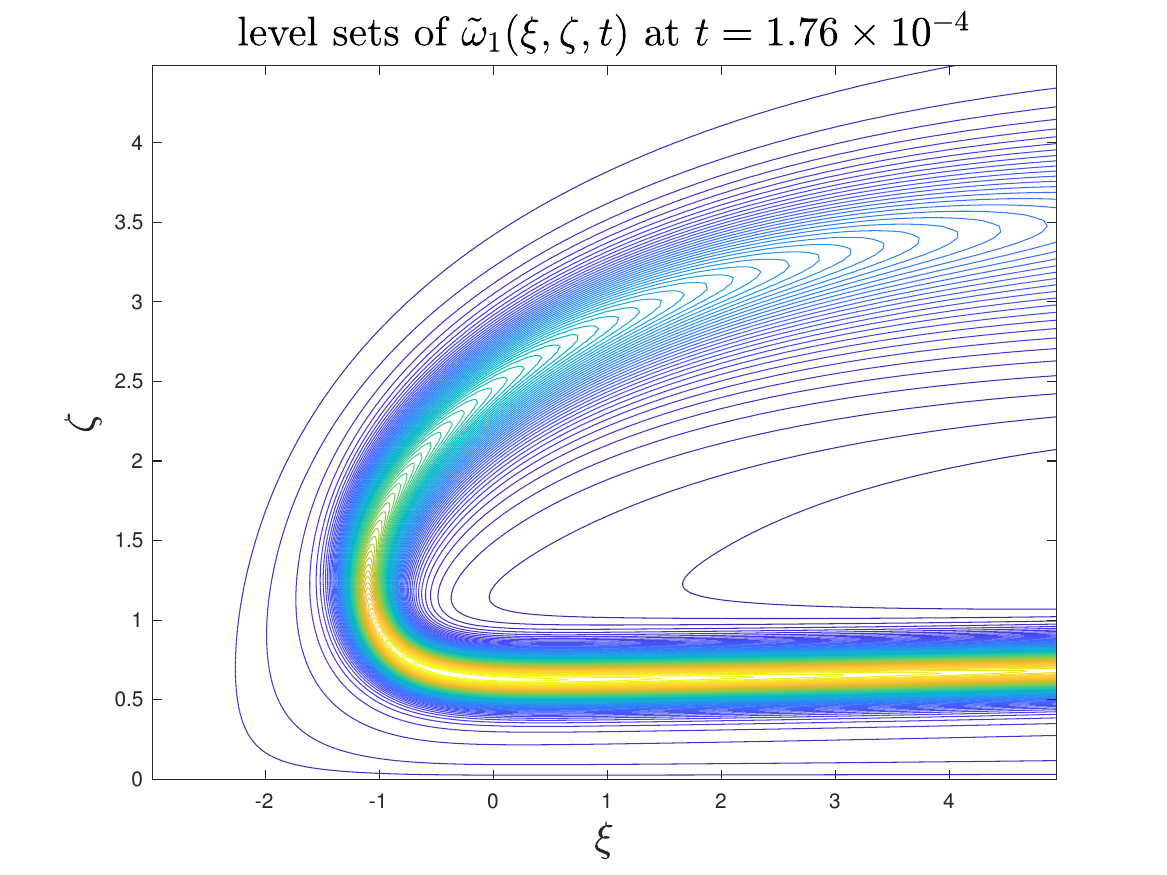} 
    \includegraphics[width=0.32\textwidth]{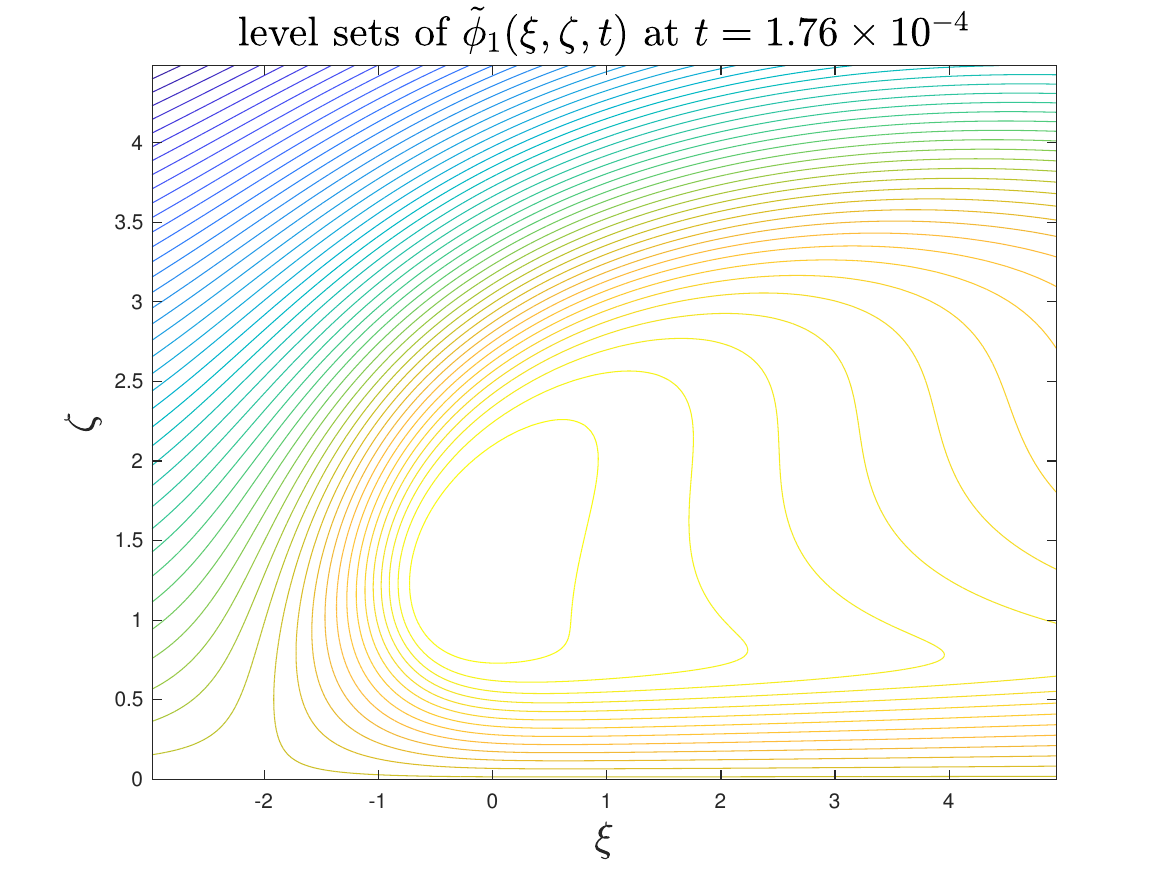} 
    \caption[Self-similar level sets]{Level sets of the stretched functions $\tilde{u}_1$ (left), $\tilde{\om}_1$ (middle) and $\tilde{\phi}_1$ (right) at $t = 1.76\times 10^{-4}$.}  
    \label{fig:selfsimilar_level_set}
       \vspace{-0.05in}
\end{figure}

Figure \ref{fig:selfsimilar_level_set} compares the level sets of the stretched and shifted functions 
\begin{align*}
\tilde{u}_1(\xi,\zeta,t) &= u_1(C_l(t)\xi + R(t),C_l(t)\zeta,t),\\
\tilde{\om}_1(\xi,\zeta,t) &= \om_1(C_l(t)\xi + R(t),C_l(t)\zeta,t),\\
\text{and}\quad\tilde{\phi}_1(\xi,\zeta,t) &= \psi_1(C_l(t)\xi + R(t),C_l(t)\zeta,t) - c_s\zeta
\end{align*}
at $t = 1.76\times 10^{-4}$ in the $\xi\zeta$-plane. These functions are the same as $U,\Omega, \Phi$ up to rescaling in magnitude. As we can see, though the contours of these functions are not exactly the same in the local neighborhood (of length scale $Z(t)$) of the point $(R(t),Z(t))$, they have surprising geometric similarities. In particular, the level sets of $\tilde{u}_1$ and $\tilde{\om}_1$ are almost parallel to each other along curved band where lies the thin structure of $\tilde{\om}_1$. We also notice that the level sets of $\tilde{\phi}_1$ are less geometrically similar to those of $\tilde{u}_1$ and $\tilde{\om}_1$, which is possibly because the rescaling constant $c_s=1/2$ is only valid in the asymptotic limit $t\rightarrow T$.

\begin{figure}[!ht]
\centering
    \includegraphics[width=0.35\textwidth]{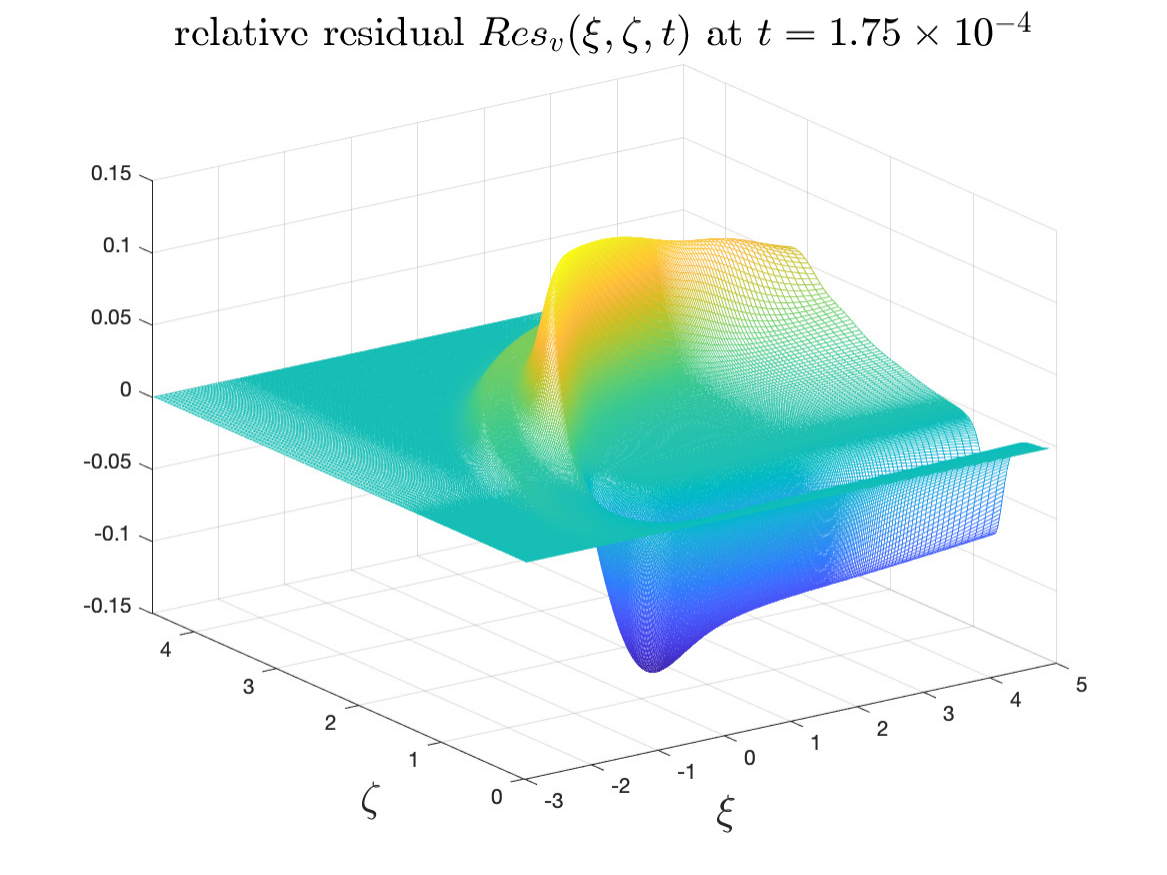}
    \includegraphics[width=0.35\textwidth]{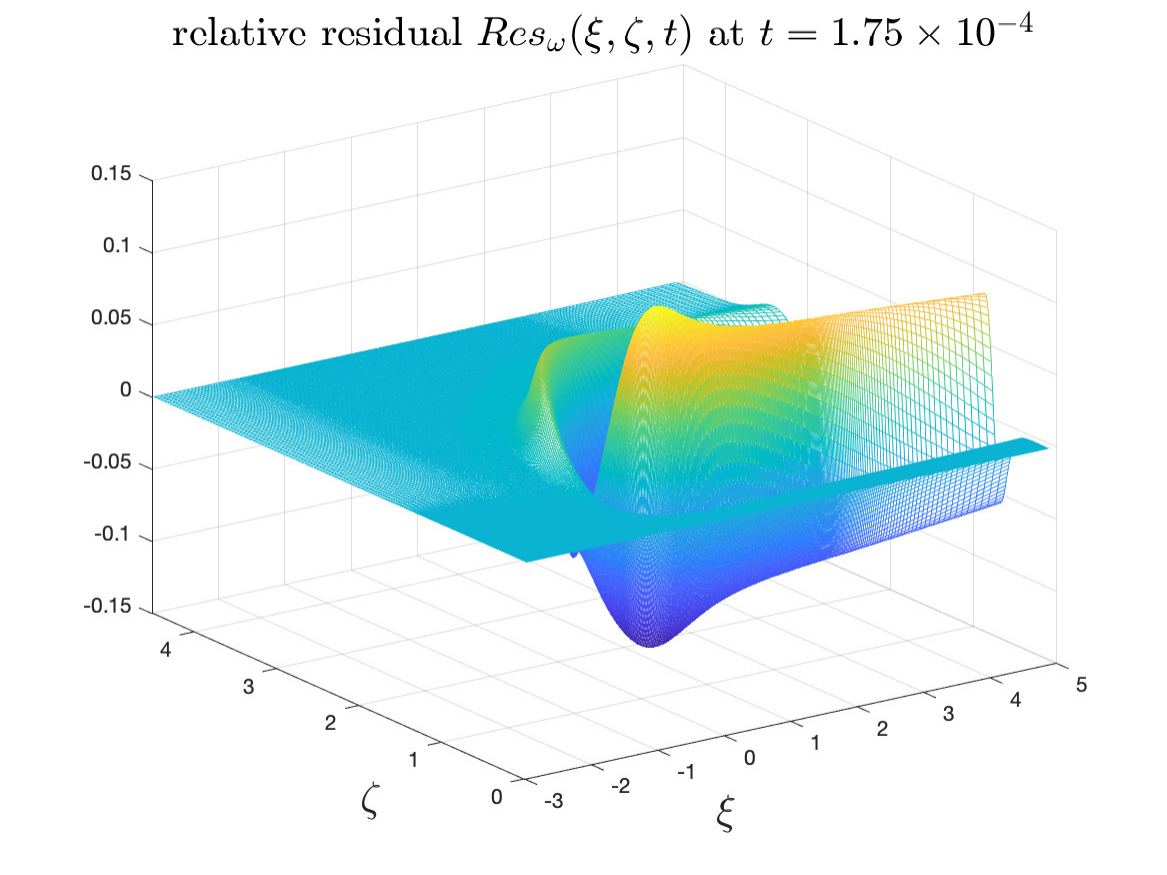}
    \includegraphics[width=0.35\textwidth]{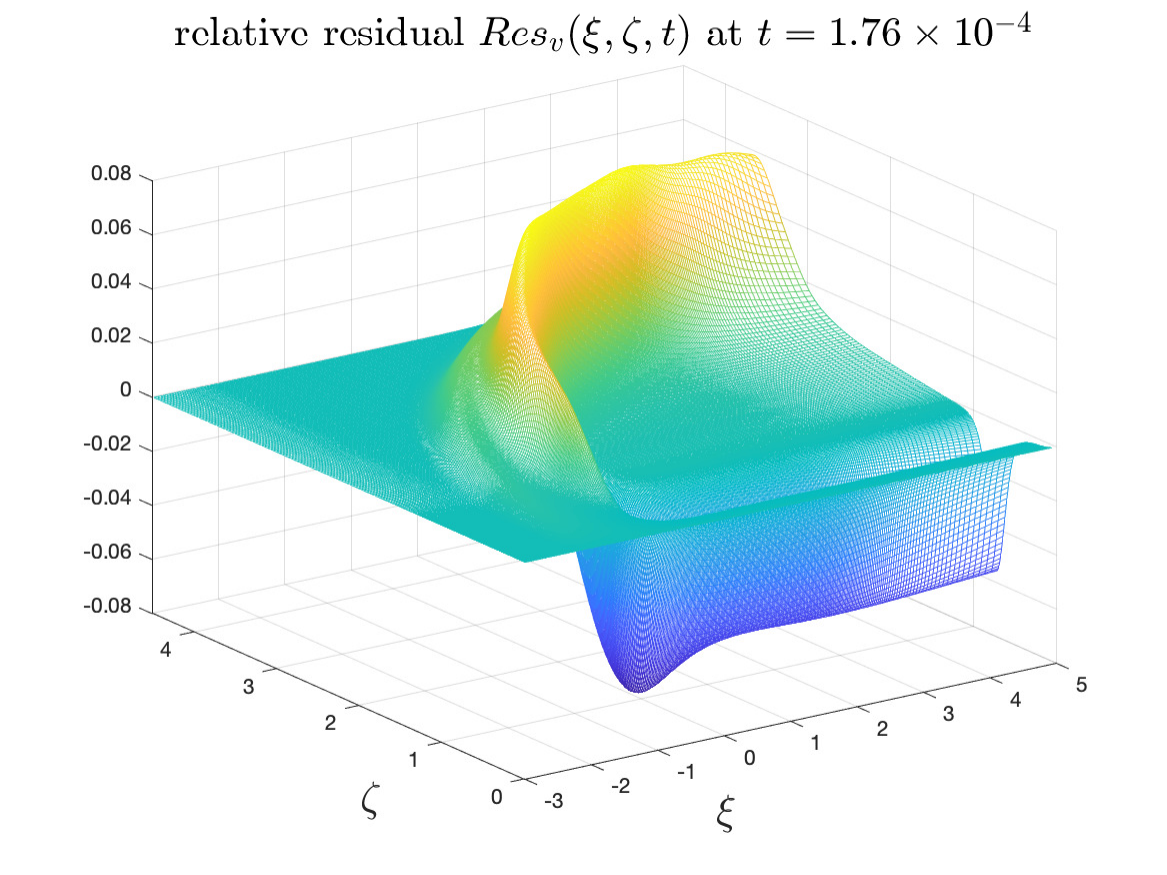}
    \includegraphics[width=0.35\textwidth]{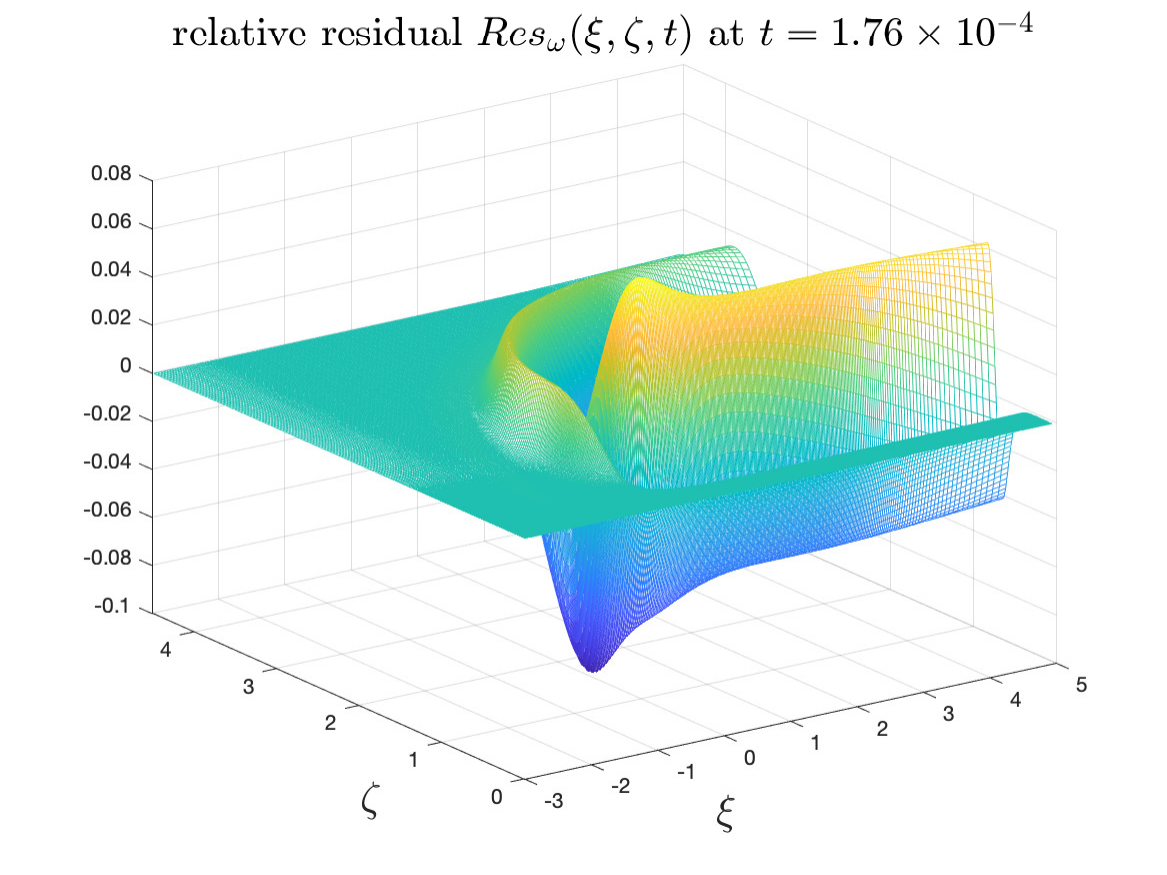}
    \caption[Relative residual]{The relative residuals $Res_{u}(\xi,\zeta,t)$ and $Res_{\om}(\xi,\zeta,t)$ at $t = 1.75\times 10^{-4}$ (upper row) and $t = 1.76\times 10^{-4}$ (lower row).}   
    \label{fig:relative_residual}
       \vspace{-0.05in}
\end{figure}

To further justify the level set condition, we investigate the (time-dependent) relative residuals of the equations in \eqref{eq:necessary_condition}. The relative residual $Res_{u}$ of the first equation is defined as
\[Res_{u}(\xi,\zeta,t) = \frac{\tilde{\phi}_{1,\xi}\tilde{u}_{1,\zeta} - \tilde{\phi}_{1,\zeta}\tilde{u}_{1,\xi}}{M_u(t)},\quad \text{where}\quad M_u(t) = \max_{\xi,\zeta}\sqrt{(\tilde{\phi}_{1,\xi}\tilde{u}_{1,\zeta})^2 + (\tilde{\phi}_{1,\zeta}\tilde{u}_{1,\xi})^2}.\]
The relative residual $Res_{\om}(\xi,\zeta,t)$ of the second equation is defined similarly. Figure \ref{fig:relative_residual} shows the profiles of $Res_{u}$ and $Res_{\om}$ at two time instants $t_1 = 1.75\times 10^{-4}$ and $t_2 = 1.76\times 10^{-4}$ in the late stage of our computation. One can see that the magnitudes of the relative residuals are reasonably small and their maximums are decreasing in time, from $(\|Res_{u}(t_1)\|_{L^\infty},\|Res_{\om}(t_1)\|_{L^\infty})=(0.1288,0.1231)$ to $(\|Res_{u}(t_2)\|_{L^\infty},\|Res_{\om}(t_2)\|_{L^\infty})=(0.0723,0.0849)$.

These numerical observations are strong evidence of the validity of the condition \eqref{eq:necessary_condition} in the critical blowup region and partial justification of our asymptotic analysis of the potential locally self-similar blowup. As mentioned above, we can also understand this interesting phenomenon from a different angle: it is because the level set condition \eqref{eq:necessary_condition_physical} is well satisfied in a local region around $(R(t),Z(t))$ that the solution can possibly develop a locally self-similar blowup in the form of \eqref{eq:self-similar_ansatz}.

\subsubsection{Local compression along the radial direction} The underlined terms in \eqref{eq:profile_equation_1} can also explain the sharp fronts in the profiles of $u_1,\om_1$ in the $r$ direction (see Figure \ref{fig:zoomin_profile}). At any time $t$ before the blowup time $T$, namely when $\tau$ is finite, the conditions \eqref{eq:profile_equation_u1_2} and \eqref{eq:profile_equation_w1_2} are not expected to be satisfied exactly. Therefore, the underlined terms in equations \eqref{eq:profile_equation_u1_1} and \eqref{eq:profile_equation_w1_1} are active convection terms that serve to shape the profiles of $U,\Omega$. The corresponding velocity field is given by 
\[U^\xi: = RC_l^{-1}(c_s - \Psi_\zeta),\quad U^\zeta = RC_l^{-1}\Psi_\xi,\]
which is formally much stronger than the velocity field corresponding to the non-underlined advection terms. This strong advection effect will transport and redistribute the profiles of $U,\Omega$ so that the conditions \eqref{eq:profile_equation_u1_2} and \eqref{eq:profile_equation_w1_2} are satisfied locally.

Let us now focus on the velocity $U^\xi$ in the $\xi$ direction (i.e. the $r$ direction). In particular, we study how the sign of $U^\xi$ changes. Suppose that the profile solution $U$ is already very close to the steady state, so that the time derivative of the maximum of $U$ is negligible. By evaluating the $U$-equation at the maximum point of $U$, namely at the point $(\xi_0,\zeta_0) = (0,Z(t)/C_l(t))$, we have 
\[(2\Psi_\zeta - c_u)U \approx - \tilde{\nu}^rU_{\xi\xi} - \tilde{\nu}^z U_{\zeta\zeta} > 0.\]
Since $ U (\xi_0,\zeta_0) > 0 $, we obtain
\[\Psi_\zeta(\xi_0,\zeta_0) \geq c_u/2 = \frac{1}{2}.\] 
As a result, the velocity $U^\xi$ is negative (traveling from right to left) around the point $(\xi_0,\zeta_0)$ since
\[U^\xi(\xi_0,\zeta_0) = RC_l^{-1}(c_s - \Psi_\zeta(\xi_0,\zeta_0)) = RC_l^{-1}\left(\frac{1}{2} - \Psi_\zeta(\xi_0,\zeta_0)\right) < 0.\]
Also, we have learned from our numerical observations that the maximum point of the cross-section $\psi_{1,z}(r,Z(t))$ (as a function of $r$) aligns with the maximum point of $u_1$. The value of $\psi_{1,z}(r,Z(t))$ decreases as $r$ goes away from $R(t)$ (see Figure \ref{fig:alignment}(a)); correspondingly, the value of $\Psi_\zeta$ (the rescaled version of $\psi_{1,z}$) will drop below $1/2$ as $|\xi|$ becomes large. As a result, the velocity $U^\xi$ is positive when $|\xi|$ is away from $0$ by a small distance. 

\begin{figure}[!ht]
\centering
    \includegraphics[width=0.4\textwidth]{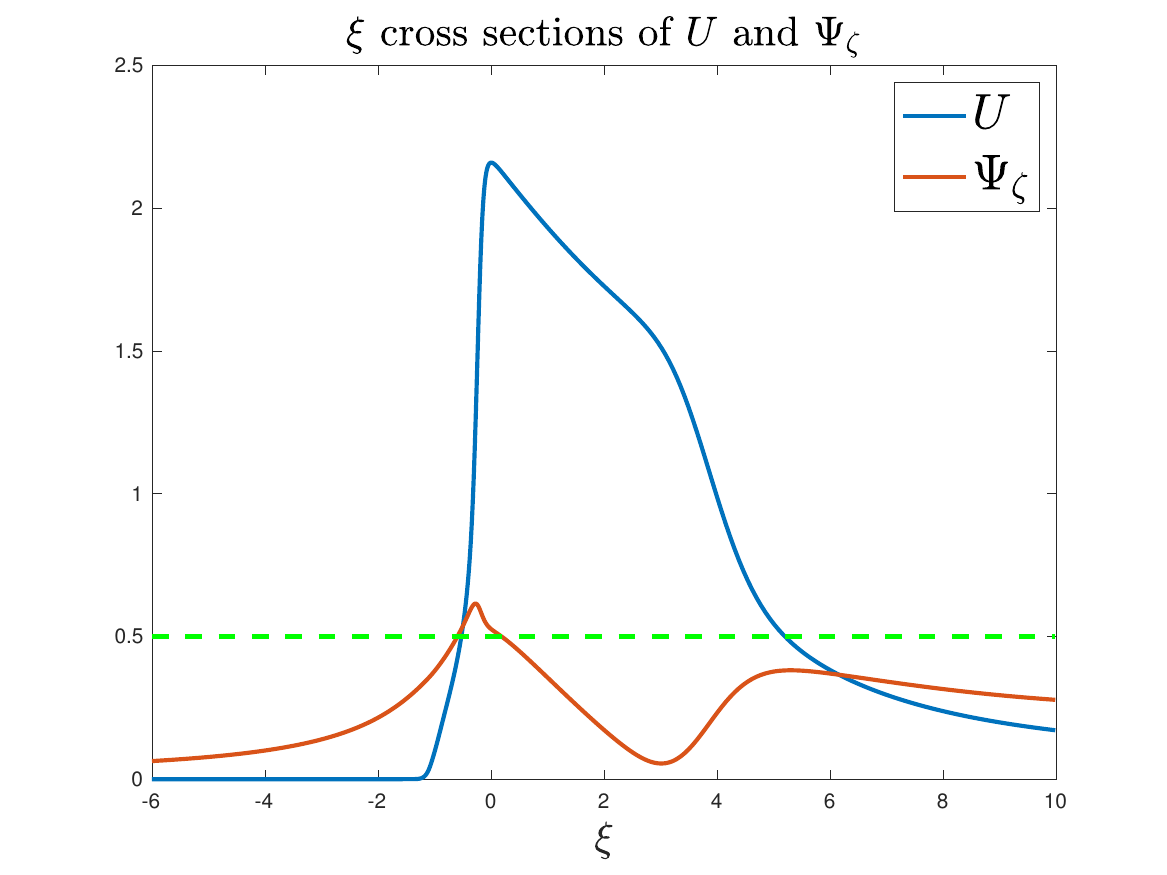}
    \includegraphics[width=0.4\textwidth]{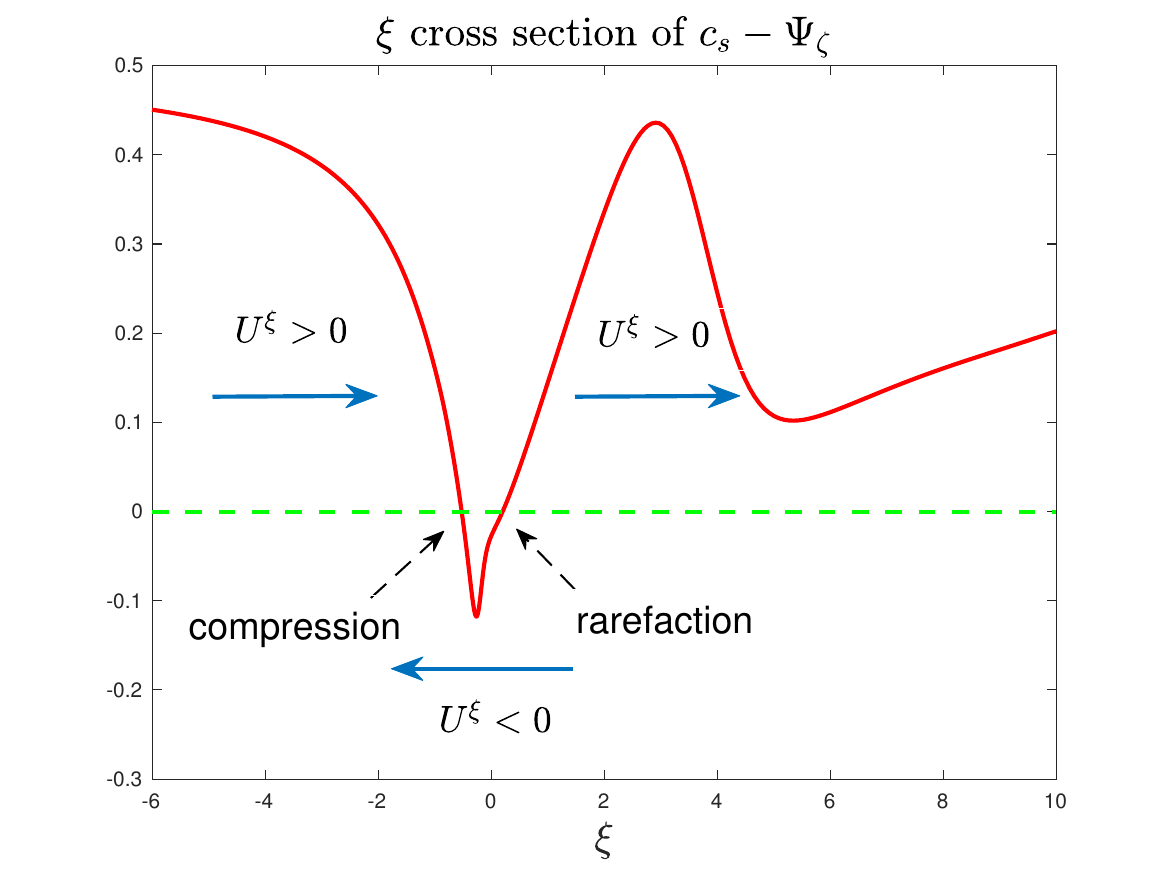} 
    \caption[Compression flow]{Left: The $\xi$ cross sections of $U$ and $\Psi_\zeta$; the green dashed line makers the level of $1/2$. Right: The $\xi$ cross section of $1/2 - \Psi_\zeta$; the green dashed line makers the level of $0$ and the blue arrows indicate the direction of the flow.} 
    \label{fig:compression_flow}
       \vspace{-0.05in} 
\end{figure}

To see this more clearly, we plot the $\xi$ cross sections of $U,\Psi_\zeta$ through the point $(\xi_0,\zeta_0)$ at physical time $t=1.7\times 10^{-4}$ in Figure \ref{fig:compression_flow}(a) and the corresponding $\xi$ cross section of $c_s-\Psi_\zeta$ in Figure \ref{fig:compression_flow}(b). The rescaled profiles $U$ and $\Psi_z$ are obtained from the numerical solution $u_1$ and $\psi_{1,z}$ via the formulas \eqref{eq:dynamic_profile} with the blowup rates given by \eqref{eq:asymptotic_powers} and the blowup time estimated by the linear fitting of $u_1$ in Section \ref{sec:growth_fitting} (i.e. $T = 1.791\times 10^{-4}$). We can see that, as we have argued, $U^\xi$ is negative near the maximum point $(0,\zeta_0)$ of $U$ and is positive for $\xi$ away from $0$. Consequently, there are two hyperbolic points in the rescaled radial velocity $U^\xi$: the left one is a compression point, while the right one is a rarefaction point. The local compressive nature of the flow along the radial direction explains why $u_1$ (or $U$) develops a sharp front to the left of its maximum point in the $r$ direction; the rarefaction wave to the right of the its maximum point explains why the local profile of $u_1$ is relatively flat on the right side of its maximum point (see Figure \ref{fig:zoomin_profile}(a)). We remark that we can observe a similar flow structure in other $\xi$ cross sections of $U^\xi$ that correspond to different values of $\zeta$.

In summary, the flow structure of the rescaled radial velocity $U^\xi$ forces the profile of $U$ to develop a two-phase structure via the underlined advection term in the equation \eqref{eq:profile_equation_u1_1}. This advection effect continues to shape the profile of $U$ until the level set condition \eqref{eq:profile_constraint_u1} is well satisfied locally. Conversely, the fact that the numerical solution $u_1$ develops a sharp front to the left of its maximum point (in the $r$ direction) and a smooth tail on the right strongly supports the validity of our preceding analysis for the self-similar blowup in a local region around $(R(t),Z(t))$.

\subsection{On the choice of viscosity coefficients}
As we have demonstrated in the previous subsection, the fact that $c_u=1$ and $c_s = 1/2$, i.e. 
\[\|u_1(t)\|_{L^\infty} \sim (T-t)^{-1},\quad R(t) \sim (T-t)^{1/2},\]
is an intrinsic property of the equations \eqref{eq:axisymmetric_NSE_1}, which does not depend on the choice of the viscosity coefficients (though we have used that $\nu^r_r\geq0$, but it is not essential). However, the result that $c_l= 1$ and its consequences (such as $c_\om=2$) rely on the particular asymptotic behavior of the degenerate viscosity coefficients:
\[\nu^r\sim\nu^z \sim R(t)^2 \sim (T-t)^1,\]
which does not seem to be essential to the potential blowup. 

We have also solved the initial-boundary value problem \eqref{eq:axisymmetric_NSE_1}--\eqref{eq:initial_data} with stronger or weaker viscosity but failed to observe convincing evidences of a sustainable (self-similar) blowup. Since the smooth viscosity coefficients are even functions of $r$ and $z$ (see Section \ref{sec:settings}), we can only choose $\nu^r,\nu^z$ to have the asymptotic behavior 
\[\nu^r,\nu^z = O(r^{2p}) + O(z^{2q})\]
for some integers $p,q$, where we have ignored the time-dependent part of $\nu^r,\nu^z$. If the viscosity is too strong ($p,q = 0$), the two-scale feature cannot survive and there is no blowup observed, as we have reported in Section \ref{sec:original_NSE}. If the viscosity is too weak or if there is no viscosity, the solution quickly becomes very unstable in the early stage of the computation before a stable self-similar blowup can be observed. Even in the case of $p = q = 1$, the constants in front of $r^2$ and $z^2$ need to be chosen carefully so that the viscosity is strong enough to control the mild oscillations in the tail region. 

It would be interesting to investigate whether the scale $c_l=1$ is intrinsic to the two-scale singularity or it is determined by the order of degeneracy $O(r^2)+O(z^2)$ in our variable viscosity coefficients. 
We plan to investigate this question by applying degenerate viscosity coefficients with different orders of degeneracy at the origin. As in the case of $p=q=1$, the constants in front of $r^{2p}$ and $z^{2q}$ need to be chosen so that the viscosity is strong enough to control the mild oscillations in the tail region. We shall leave this question to our future work.

\section{Concluding Remarks}
\label{sec:conclusion}

In this paper, we presented strong numerical evidences that the axisymmetric Euler equations with degenerate variable viscosity coefficients develop a finite-time singularity at the origin.
An important feature of this potential singularity is that the solution develops a two-scale traveling wave solution that travels towards the origin. The antisymmetric vortex dipole and the odd symmetry in the initial data generate a strong shear flow that pushes the solution travels towards the symmetry plane $z=0$ rapidly. The flow is then transported towards the symmetry axis $r=0$ by the strong negative radial velocity induced by the vortex dipole. 
The hyperbolic flow structure near the center of the traveling wave generates a no-spinning region near the symmetry axis within which the angular velocity is almost zero.
The special design of our initial data and the dynamic formation of this no-spinning region generate a positive feedback loop that enforces strong nonlinear alignment in vortex stretching, leading to a potential locally self-similar blowup at the origin. We performed resolution study and asymptotic scaling analysis to provide further support of the potential locally self-similar blowup.

The degeneracy of the variable viscosity coefficients at the origin plays an essential role in stabilizing this potential singularity formation for the incompressible axisymmetric Euler equations. We have also studied the incompressible Navier--Stokes equations with constant viscosity coefficient using the same initial data. Our numerical study revealed that the constant viscosity regularizes the smaller scale of the two-scale traveling wave solution and destabilizes the nonlinear alignment in the vortex stretching term. The solution of the Navier--Stokes equation with constant viscosity coefficient behaves completely differently. We did not observe the finite-time singularity formation that we observed for the Euler equations with degenerate variable viscosity coefficients. 

We also performed some preliminary study of the $3$D Euler equations using the same initial data. Our study showed that the solution of the Euler equations grows even faster than the solution of the $3$D Euler equations with degenerate viscosity coefficients during the warm-up phase. However, the Euler solution quickly developed a very thin structure near the sharp front before the nonlinear vortex stretching had a chance to develop a strong alignment. Without viscous regularization, the thickness of the sharp front collapses to zero faster than $Z(t)$. Thus, the solution of the $3$D Euler equations seems to develop a $3$-scale structure, which is extremely difficult to resolve numerically. We presented some preliminary numerical evidence that seems to indicate that the $3$D Euler equations may develop a potentially singular behavior in a way similar to that of the $3$D Euler equations with degenerate viscosity coefficients. However, without viscous regularization, we were not able to produce convincing numerical evidence for the potential blowup of the $3$ Euler equations. On the other hand, by applying a first order numerical viscosity plus a time-dependent vanishing viscosity of order $O(R(t)^2)+O(Z(t)^2)$, we obtained strong numerical evidence that the $3$D Euler equations develop a finite-time singularity with scaling properties similar to those of the Euler equations with degenerate viscosity coefficients, see \cite{houhuang2022} for more discussion.

Our current computation still suffers from two limitations. The first one is that our adaptive mesh strategy would break down in the limit of $R(t)/Z(t) \rightarrow \infty$. We need to develop a more effective adaptive mesh strategy to overcome this difficulty. The second one is that as we approach the singularity time, the shearing induced oscillations in the tail region become more severe. We need to apply stronger filtering in the tail region to control these oscillations in the tail region, which compromises the accuracy of our computation unless we use a very fine mesh. We hope to address these limitations in our future work.

\vspace{0.2in}
{\bf Acknowledgments.} This research was in part supported by NSF Grants DMS-1907977 and DMS-1912654. DH is supported by the National Key R\&D Program of China under the grant 2021YFA1001500, and he gratefully acknowledges the supports from the Choi Family Postdoc Gift Fund. We have benefited a lot from the AIM SQarRE ``Towards a $3$D Euler singularity'', which has generated many stimulating discussions related to the $3$D Euler singularity.

\appendix

\section{The Numerical Methods}\label{apdx:numerical_methods}
In this appendix, we describe the numerical method that we use to compute the potential finite-time singularity of the initial-boundary value problem \eqref{eq:axisymmetric_NSE_1}--\eqref{eq:initial_data}. In our scenario, the solution profiles shrink in space and quickly develop complex geometric structures, which makes it extremely challenging to compute the numerical solution accurately. In order to overcome this difficulty, we have designed in Appendix \ref{apdx:adaptive_mesh} a special adaptive mesh strategy to resolve the singularity formation near the origin $(r,z) = (0,0)$. 
In Section \ref{apdx:regularization}, we introduce the numerical regularization techniques that we use to control the mild oscillations in the solution, especially in the late stage of our computation. The overall numerical algorithm is outlined in Appendix \ref{apdx:overall_algorithm}.

\subsection{The adaptive mesh algorithm}\label{apdx:adaptive_mesh} 
To effectively and accurately compute the potential blowup, we have designed a special meshing strategy that is dynamically adaptive to the more and more singular structure of the solution. The adaptive mesh covering the half-period computational domain $\mathcal{D}_1 = \{(r,z):0\leq r\leq 1,0\leq z\leq 1/2\}$ is characterized by a pair of analytic mesh mapping functions
\[r = r(\rho),\quad \rho\in [0,1];\quad z = z(\eta),\quad \eta\in[0,1].\]
These mesh mapping functions are both monotonically increasing and infinitely differentiable on $[0,1]$, and satisfy 
\[r(0) = 0,\quad r(1) = 1,\quad z(0) = 0,\quad z(1) = 1/2.\]
In particular, we construct these mapping functions by designing their derivatives
\[r_\rho = r'(\rho), \quad z_\eta = z'(\eta),\] 
using analytic functions that are even at $0$. The even symmetries ensure that the resulting mesh can be smoothly extended to the full-period cylinder $ \{(r,z):0\leq r\leq 1,-1/2\leq z\leq 1/2\}$. The mapping functions contain a small number of parameters, which are dynamically adjusted to the solution. Once the mesh mapping functions are constructed, the computational domain is covered with a tensor-product mesh:
\begin{equation}\label{eq:mesh}
\mathcal{G} = \{(r_i,z_j): 0\leq i\leq n,\ 0\leq j\leq m\},
\end{equation}
where $r_i = r(ih_\rho),\quad h_\rho = 1/n;\quad z_j = z(jh_\eta),\quad h_\eta = 1/m.$
The precise definition and construction of the mesh mapping functions are described in Appendix \ref{apdx:adaptive_mesh_construction}. 

Figure \ref{fig:map_density} gives an example of the densities $r_\rho,z_\eta$ (in log scale) we use in the computation. As we will see in Section \ref{sec:profile_evolution}, the blowup solution of the initial-boundary value problem \eqref{eq:axisymmetric_NSE_1}--\eqref{eq:initial_data} develops a sharp front in the $r$ direction and a thin profile of two spatial scales: (i) the smaller scale corresponds to the inner profile near the sharp front, and (ii) the larger scale corresponds to the outer profile of the traveling wave. Correspondingly, we design the densities $r_\rho,z_\eta$ to have three phases: 
\begin{itemize}
\item Phase $1$ covers the inner profile of the smaller scale near the sharp front;
\item Phase $2$ covers the outer profile of the larger scale of the traveling wave;
\item Phase $3$ covers the (far-field) solution away from the symmetry axis $r=0$.
\end{itemize} 
Moreover, we have observed that there is a no-spinning region between the sharp front and the symmetry axis $r=0$ where the solutions $u_1,\om_1$ are smooth and extremely small. Hence, we add a phase $0$ in the density $r_\rho$ of the $r$ direction to cover this region. In our computation, the number (percentage) of mesh points in each phase are fixed, but the physical location of each phase will change in time, dynamically adaptive to the structure of the solution. Between every two neighboring phases, there is also a smooth transition region that occupies a fixed percentage of mesh points. 

\begin{figure}[!ht]
\centering
    \includegraphics[width=0.40\textwidth]{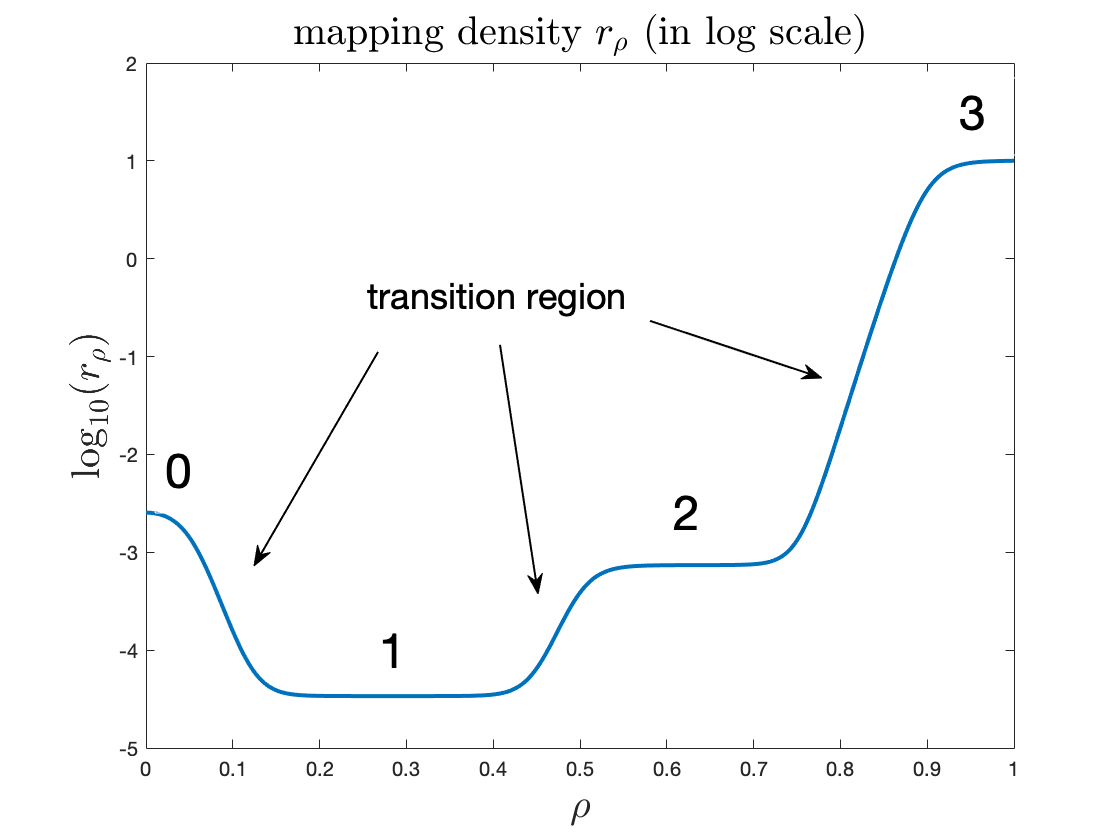}
    \includegraphics[width=0.40\textwidth]{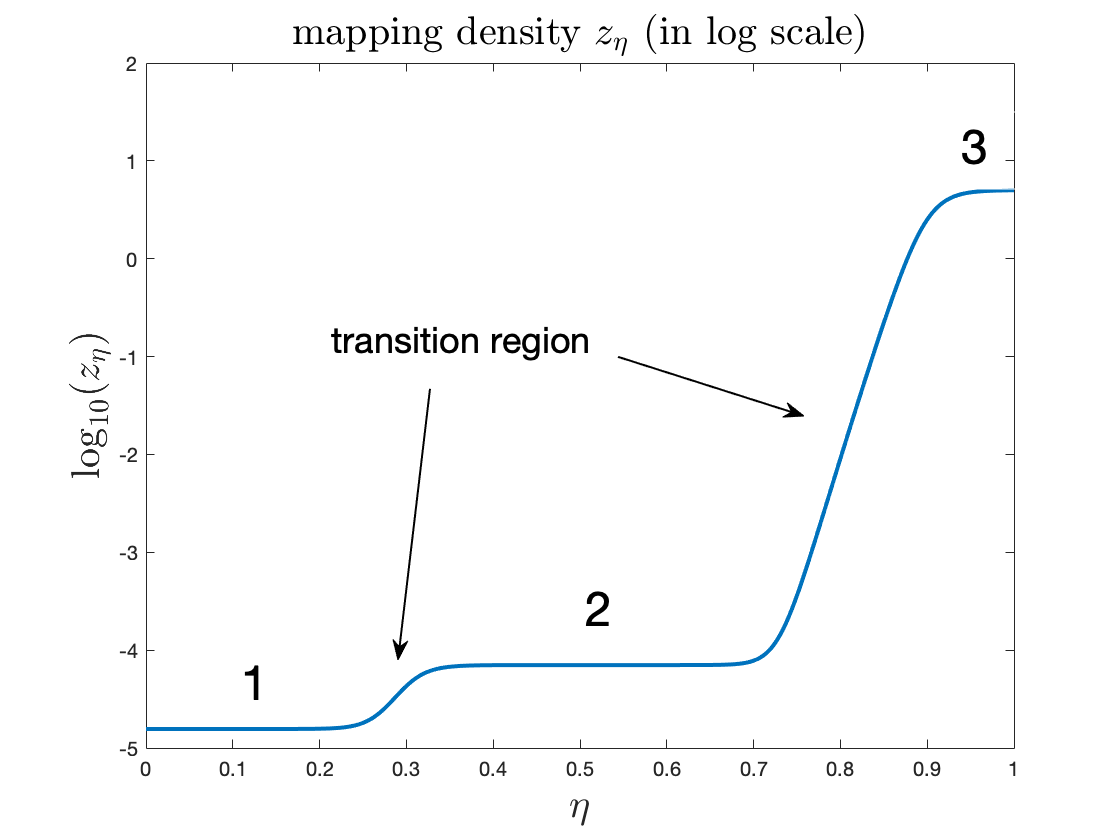}
    \caption[Map densities]{The mapping densities $r_\rho$ (left) and $z_\eta$ (right) with phase numbers labeled.}  
    \label{fig:map_density}
       \vspace{-0.05in}
\end{figure}

The adaptive mesh is evolved using the following procedure. Given a reference time $t_0$, we first use the solution at $t_0$ to construct an adaptive mesh $\mathcal{G}_0$. To do so, we need to construct the mapping densities $r_\rho,z_\eta$ that are adaptive to the reference solution. We find the maximum point $(R(t_0),Z(t_0))$ of $u_1(r,z,t_0)$ and find the coordinate $R_r(t_0)$ such that 
\[u_{1,r}(R_r(t_0),Z(t_0),t_0) =\max_{0\leq r\leq 1}u_{1,r}(r,Z(t_0),t_0).\]
We use $d_r(t_0) := R(t_0) - R_r(t_0)$ as the feature size of the smaller scale ($R(t_0)$ is always greater than $R_r(t_0)$ in our computation). Let $[\rho_1,\rho_2]\subset [0,1]$ denote the interval of phase $1$ in $r_\rho$, and let $\rho_3$ denote the right end point of phase $2$. We then choose the parameters in $r_\rho$ such that 
\begin{itemize}
\item the most singular region $[R_r(t_0) - d_r(t_0), R(t_0) + d_r(t_0)]$ of the solution in the $r$ direction is contained in the interval $[r(\rho_1),r(\rho_2)]$, and 
\item the larger scale $[0,3R(t_0)]$ is within the outer boundary of phase $2$, i.e. $3R(t_0)\leq r(\rho_3)$.
\end{itemize}
Let $[0,\eta_2]\subset [0,1]$ denote the interval of phase $1$ in $z_\eta$, and let $\eta_3$ denote the right end point of phase $2$. Similar, we choose the parameters in $z_\eta$ such that 
\begin{itemize}
\item the most singular region $[0, (3/2)Z(t_0)]$ of the solution in the $z$ direction is contained in the interval $[0,z(\eta_2)]$, and 
\item the larger scale $[0,15Z(t_0)]$ is within the outer boundary of phase $2$, i.e. $15Z(t_0)\leq z(\eta_3)$.
\end{itemize}
With the adaptive mesh $\mathcal{G}_0$ constructed using the reference solution at time $t_0$, we compute the solution of equations \eqref{eq:axisymmetric_NSE_1} forward in time, until at some time $t_1$ at least one of the following happens:
\begin{itemize}
\item[(i)] the region $[R_r(t_1) - s_r(t_1), R(t_1) + s_r(t_1)]$ is no longer contained in the interval $[r(\rho_1),r(\rho_2)]$; 
\item[(ii)] the region $[0, (3/2)Z(t_1)]$ is no longer contained in the interval $[0,z(\eta_2)]$;
\item[(iii)] the number of mesh points in the region $[R_r(t_1), R(t_1)]$ drops below a certain value;
\item[(iv)] the number of mesh points in the region $[0, Z(t_1)]$ drops below a certain value.
\end{itemize}
Then we use $t_1$ as the new reference time and construct a new adaptive mesh $G_1$ using the solution at $t_1$ with the same strategy. The solution is interpolated from $G_0$ to $G_1$ in the $\rho\eta$-space using a $4$th-order piecewise polynomial interpolation in both directions. We denote this interpolation operation by 
\begin{equation}\label{eq:interpolation}
\hat{f} = \text{IP}_4(f;\mathcal{G}_0,\mathcal{G}_1),
\end{equation}
where $f$ is a discrete function on $\mathcal{G}_0$ while $\hat{f}$ is a discrete function on $\mathcal{G}_1$. In what follows, we will always use the $4$th-order piecewise polynomial interpolation in the $\rho\eta$-space when we need to interpolate the solution from one mesh to another. In summary, the computation is carried out on a mesh that is constructed with respect to the reference solution at the latest reference time, and the reference time is replaced by the current time when one of the mesh update criteria is satisfied.

We remark that the mesh update criteria described above are designed to ensure that the most singular part of the solution is covered and well resolved in the finest phase $1$ region of the adaptive mesh. In our computation, we choose the parameters in $r_\rho,z_\eta$ such that the partitioning of different phases is given by (as shown in Figure \ref{fig:map_density})
\begin{align*}
r_\rho&:\quad \text{phase $0$}: [0,0.1],\quad \text{phase $1$}: [0.1,0.5],\quad \text{phase $2$}: [0.5,0.85],\quad \text{phase $3$}: [0.85,1],\\
z_\eta&:\quad \text{phase $1$}: [0,0.3],\quad \text{phase $2$}: [0.3,0.85],\quad \text{phase $3$}: [0.85,1].
\end{align*}
Therefore, the most singular part of the solution is covered in $[0.05,0.5]\times [0,0.3]$ in the $\rho\eta$-space.

The mesh adaptation strategy described above has several advantages compared with the conventional moving mesh method. First of all, it can automatically resolve a blowup solution regardless of its spatial scalings. This is crucial to the success of our computation, because the solution evolves rapidly in time. Secondly, the method always places enough mesh points outside the most singular region, ensuring a well-behaved and stable mesh. Thirdly, the analytic representation of the mapping densities ensures accurate approximations of space derivatives, greatly improving the quality of the computation.

\subsection{The B-Spline based Galerkin Poisson solver}\label{apdx:poisson_solver}
One crucial step in our computation is to solve the Poisson problem \eqref{eq:as_NSE_1_c} accurately. The Poisson solver we use should be compatible to the our adaptive mesh setting. Moreover, the finite dimensional solver needs to be easy to construct from the mesh, as the mesh is updated frequently in our computation. For these reasons, we choose to implement the Galerkin finite element method based on a tensorization of B-spline functions. In fact, this framework has been used by Luo--Hou \cite{luo2014toward} to compute the potential singularity formation of the $3$D axisymmetric Euler equations on the solid boundary. The equations they computed are exactly the equations \eqref{eq:axisymmetric_NSE_1} with zero viscosity ($\nu=0$). In their computation, an adaptive mesh strategy was also adopted for resolving a focusing blowup solution, and the B-Spline based Galerkin Poisson solver was found to be very effective and stable.

We will follow the framework and notations in \cite{luo2014toward}. The Poisson equation \eqref{eq:as_NSE_1_c} is solved in the $\rho\eta$-space in the following procedure. First, the equation is rewritten in the $\rho\eta$-coordinates:
\[-\frac{1}{r^3r_\rho}\left(r^3\frac{\psi_\rho}{r_\rho}\right)_\rho - \frac{1}{z_\eta}\left(\frac{\psi_\eta}{z_\eta}\right)_\eta = \om,\quad (\rho,\eta)\in [0,1]^2,\]
where for clarity we have written $\psi$ for $\psi_1$ and $\om$ for $\om_1$ (only in this section). Next, the equation is multiplied by $r^3r_\rho z_\eta\phi$ for some suitable test function $\phi\in V$ (to be defined below) and is integrated over the domain $[0,1]^2$. After integration by part, this yields the weak formulation of the Poisson equation \eqref{eq:as_NSE_1_c}: find $\psi\in V$ such that 
\begin{align*}
a(\psi,\phi) :&= \int_{[0,1]^2} \left(\frac{\psi_\rho}{r_\rho}\frac{\psi_\rho}{r_\rho} + \frac{\psi_\eta}{z_\eta}\frac{\phi_\eta}{z_\eta}\right)r^3r_\rho z_\eta \idiff \rho \idiff \eta \\
&= \int_{[0,1]^2} \om \phi r^3 r_\rho z_\eta \idiff \rho \idiff \eta =: f(\phi),\quad \forall \phi\in V,
\end{align*}
where the function space $V$ is given by (recall the symmetry properties of $\psi$) 
\begin{equation}\label{eq:function_space}
\begin{split}
V = \text{span}\Big\{ \phi\in H^1([0,1]^2): &\phi(-\rho,\eta) = \phi(\rho,\eta),\ \phi(1,\eta) = 0,\\
&\phi(\rho,-\eta) = -\phi(\rho,\eta),\ \phi(\rho,1-\eta) = - \phi(\rho,1+\eta)\Big\}.
\end{split}
\end{equation}

To introduce the Galerkin approximation, we define the finite-dimensional subspace of weighted uniform B-splines (H\"ollig \cite{hollig2003finite}) of even order $k$:
\[V_h := V_{w,h}^k = \text{span}\left\{w(\rho)b_{i,h_\rho}^k(\rho)b_{j,h_\eta}^k(\eta)\right\}\cap V,\]
where 
$w(\rho) = 1-\rho^2,\rho\in [0,1]$
is a nonnegative weight function that is vanishing at $\rho=1$, and 
\[b_{l,h}^k(s) = b^k\big((s/h) - (l-k/2)\big),\quad s\in[0,1],\ l\in \mathbb{Z},\]
is the shifted and rescaled uniform B-spline of order $k$. Here $b^k$ are the standard B-spline functions that satisfy the recursion (\cite{hollig2003finite})
\begin{equation*}
b^1(x) = \left\{\begin{array}{ll}
1, & x\in[0,1]\\
0, & x\notin[0,1]
\end{array}\right.; \quad b^k(x) = \int_{x-1}^xb^{k-1}(y)\idiff y,\quad k\geq 2.
\end{equation*}
The Galerkin finite element method then finds $\psi_h\in V_h$ such that 
\begin{equation}\label{eq:variational_form}
a(\psi_h,\phi_h) = f(\phi_h),\quad \forall \phi_h\in V_h.
\end{equation}
With suitably chosen basis functions of $V_h$, this leads to to a symmetric, positive definite linear system $Ax = b$ which can be solved to yield the Galerkin solution $\psi_h$. In particular, we choose the finite element basis of $V_h$ to be
\[B^w_{i,j}(\rho,\eta) := w(\rho)B_i(\rho)B_j(\eta),\quad 0\leq i\leq n + k/2 -1,\quad 1\leq j\leq m-1,\]
where 
\[B_i(\rho) = \frac{b^k_{i,h_\rho}(\rho)+ b^k_{i,h_\rho}(-\rho)}{1+\delta_{i0}},\quad B_j(\eta) = b^k_{j,h_\eta}(\eta) - b^k_{j,h_\eta}(-\eta) - b^k_{j,h_\eta}(2 - \eta). \]
If we write the finite element solution as 
\[\psi_h = \sum_{i,j}c_{ij}B^w_{ij}(\rho,\eta),\]
the variational problem \eqref{eq:variational_form} is transformed to an equivalent linear system of the coefficient $\{c_{ij}\}$:
\[\sum_{ij}a(B^w_{i'j'},B^w_{ij})c_{ij} = f(B^w_{i'j'}),\quad \text{for all $0\leq i'\leq n + k/2 -1, 1\leq j'\leq m-1$}.\]
In our implementation, the entries $a(B^w_{i'j'},B^w_{ij})$ of the stiffness matrix and the entries $f(B^w_{i'j'})$ of the load vector are approximated using the composite 6-point Gauss--Legendre quadrature rule. The resulting sparse linear system is then solved using the built-in linear solver of MATLAB R2019a. 

In our computation, we choose $k=2$. It is observed that the corresponding Poisson solver is at least $2$nd-order convergent for the stream function $\psi_1$. We will explain why we choose to use a $2$nd-order method in Appendix \ref{apdx:regularization}. For more detailed error analysis of the B-Spline based Galerkin Poisson solver, one may refer to \cite{hollig2003finite,luo2014toward}.

\subsection{Spatial derivatives}\label{apdx:derivative} To solve the equations \eqref{eq:axisymmetric_NSE_1} in time, we also need to approximate the spatial derivatives of the solution. For the reason to be explained in Appendix \ref{apdx:regularization}, we choose to use $2$nd-order centered difference schemes with respect to the $\rho\eta$-coordinates. Let $v$ be some solution variable. For example, the first and the second derivatives of $v$ in $r$ are approximated as
\begin{subequations}\label{eq:centered_difference}
\begin{align}
(v_r)_{i,j} = \frac{(v_\rho)_{i,j}}{(r_\rho)_i} &\approx \frac{1}{(r_\rho)_i}\cdot\frac{v_{i+1,j}-v_{i-1,j}}{2h_\rho},\\
(v_{rr})_{i,j} = \frac{(v_{\rho\rho})_{i,j}}{(r_\rho)_i^2} - \frac{(r_{\rho\rho})_{i,j}(v_\rho)_{i,j}}{(r_\rho)_i^3} &\approx \frac{1}{(r_\rho)_i^2}\cdot\frac{v_{i+1,j}-2v_{i,j} + v_{i-1,j}}{h_\rho^2} - \frac{(r_{\rho\rho})_{i,j}}{(r_\rho)_i^3}\cdot \frac{v_{i+1,j}-v_{i-1,j}}{2h_\rho}.
\end{align}
\end{subequations} 
Note that $r_{\rho\rho}$ can be computed directly from our analytic formula of $r_\rho$ without approximation. This is crucial for the accurate evaluation of $v_r,v_{rr}$, especially in the smaller scale regions where the mesh density $r_\rho$ is close to $0$. 

The centered difference formulas described above need to be supplemented by numerical boundary conditions near $\rho,\eta=0,1$. In particular, we can make use of the symmetry properties \eqref{eq:BC} of the solution. Suppose that $v=u_1$. In the $\eta$ direction, the odd symmetry conditions 
\[v_{i,-1} = -v_{i,1},\quad v_{i,m+1} = -v_{i,m-1},\quad 0\leq i\leq n,\]
are used at $\eta = 0$ and $\eta = 1$. In the $\rho$ direction, the even symmetry 
\[v_{-1,j} = v_{1,j},\quad 0\leq j\leq m,\]
is used at $\rho=0$, and the extrapolation condition 
\[v_{n+1,j} = 3v_{n,j} - 3v_{n-1,j} + v_{n-2,j},\quad 0\leq j\leq m,\]
is applied near the solid boundary $\rho=1$ when necessary. The extrapolation condition is known to be GKS stable for linear hyperbolic problems (Gustafsson et al., \cite{gustafsson1995time}), and is expected to remain stable when applied to the NSE equations as long as the underlying solution is sufficiently smooth. Moreover, the no-slip boundary condition \eqref{eq:no-slip1} will also be enforced: 
\[(u_1)_{n,j} = 0,\quad (\om_1)_{n,j} = -(\psi_{1,rr})_{n,j},\quad 0\leq j\leq m.\]

We remark that the function $\psi_1$ is very smooth in the $\rho\eta$-coordinates since it is recovered from $\om_1$ through the Poisson equation \eqref{eq:as_NSE_1_c}. Correspondingly, the numerical solution of $\psi_1$ obtained from the $2$nd-order finite element method in Appendix \ref{apdx:poisson_solver} will have a smooth error expansion. The expected order of accuracy thus will not be compromised when computing the derivatives of $\psi_1$ using the $2$nd-order centered difference schemes \eqref{eq:centered_difference}. In particular, the computation of the velocities $u^r,u^z$ (via formulas \eqref{eq:as_NSE_1_d}) should also be $2$nd-order accurate, which has been verified in Section \ref{sec:resolution_study}.

\subsection{Numerical regularization}\label{apdx:regularization}
The potential blowup solution we compute develops a long thin tail structure, stretching from the sharp front to the far field. This tail structure develops some mild oscillations in the late stage of the computation, and may affect the blowup mechanism. Therefore, we have chosen to apply some numerical regularization techniques to stabilize the solution, especially its tail part. 

Consider the computational mesh $\mathcal{G}$ \eqref{eq:mesh} that is constructed with respect to the latest reference time. Let $f = (f_{i,j})$ be the discretization of some solution variable. The numerical regularization technique we use is a local smoothing operation (denoted by $\text{LPF}_c$):
\begin{equation}\label{eq:LPF}
\begin{split}
\tilde{f}_{i,j} &= f_{i,j} + \frac{c(\rho_i,\eta_j)}{4}\big(f_{i-1,j} + f_{i+1,j} - 2 f_{i,j}\big),\\
\big(\text{LPF}_c(f)\big)_{i,j} &= \tilde{f}_{i,j} + \frac{c(\rho_i,\eta_j)}{4}\big(\tilde{f}_{i,j-1} + \tilde{f}_{i,j+1} - 2 \tilde{f}_{i,j}\big),
\end{split}
\end{equation}
where $c(\rho,\eta):[0,1]^2\mapsto [0,1]$ is a nonnegative weight functions. To define $\big(\text{LPF}_c(f)\big)_{i,j}$ near the boundaries $\rho,\eta=0,1$, we need to use numerical boundary conditions discussed in Appendix \ref{apdx:adaptive_mesh}. This smoothing operation can be seen as a low-pass filtering in the $\rho\eta$-space. It can also be viewed as solving a heat equation for one step in the $\rho\eta$-space. The weight function $c(\rho,\eta)$ measures how strong the smoothing effect is; the smoothing effect can have different strength at different part of the solution, depending on the value of $c(\rho,\eta)$. We will use $\text{LPF}_c^k(\cdot)$ to denote the operation of applying $k$ times of $\text{LPF}_c$. Note that this low-pass filtering will introduce a $2$nd-order error of size $O(h_\rho^2 + h_\eta^2)$.

\subsection{The overall algorithm}\label{apdx:overall_algorithm}
Given an adaptive mesh $\mathcal{G}$ and the data $(u_1,\om_1)$ defined on it, the numerical solution of equations \eqref{eq:axisymmetric_NSE_1} is advanced in time via the following procedure. 
\begin{enumerate}
\item The Poisson equation \eqref{eq:as_NSE_1_c} is solved for $\psi_1$ using the $2$nd-order B-spline based Galerkin method introduced by Luo--Hou in \cite{luo2014toward}.
\item The spacial derivatives are computed using the $2$nd-order centered difference schemes 
and $(u^r,u^z)$ is evaluated at the grid points using \eqref{eq:as_NSE_1_d}.  
\item An adaptive time stepping $\delta t$ is computed on $\mathcal{G}$ so that the CFL condition is satisfied with a suitably small CFL number (e.g. 0.1):
\begin{align*}
\delta t_1 &= 0.1\min\left\{\min_{\rho,\eta}\frac{h_\rho r_\rho}{u^r}\, ,\, \min_{\rho,\eta}\frac{h_\eta z_\eta}{u^z} \right\},\quad (\text{stability for convection})\\
\delta t_2 &= 0.1\min\left\{\min_{\rho,\eta}\frac{(h_\rho r_\rho)^2}{\nu^r}\, ,\, \min_{\rho,\eta}\frac{(h_\eta z_\eta)^2}{\nu^z}\right\},\quad (\text{stability for viscosity})\\
\delta t &= \min\{\delta t_1\, ,\,\delta t_2\}.
\end{align*}
We remark that with this choice of $\delta t$, the relative growth of the maximum value of the solution in one step is observed to remain below $1\%$.
\item The solution $(u_1,\om_1)$ is advanced in time by $\delta t$ using a $2$-stage, $2$nd-order explicit Runge--Kutta method associated with the Butcher tableau:
\[
\begin{tabular}{c|c c}
0 & 0 & 0\\
1 & 1 & 0 \\ \hline
  & 1/2 & 1/2
\end{tabular}
\]

\item The mesh $\mathcal{G}$ is updated if necessary (according to the strategy in Appendix \ref{apdx:adaptive_mesh}).
\end{enumerate}

It can be shown that the overall algorithm without numerical regularization is formally $2$nd-order accurate in space and in time. 

To control the mild oscillations in the tail region in the late stage of the computation, we will apply numerical regularizations in our computations as described in Appendix \ref{apdx:regularization}. In particular, we will do the following each time when we update {\it the time derivatives} of $u_1,\om_1$ in equations \eqref{eq:axisymmetric_NSE_1}. At the beginning of each iteration, we apply to the solution variables $u_1,\om_1$ one time low-pass filtering \eqref{eq:LPF} with uniform strength $c(\rho,\eta)=0.1$ and then one time low-pass filtering with strength function $c(\rho,\eta) = c_L(\rho,\eta)$:
\[\tilde{u}_1 = \text{LPF}_{c_L}(\text{LPF}_{0.1}(u_1)),\quad \tilde{\om}_1 = \text{LPF}_{c_L}(\text{LPF}_{0.1}(\om_1)).\]
The strength function $c_L(\rho,\eta)$ is given by 
\begin{equation}\label{eq:regularization_strength}
c_L(\rho,\eta) = \bar{f}_{sc}(\rho,0.8,0.05)\cdot\bar{f}_{sc}(\eta,0.8,0.05)\cdot \big(1 - \bar{f}_{sc}(\rho,0.45,0.05)\cdot\bar{f}_{sc}(\eta,0.45,0.05)\big),
\end{equation}
where $\bar{f}_{sc} = 1-f_{sc}$ and the soft-cutoff function $f_{sc}$ is given by \eqref{eq:cutoff}. 
The point of using this ``L'' shape function $c(\rho,\eta)$ is to only suppress the oscillations occurring at the tail part of the solution without affecting the critical region near the front. We will use $\tilde{u}_1,\tilde{\om}_1$ to compute the stream function $\psi_1$ and all the spatial derivatives. But when advancing the solution using the $2$nd-order Runge--Kutta method, we still use the original $u_1,\om_1$ as the one-step initial values. That is, 
\[(u_1(t+\delta_t),\om_1(t+\delta_t)) = (u_1(t),\om_1(t)) + F(\tilde{u}_1(t),\tilde{\om}_1(t),\delta_t),\]
where $F(\tilde{u}_1(t),\tilde{\om}_1(t),\delta_t)$ is the update of the solution in one time step. 

Moreover, after solving the Poisson equation for $\psi_1$ (using $\tilde{\om}_1)$, we regularize $\psi_1$ by doing one time of low-pass filtering with uniform strength $c(\rho,\eta)=1$:
\[\psi_1 \ \longleftarrow\ \text{LPF}_1(\psi_1).\]
Then we compute the spatial derivatives $\psi_{1,r},\psi_{1,z}$ and regularize them by doing
\begin{align*}
\psi_{1,r}\ \longleftarrow\ \text{LPF}_1(\psi_{1,r}),\quad \psi_{1,z}\ \longleftarrow\ \text{LPF}_1(\psi_{1,z}).
\end{align*}
This step is meant to regularize the velocity field in the tail part region of the solution, as the velocity is computed from $\psi,\psi_r,\psi_z$ via \eqref{eq:as_NSE_1_d}. 

\section{Construction of the Adaptive Mesh}
\label{apdx:adaptive_mesh_construction}

In this appendix, we provide the details of the construction of our adaptive meshes. As mentioned in Appendix \ref{apdx:adaptive_mesh}, we will construct the mapping densities $r_\rho,z_\eta$ using analytic functions. Then the mapping functions are given by 
\[r(\rho) = \int_0^\rho r_\rho(s)\idiff s,\quad z(\eta) = \int_0^\eta z_\eta(s)\idiff s,\]
which can be computed explicitly using the analytic formulas of $r_\rho,z_\eta$. In particular, we use the following parametrized formulas for the mapping densities
\[r_\rho(s) = p(s;\mathbf{a}^\rho,\mathbf{s}^\rho),\quad z_\eta(s) = p(s;\mathbf{a}^\eta,\mathbf{s}^\eta),\]
where $\mathbf{a}^\rho = \{a^\rho_i\}_{i=0}^3,\mathbf{a}^\eta = \{a^\eta_i\}_{i=0}^3, \mathbf{s}^\rho = \{s^\rho_j\}_{j=1}^3,\mathbf{s}^\eta = \{s^\eta_j\}_{j=1}^3$ are parameters, 
\begin{equation}\label{eq:general_density}
p(s; \mathbf{a}, \mathbf{s}) = a_1 + a_2 q_b(s-s_2) + a_3 q_b(s-s_3) + a_0(q_b(s_1-s) +q_b(s_1+s)-1) ,\quad s\in [0,1],
\end{equation}
and 
\[q_b(x) = \frac{(1+x)^b}{1+(1+x)^b},\quad x\in [-1,1],\]
for some even integer $b$. We will determine the parameters by adjusting $r_\rho,z_\eta$ to a reference solution. In what follows, we explain how to construct $r_\rho(s) = p(s;\mathbf{a},\mathbf{s})$. The construction of $z_\eta$ is similar. 

It is not hard to see that, for a large integer $b$, 
\[q_b(x)\approx \mathbf{1}_{x\geq0} := \left\{\begin{array}{ll}
1,& x\in[0,1],\\
0,& x\in[-1,0).
\end{array}\right.\]
Therefore, we can choose $b$ to be a large number and treat $q_b(x)$ ideally as the indicator function $\mathbf{1}_{x\geq0}$. With this approximation, we have
\begin{equation}\label{eq:density_approximation}
p(s; \mathbf{a}, \mathbf{s}) \approx \hat{p}(s;\mathbf{a}, \mathbf{s}) := \left\{\begin{array}{ll}
a_0+a_1,& s\in[0,s_1),\\
a_1,&  s\in[s_1,s_2),\\
a_1+a_2, & s\in[s_2,s_3),\\
a_1+a_2+a_3, & s\in(s_3,1].
\end{array}\right.
\end{equation}
Note that $\hat{p}(s;\mathbf{a}, \mathbf{s})$ is formally an even function of $s$ at $s=0$ and at $s=1$. In fact, when $b$ is large enough, $p(s; \mathbf{a}, \mathbf{s})$ can be viewed as an even function of $s$ at $s=0$ and at $s=1$ up to machine precision. In our computation, we choose $b=60$. 

The intervals $[0,s_1),[s_1,s_2),[s_2,s_3),[s_3,1]$ stand for the $4$ phases in $r_\rho(s)$. Therefore, the numbers $s_1,s_2-s_1,s_3-s_2,1-s_3$ are the percentages of mesh points that are allocated the $4$ phases respectively. Throughout our computation, we choose $s_1,s_2,s_3$ to be some constants independent of the reference solution:
\begin{align*}
s^\rho_1 &= 0.1, \quad s^\rho_2 = 0.5,\quad s^\rho_3 = 0.85,\\
s^\eta_1 &= 0, \quad s^\eta_2 = 0.3,\quad s^\eta_3 = 0.85.
\end{align*}
The case of $s^\eta_1=0$ corresponds to $z_\eta$ having only three phases. 

To determine the parameters $\mathbf{a}$, we need to use the information of the reference solution so that the mapping density $p(s;\mathbf{a},\mathbf{s})$ is adaptive to it. This is done by determining the images of the $4$ phases in the physical domain. Let $P(x; \mathbf{a}, \mathbf{s})$ be the mapping function corresponding to the density $p(s; \mathbf{a}, \mathbf{s})$, 
\begin{equation}\label{eq:mapping_function}
P(x; \mathbf{a}, \mathbf{s}) = \int_0^x p(s; \mathbf{a}, \mathbf{s}) \idiff s, \quad x\in[0,1].
\end{equation}
Let $y$ denote a physical coordinate ($y=r$ or $y=z$). We will use the reference solution to determine some coordinates $y_1,y_2,y_3$ such that we want to make
\begin{equation}\label{eq:density_constraint}
P(s_1; \mathbf{a}, \mathbf{s}) = y_1,\quad P(s_2; \mathbf{a}, \mathbf{s}) = y_2,\quad P(s_3; \mathbf{a}, \mathbf{s}) = y_3,\quad P(1; \mathbf{a}, \mathbf{s}) = L.
\end{equation}
Here $L=1$ for $y=r$ and $L=1/2$ for $y=z$. That is to say, the intervals $[0,y_1),[y_1,y_2),[y_2,y_3),[y_3,L]$ correspond to the $4$ phases in the physical domain $[0,L]$. Using the approximation \eqref{eq:density_approximation}, the constraints \eqref{eq:density_constraint} can be approximated by 
\begin{align*}
(a_0 + a_1)s_1 &= y_1,\\
a_1(s_2-s_1) &= y_2 - y_1,\\
(a_1+a_2)(s_3-s_2) &= y_3 - y_2,\\
(a_1+a_2+a_3)(1-s_3) &= L - y_3. 
\end{align*}  
This gives 
\begin{equation}\label{eq:solution_a_1}
\begin{split}
a_1 &= \frac{y_2-y_1}{s_2-s_1},\quad a_0 = \frac{y_1}{s_1} - \frac{y_2-y_1}{s_2-s_1}\\
a_2 &= \frac{y_3-y_2}{s_3-s_2} - \frac{y_2-y_1}{s_2-s_1},\quad a_3 = \frac{L-y_3}{1-s_3} - \frac{y_3-y_2}{s_3-s_2}.
\end{split}
\end{equation}
In the case $s_1 = 0$, we will choose $y_1=0$ and $a_0 = 0$. Then we have
\begin{equation}\label{eq:solution_a_2}
a_1=\frac{y_2}{s_2},\quad a_2 = \frac{y_3-y_2}{s_3-s_2} - \frac{y_2}{s_2},\quad a_3 = \frac{L-y_3}{1-s_3} - \frac{y_3-y_2}{s_3-s_2}.
\end{equation}
Since we want $p(s; \mathbf{a}, \mathbf{s})$ to be a density function, we need $a_0,a_1,a_2,a_3$ to be non-negative. This yields the constraints on $y_1,y_2,y_3$:
\[0\leq y_1\leq y_2\leq \frac{s_2}{s_1}y_1,\quad y_2 + \frac{s_3-s_2}{s_2-s_1}(y_2-y_1)\leq y_3\leq \frac{s_3-s_2}{1-s_2}L + \frac{1-s_3}{1-s_2}y_2.\]
In our computation, since the blowup is focusing near the origin, the most interesting region is far away from the boundary $L$. Therefore, we always have $y_3\ll L$, and we do not need to worry about the last constraint $y_3\leq L(s_3-s_2)/(1-s_2) + y_2(1-s_3)(1-s_2)$. 

After $a_0,a_1,a_2$ and $a_3$ have been computed, we then construct the mapping density $p(s; \mathbf{a}, \mathbf{s})$ using \eqref{eq:general_density} and the mapping function $P(x;\mathbf{a}, \mathbf{s})$ using \eqref{eq:mapping_function}. However, since the approximation \eqref{eq:density_approximation} is not exact, the right boundary condition of $P(x;\mathbf{a}, \mathbf{s})$ is only approximately satisfied: $P(1;\mathbf{a}, \mathbf{s})\approx L$. To make sure that $P(1) = L$ exactly, we need to rescale $p(s; \mathbf{a}, \mathbf{s})$ by a factor $L/P(1;\mathbf{a}, \mathbf{s})\approx1$, i.e. 
\[\tilde{a}_i = \frac{L}{P(1;\mathbf{a}, \mathbf{s})}a_i, \quad i=0,1,2,3,\]
and then use $p(s; \tilde{\mathbf{a}}, \mathbf{s})$ and $P(x; \tilde{\mathbf{a}}, \mathbf{s})$ in the real computation.

Next, we explain how to choose $y_1,y_2,y_3$ for $r_\rho$ and $z_\eta$. Again, let $(R,Z)$ denote the maximum point of the reference solution $u_1(r,z)$. This coordinate locates the most singular part of the solution. We first consider $r_\rho$, i.e. $y = r$. Since we want the most singular part of the solution to be covered in phase $1$, we should choose $r_1,r_2$ such that $[r_1,r_2]$ is a local neighborhood of the smaller scale around $R(t)$. In particular, we choose
\[r_2 = R + 3d_r,\quad r_1 = \max\big\{R_r - 4d_r\,,\,\frac{s^\rho_1}{s^\rho_2}r_2\big\},\]
where $d_r = R(t) - R_r(t)$ is the feature size of the smaller scale, and
\[R_r = \arg\max_{0\leq r\leq 1}u_{1,r}(r,Z,t)\] 
is the point where the $r$ derivative $u_{1,r}$ achieves its maximum on the cross section $\{(r,Z):0\leq r\leq 1\}$. Moreover, we impose the constraint that the outer part of solution of the larger scale does not go beyond phase $2$. Note that the larger scale of the solution is characterized by the coordinate $R$, and the tail part of the solution mainly lies in $[0,3R]$. Therefore, we choose 
\[r_3 = \max\big\{3R\,,\,r_2 + \frac{s^\rho_3-s^\rho_2}{s^\rho_2-s^\rho_1}(r_2-r_1)\big\}.\]

In the $z$ direction, the smaller scale is characterized by $Z$. Then in order for the most singular region of the smaller scale to be covered in phase $1$, we choose 
$z_1 = 3Z/2$.
We have observed that the outer part of the solution in the $z$ direction mainly lie in the range $[0,10Z]$. Therefore, for the outer part to be covered in phase $2$, we choose 
$z_2 = 15Z.$
This completes the construction of the mapping functions $r(\rho)$ and $z(\eta)$.

\section{Beyond the stable phase}\label{sec:beyond}
In this subsection, we study the solution behavior beyond the stable phase, namely for $t> 1.75\times 10^{-4}$. In particular, we investigate the mild oscillations that occur at the tail part of the solution in the very late stage of our computation, which seem to be induced by under-resolution and amplified by the vortex shedding from the strong shearing flow in the near field. Such oscillations do not cause or contribute to the potential singularity of interest; instead, they may obstruct the blowup process by compromising the nice coupling phenomena of the solution that we discussed in Section~\ref{sec:mechanism}. As an evidence of this, the important nonlinear alignment between $\psi_{1,z}$ and $u_1$ starts to trail off beyond the stable phase. 

\subsection{Oscillations in the tail region} We have mentioned earlier that we could not extend the computation all the way to the potential singularity time in Case $1$. One of the reasons is that our mesh strategy with an affordable mesh size will lose effectiveness when the contrast between the two separate scales becomes extremely huge. Another reason is that the solution triggers some mild oscillations in the late stage of our computation, particularly when $t\geq 1.76\times 10^{-4}$. According to our observations, the oscillations mainly occur along the long tail part of the solution where a strong velocity shear appears. 

\begin{figure}[!ht]
\centering
    \includegraphics[width=0.40\textwidth]{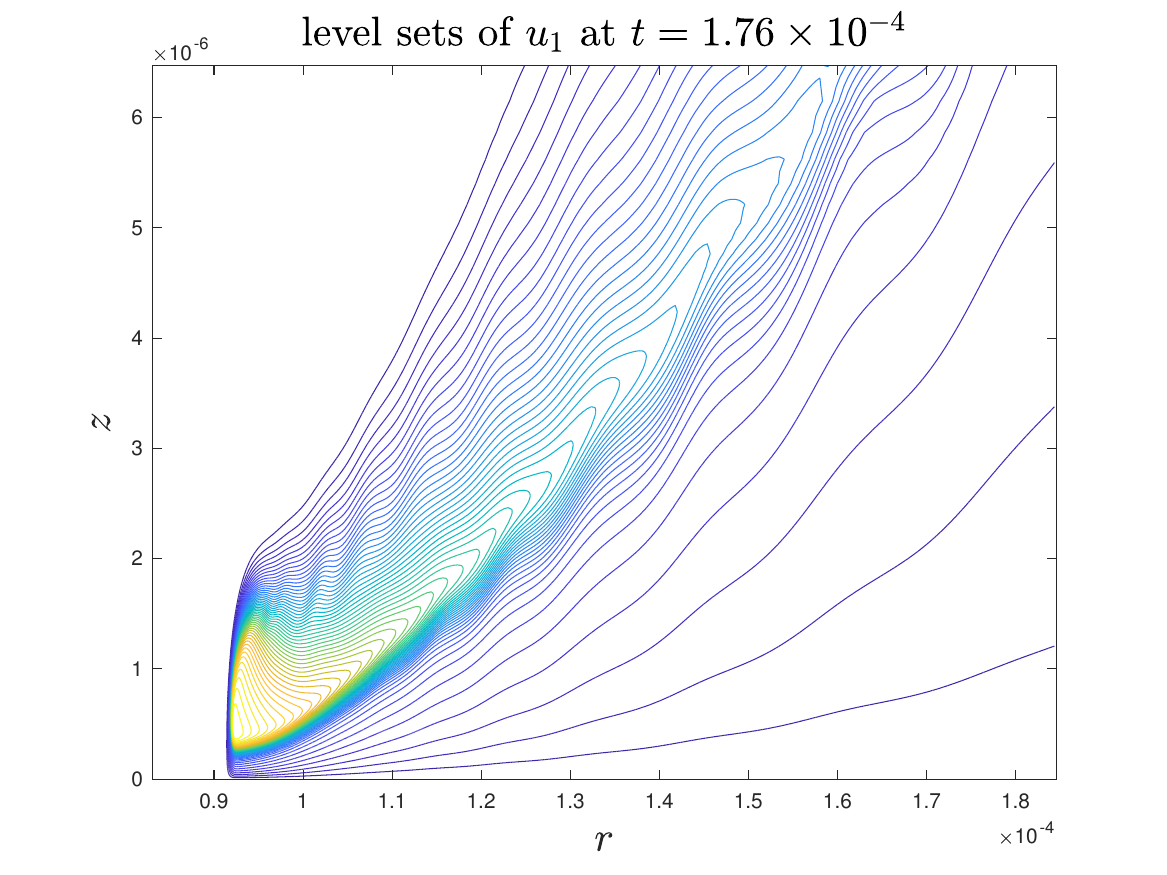}
    \includegraphics[width=0.40\textwidth]{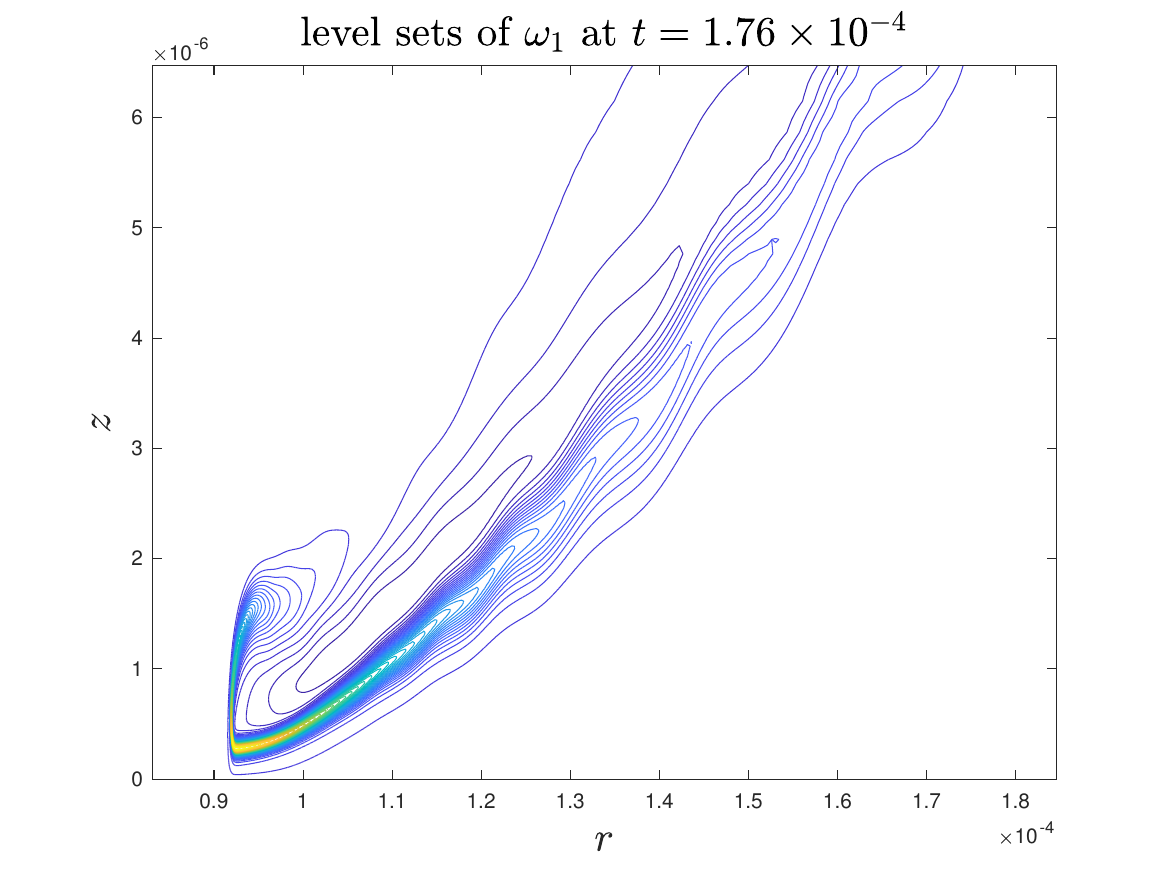} 
    \includegraphics[width=0.40\textwidth]{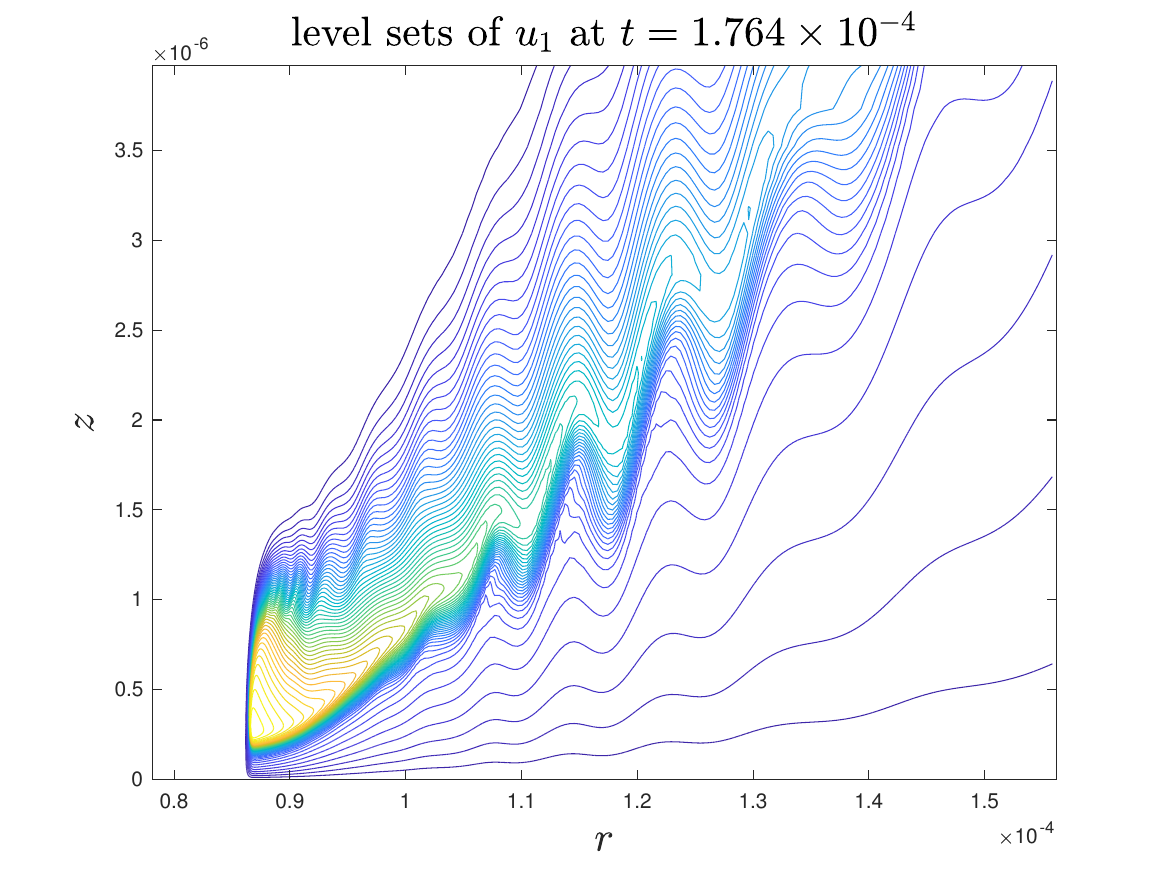}
    \includegraphics[width=0.40\textwidth]{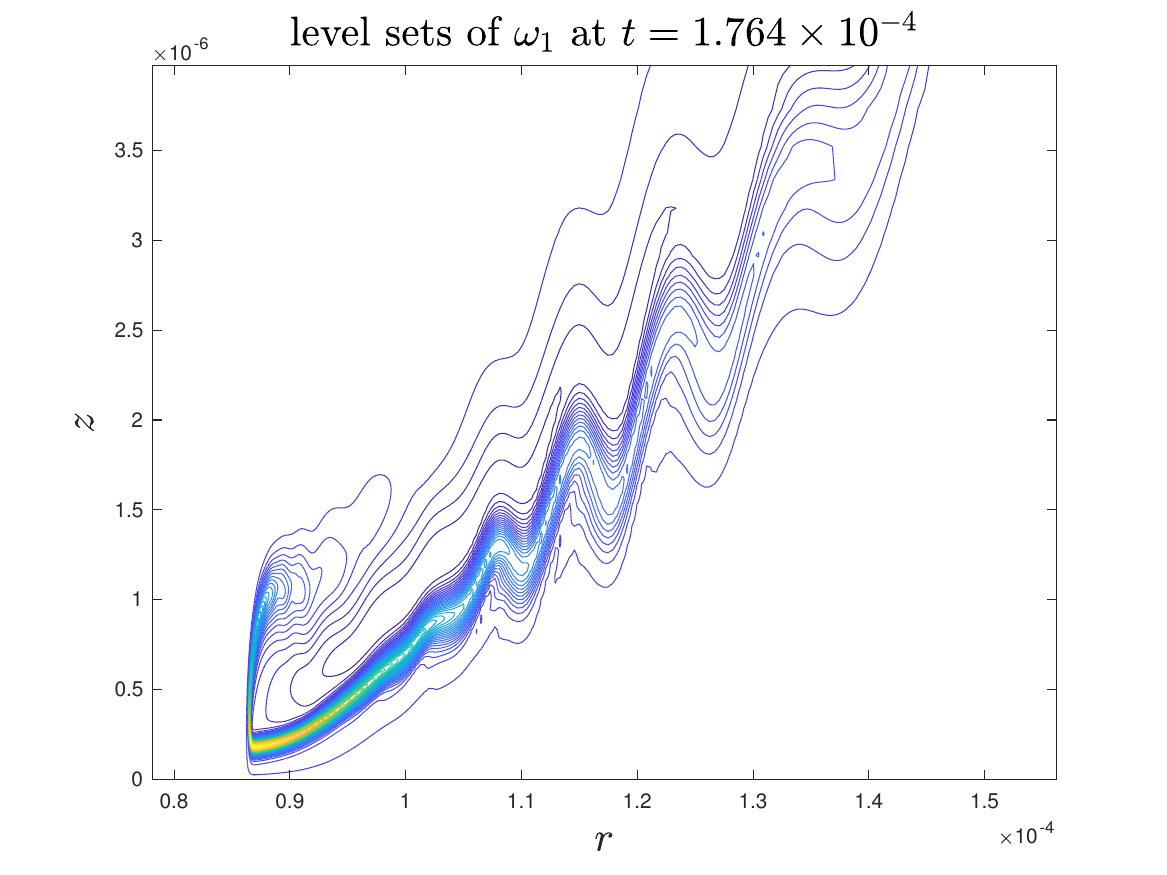} 
    \caption[Oscillation]{The level sets of $u_1,\om_1$ in Case $1$ at $t=1.76\times 10^{-4}$ (first row) and $t=1.764\times 10^{-4}$ (second row). There oscillations begin to appear at $t=1.76\times 10^{-4}$ and become severer at $t=1.764\times 10^{-4}$.}  
    \label{fig:tail_oscillation}
       \vspace{-0.05in}
\end{figure}

The first row of Figure~\ref{fig:tail_oscillation} shows the level sets of $u_1,\om_1$ (computed on mesh of size $1024\times 512$ in Case $1$) at $t=1.76\times 10^{-4}$ in a local region where the oscillations begin to appear. We can see that the tails of $u_1$ and $\om_1$ have already had mild oscillations that propagate to the far field. The second row of Figure~\ref{fig:tail_oscillation} shows that the oscillations become stronger at $t=1.764\times 10^{-4}$. We remark that these oscillations we observe are not the same as the grid-scale high frequency instability that we typically observe for the Kelvin--Helmholtz instability since the degenerate viscosity has regularized the high frequency components of the solution. In fact, the wavelengths of these oscillations remain the same as we increase the mesh size from $(n,m) = (1024,512)$ to $(n,m) = (2048,1024)$. 

\subsubsection{Source of the tail oscillations}
\begin{figure}[!ht]
\centering
	\begin{subfigure}[b]{0.32\textwidth}
    \includegraphics[width=1\textwidth]{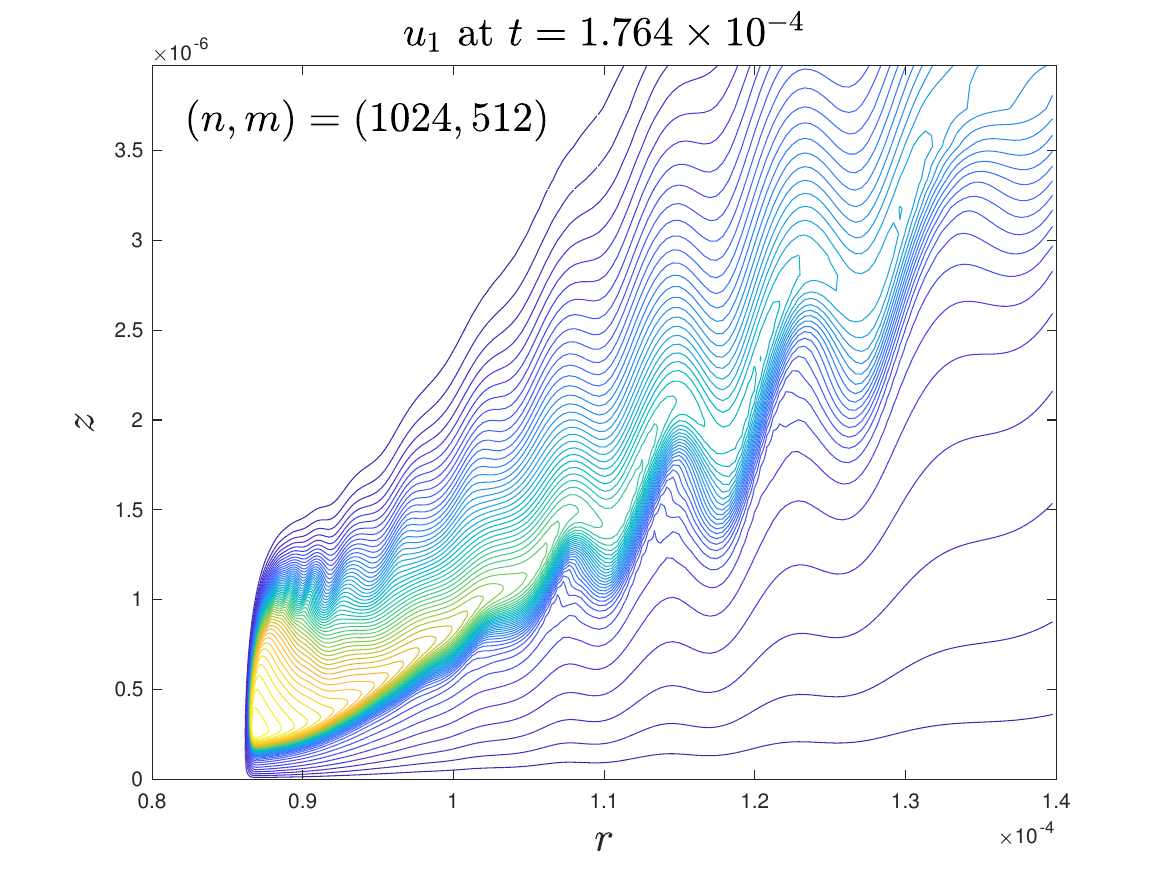}
    \caption{$(n,m)=(1024,512)$}
    \end{subfigure}
  	\begin{subfigure}[b]{0.32\textwidth}
    \includegraphics[width=1\textwidth]{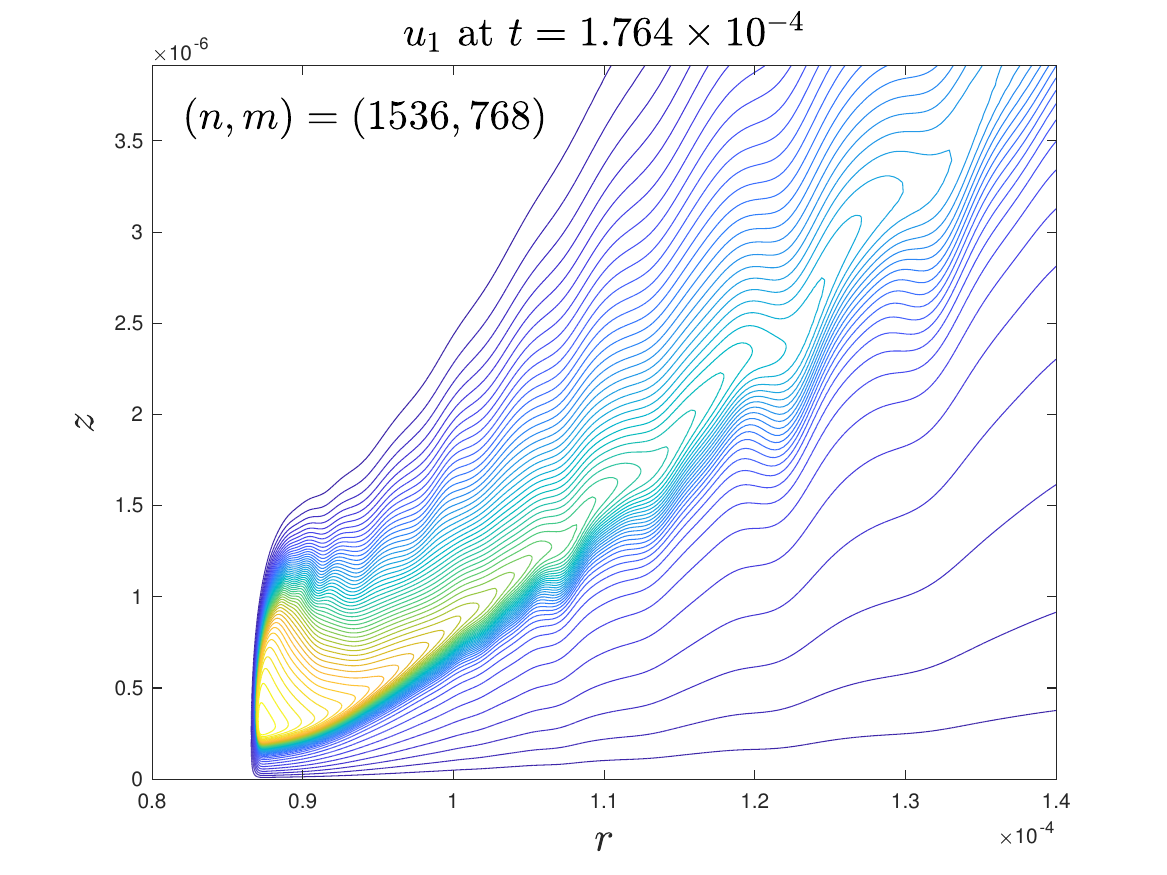}
    \caption{$(n,m)=(1536,768)$}
    \end{subfigure} 
	\begin{subfigure}[b]{0.32\textwidth}
    \includegraphics[width=1\textwidth]{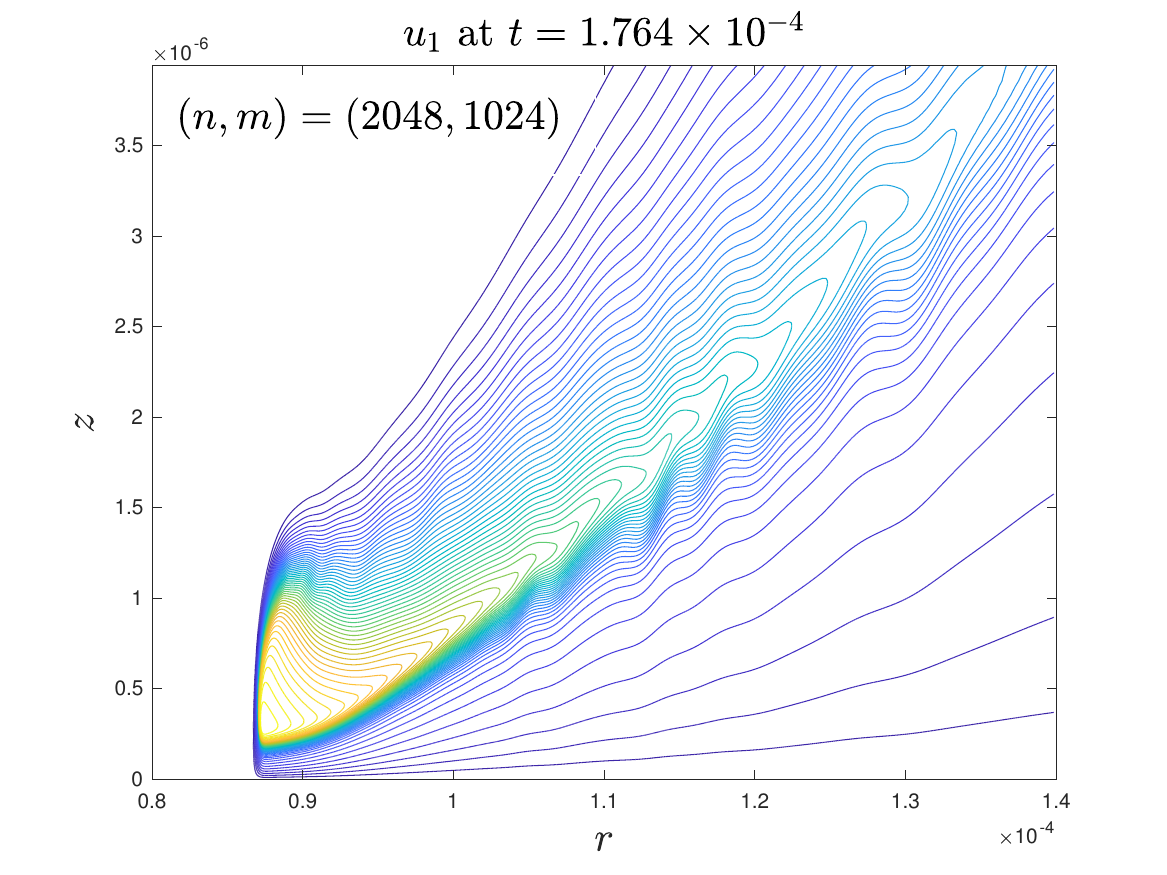}
    \caption{$(n,m)=(2048,1024)$}
    \end{subfigure}
    \caption[Oscillation vs resolution]{Contours of $u_1$ at $1.764\times 10^{-4}$ with different resolutions. The tail oscillations become weaker as the resolution of the numerical solution increases.} 
    \label{fig:Oscillation_vs_resolution}
\end{figure}

To understand the cause of the oscillations in the tail region, we first perform a resolution study to see if the solution is less oscillatory with a finer resolution. In Figure \ref{fig:Oscillation_vs_resolution} we plot the contours of $u_1$ computed with different mesh sizes at the same time $t = 1.764\times 10^{-4}$. It is clear that the oscillations become weaker as the resolution increases. This observation suggests that the tail oscillation is a consequence of under-resolution of the numerical solution in the late stage of our computation. As our adaptive mesh strategy becomes less effective in the late stage of the computation, the relative errors due to approximation and interpolation become larger and thus induce the visible oscillations.  

\subsubsection{Amplification of the oscillations}
\begin{figure}[!ht]
\centering
    \includegraphics[width=0.5\textwidth]{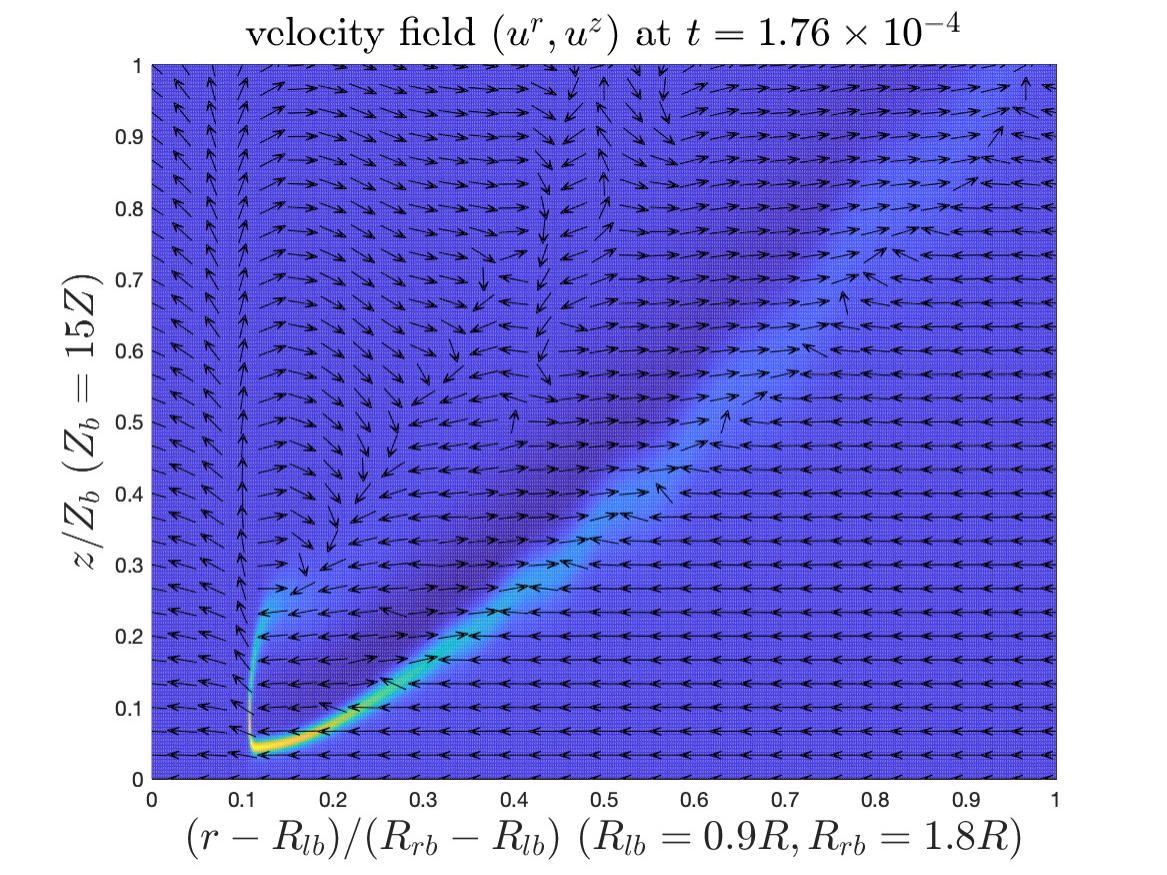}
    \caption[Shearing]{Directions of the velocity $(u^r,u^z)$ (vectors) and the profile of $\om_1$ (background color map) in rescaled coordinates at $t=1.76\times 10^{-4}$. The velocity field has a strong shearing along the tail of $\om_1$.}  
    \label{fig:tail_shear}
\end{figure}

Though the oscillations are induced by under-resolution of the numerical solution, it is likely that they are amplified by a different source of instability, as they are mainly orthogonal to the direction of the long tail. To under this particular pattern of the oscillations, we investigate the velocity field near the tail part of the solution that is shown in Figure~\ref{fig:tail_shear} (in the rescaled coordinates). The arrows represent the direction of the velocity vector field, and the background color represents the profile (magnitude) of $\om_1$. One can see that the flows on the two sides of the tail have opposite directions, meaning that there is strong a velocity shearing along the tail. This shearing can trigger unstable modes and roll into oscillations under slight perturbations (generated by under-resolution). In fact, this long thin tail structure of $\om_1$ mimics the behavior of a regularized vortex sheet. This explains the occurrence of the visible instability in the late stage of our computation when the structure of $\om_1$ becomes very thin. 

In general, viscosity shall prevent the occurrence of shearing-type instability in fluids. However, the variable viscosity coefficients \eqref{eq:viscosity_coefficient} in our equations \eqref{eq:axisymmetric_NSE_1} are degenerate in both space and time. The space-dependent part of the viscosity coefficients is degenerate at $r=0$ and $z=0$, and the time-dependent part decreases rapidly in time since $\max \om^\theta$ grows rapidly. Therefore, as the main profiles of the solution grow and travel towards the origin $(r,z)=(0,0)$, the effective viscosity coefficients in the region of interest become smaller and smaller. The viscosity in the tail part of the solution is then too weak to control the shearing instability in the late stage of the computation.

\subsection{Weakening of the nonlinear alignment} 
\begin{figure}[!ht]
\centering
	\begin{subfigure}[b]{0.32\textwidth}
    \includegraphics[width=1\textwidth]{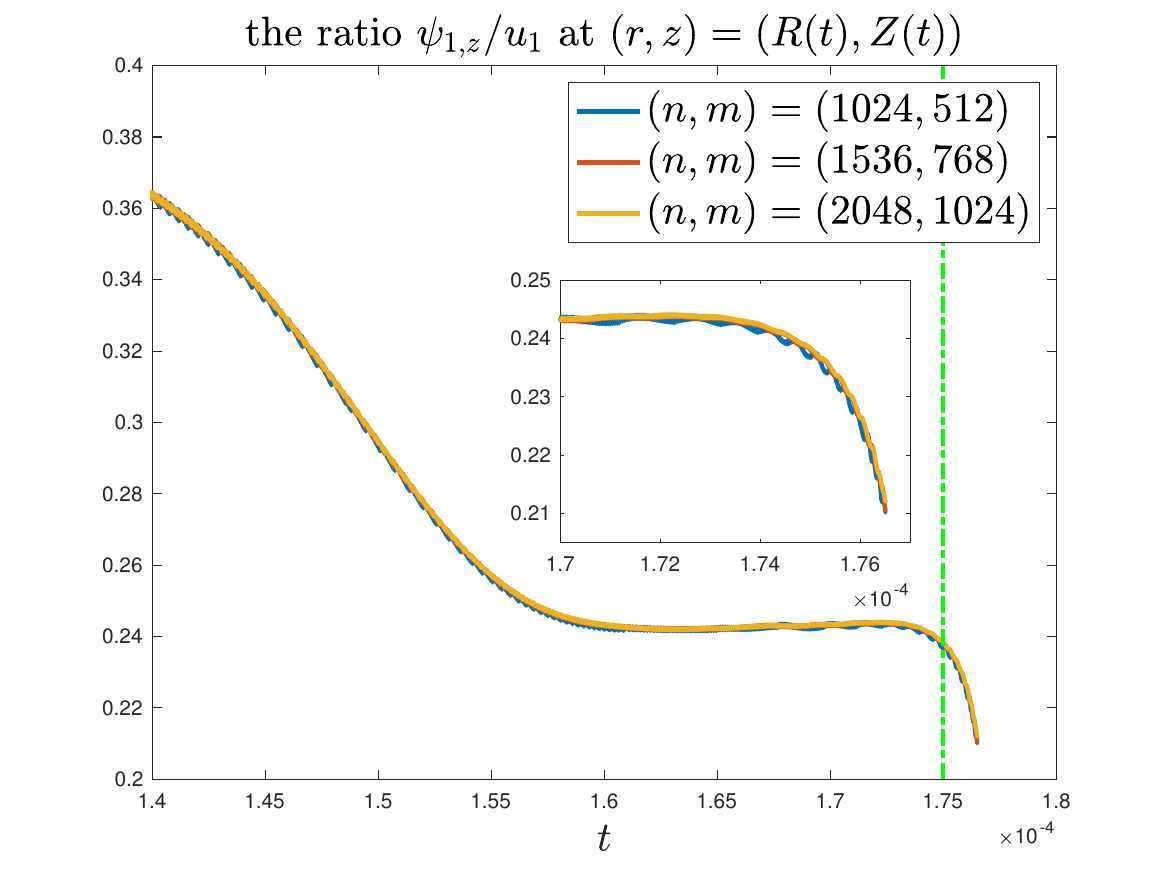}  
    \caption{$\psi_{1,z}/u_1$ till $t=1.765\times 10^{-4}$}
    \end{subfigure}
    \begin{subfigure}[b]{0.32\textwidth}
    \includegraphics[width=1\textwidth]{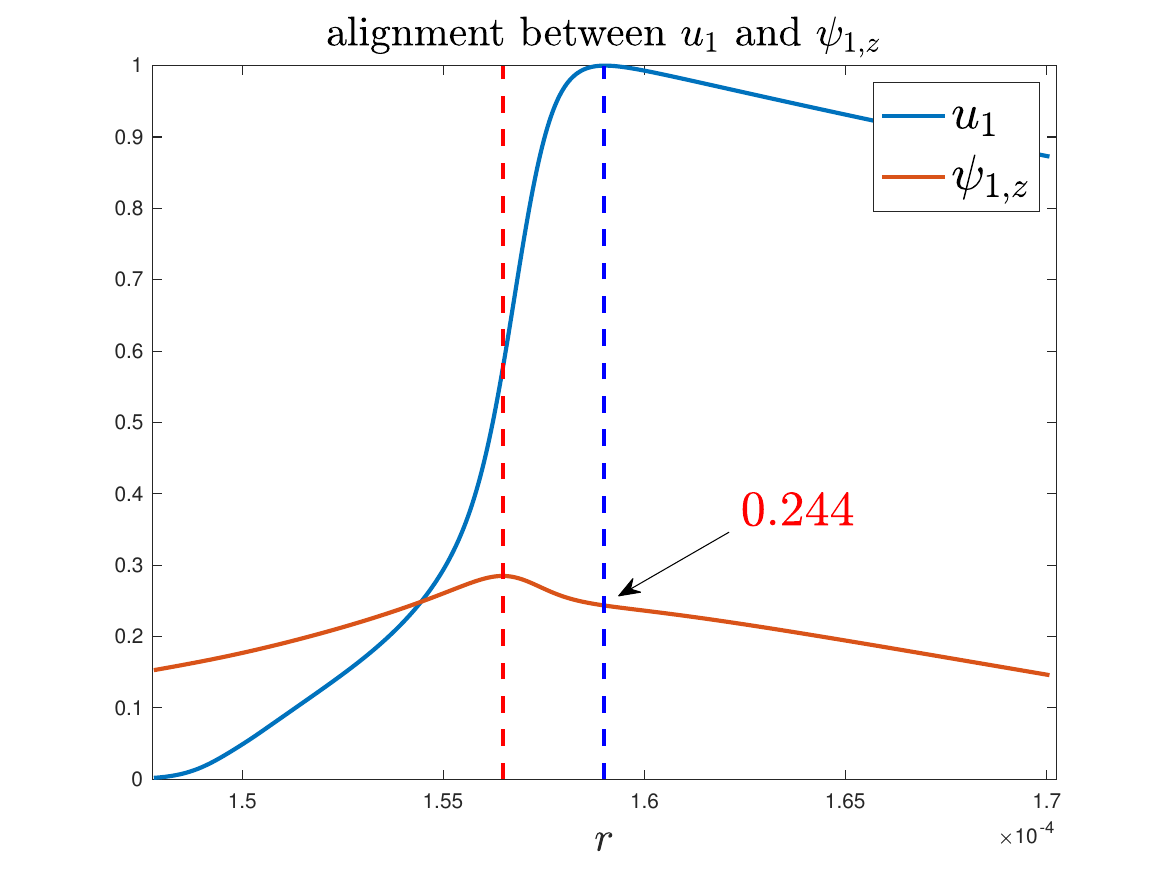}
    \caption{$t=1.7\times10^{-4}$}
    \end{subfigure}
  	\begin{subfigure}[b]{0.32\textwidth}
    \includegraphics[width=1\textwidth]{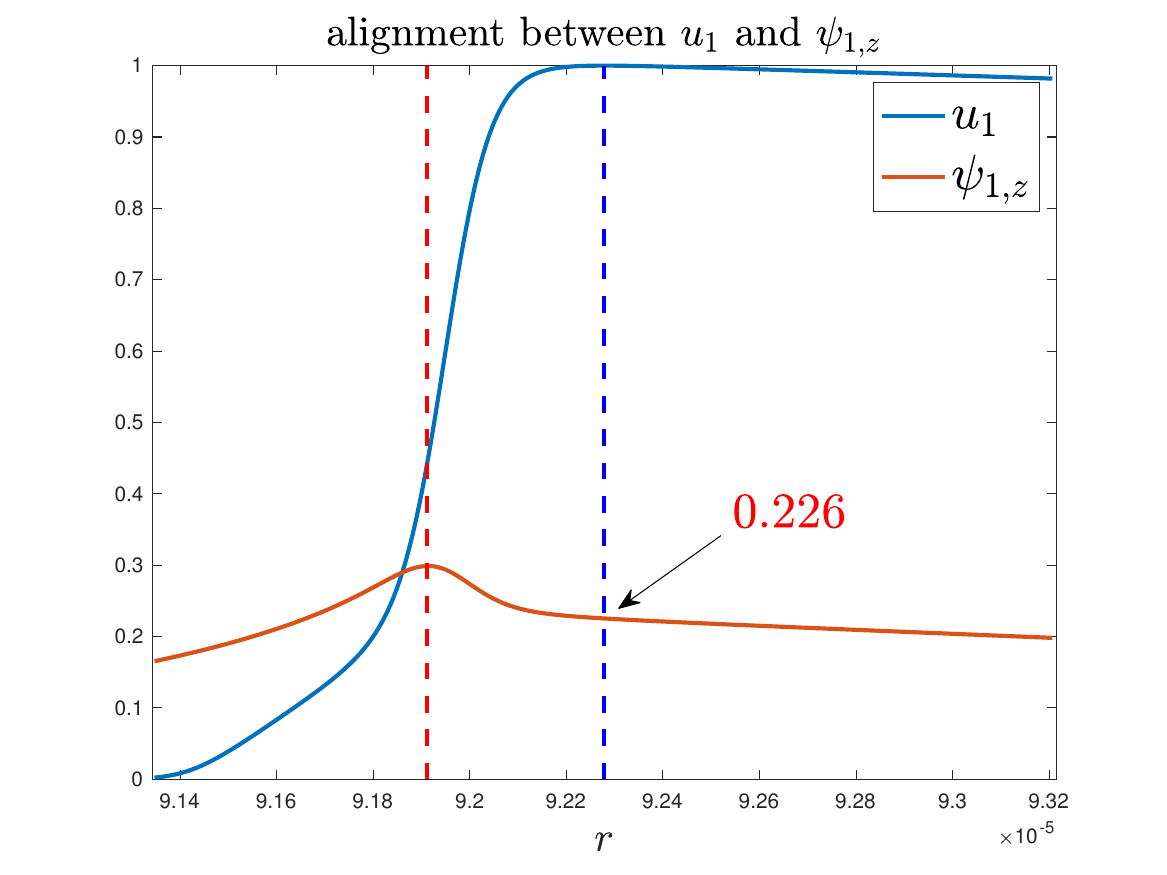}
    \caption{$t=1.76\times10^{-4}$}
    \end{subfigure}
    \caption[Alignment continued]{(a): The ratio $\psi_{1,z}/u_1$ at $(R(t),Z(t)$ as a function of time. The last time instant in the plot is $1.765\times 10^{-4}$. (b) and (c): The $r$ cross sections of $\psi_{1,z}$ and $u_1$ through the point $(R(t),Z(t)$ at $t = 1.7\times 10^{-4}$ and $t = 1.76\times 10^{-4}$, respectively. The vertical dashed lines indicate respectively the maximum points of $\psi_{1,z}$ and $u_1$ on the particular cross section.}  
     \label{fig:alignment_continued}
        \vspace{-0.05in}
\end{figure}

Another feature beyond the stale phase is the weakening of the nonlinear alignment between $\psi_{1,z}$ and $u_1$ at the maximum point of $u_1$. In Figure \ref{fig:alignment_continued}(a), we plot the ratio $\psi_{1,z}/u_1$ at the point $(R(t),Z(t))$ up to time $t=1.765\times 10^{-4}$ (which is a continuation of Figure \ref{fig:alignment}(c)). One can see that the ratio $\psi_{1,z}(R(t),Z(t))/u_1(R(t),Z(t))$ drops quickly beyond the stable phase. Increasing the resolution does not affect the decrease of the nonlinear alignment. In fact, the three curves computed with different mesh sizes almost coincide with each other, implying that this phenomenon is not a consequence of under-resolution. As we have argued in the previous subsection that the alignment between $\psi_{1,z}$ and $u_1$ is critical for the blowup mechanism. The weakening trend in the alignment seems to imply that the potential blowup is unstable beyond $t=1.75\times 10^{-4}$. 

To understand this trend, we compare the $r$ cross sections of $\psi_{1,z}$ and $u_1$ at two different time instants, one in the stable phase and one beyond the stable phase. In Figure \ref{fig:alignment_continued}(a) and (b), we plot the cross sections through the point $(R(t),Z(t))$ at $t = 1.7\times 10^{-4}$ and $t=1.76\times 10^{-4}$, respectively, and in a local region comparable to the smaller scale. The plots are rescaled in magnitude so that the maximum of $u_1$ is equal to $1$. One can see that the local profile of $u_1$ to the right of the sharp front becomes much flatter at a later time. As a result, the maximum point of $u_1$ (the blue dashed line) is farther (relative to the smaller scale) from the maximum point of $\psi_{1,z}$ (the red dashed line) on the particular cross section. This explains the weakening of the nonlinear alignment between $\psi_{1,z}$ and $u_1$ at $(R(t),Z(t))$. In fact, it can be numerically unstable to track the information at the maximum point of $u_1$ when its local profile is relatively flat. 

Therefore, the weakening of the nonlinear alignment may not be a sign of non-blowup. It is plausible that the alignment may settle down to another level after a transient stage of weakening of the alignment, as we observed beyond the warm-up phase. Such transitional stage may be due to the adjustment of the local profile of the potential blowup solution. However, we currently do not have the computational capacity that would allow us to obtain a well-resolved solution beyond the stable phase to draw a convincing conclusion regarding how the solution may develop beyond the stable phase.

\bibliographystyle{myalpha}
\bibliography{reference}

\end{document}